\documentclass[11pt]{iopart}
\usepackage{pifont}
\usepackage{amsfonts,amssymb,theorem,graphicx,listings}
\usepackage{dsfont}
\usepackage{graphics}
\usepackage{caption}
\usepackage{flafter}
\usepackage{indentfirst}
\usepackage{epsfig}
\usepackage{mathrsfs}
\usepackage{subfigure}
\usepackage{cite}
\usepackage[colorlinks=true]{hyperref}
\usepackage{cases}
\usepackage{algorithmic}
\usepackage{algorithm}
\usepackage{graphicx}
\usepackage{times}
\usepackage{color}

\usepackage{iopams}

\begin{document}

\title[Dynamic SPECT reconstruction ]{Dynamic SPECT reconstruction with temporal edge correlation}
\author{ Qiaoqiao Ding$^{a}$, Martin Burger$^{b}$, Xiaoqun Zhang$^{a,c}$}
\address{$^{a}$ School of Mathematical Sciences,  Shanghai Jiao Tong University, 800 Dongchuan Road, 200240 Shanghai, China}
\vspace{5pt}
\address{$^{b}$Institute for Computational and Applied Mathematics, University of M\"unster, Einsteinstra\ss e 62, 48149 M\"unster, Germany}
\vspace{5pt}
\address{$^{c}$ Institute of Natural Sciences, Shanghai Jiao Tong University, 800 Dongchuan Road, 200240 Shanghai, China}
\vspace{5pt}
\ead{\url{martin.burger@wwu.de}~~~\url{ xqzhang@sjtu.edu.cn}}
\begin{indented}
\item[]May 2017
\end{indented}

\begin{abstract}
 In dynamic imaging, a key challenge is  to reconstruct image sequences with high temporal resolution from strong undersampling  projections due to a
 relatively slow data acquisition speed. In this paper, we propose a variational model using the infimal convolution of Bregman distance with respect
 to total variation to  model edge dependence of sequential frames. The proposed model is solved via an alternating  iterative scheme, for which each
 subproblem is convex and can be solved by existing  algorithms. The proposed model is formulated under  both Gaussian and Poisson noise assumption
 and the simulation on  two sets of  dynamic images shows the advantage of the proposed method compared to previous methods.
 \\[2mm]
\end{abstract}

\clearpage
\section{Introduction}
\label{sec:Introduction}
Single Photon Emission Computed Tomography (SPECT) and Positron emission tomography
(PET)\cite{vardi1985statistical,kim2006minimally,ferl2007estimation,shoghi2007hybrid} are nuclear medical imaging modalities that detect the trace
concentrations of radioactively labeled pharmaceutical injected in the body within chosen volumes of interest.
 After an isotope tagged to a biochemical compound is injected into a patient¡¯s vein, the biochemical compound  travels to body organs (liver, kidney,
 brain, heart and the peripheral vascular system) through the blood stream, and is absorbed by these organs according to their affinity for the
 particular compound \cite{gullberg2010dynamic,reader20144d}. The SPECT system, usually consisting of one, two or three detector
 head(s)\cite{vanzi2004kinetic,sugihara2001estimation,marshall2001kinetic},  can record radiopharmaceutical exchange between biological compartments and
 isotope decay in the patient¡¯s body as the detector(s) rotate around the body.

As very few views can be obtained in one time interval, dynamic SPECT reconstruction is an ill-posed inverse problem with incomplete noisy data.
On assuming that  motion and deformation are negligible during the data acquiring procedure, we aim to reconstruct the dynamic radioiostrope
distribution with high temporal resolution. In fact, it is crucial to extract the actual decay of the isotope, i.e. time activity curves (TACs) of
different organ compartments \cite{carson1986maximum} \cite{formiconi1993least}, either from  projections or reconstructed images
\cite{huang1986principles,formiconi1993least,gullberg2010dynamic,zeng1995using,zan2013fast,humphries2011slow,sitek2001reconstruction,farncombe2001incorporation}.

Besides methods that reconstruct each frame of dynamic sequence independently  \cite{farncombe2001fully,feng2006simultaneous},  many approaches
\cite{vanzi2004kinetic,iida1998quantitative,komatani2004development,smith1996experimental} have been proposed to monitor the tracer concentrations over
time, under the assumption of static radioactivity concentration during the acquisition period. However, as the measuring procedure usually takes a
considerable amount of time, physiological processes in the body are  dynamic and some organs (kidney, heart) show a significant change of activity.
Hence, ignoring the dynamics of radioisotopes over the acquisition and applying  the conventional tomographic reconstruction method (such as filtered
back projection method (FBP)) yields inaccurate reconstructions with serious artifacts. Joint reconstruction approaches for dynamic imaging were
proposed in  \cite{jin2005reconstruction,jin2007dynamic}, where
dynamic images of the different time are treated collectively through motion compensation and temporal basis function.
Then, the dynamic processes of blood flow, tissue perfusion and  metabolism can be described by tracer kinetic modeling \cite{gullberg2010dynamic}.
For example, a spatial segmentation and temporal B-Splines were combined to reconstruct spatiotemporal
distributions from projection data \cite{reutter2000direct,zan2013fast}. 
In a different form,  dynamic SPECT
images are  represented by low-rank factorization models in \cite{ding2015dynamic,burger2016simultaneous}, and further constraints are enforced for the
representation coefficients and basis. In some recent work, motion of the organs has been taken into account for the dynamic CT reconstruction in
\cite{gravier2007tomographic,lamare2006respiratory,gigengack2012motion,katsevich2010accurate,roux2004exact,isola2008motion,Hahn2014efficient}.

In this paper, we propose a new variational model in which we take local and global coherence in temporal-spatial domain  into consideration.
{{The key idea is that dynamic image sequences  possess similar structures of radioactivity concentrations. In fact, the boundary of organs, which are
the locations with large gradient, are preserved or changed mildly along time.  Inspired by color Bregman TV for color image reconstruction
\cite{moeller2013color} and PET-MRI joint reconstruction  \cite{rasch2017joint}, we introduce  the infimal convolution of Bregman distance with respect
to total variation as a regularization, to obtain the dependence of edges of sequential images.}} Furthermore, based on our previous work
\cite{ding2015dynamic,burger2016simultaneous}, low rank matrix approximation  $U=\alpha B^T$ is demonstrated to be  a robust representation for dynamic
images, especially with proper regularization on the coefficient $\alpha$ and the basis  $B$ that corresponds to the TACs of different compartments.
Specially, the group sparsity on the coefficient $\alpha$ is enforced  as the concentration distribution is mixed from few basis elements for each
voxel. Finally, the proposed variational model is composed of a data fidelity term and several regularization elements, to overcome the incompleteness
and ill-posedness of the reconstruction problem.

The proposed model is solved alternatingly, with each subproblems can be solved with popular operator splitting methods at ease. In particular, the
primal-dual hybrid gradient (PDHG) algorithm  \cite{zhu08anefficient,esser2010a,chambolle2011a} is applied to solve the subproblem for the images $U$
and Proximal Forward-Backward Spliting (PFBS) for the coefficients $\alpha$ and the basis $B$. Our numerical experiments on simulated phantom  shows the
feasibility of the proposed model for reconstruction from highly undersampled data with noise, compared to conventional methods such as FBP method and
least square methods (or EM). Monte-Carlo simulation of noisy data is also performed to demonstrate the robustness of the proposed model on two phantom
images.

The paper is organized as followed:
 Section \ref{sec:Model} presents the proposed model, where section \ref{sec:Correlation}   briefly introduces the concept of infimal convolution of
 Bregman distance and models the  edge alignments of images by infimal convolution and section \ref{sec:Spatuotemporal} models the   spatiotemporal
 property of dynamic image.
 Section \ref{sec:NumericalSolution} describe the numerical algorithm.
Finally, Section \ref{sec:ComputationalResults} demonstrates  numerical results on simulated dynamic images, in the settings  of Gaussian, Poisson noise
and Monte Carlo simulation.




\section{Model}
\label{sec:Model}
The regularization of our proposed variational model is based on some properties of dynamic SPECT images. First of all, dynamic image sequences are
originated from the radioactivity concentration of few compartments in the field of view. Thus we can naturally use a low rank matrix factorization
representation for the image, with the basis related to the time activity curves, which are usually smooth,  and the coefficients are related to the
locations of the organs that are piecewise constant.
For the scenario of large motion, such as the case of  cardiac and respiratory dynamic imaging,  a proper modeling of motion correction is necessary for a reconstruction with high accuracy.  In this paper, we assume that the body movement is minor which means the boundaries of organs of image
sequence are  almost static and aim to capture the activity decay in each region. One can incorporate a registration or motion correction process in the model for the extension to the case of large motion.
 In other words, the edges of the image sequence share the similar locations. We will then use the tool of infimal convolution
of the Bregman distance with respect to total variation to enforce edge alignments. In the following, we will present the regularization terms in
details.

\subsection{Observation model}

In  dynamic SPECT, the goal is to reconstruct  a spatiotemporal radioiostrope distribution $u_t(x)$ for $x\in\Omega\subset\mathcal{R}^2$ in a given time interval. Given a sequence of projection data $f_{1}, f_{2} ,\cdots f_{T}$ with different view angles, we aim to reconstruct the samples of continuous image $u_t(x)$ at $t$-th time interval, i.e. the image sequence $u_{1}(x), u_{2}(x) ,\cdots u_{T}(x)$. If we  denote $A_{1}, A_{2},\cdots
A_{T}$  as $T$ corresponding projection matrices,  the observation model can be described as
 $$A_{t}u_{t}=f_{t},\qquad  t=1,2,\cdots,T.$$
 
 For ease of notations, we present the discrete form of  $u_t(x)$ with  $M$ pixel/voxel at each frame, and the dynamic image is represented as
 $U\in\mathbb{R}^{M
 \times T}$.  The sequence of projections are formulated in linear form:

 \begin{equation}
\mathcal{A}U=f
\label{eq:pro}
\end{equation}
where  $\mathcal{A} U=(A_{1}u_{1}, A_{2}u_{2} ,\cdots A_{T}u_{T})$, and $f=(f_{1}, f_{2} ,\cdots f_{T})$.
 In practice, the observed projection data  often inevitably accompany with noise. If white  Gaussian noise  is considered, i.e.
\begin{equation}
\mathcal{A}U+  N(0,\sigma^2)=f.\nonumber
\end{equation}
 The negative log likelihood functional leads to
\begin{equation}
H(f,\mathcal{A}U)=\frac{1}{2}\sum_{i=1}^T\|A_iu_i-f_i\|_2^2\triangleq\frac{1}{2}\|\mathcal{A}U-f\|_F^2.
\end{equation}
Similarily, if Poisson noise is considered, this term can be replaced by the log likelihood
\begin{eqnarray}\label{eq:kl}
H(f,\mathcal{A}U)=D_{\mathrm{KL}}(f,\mathcal{A}U)\nonumber\\
=\sum_{i=1}^T \bigg{(}\langle A_iu_i,\mathrm{1}\rangle-\langle f_i, \log{(A_iu_i)} \rangle\bigg{)}\triangleq\langle \mathcal{A}U,
\mathrm{1}\rangle-\langle f,\log{(\mathcal{A}U)}\rangle.
\end{eqnarray}


\subsection{Edge alignments}
\label{sec:Correlation}

For any two nonzero vectors $\vec{p},\vec{q} \in \mathcal{R}^2 $, we first define a relative  distance measure
\begin{equation}
\label{eq:d}
d(\vec{p}, \vec{q}):=\frac{\|\vec{p}\|\|\vec{q}\|-\vec{p}\cdot\vec{q}}{\|
\vec{q}\|}=\|\vec{p}\|(1-\frac{\vec{p}}{\|\vec{p}\|}\cdot\frac{\vec{q}}{\|q\|})
\end{equation}
where $\cdot$ is the standard dot product  on $\mathcal{R}^d$,  $\|\vec{p}\|$
denotes Euclidian norm. It is easy to see that if and only if $\vec{p} $ and $\vec{q} $ are parallel and point to the same direction,
$d(\vec{p},\vec{q})=0$.
In order to avoid penalizing the opposite direction,  one can define a "symmetric" distance as
\begin{equation}
\label{eq:d_sym}
\mathcal{D}(\vec{p}, \vec{q}):=\frac{\|\vec{p}\|\|\vec{q}\|-|\vec{p}\cdot \vec{q}
|}{\|\vec{q}\|}=\|\vec{p}\|(1-|\frac{\vec{p}}{\|\vec{p}\|}\cdot\frac{\vec{q}}{\|\vec{q}\|}|).
\end{equation}


%

For  two given images $u$ and $v$, we are interested in measuring the degree of parallelism of the gradients  at each pixel  as a correlation criterion
of two images. We first consider the Bregman distance with respect to the total variation of two images (in a continuous setting and assume that
$\frac{\nabla u}{|\nabla u|}$ and $\frac{\nabla v}{|\nabla v|}$ are well defined a.e.):
$$D_{\mathrm{TV}}^p(u,v)=\int_\Omega |\nabla u|(1-\frac{\nabla u}{|\nabla u|}\cdot \frac{\nabla v}{|\nabla v|})= \int_\Omega d(\nabla u, \nabla v)$$
where $p=\nabla^T\frac{\nabla v}{\|\nabla v\|}\in \partial J(v)$ and the Bregman distance is defined as
\begin{equation}\label{def:breg}
D_J^p (u,v) = J(u) - J(v)-\langle p, u-v\rangle
\end{equation}
for a convex and nonnegative functional. We can see that there is no penalty for aligned image gradients if the angle between  $\frac{\nabla u}{|\nabla
u|}$ and $\frac{\nabla v}{|\nabla v|}$ is zero, independent of the magnitude of the jump in $u$ and $v$.

To define a symmetric distance as in (\ref{eq:d_sym}), the infimal convolution of the Bregman distance in \cite{moeller2013color} is introduced to gain
independence of the direction of the gradient vector. The infimal convolution\cite{Heinz2011convex})
of two convex function $I,J :\mathcal{X} \rightarrow (-\infty, \infty]$ is defined as
$$I\Box J(x) :=\inf_{z\in \mathcal{X}}\{ I(x-z)+J(z)\}.$$

For example, if  $p,q\in \mathcal{R}^n$, the infimal convolution between the $\ell^1$ Bregman distances $D^s_{\|\cdot\|_{1}}(p, q) $ and
$D^{-s}_{\|\cdot\|_{1}}(p, -q) $ is given by
$$D^s_{\|\cdot\|_1}(p, q)\Box D^{-s}_{\|\cdot\|_{1}}(p, -q) =\sum_{i=1}^n|p_i|-|s_ip_i|$$
We can see that $D^s_{\|\cdot\|_1}(p, q)\Box D^{-s}_{\|\cdot\|_{1}}(p, -q)$ is zero if $|s(i)|=1$ for $|p_i|\neq 0$, thus the support of $p$ is
contained in the support of $q$, no matter the sign of $p(i)$ and $q(i)$. Inspired of this,  the regularization in the form of infimal convolution  of
Bregman distance with respect to total variation is considered as a measure of parallelism of the edges between two images:
 \begin{eqnarray}
R(u,v):=D_{TV}^p(u,v)\Box D_{TV}^{-p}(u,-v)
\end{eqnarray}
This formulation is  originally proposed as color Bregman total variation in \cite{moeller2013color} to couple different channels of color images.
Ehrhardt et al. \cite{ehrhardt2015joint}\cite{JulianMas}\cite{rasch2017joint} used the above derivation and the resulting measure for PET-MRI joint
reconstruction. More rigourous definition in bounded variation space and geometric interpretation, one can refer to
\cite{moeller2013color}\cite{rasch2017joint}.

For the dynamic SPECT image,  we consider the infimal convolution of the Bregman distance of the total variation to  enforce the alignment of the edge
sets of sequential frames. Specifically, let $u_i^n$ be the estimate of frame $i$ at iteration $n$,  we consider an average of the deviation of next
estimate to this image and the other frames  $u_1^n, u_2^n,\cdots u_{i-1}^n , u_{i+1}^n,\cdots u_T^n$. That is,
\begin{equation}\label{eq:infcon}
R(u_i)=w_{i,i}D_{\mathrm{TV}}^{p_i^n}(u_i,u^n_i) +\sum_{j=1,~j\neq i}^{T} w_{i,j}D_{\mathrm{TV}}^{p_{j}^n}(u_i,u^n _{j})\square
D_{\mathrm{TV}}^{-p_{j}^n}(u_i,-u^n_{j})
\end{equation}
where  $\sum_{j=1}^{T} w_{i,j}=1.$

\subsection{Low rank and sparse approximation}
\label{sec:Spatuotemporal}

The compartment model is often used to describe the concentration change of the tracer in the dynamic image  \cite{phelps1985positron}. It is assumed
that the transportation and mixture take places between the different physical compartments, such as  organ and tissue.
We impose the low-rank structure of the dynamic images by assuming that the unknown concentration distribution is a sparse linear combination of a few temporal basis functions which represent the TACs of different compartments.
  In other words, it assumes that concentration distribution of the radioisotope $u_t(x)$, for each pixel/voxel $x\in\Omega$ at time $t$  can be
  approximated as  a linear combination of  some basis TACs:
\begin{equation}
u_t(x)=\sum_{k=1}^K \alpha_k (x) B_k(t),
\end{equation}
where $B_k(t)$ denotes the TAC for $k$-th compartment at time $t$, and $\alpha_k(x)$ denotes the mixed coefficients.  This can be written in matrix form
as
 $$U=\alpha B^T,$$
 where $\alpha\in\mathbb{R}^{M\times K}$ and $B\in\mathbb{R}^{T\times K}$ with $K$ as  the number of compartments. As in general $K$ is a small number
 compared to the number of time intervals $T$, we naturally obtain a low rank matrix representation for $U$.

 Furthermore, $\alpha_{k}$ is the $k$-th column of the coefficient $\alpha$, and each element of $\alpha_k$ represents the contribution of the $k$-th
 basis to the current pixel.  For the image only have few compartments, the nonzero coefficients in  $\alpha$ is sparse. In temporal direction, we want
 to use the least number of bases, that is, as many columns of $\alpha$ as possible are entirely zero. We can use $\ell_{1,\infty}$ to describe this
 quantity
\begin{equation}
 \|\alpha\|_{1,\infty}=\sum_{j=1}^{K}\max_{i}|\alpha_{i,j}|. \nonumber
 \end{equation}
This  term is designed to select few number of  basis to represent the images $u_t$. This type of column sparsity was previously studied in
\cite{tropp2006algorithms} and in  \cite{esser2012convex} for hyperspectral image classification.

 Finally, to enfore the smoothness of the decay of radioactive distribution, we also use
$$\|\partial_t U\|^2_2=\sum_{m=1}^M\sum_{t=1}^T (u_{t+1}(m)-u_t(m))^2$$ as another regularization.

Now, we summarize the model that we propose for the reconstruction of dynamic images. Given an estimate of $U^n$ at step $n$, we propose to solve the
following reconstruction model
\begin{equation}
\label{eq:model}
\min_{U\geq0, \alpha, B}H(f,\mathcal{A}U)+\frac{\gamma}{2}\|U-\alpha B^{T}\|_2^2+\beta\|\alpha\|_{1,\infty}
+\frac{\eta}{2}\|\partial_t U\|_2^2+\lambda\sum_{i=1}^TR(u_i),
\end{equation}
where $\sum_{j=1}^{T} w_{i,j}=1, i=1,2\cdots T$, $\gamma,\beta,\eta, \lambda>0$ and $R(u_i)$ is defined in (\ref{eq:infcon}).

The combination of multi regularization aims to take into account the spatial-temporal factorization with the constraints on the basis and the sparsity
of representation coefficients, and  the alignments of edges. We will show that the proposed model is robust for overcoming  the incompleteness of
projection data and noise.
\section{Numerical algorithms}
\label{sec:NumericalSolution}
There are  three variables $U,\alpha$ and $B$ in the proposed model (\ref{eq:model}), and non-smooth and complex regularization terms are involved.  Thus it is a rather  complex  problem to solve  directly for the three  variables. We propose to solve the  nonconvex optimization problem with alternating scheme on updating the image $U$, the coefficient $\alpha$ and the basis  $B$.

In the following, we present the alternating algorithm that solves (\ref{eq:model}) with $H(f, \mathcal{A}U)$ is defined as (\ref{eq:kl}).
\begin{subequations}
 \begin{numcases}
 {}U^{n+1}=\arg\min_{U\geq0} H(f,\mathcal{A}U)+\frac{\gamma}{2}\|U-\alpha^{n} ({B^n})^T\|_2^2+\frac{\eta}{2}\|\partial_t U\|_2^2+\lambda\sum_{i=1}^TR(u_i),  \nonumber\\\label{suba}\\
\alpha^{n+1}=\arg\min_\alpha\beta\|\alpha\|_{1,\infty}+\frac{\gamma}{2}\|U^n-\alpha (B^n)^{T}\|_2^2,\label{subb}\\
B^{n+1}=\arg\min_B\frac{\gamma}{2}\|U^{n+1}-\alpha^{n+1} B^T\|_2^2.\label{subc}
\end{numcases}
\end{subequations}
\begin{itemize}
\item  In the subproblem (\ref{suba}),  given $u_i^n$ and $p_i^n\in\partial J_{TV}(u_i^n)$, $R(u_i)$ can be rewritten as followed:
\begin{eqnarray}
R(u_i)&=w_{i,i} \bigg{(}\mathrm{TV}(u_i)- \langle p_i^n, u_i\rangle \bigg{)}+\sum_{ j=1,j\neq i}^{T}\inf_{z_{ij}} w_{i,j}\bigg{(}\mathrm{TV}(u_i-z_{ij})\nonumber\\
&\quad-\langle p_{j}^n,( u_i-z_{ij})\rangle+\mathrm{TV}(z_{ij})+\langle p_{j}^n, z_{ij}\rangle\bigg{)}\nonumber\\
&=w_{i,i} \bigg{(}\mathrm{TV}(u_i)- \langle q_i^n,\nabla u_i\rangle \bigg{)}+\sum_{ j=1,j\neq i}^{T}\inf_{z_{ij}} w_{i,j}\bigg{(}\mathrm{TV}(u_i-z_{ij})\nonumber\\
&\quad-\langle q_{j}^n,\nabla( u_i-z_{ij})\rangle+\mathrm{TV}(z_{ij})+\langle q_{j}^n, \nabla z_{ij}\rangle\bigg{)}\\
&\triangleq\inf_{\{z_{ij}\}_{j=1}^T, j\neq i}\tilde{R}(u_i,z_{i,\cdot}),\nonumber
\end{eqnarray}
where $\nabla^T q_i^n=p_i^n, i=1,\cdots, T.$ Then, the subproblem (\ref{suba}) is reformulated as
 \begin{eqnarray}
\min_{U\geq0, Z} H(f,\mathcal{A}U)+\frac{\gamma}{2}\|U-\alpha^{n} ({B^n})^T\|_2^2+\frac{\eta}{2}\|\nabla_t U\|_2^2+\lambda\sum_{i=1}^T\tilde{R}(u_i,z_{i,\cdot}).\label{equalsuba}
 \end{eqnarray}
We solve problem (\ref{equalsuba})  by primal-dual hybrid gradient (PDHG)   Algorithm \cite{esser2010a,chambolle2011a}.
The problem formulation is as follows,
 \begin{eqnarray}
  \min_X\max_Y -F^\ast(Y)+G(X)+\langle \mathcal{K}X , Y\rangle,
\end{eqnarray}
where $\mathcal{K}: \mathcal{U} \rightarrow \mathcal{V}$ is a linear and continuous operator between two finite dimensional vector spaces $\mathcal{U}$ and $\mathcal{V}$. $F :\mathcal{V}\rightarrow [0,+\infty] $ and $G: \mathcal{U} \rightarrow [0,+\infty] $ are proper, convex and lower semi-continuous functions. $F^\ast$ is the conjugate of $F$ and $\mathcal{W}^\ast=\mathcal{V}$.
The PDHG iteration is,
\begin{subequations}
\label{eq:pdhg}
 \begin{numcases}
   {}Y^{k+1}=(I+\sigma\partial F^{\ast})^{-1}(Y^k+\sigma  \mathcal{K}\bar{X}^k)\label{eq:pdhg1}\\
   X^{k+1}=(I+\tau\partial G)^{-1}(X^k-\tau \mathcal{K}^{\ast}Y^{k+1})\label{eq:pdhg2}\\
   \bar{X}^{k+1}=X^{k+1}+\theta(X^{k+1}-X^{k}).\label{eq:pdhg3}
 \end{numcases}
 \end{subequations}
For our model,  we have
\begin{eqnarray}
F(g,b,d_{ii}, d_{ij}^+, d_{ij}^-)&:= H(f,g)+\frac{\eta}{2}\|b\|_2^2 +\lambda\sum_{i=1}^T \Bigg{(} w_{i,i} \bigg{(}\|d_{ii}\|_1- \langle q_i^n,d_{ii}\rangle \bigg{)}\nonumber\\
 &+\sum_{ j=0,j\neq i}^{T}w_{i,j}\bigg{(}\|d^+_{ij}\|_1-\langle q_{j}^n,d^+_{ij}\rangle+\|d^-_{ij}\|_1+\langle q_{j}^n, d^-_{ij}\rangle\bigg{)}\Bigg{)}, \nonumber
  \end{eqnarray}
where $g=\mathcal{A}U; b=\partial_t U; d_{ii}=\nabla u_i; d_{ij}^+=u_i-z_{ij}; d_{ij}^-=\nabla z_{ij}$ and
\begin{eqnarray}
G(U):=&\iota_{\{U\geq 0\}}(U)+\frac{\gamma}{2}\|U&-\alpha^n( B^n)^T\|_2^2.\nonumber
\end{eqnarray}
where $\iota_{\{U\geq 0\}}(U)$ is the characteristic function of the set $\{U\geq 0\}$.
According to (\ref{eq:pdhg}), we obtain the iterative scheme for each subproblem.
\begin{itemize}
\item For the subproblem (\ref{eq:pdhg1}), the  dual variables can be updated as
  \begin{subequations}
 \begin{numcases}
   {}g^{k+1} = \frac{\tilde{g}+1}{2}-\sqrt{(\frac{\tilde{g}+1}{2})^2+\sigma f-\tilde{g}}, ~~\tilde{g}=g^k+\sigma A\bar{U}^k,\nonumber\\
  b^{k+1}= \frac{\eta}{\eta+\sigma}\tilde{b},~~ \tilde{b}=b^k+\sigma\partial_t\bar{U}^k,\nonumber\\
  (d_{ii})^{k+1}=w_{i,i}\Pi_{B^{\infty}_{(1)}}(\tilde{d}_{ii}/w_{i,i}-q^n_{i})  +w_{i,i}q^n_{i},  ~~ \tilde{d}_{ii}=(d_{ii})^k+\sigma\nabla \bar{u}_i^k,\nonumber\\
  (d^+_{ij})^{k+1}=w_{i,j}\Pi_{B^{\infty}_{(1)}}(\tilde{d}^+_{ij}/w_{i,j}-q^n_{j})  +w_{i,j}q^n_{j}, ~~ \tilde{d}^+_{ij}=(d^+_{ij})^k+\sigma\nabla (\bar{u}_i^k-\bar{z}^k_{ij}),\nonumber\\
  (d^-_{ij})^{k+1}=w_{i,j}\Pi_{B^{\infty}_{(1)}}(\tilde{d}^-_{ij}/w_{i,j}+q^n_{j})  -w_{i,j}q^n_{j}, ~~ \tilde{d}^-_{ij}=(d^-_{ij})^k+\sigma\nabla \bar{z}^k_{ij},\nonumber
 \end{numcases}
 \end{subequations}
where $i=1\cdots T$, $j=1,\cdots,i-1,i+1,\cdots, T$ and $\Pi_{B^{\infty}_{(1)}}(z)=\frac{z}{\max{(1,z)}}$.
 \item For the subproblem (\ref{eq:pdhg2}), the primal variables are updated as
 \begin{subequations}
 \begin{numcases}
   {}U^{k+1}=\mathrm{argmin}_{U}\iota_{\{U\geq 0\}}(U)  +\frac{\gamma}{2}\|U-\alpha^n(B^n)^T\|_2^2+\frac{1}{2\tau}\|U-U'\|_2^2\nonumber\\
(z_{ij})^{k+1}= (z_{ij})^{k}-\tau \nabla^T\bigg{(}(d^-_{ij})^{k+1}-(d^+_{ij})^{k+1}\bigg{)},\nonumber
 \end{numcases}
 \end{subequations}
 where $U'=U^k-\tau\Bigg{(}A^T g^{k+1}+\partial_t^T b^{k+1}+\nabla^T\Big{(} d^{k+1}+\sum_i^T (d^+_{i\cdot })^{k+1}\Big{)}\Bigg{)}$ .
 It is easy to see that $U^{k+1}$ is the nonnegative projection of $(\frac{\gamma \alpha^n (B^T)^n+\frac{U'}{\tau}}{\gamma+\frac{1}{\tau}}).$
 \item The update on the subproblem (\ref{eq:pdhg3}) is straightforward.
\end{itemize}

We note that  the variables $q_i^n$ can be updated by the following property
 \begin{eqnarray} \mathcal{K}X &\in \partial F^\ast( Y)\Longrightarrow  Y\in \partial F(\mathcal{K} X),\nonumber\\
\Longrightarrow&d_{ii}\in \partial_{\nabla u_i} w_{i,i}\big{(}\|\nabla u_i\|_1-\langle q_i^n,\nabla u_i\rangle\big{)},\nonumber\\
\Longrightarrow& d_{ii}+ w_{i,i}  q^n_i \in w_{i,i} \partial_{\nabla u_i} \|\nabla u_i\|_1.\nonumber
\end{eqnarray}
We update $q_i^{n+1}=\frac{d^n_{ii}}{w_{i,i}}+ q^n_i$.

The overall algorithm for the $U$ subproblem is summarized in  Algorithm \label{alg:CP}.

 \begin{algorithm}[H]
 \caption{ PDHG Algorithm for $U$ }
\begin{algorithmic}
\STATE  Input: $\sigma, \tau $, $\lambda, w_{i,j},  i,j=1,2\cdots T$,
\STATE  Initial: $U^0=U^n$, $z_{ij}^0$,  $g^0$, $b^0$, $d_{ii}^0$, $(d_{ij}^+ )^0$, $(d_{ij}^-)^0$,   $i=1,2\cdots T$, $j=1,2\cdots T, j\neq i $.
\WHILE { not satisfy stopping conditions}
\STATE  \bf{Dual update}:
\STATE $g^{k+1} = \frac{g^k+\sigma A\bar{U}^k+1}{2}-\sqrt{(\frac{g^k+\sigma A\bar{U}^k+1}{2})^2+\sigma f-g^k-\sigma A\bar{U}^k}.$
\STATE $b^{k+1}= \frac{\eta}{\eta+\sigma}(b^k+\sigma\partial_t\bar{U}^k).$
\FOR { $ i=1,2 \cdots T $ }
\STATE $(d_{ii})^{k+1}=w_{i,i}\Pi_{B^{\infty}_{(1)}}(\frac{(d_{ii})^k+\sigma\nabla \bar{u}_i^k}{w_{i,i}}-q^n_{i})  +w_{i,i}q^n_{i},$
\STATE $(d^+_{ij})^{k+1}=w_{i,j}\Pi_{B^{\infty}_{(1)}}(\frac{(d^+_{ij})^k+\sigma\nabla (\bar{u}_i^k-\bar{z}^k_{ij})}{w_{i,j}}-q^n_{j})  +w_{i,j}q^n_{j},$
\STATE $(d^-_{ij})^{k+1}=w_{i,j}\Pi_{B^{\infty}_{(1)}}(\frac{(d^-_{ij})^k+\sigma\nabla \bar{z}^k_{ij}}{w_{i,j}}+q^n_{j})  -w_{i,j}q^n_{j}, $
\ENDFOR
\STATE  \bf{Primal update}:
\STATE $U'=U^k-\tau\Bigg{(}A^T g^{k+1}+\partial_t^T b^{k+1}+\nabla^T\Big{(} d^{k+1}+\sum_i^T (d^+_{i\cdot })^{k+1}\Big{)}\Bigg{)},$
\STATE $U^{k+1}= \Pi_{\{U\geq 0\}}(U),$ where $U=\frac{\frac{U'}{\tau}+\gamma \alpha^n (B^T)^n}{\tau+\frac{1}{\gamma}}.$
\STATE  $ (z_{ij})^{k+1}= (z_{ij})^{k}-\tau \nabla^T\bigg{(}(d^-_{ij})^{k+1}-(d^+_{ij})^{k+1}\bigg{)},$ $i,j=1,2\cdots T, j\neq i.$
 \STATE  \bf{Relaxation}
 \STATE  $\bar{U}^{k+1}=U^{k+1}+\theta(U^{k+1}-U^{k})$
 \STATE  $(\bar{z}_{ij})^{k+1}=({z}_{ij})^{k+1}+\theta(({z}_{ij})^{k+1}-({z}_{ij})^{k+1} ),$ $i,j=1,2\cdots T, j\neq i. $
 \ENDWHILE
 \STATE $U^{n+1}=U^{k+1}.$
\end{algorithmic}
\label{alg:CP}
\end{algorithm}

\item The subproblem  (\ref{subb}) can be solved by PFBS (Proximal Forward Backward Spliting)  
\begin{subequations}
 \begin{numcases}
 {}\alpha^{k+\frac{1}{2}} = \alpha^k-\tau(U^n-\alpha^k (B^n)^T)B^n, \\
 \alpha^{k+1}=\arg\min_\alpha\beta\|\alpha\|_{1,\infty}+\frac{1}{2\tau}\| \alpha -\alpha^{k+\frac{1}{2}}\|_2^2.\label{PFBSb}
  \end{numcases}
\end{subequations}
(\ref{PFBSb}) can be rewritten as
\begin{eqnarray}
\alpha^{k+1}=\arg\min_{\alpha}\beta\sum_{j=1}^{K}\max_{i}|\alpha_{i,j}|+\frac{1}{2\tau}\| \alpha-\alpha^{k+\frac{1}{2}}\|_{2}^{2}.\nonumber
 \end{eqnarray}

 We can see that $\alpha$ is separable in column from the formulation above. Thus, we can solve   for each column $\alpha_{j}$:
 \begin{eqnarray}
\alpha_{j}^{k+1}=\arg\min_{\alpha_{j}}\beta\max_{i}|\alpha_{i,j}|+\frac{1}{2\tau}\| \alpha_{j}-\alpha^{k+\frac{1}{2}}_{j}\|_{2}^{2}.
 \end{eqnarray}
 The solution of this problem is given by Moreau decomposition \cite{moreau1965proximite,rockafellar1996convex,ekeland1999convex}.
 The Moreau decomposition of  a convex function $J$ in $\mathbb{R}^{n}$ is defined as
 \begin{eqnarray}
x=\arg\min_{u\in\mathbb{R}^{n}}J(u)+ \frac{1}{2\sigma}\| u-x\|_{2}^{2}
+\sigma\arg\min_{p\in\mathbb{R}^{n} }J^{\ast}(p)+\frac{\sigma}{2}\|p-\frac{x}{\sigma} \|_{2}^{2},\nonumber
 \end{eqnarray}
 where  $J^{\ast}(p)$ is  the conjugate function of $J$  \cite{rockafellar1996convex,arnol1989mathematical}. Let $J(\alpha)\triangleq\beta\max_{i}|\alpha_{i,j}|$, then the conjugate function
 \begin{eqnarray}
 J^{\ast}(p)&=&\sup_{\alpha_{j}}\langle p,\alpha_{j}\rangle-\beta\max_{i}|\alpha_{i,j}|\nonumber \\
 &=&\sup_{\alpha_{j}}(\max_{i}|\alpha_{i,j}|)(\| \max(p,0)\|_{1}-\beta) \nonumber\\
 &=& \left \{
 \begin{array}{ll}
  \infty {} & \textrm{otherwise}\\
 0  &\textrm{if $\| \max(p,0) \|_{1}\leq \beta$}.\\
 \end{array}
 \right.\label{MD}
 \end{eqnarray}
 From equation (\ref{MD}), we can see that $J^{\ast}(p)$ is the set characteristic  function of the convex set $C_{\beta}=\{p\in \mathbb{R}^{n}, \| \max(p,0)\|_{1}\leq \beta\}$.  By Moreau decomposition,  $\alpha^{k+1}$ is then given by
 \begin{eqnarray}
 \alpha^{k+1}={\alpha}^{k+\frac{1}{2}}-\Pi_{C_{\delta}}({\alpha}^{k+\frac{1}{2}}),
 \end{eqnarray}
 where $\Pi_{C_{\delta}}({\alpha}^{k+\frac{1}{2}})$ is the projection of each column of $ {\alpha}^{k+\frac{1}{2}} $ onto $ C_{\delta}$ and $\delta =\tau\beta$.
 \begin{algorithm}[H]
\caption{PFBS for $\alpha$}
\begin{algorithmic}
\STATE  Input: $\tau>0, \beta>0, \delta=\tau\beta$,
\STATE  Initial: $\alpha^{0}=\alpha^{n}$,
\FOR { $ k=0,1,2 \cdots K $ }
\STATE  $\alpha^{k+\frac{1}{2}} = \alpha^k-\tau (U^{n+1}-\alpha^k ({B^n})^T)B^n, $
\STATE  $\alpha^{k+1}={\alpha}^{k+\frac{1}{2}}-\Pi_{C_{\delta}}({\alpha}^{k+\frac{1}{2}}),$
\ENDFOR
\STATE $\alpha^{n+1}=\alpha^{k+1}.$
\end{algorithmic}
\label{alg:PFBS}
\end{algorithm}
\item The subproblem (\ref{subc}) is a quadratic problem, which is equal to solve the following problem:
\begin{equation}
(\alpha^{n+1})^T\alpha^{n+1} B^T=(\alpha^{n+1})^T U^{n+1}.
\label{LS}
\end{equation}
Problem (\ref{LS}) is a linear system that can be solved  by an iterative method, such as conjugate gradient.
 If $(\alpha^{n+1})^T\alpha^{n+1}$ is invertible, one can also compute the inverse of the matrix directly as its size is $K\times K $,.
\end{itemize}

Overall, the  entire algorithm is summarized as follows
\begin{algorithm}[H]
\caption{ }
\label{alg:All}
\begin{algorithmic}
\STATE Input: $ \beta,  \gamma, \eta, \varepsilon$, $w_{i,j}, i,j=1,2 \cdots T.$
\STATE Initial: $U^0, \alpha^{0}, B^{0}.$
\FOR { $ n=0,1,2 \cdots $ }
\STATE Solve for $U^{n+1}$ by algorithm \ref{alg:CP},
\STATE Solve for $\alpha^{n+1}$ by algorithm \ref{alg:PFBS},
\STATE Solve linear system (\ref{LS}), 
\STATE Update $q^{n+1}_i=\frac{d^n_{ii}}{w_{i,i}}+  q^n_i, i=1,2\cdots T$,
\STATE $\gamma=\gamma\varepsilon$.
\ENDFOR
\end{algorithmic}
\end{algorithm}

\section{Simulation results}
\label{sec:ComputationalResults}
We present the simulation results to validate the proposed model and algorithm. First of all, for the edge correlation regularisation term¡¡ $R(u_i)$, in practice, it is computational expensive and also unnecessary to consider the correlation of each frame to all the other frames. In our computational results, we only consider the edge correlation to the  last iterate, and the former two images and the later two images, that is $w_{i,j}\neq 0$ when $j=i-2,i-1,i,i+1,i+2$. For the boundary frames, we have the following setting, when $i=1$, $j=1,2,3$;  $i=2$, $j=1,2,3,4$; $i=T-1$, $j=T-3,T-2,T-1,T$  and $i=T$, $j=T-2,T-1,T$, $w_{i,j}\neq 0$.

\begin{figure}[ht]
\begin{center}
\subfigure[Simulated image phantom]{
\includegraphics[width=.3\linewidth,height=.3\linewidth]{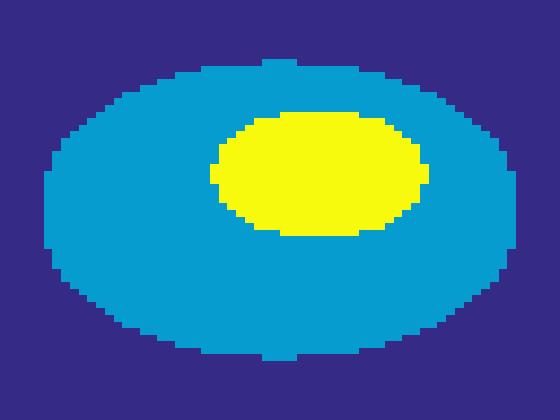}}
\subfigure[Simulated concentration curves]{
\includegraphics[width=.45\linewidth,height=.3\linewidth]{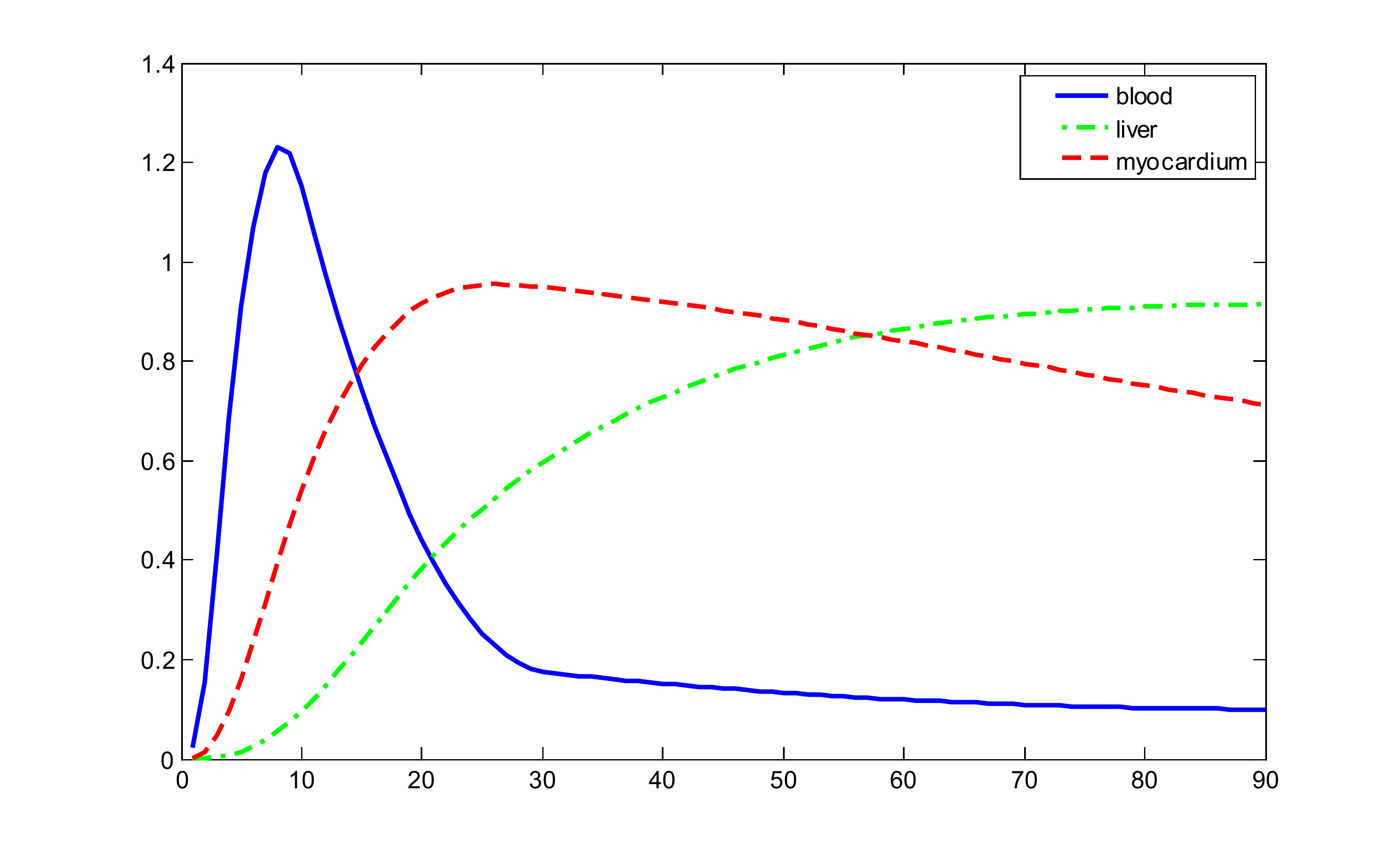}}
\end{center}
\caption{Phantom of tested data and  temporal concentration curves in all subregions.}
\label{fig:Ellreal}
\end{figure}

The proposed method is tested on numerical phantoms for a proof of concept study.
We simulate $90$  image frames of size $64\times 64$ and $2$ projections per frame.
Three time  activity curves (TAC) for blood, liver and myocardium, previously used in \cite{zan2013fast} (see Figure \ref{fig:Ellreal}), are used to simulate the dynamic images. The first simulated dynamic phantom is composed of two ellipses. In temporal direction, the positions of the  two ellipses  are stationary while the intensity in 90 frames  within the region of each ellipse is generated according to the  TAC of blood or liver. The projections are generated by using Radon transform sequentially performed for each frame.

\begin{figure}[ht]
\begin{center}
\subfigure[Simulated image phantom]{
\includegraphics[width=.3\linewidth,height=.3\linewidth]{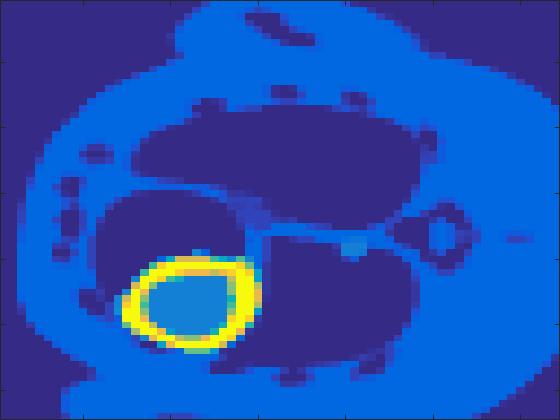}}
\subfigure[Simulated concentration curves]{
\includegraphics[width=.45\linewidth,height=.3\linewidth]{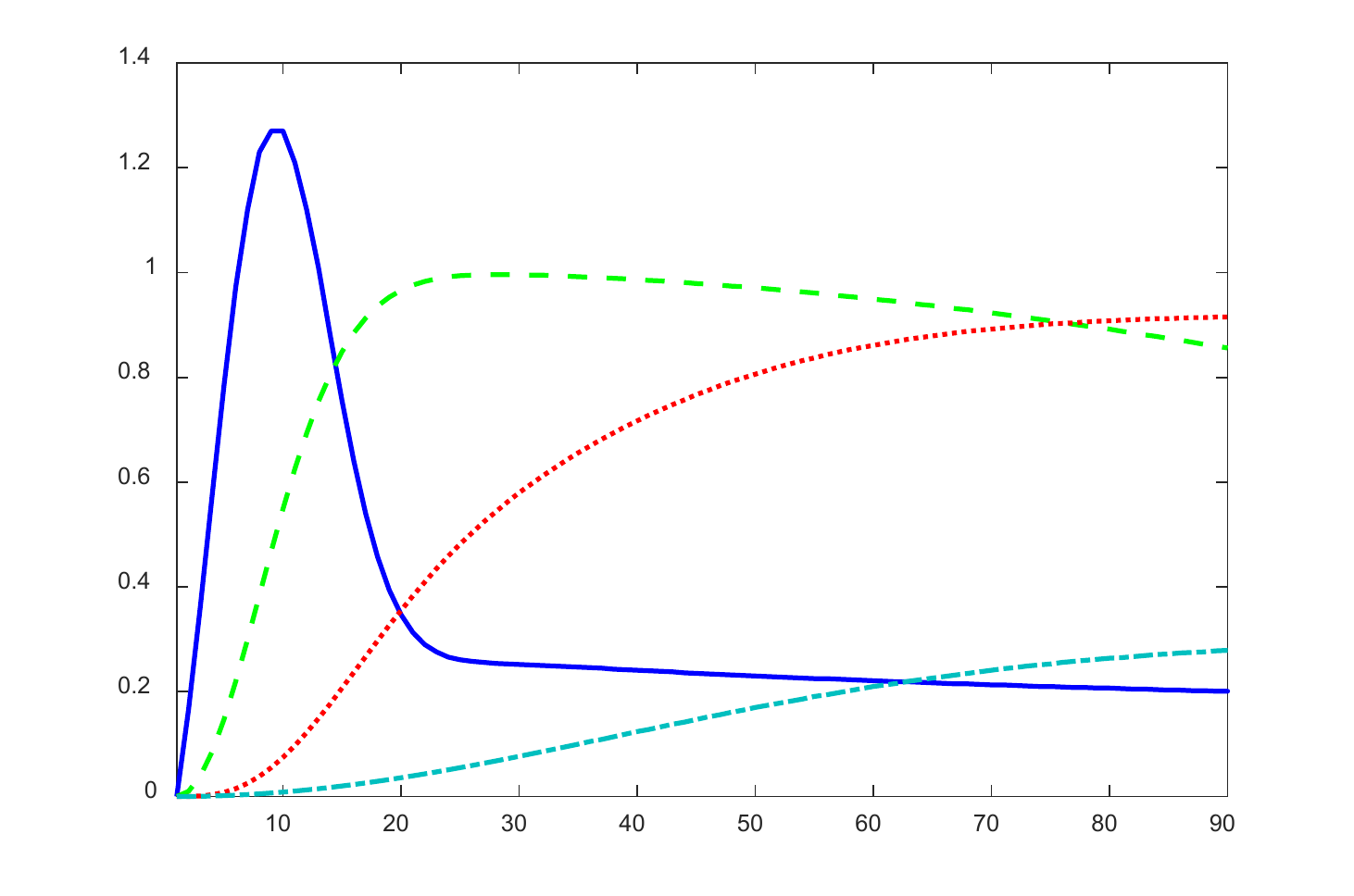}}
\end{center}
\caption{Phantom of  rat's abdomen  and  temporal concentration curves in all subregions.}
\label{fig:Liverreal}
\end{figure}

The second numerical experiment is performed on a synthetic image simulating rat's
abdomen, where the bright region represents the heart of a rat. We use the TAC in Figure \ref{fig:Liverreal} to simulate the dynamic images.

\subsection{Gaussian noise}
We compare our method with the  filtered back projection (FBP) method, the results by alternatingly solving  least square model $\arg\min_{\alpha, B}\|A \alpha B^{\top}-f\|_{F}^{2}$, our previous model, sparsity enforced matrix factorization(SEMF) proposed in \cite{ding2015dynamic}. As for the initial value of $U$, $\alpha$ and $B$,  we use uniformed  B-spline, $B\in\mathbb{R}^{90\times20}$ as initial basis to solve $\arg\min_{\alpha,B}\|A\alpha B^{\top}-f\|_{F}^{2}$  for $\alpha$ and $B $.  Then, the  same $\alpha$, $B$, and  $U=\alpha B^T$ are used as the initialization of $U$,$\alpha$ and $B$ to solve our proposed model.

\begin{figure}
\begin{tabular}{c@{\hspace{2pt}}c@{\hspace{2pt}}c@{\hspace{2pt}}c@{\hspace{2pt}}c@{\hspace{2pt}}c@{\hspace{2pt}}c@{\hspace{2pt}}c@{\hspace{2pt}}c@{\hspace{2pt}}c}
\includegraphics[width=.1\linewidth,height=.1\linewidth]{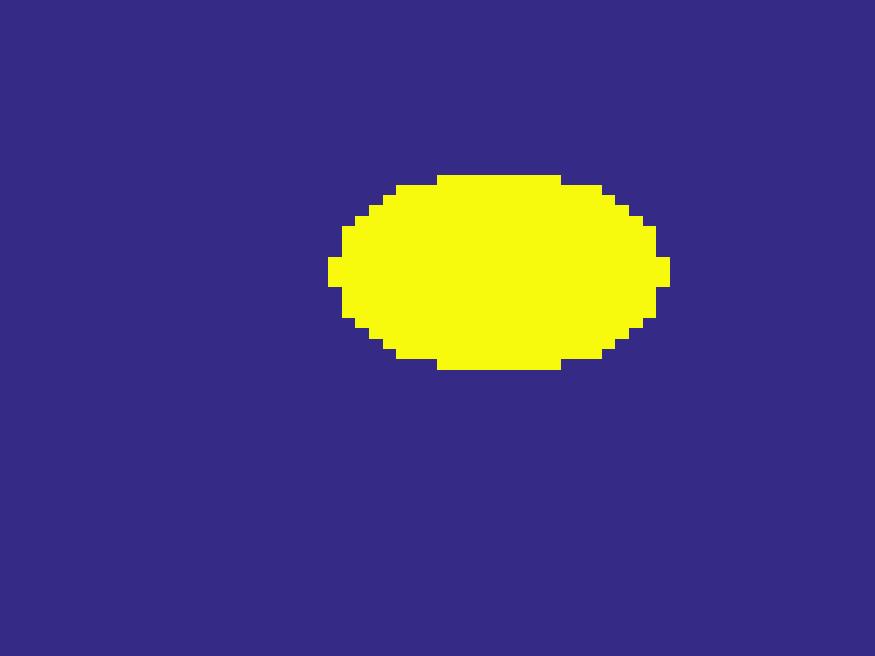}&
\includegraphics[width=.1\linewidth,height=.1\linewidth]{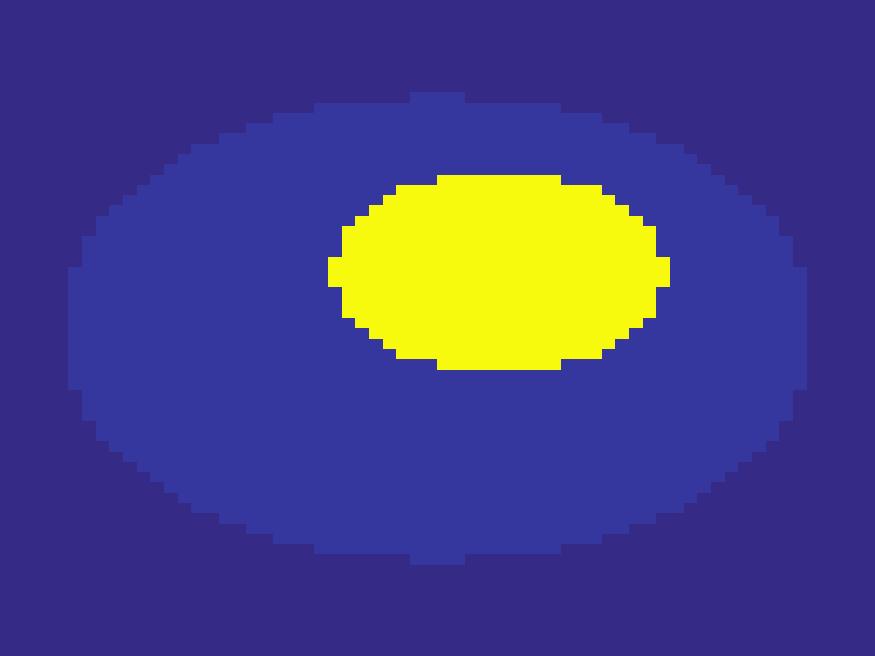}&
\includegraphics[width=.1\linewidth,height=.1\linewidth]{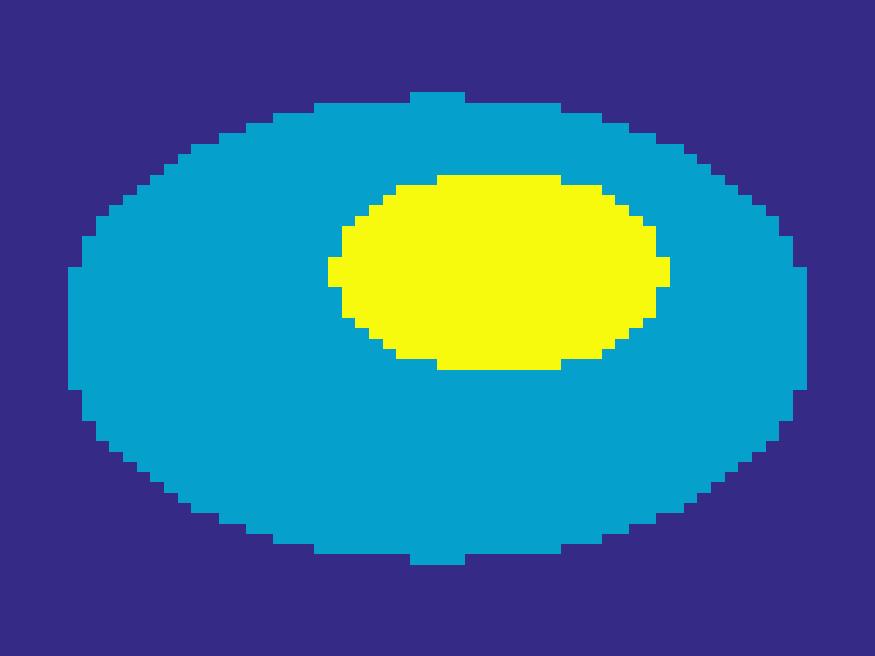}&
\includegraphics[width=.1\linewidth,height=.1\linewidth]{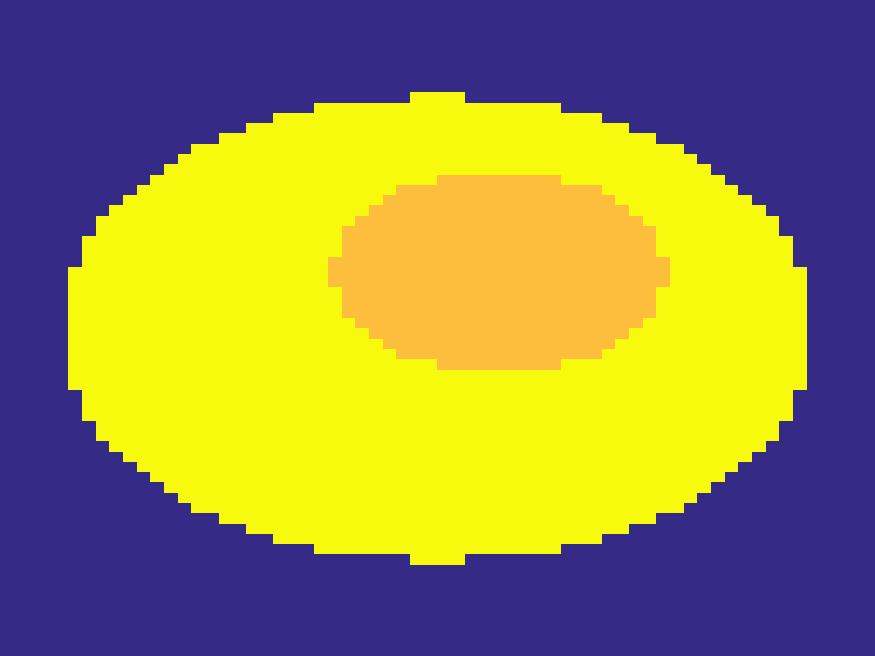}&
\includegraphics[width=.1\linewidth,height=.1\linewidth]{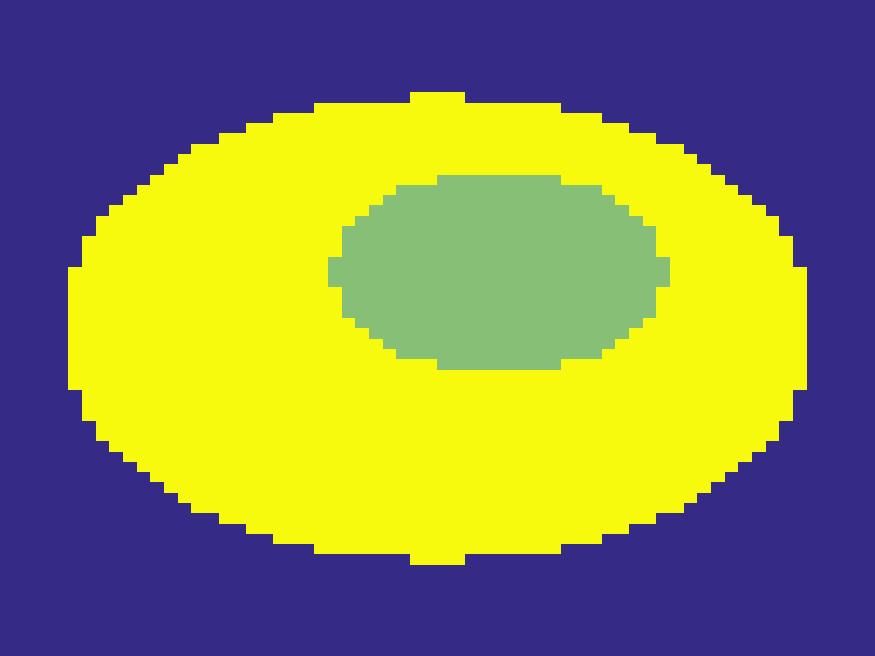}&
\includegraphics[width=.1\linewidth,height=.1\linewidth]{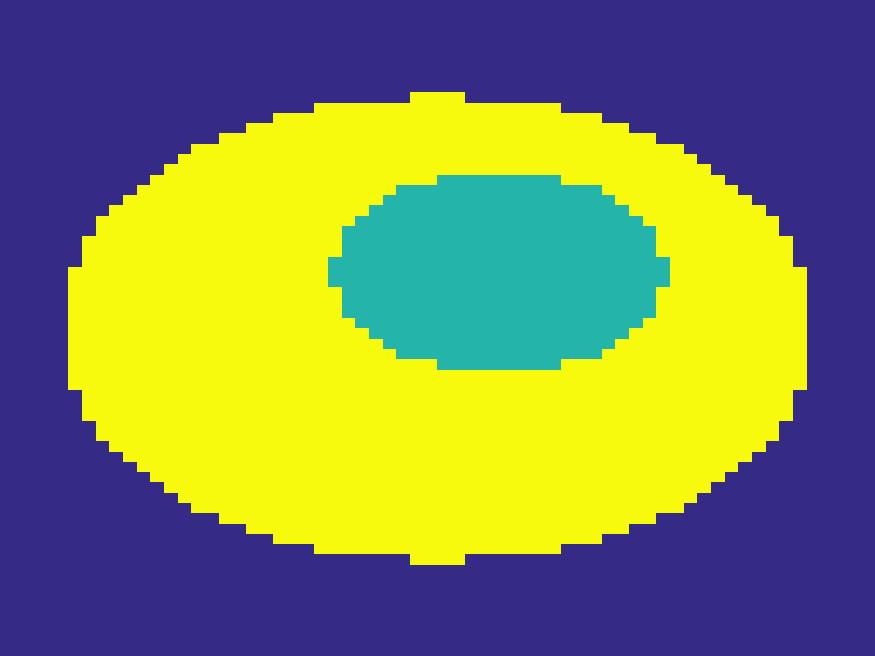}&
\includegraphics[width=.1\linewidth,height=.1\linewidth]{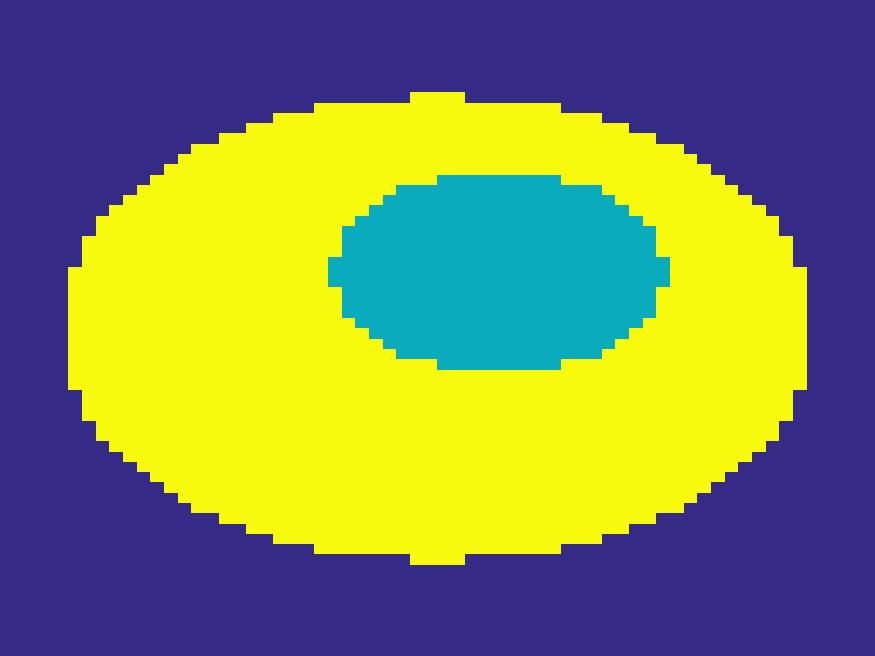}&
\includegraphics[width=.1\linewidth,height=.1\linewidth]{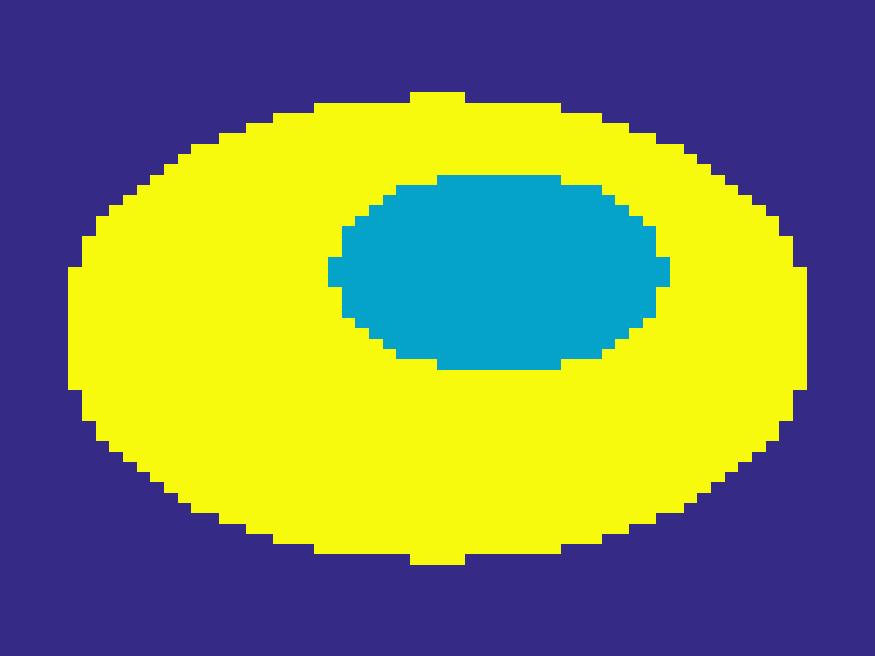}&
\includegraphics[width=.1\linewidth,height=.1\linewidth]{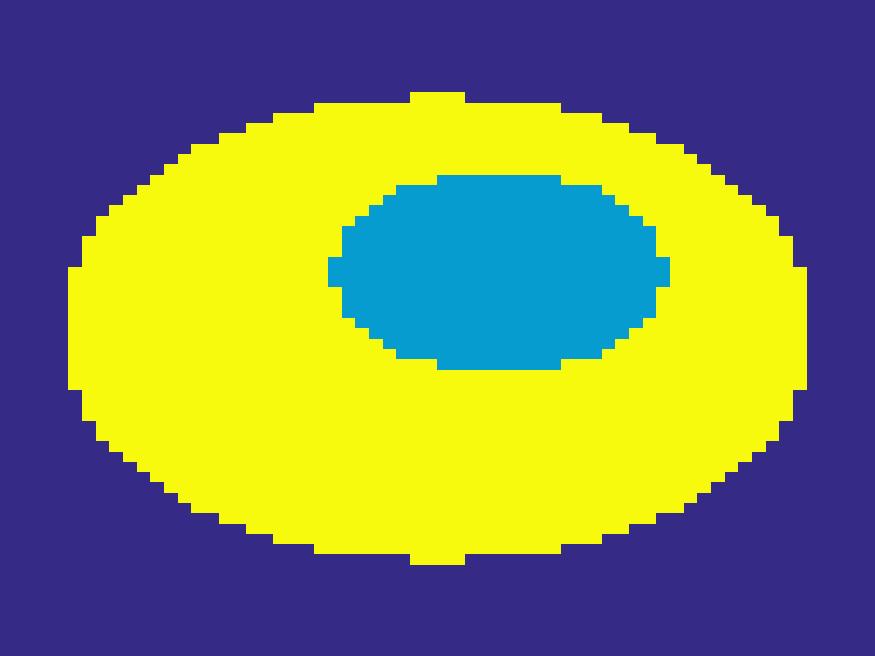}\\
\includegraphics[width=.1\linewidth,height=.1\linewidth]{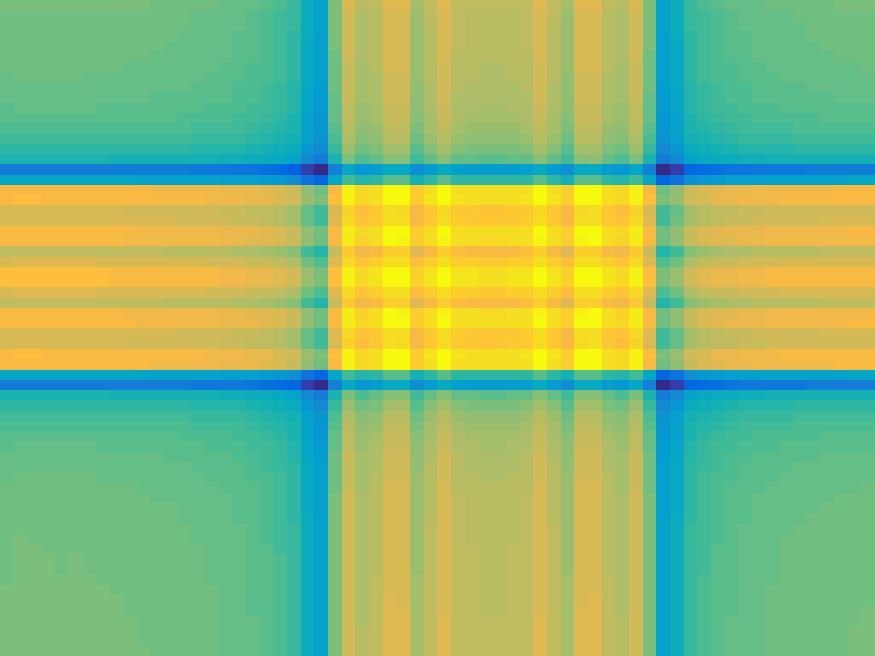}&
\includegraphics[width=.1\linewidth,height=.1\linewidth]{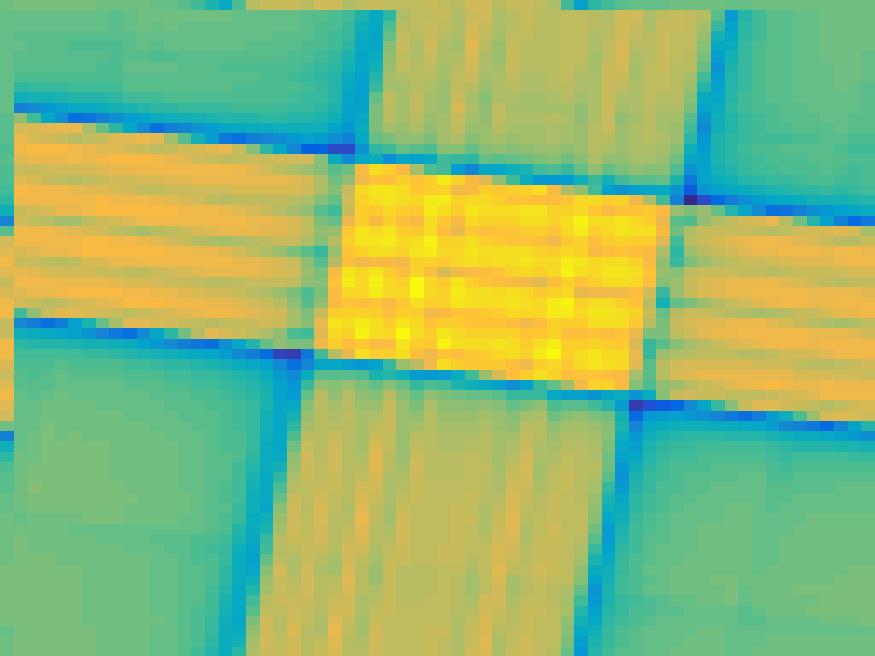}&
\includegraphics[width=.1\linewidth,height=.1\linewidth]{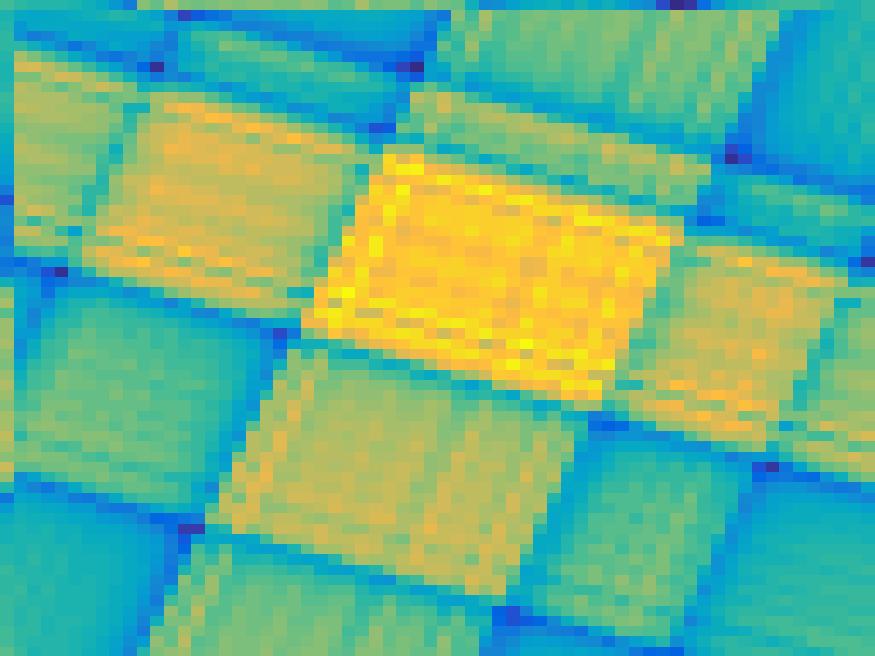}&
\includegraphics[width=.1\linewidth,height=.1\linewidth]{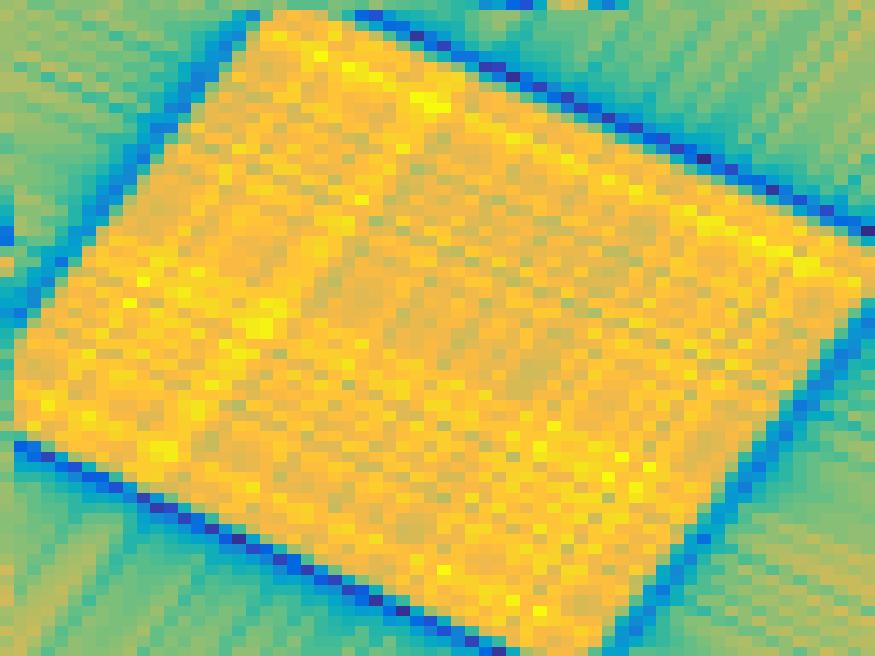}&
\includegraphics[width=.1\linewidth,height=.1\linewidth]{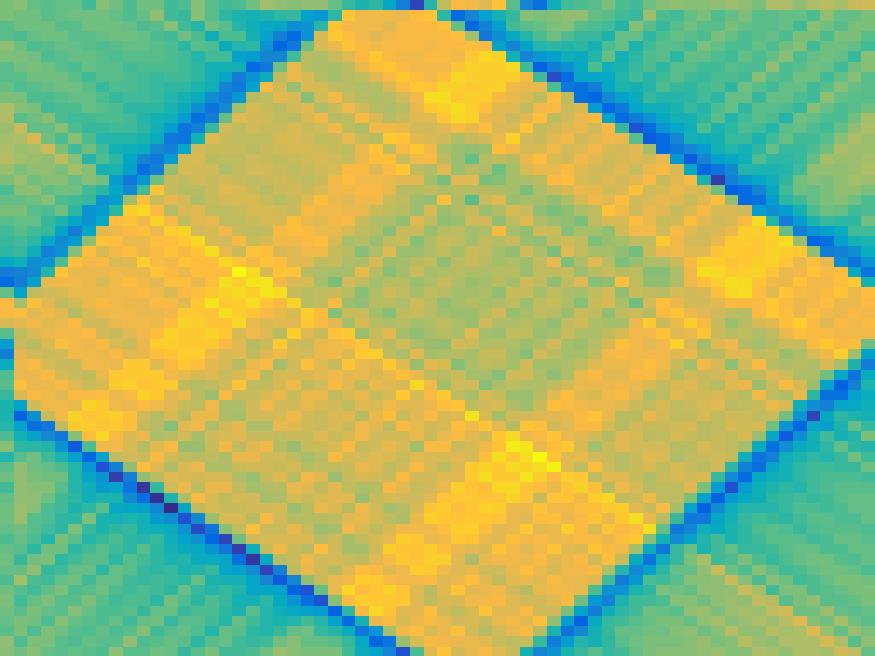}&
\includegraphics[width=.1\linewidth,height=.1\linewidth]{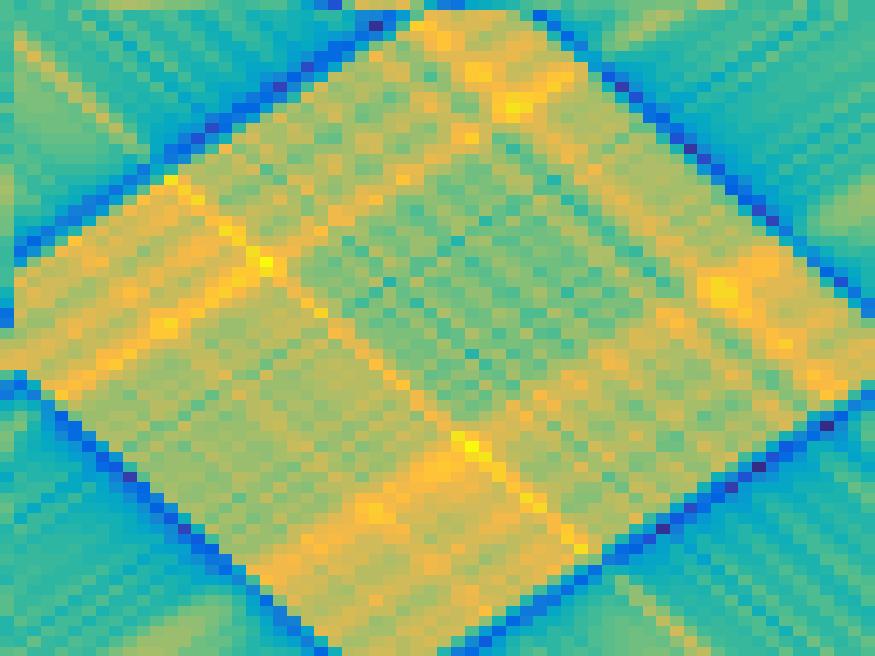}&
\includegraphics[width=.1\linewidth,height=.1\linewidth]{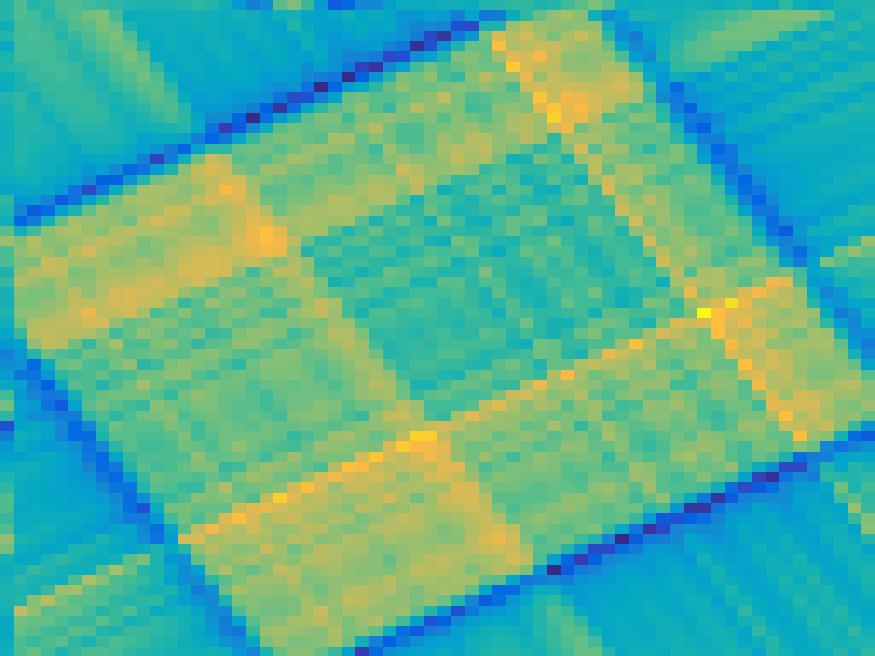}&
\includegraphics[width=.1\linewidth,height=.1\linewidth]{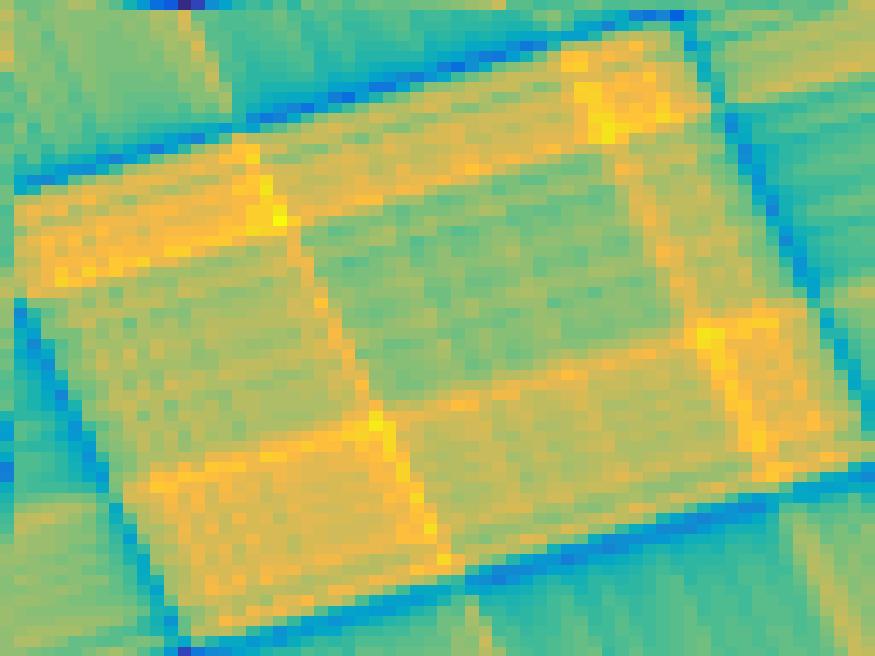}&
\includegraphics[width=.1\linewidth,height=.1\linewidth]{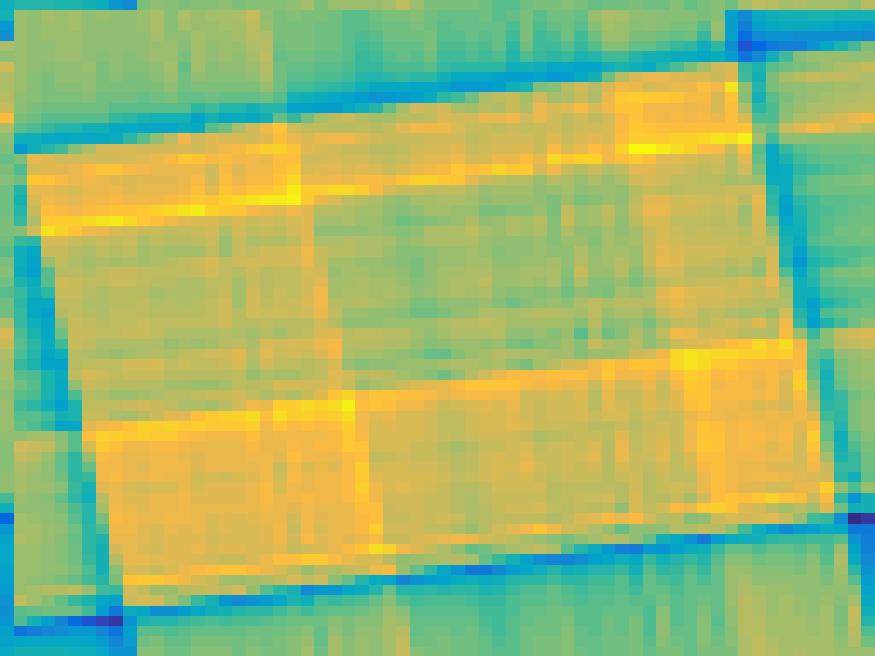}\\
\includegraphics[width=.1\linewidth,height=.1\linewidth]{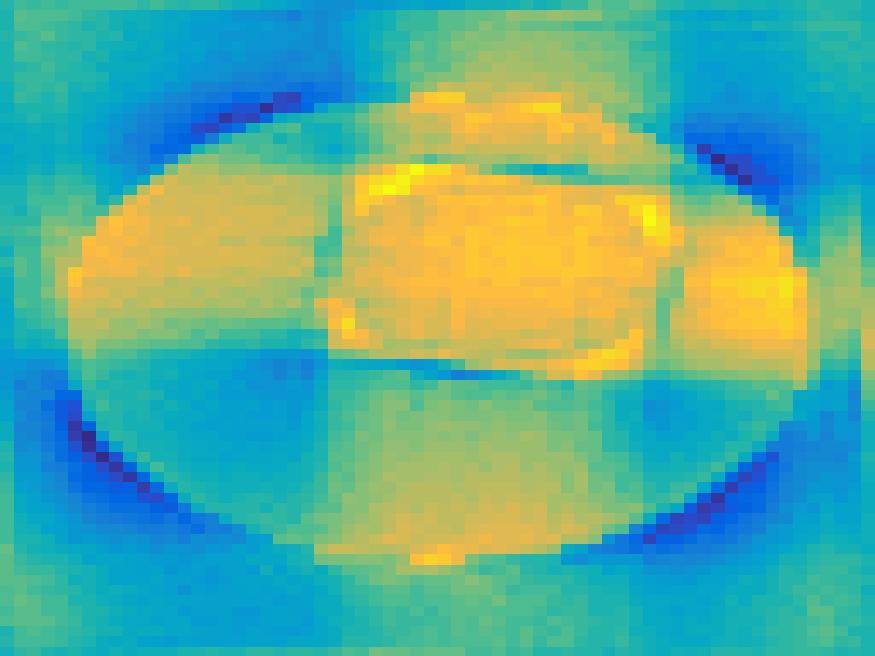}&
\includegraphics[width=.1\linewidth,height=.1\linewidth]{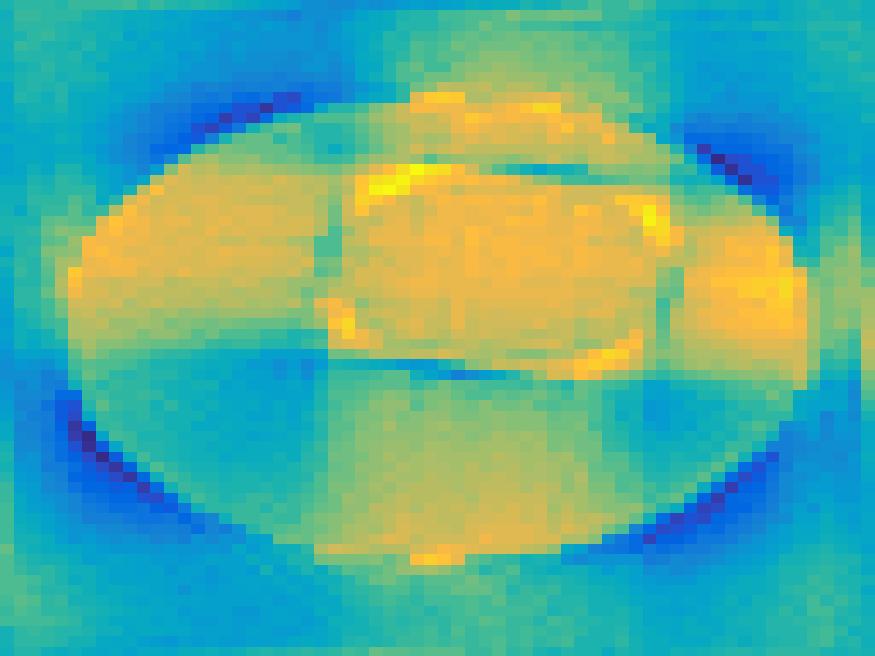}&
\includegraphics[width=.1\linewidth,height=.1\linewidth]{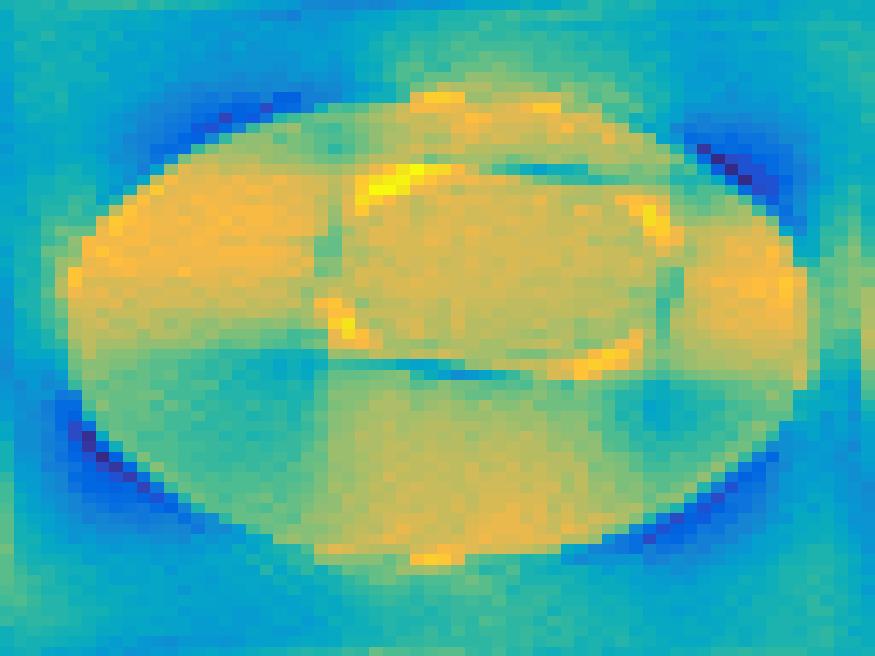}&
\includegraphics[width=.1\linewidth,height=.1\linewidth]{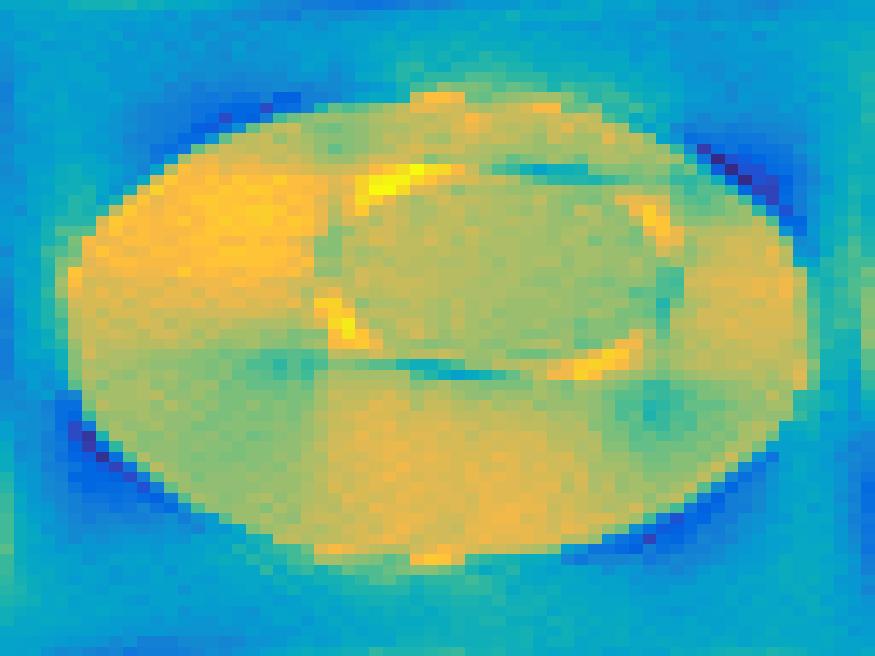}&
\includegraphics[width=.1\linewidth,height=.1\linewidth]{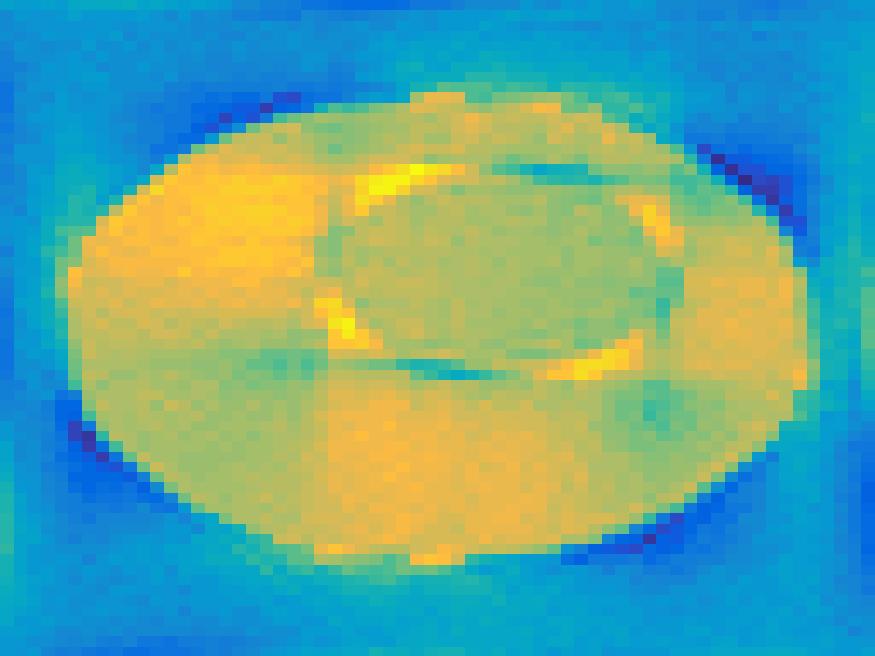}&
\includegraphics[width=.1\linewidth,height=.1\linewidth]{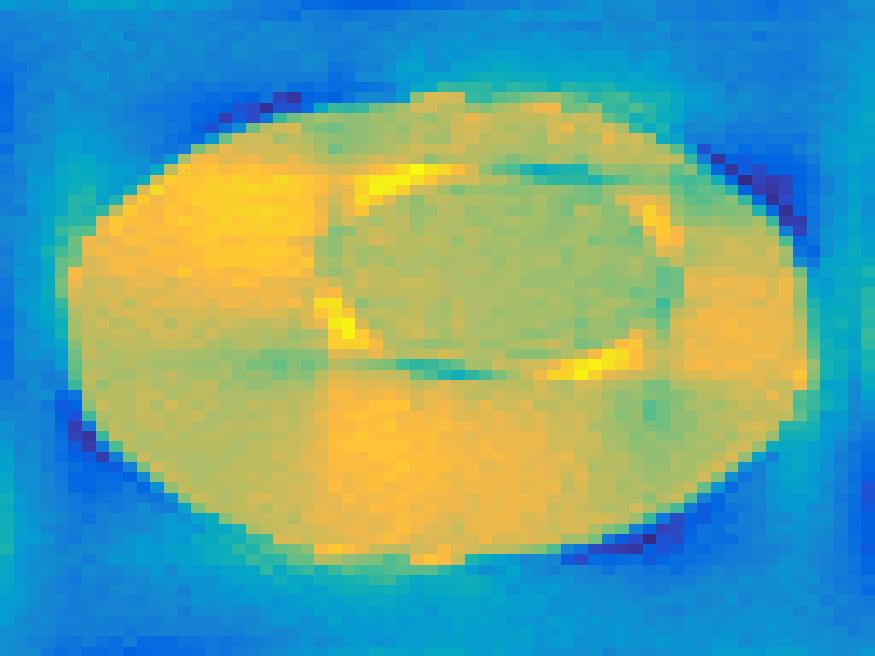}&
\includegraphics[width=.1\linewidth,height=.1\linewidth]{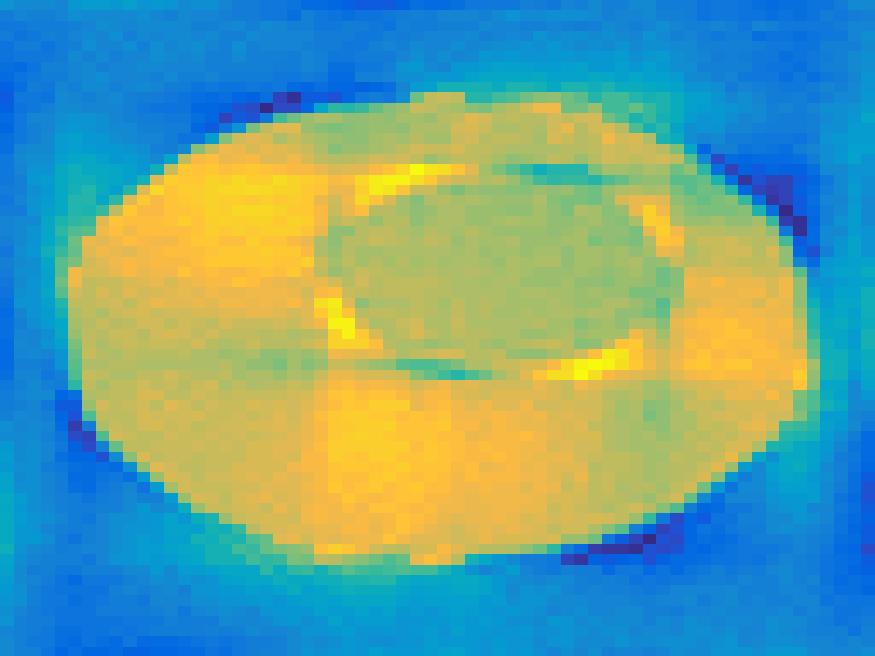}&
\includegraphics[width=.1\linewidth,height=.1\linewidth]{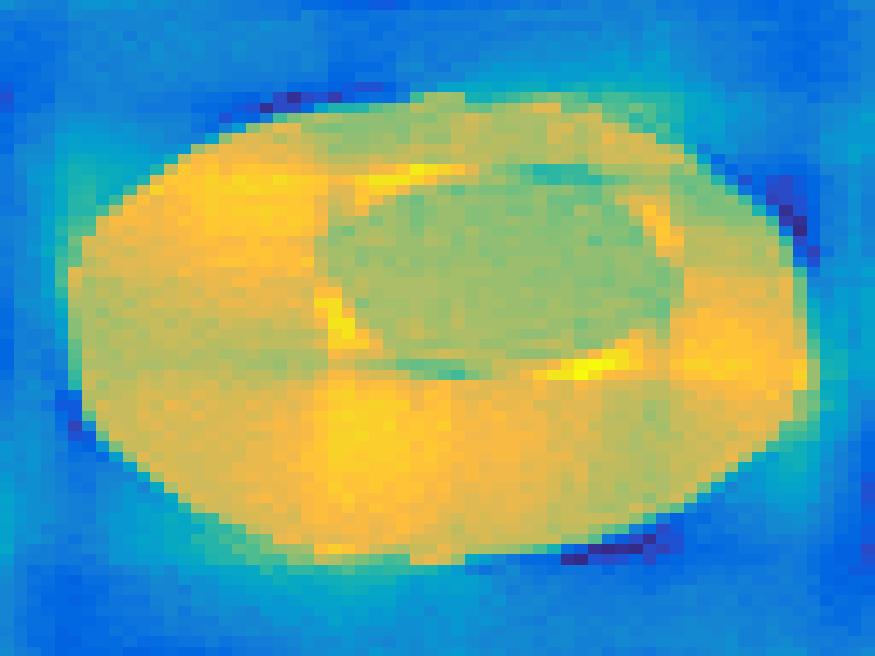}&
\includegraphics[width=.1\linewidth,height=.1\linewidth]{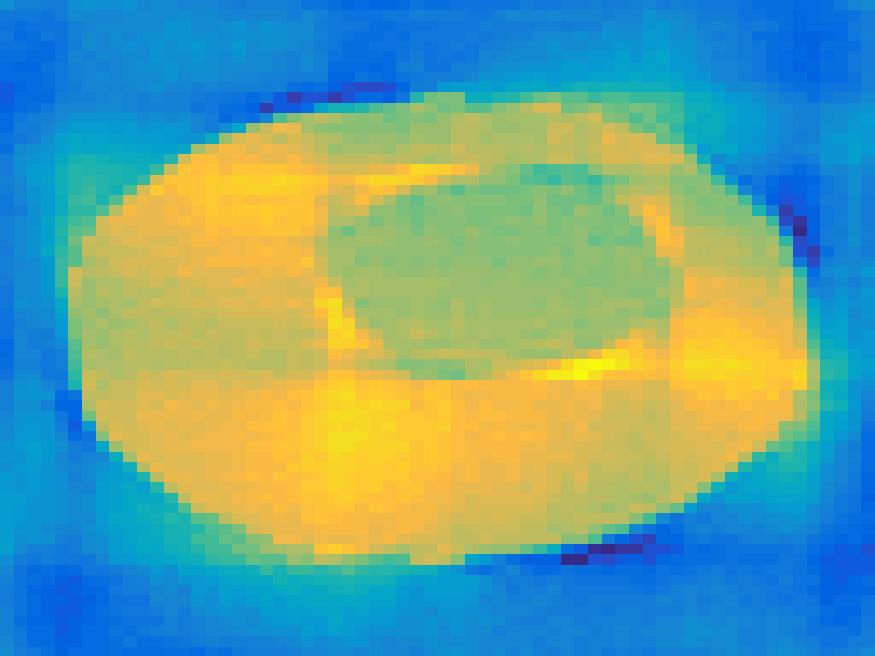}\\
\includegraphics[width=.1\linewidth,height=.1\linewidth]{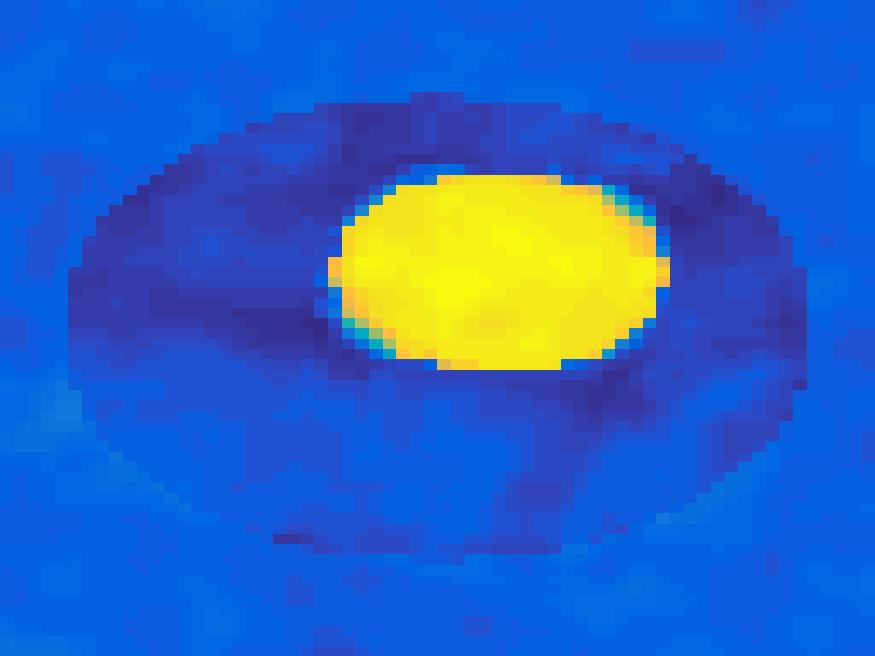}&
\includegraphics[width=.1\linewidth,height=.1\linewidth]{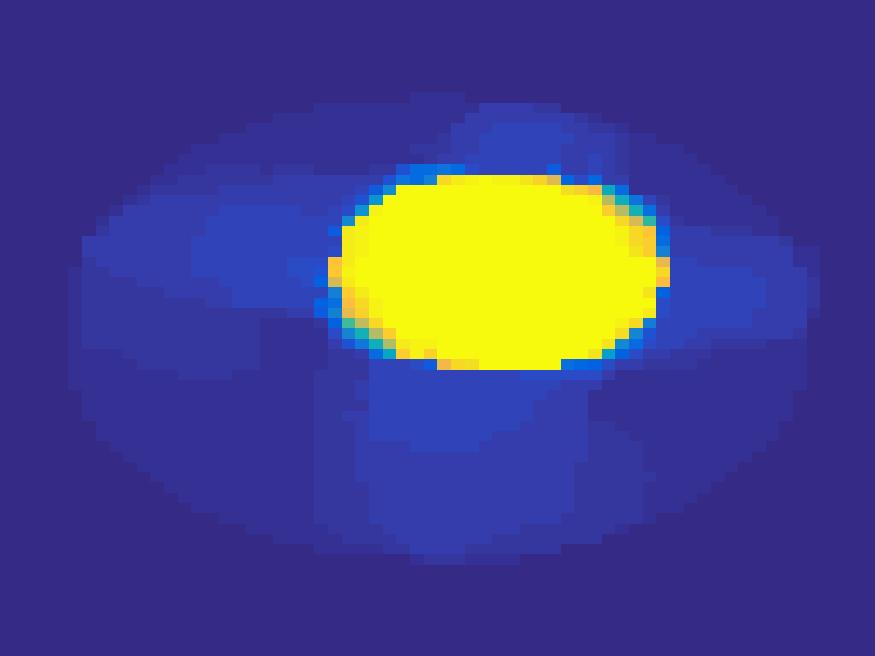}&
\includegraphics[width=.1\linewidth,height=.1\linewidth]{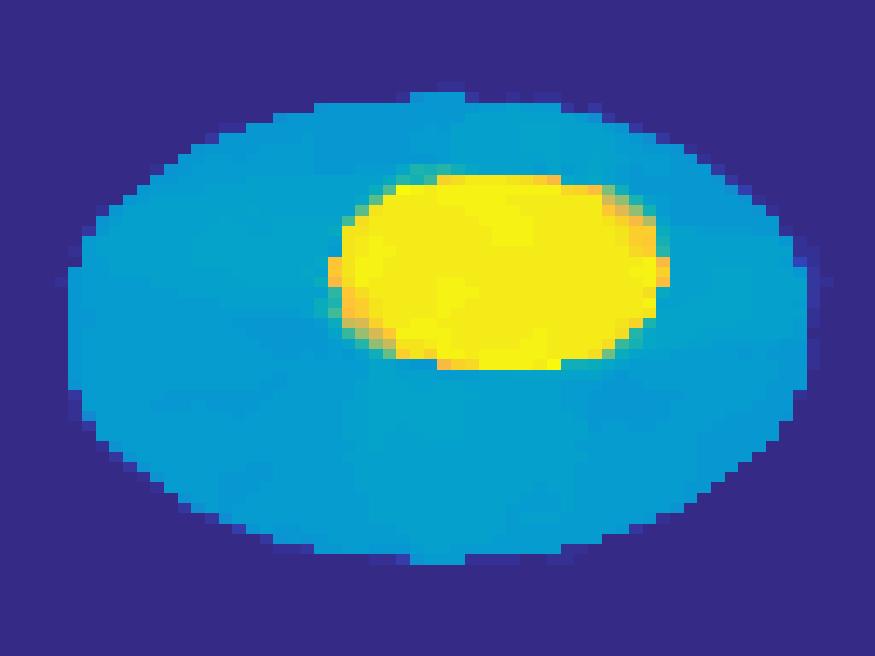}&
\includegraphics[width=.1\linewidth,height=.1\linewidth]{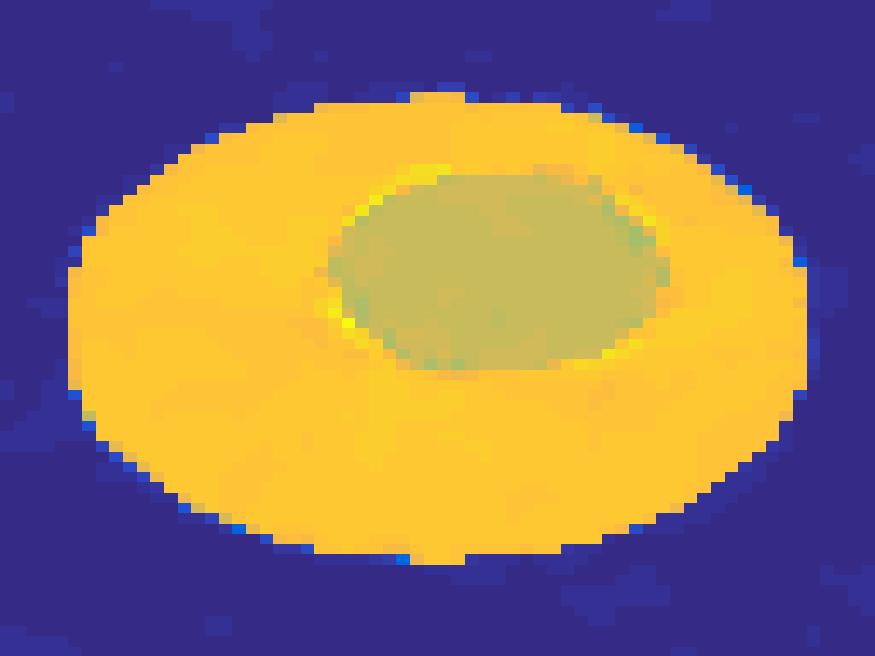}&
\includegraphics[width=.1\linewidth,height=.1\linewidth]{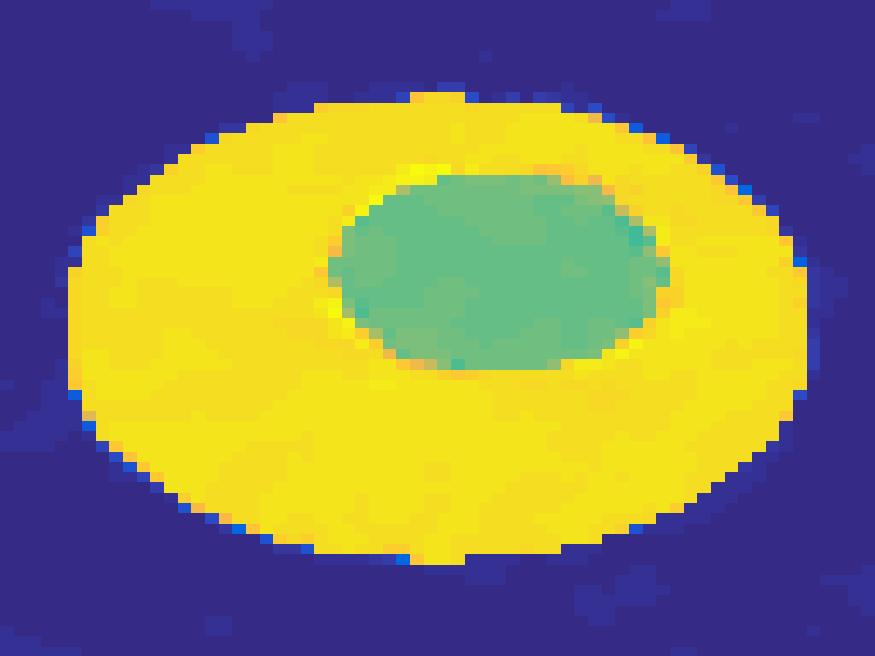}&
\includegraphics[width=.1\linewidth,height=.1\linewidth]{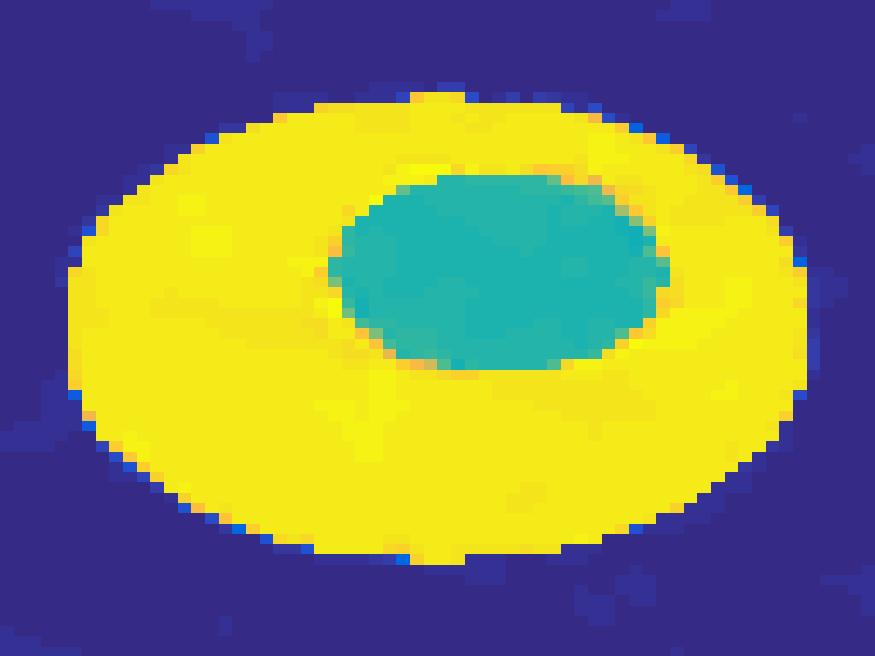}&
\includegraphics[width=.1\linewidth,height=.1\linewidth]{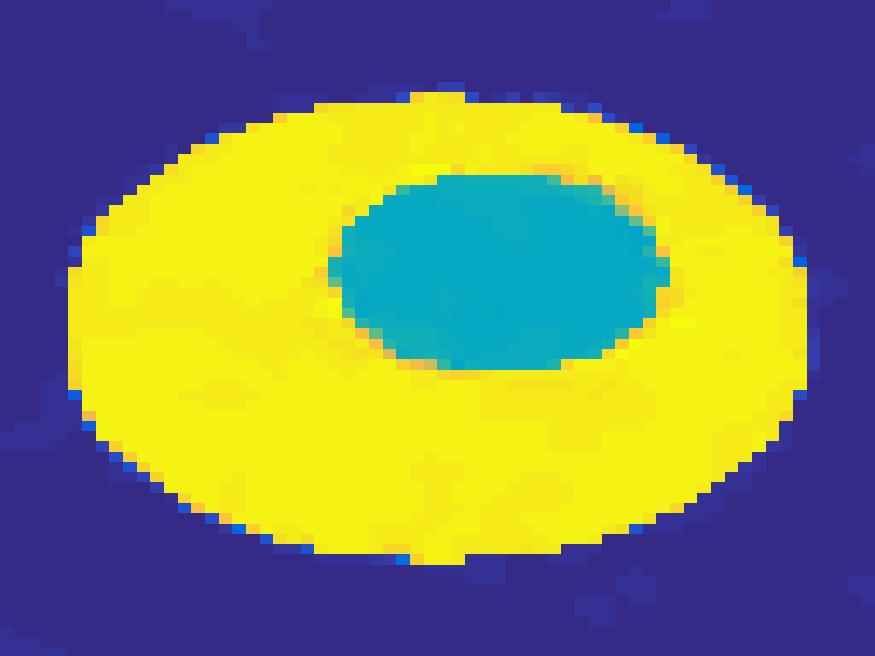}&
\includegraphics[width=.1\linewidth,height=.1\linewidth]{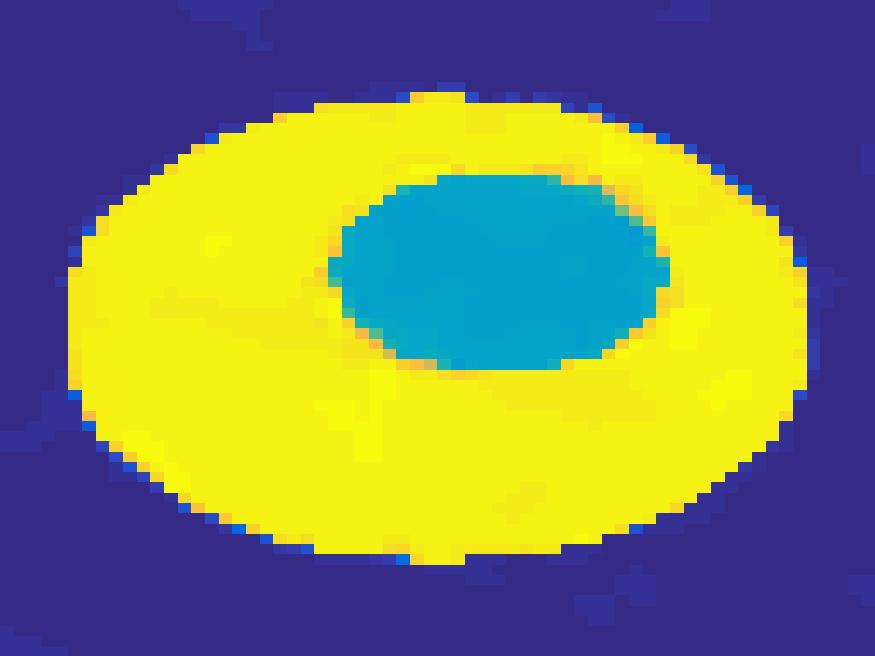}&
\includegraphics[width=.1\linewidth,height=.1\linewidth]{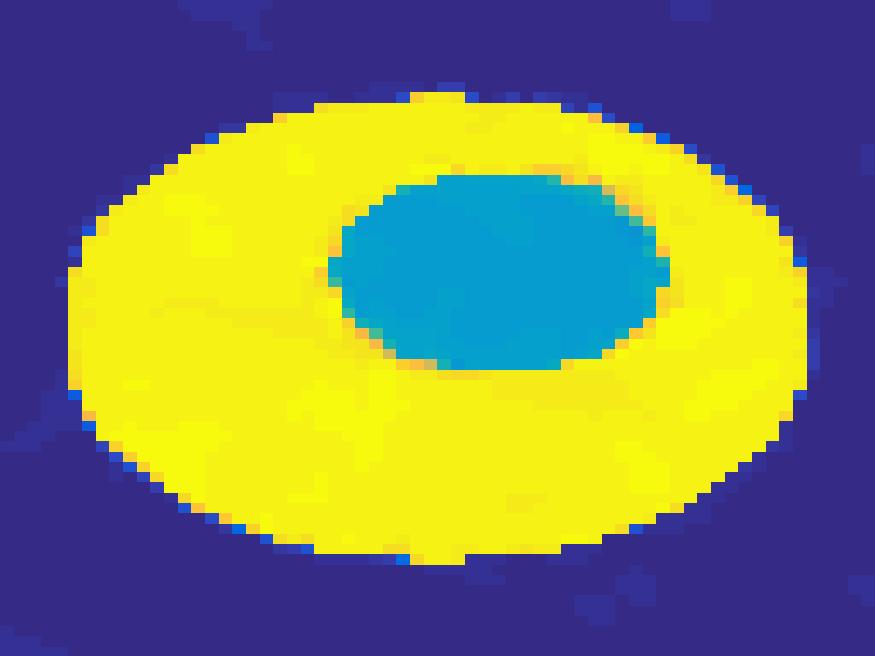}\\
\includegraphics[width=.1\linewidth,height=.1\linewidth]{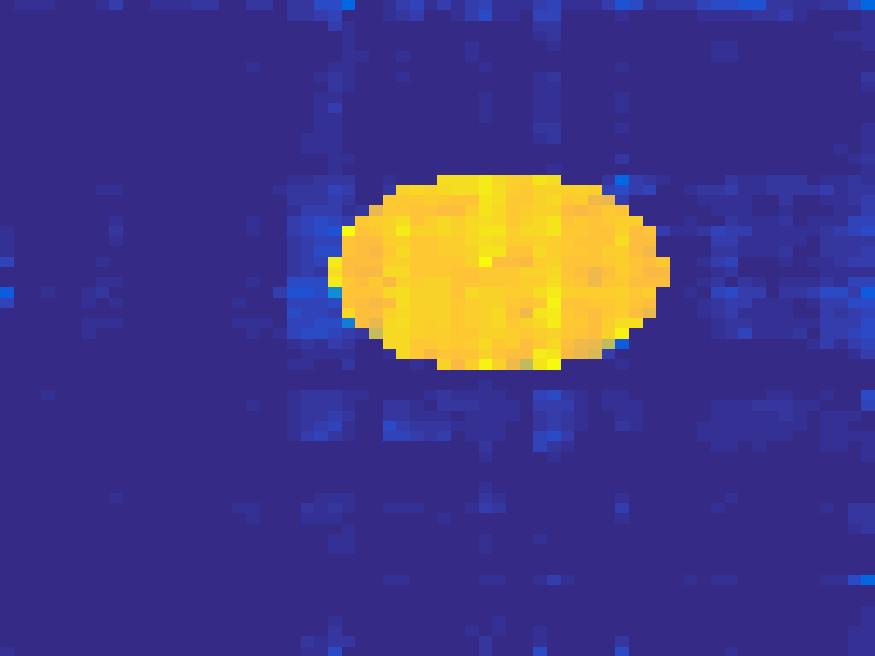}&
\includegraphics[width=.1\linewidth,height=.1\linewidth]{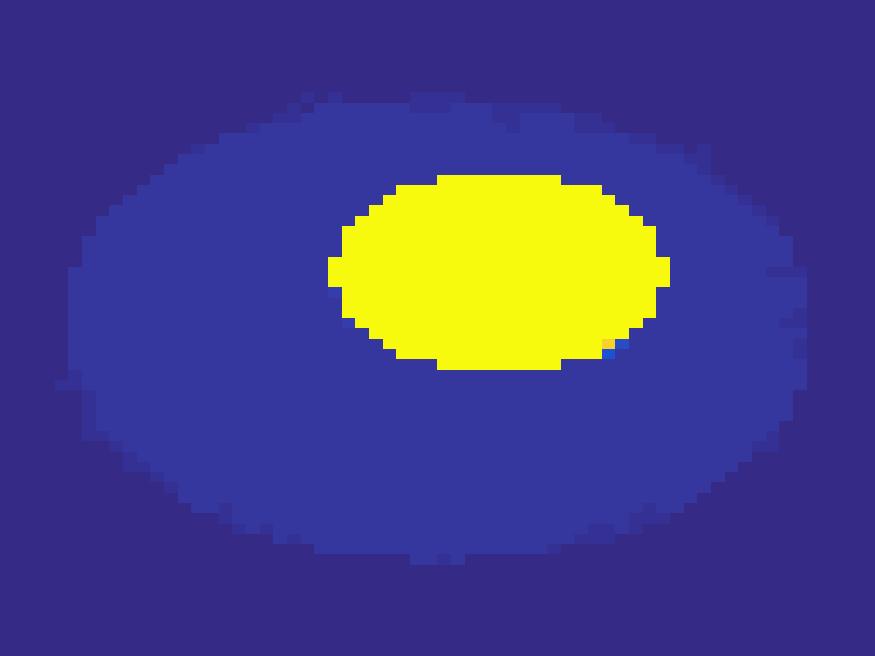}&
\includegraphics[width=.1\linewidth,height=.1\linewidth]{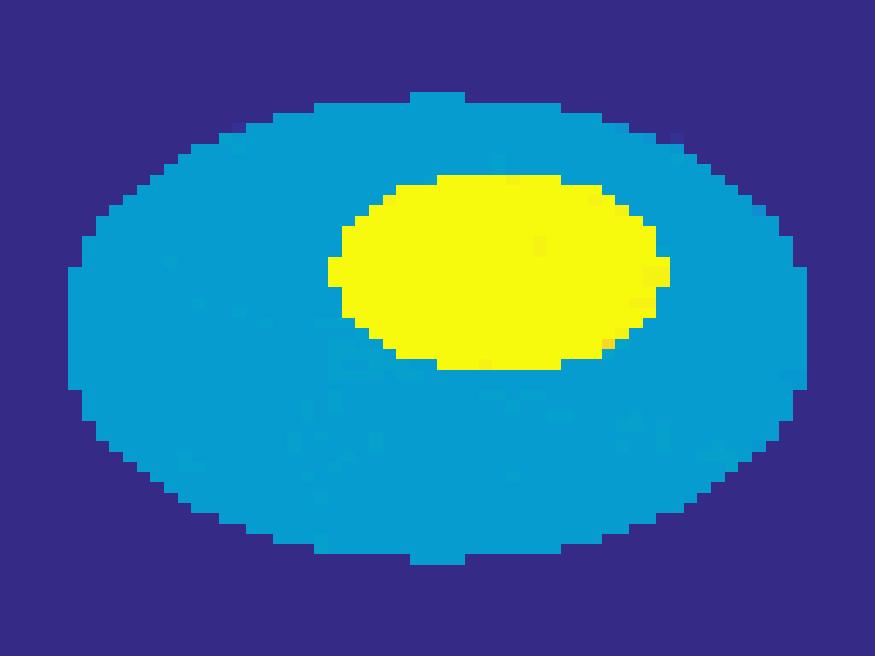}&
\includegraphics[width=.1\linewidth,height=.1\linewidth]{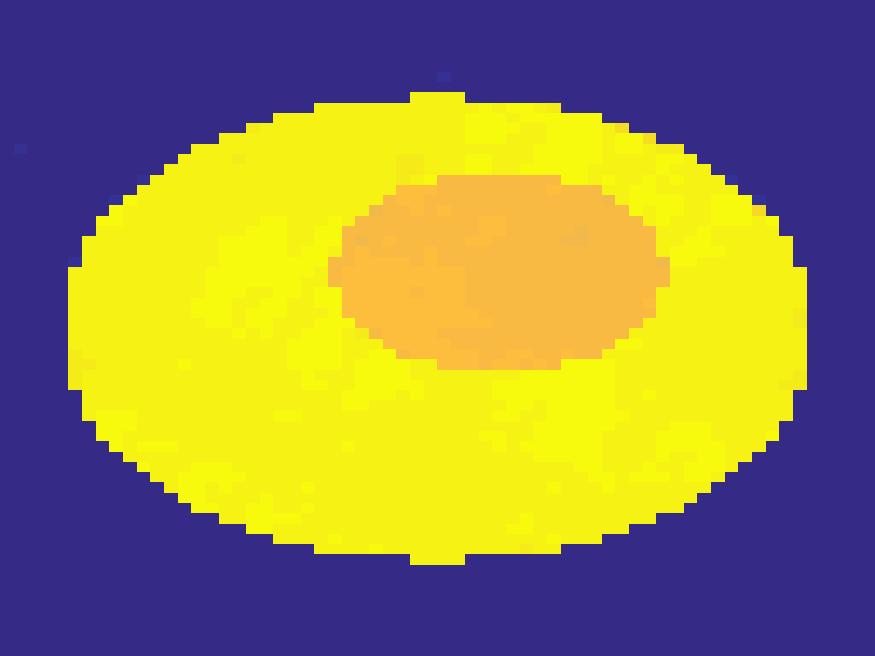}&
\includegraphics[width=.1\linewidth,height=.1\linewidth]{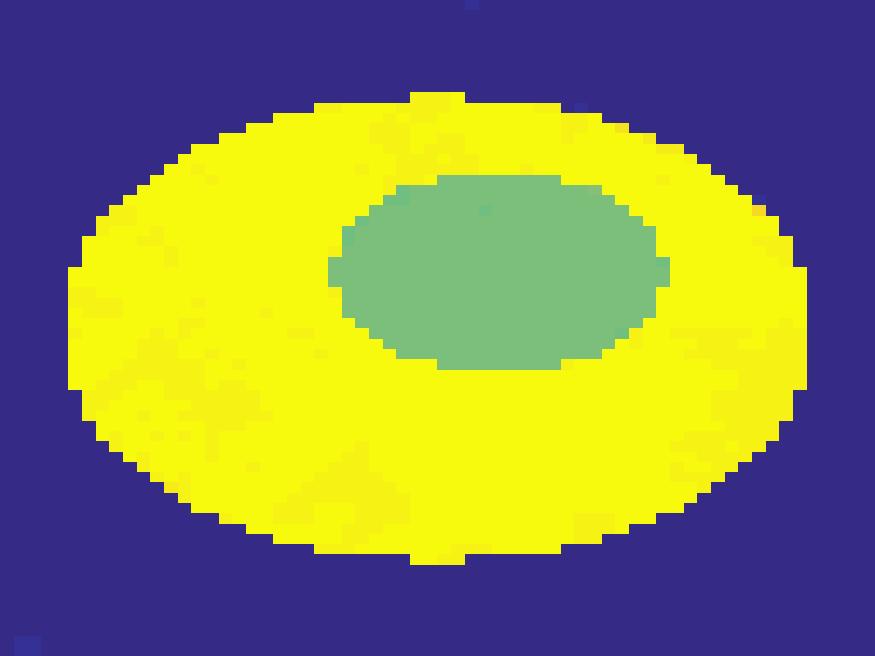}&
\includegraphics[width=.1\linewidth,height=.1\linewidth]{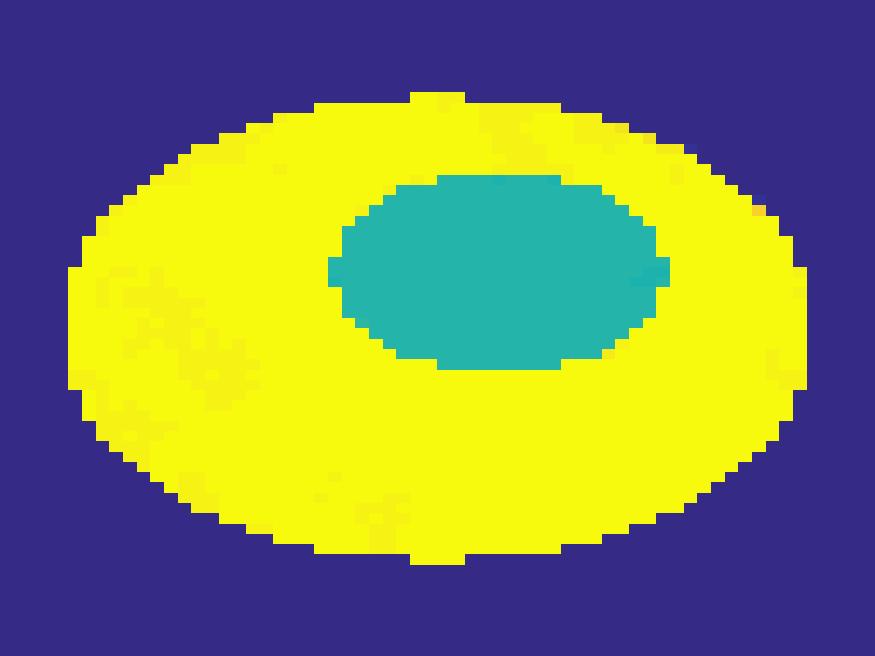}&
\includegraphics[width=.1\linewidth,height=.1\linewidth]{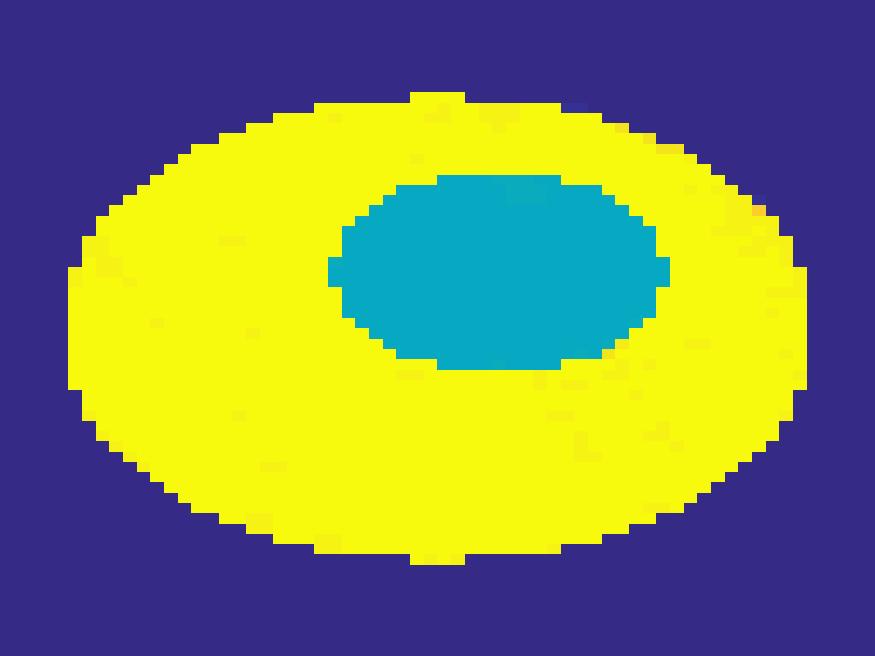}&
\includegraphics[width=.1\linewidth,height=.1\linewidth]{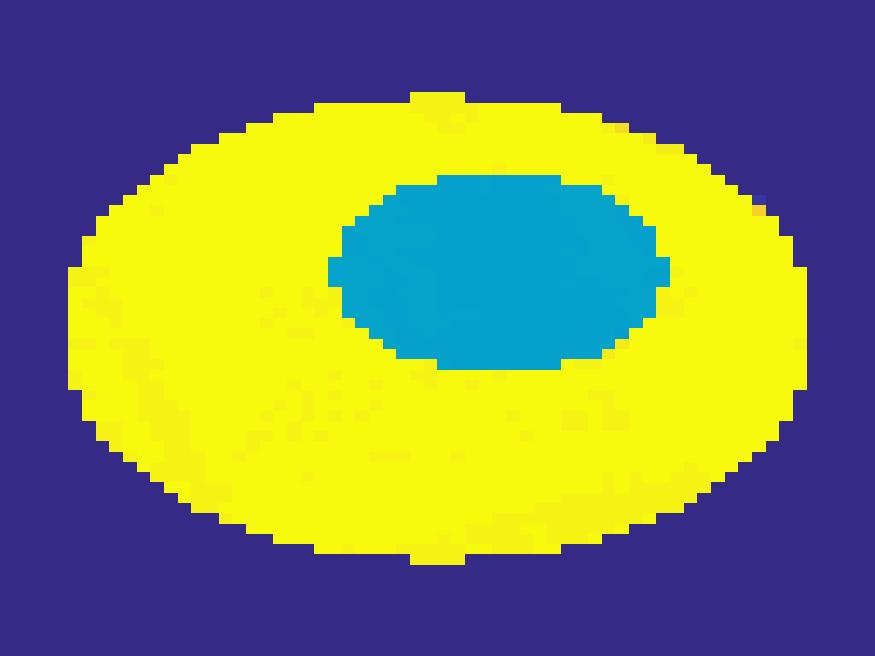}&
\includegraphics[width=.1\linewidth,height=.1\linewidth]{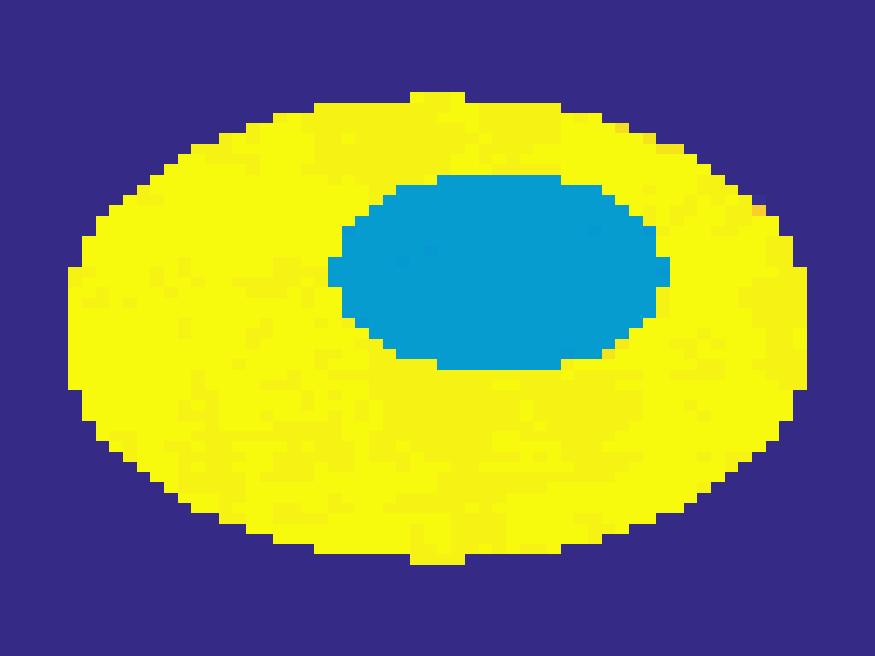}\\
{\footnotesize Frame 1}&
{\footnotesize Frame 11}&
{\footnotesize Frame 21}&
{\footnotesize Frame 31}&
{\footnotesize Frame 41}&
{\footnotesize Frame 51}&
{\footnotesize Frame 61}&
{\footnotesize Frame 71}&
{\footnotesize Frame 81}\\
\end{tabular}
\caption {First row: Ground truth; Second row: FBP; Third  row: least square method;  Forth row: SEMF \cite{ding2015dynamic}; Fifth row: Proposed model.}
\label{fig:ICgaussEll}
\end{figure}

In the tests,  projections at two  orthogonal angles are simulated for every frame to mimic 2-head camera data collection. The projection angles increase sequentially by $1^\circ$ along temporal direction.  For example, at frame 1,  projections are simulated at angle $1^ \circ$ and $91^\circ$, and at frame 2,  angle $2^\circ$ and $92^\circ$, etc.
Finally, $10\%$ white Gaussian noise is added to the projection data. Reconstruction results with different methods  are shown in Figure  \ref{fig:ICgaussEll}. Since the number of projections is very limited for each frame, the traditional FBP and least square methods cannot reconstruct the images satisfactorily, while the proposed method is capable to reconstruct the images effectively. Compared with SEMF model,
when the edge of images jump (see frame 21 -frame 31 in Figure  \ref{fig:ICgaussEll}), the proposed model can better capture the change of the tendency of  TAC.

\begin{figure}[ht]
\begin{center}
\subfigure{
\includegraphics[width=.45\linewidth,height=.3\linewidth]{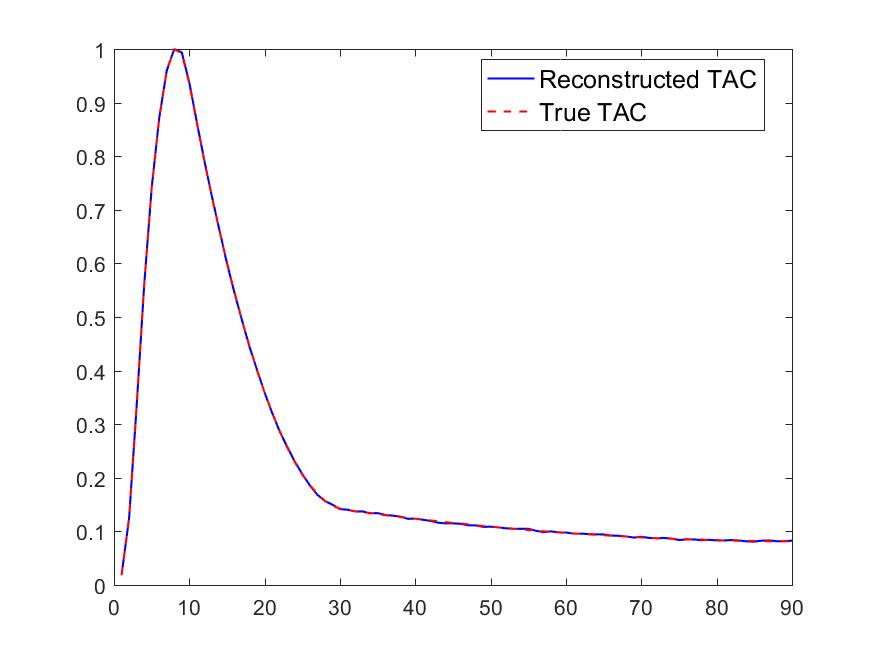}}
\subfigure{
\includegraphics[width=.45\linewidth,height=.3\linewidth]{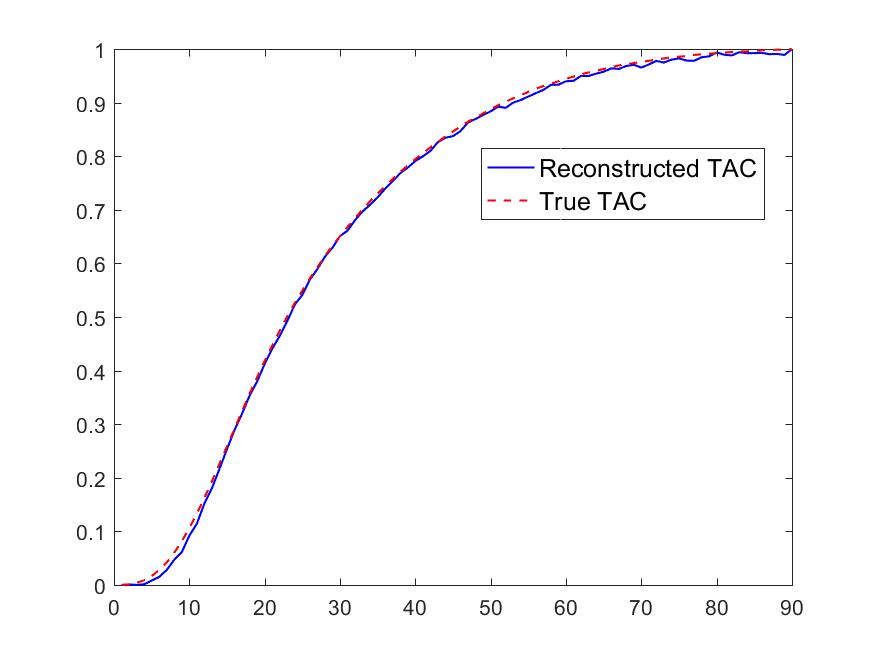}}
\end{center}
\caption{Time activity curve of two regions.}
\label{fig:ICGauEllTAC}
\end{figure}

Figure \ref{fig:ICGauEllTAC} illustrates the comparison of the TACs of blood and liver. The dash  lines are
the normalized true TACs and the solid lines are the normalized one extracted from the reconstruction images by our
method. Even with high level noise and fast change of radioisotope,  the reconstructed one  fit closely to the true one.


\begin{figure}
\begin{tabular}{c@{\hspace{2pt}}c@{\hspace{2pt}}c@{\hspace{2pt}}c@{\hspace{2pt}}c@{\hspace{2pt}}c@{\hspace{2pt}}c@{\hspace{2pt}}c@{\hspace{2pt}}c@{\hspace{2pt}}c}
\includegraphics[width=.1\linewidth,height=.1\linewidth]{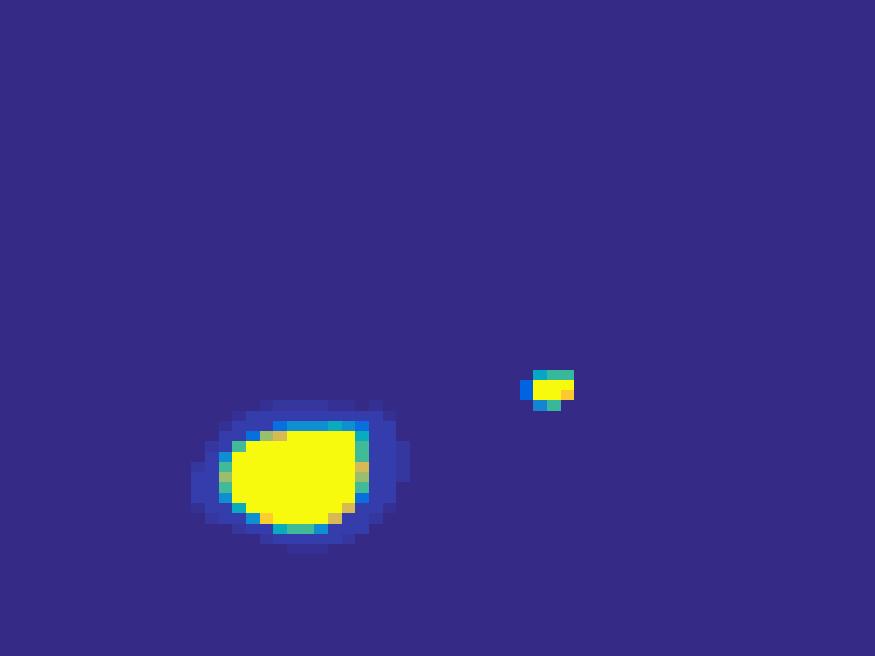}&
\includegraphics[width=.1\linewidth,height=.1\linewidth]{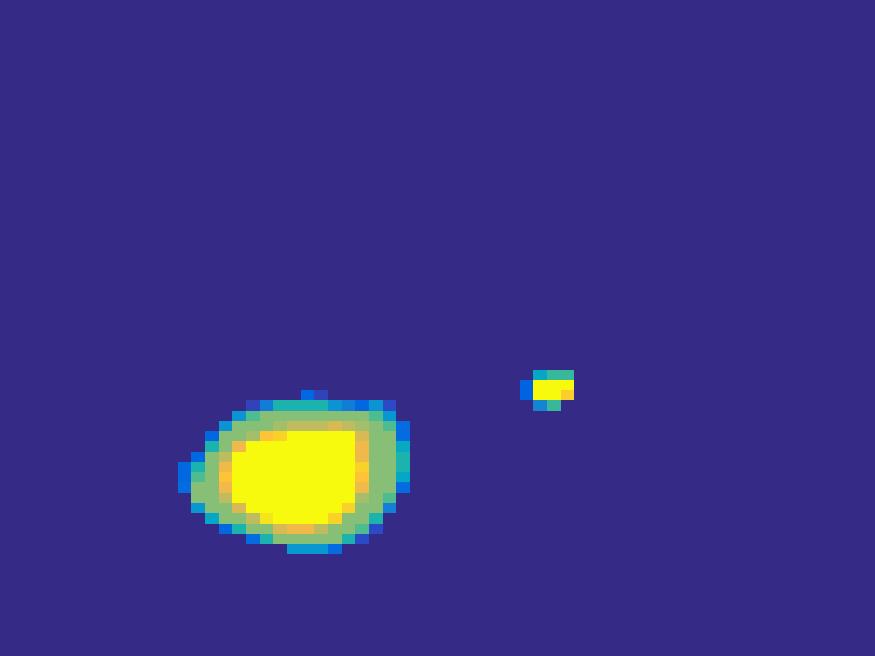}&
\includegraphics[width=.1\linewidth,height=.1\linewidth]{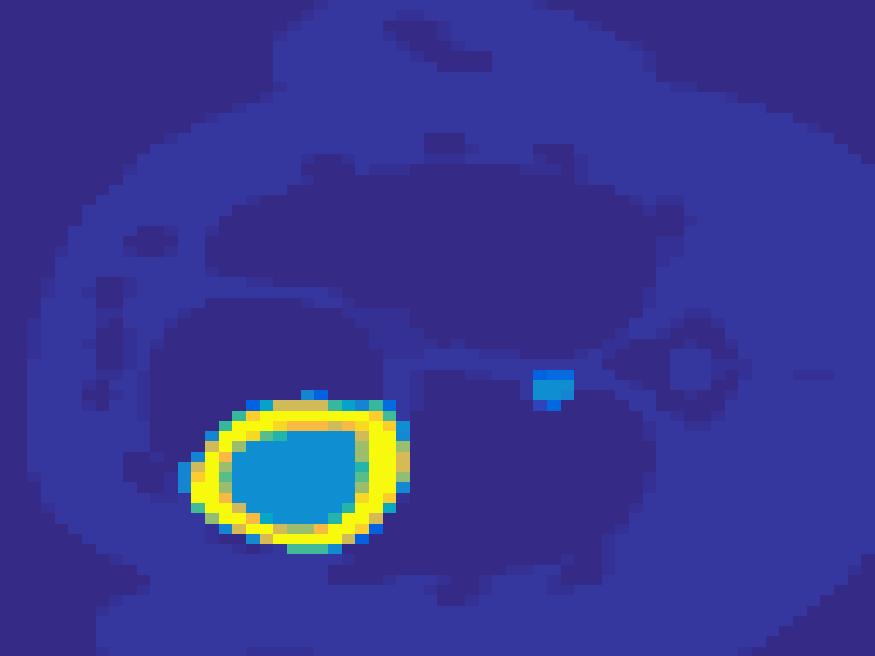}&
\includegraphics[width=.1\linewidth,height=.1\linewidth]{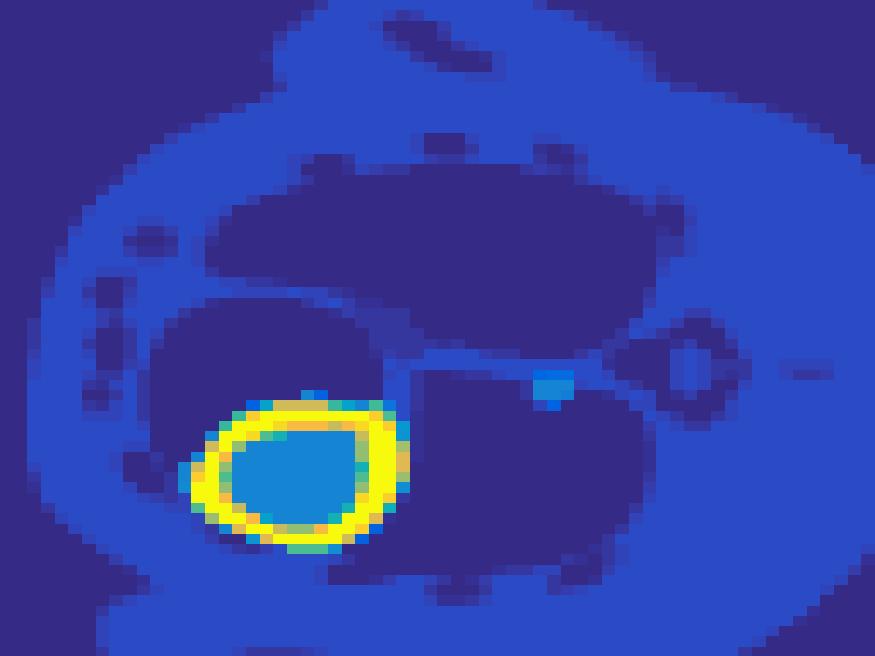}&
\includegraphics[width=.1\linewidth,height=.1\linewidth]{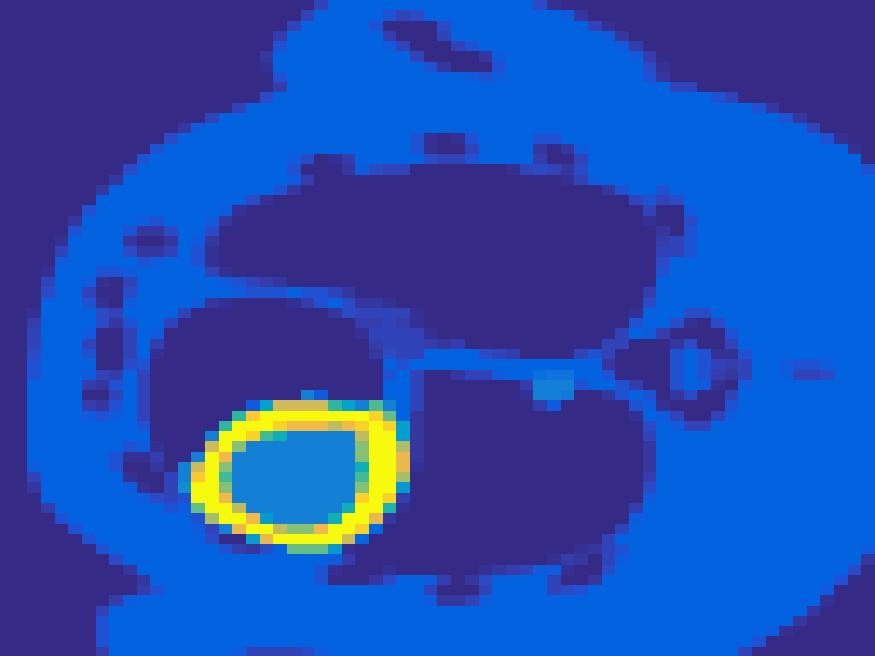}&
\includegraphics[width=.1\linewidth,height=.1\linewidth]{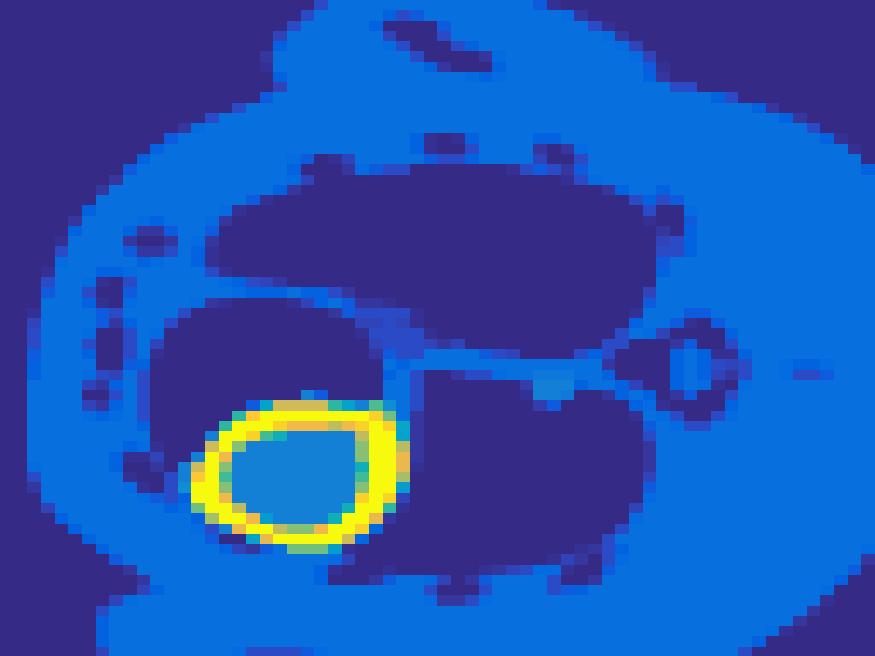}&
\includegraphics[width=.1\linewidth,height=.1\linewidth]{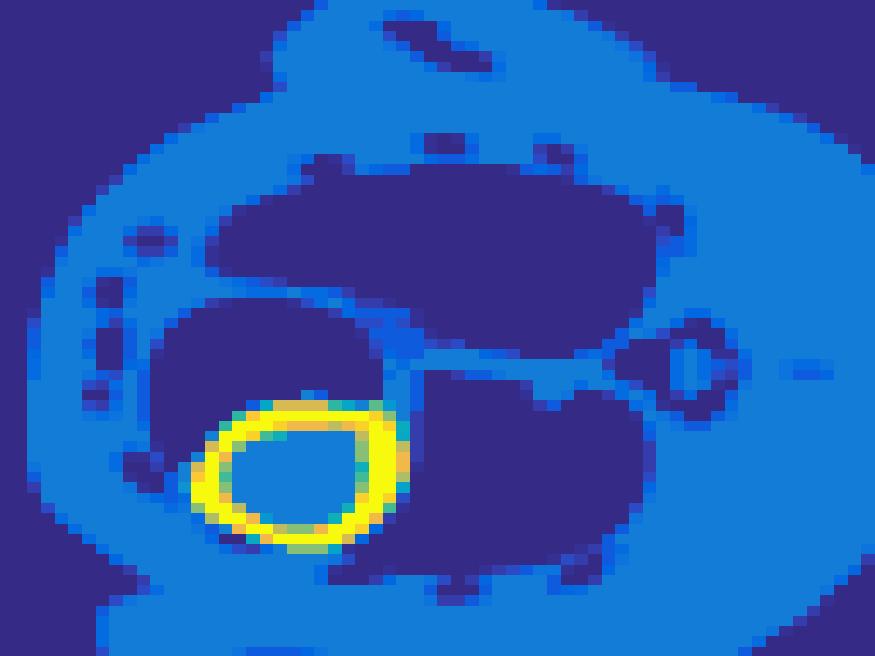}&
\includegraphics[width=.1\linewidth,height=.1\linewidth]{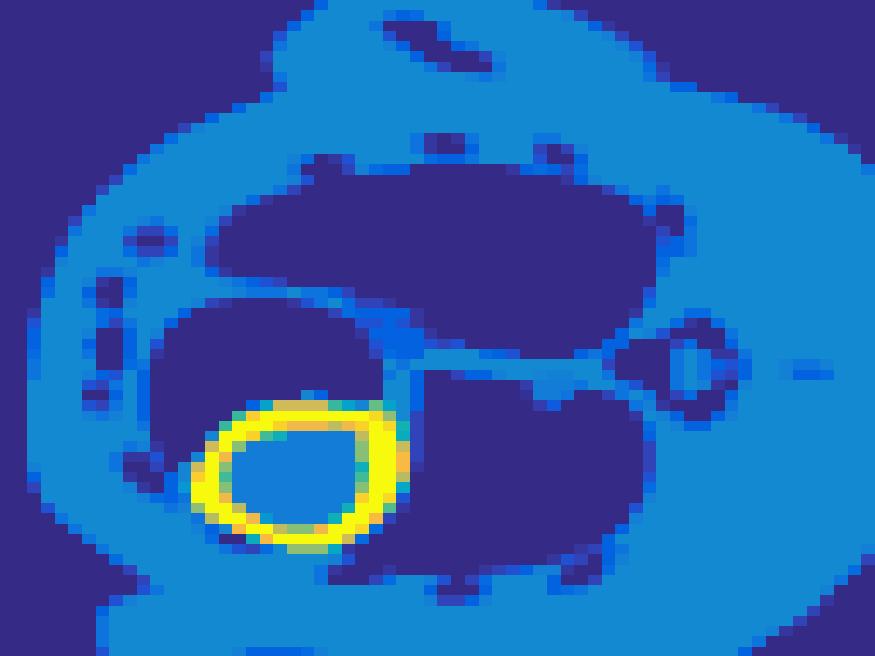}&
\includegraphics[width=.1\linewidth,height=.1\linewidth]{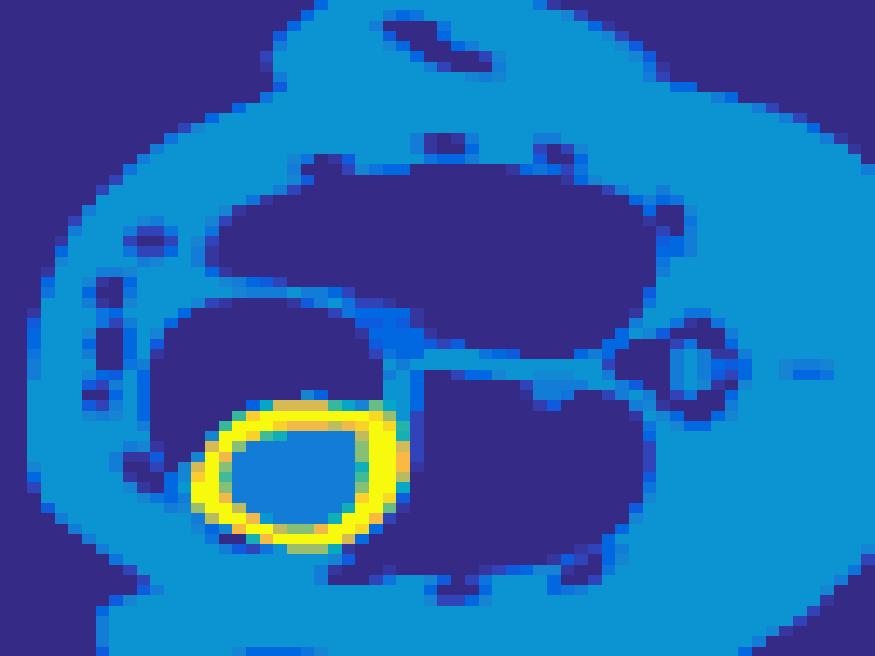}\\
\includegraphics[width=.1\linewidth,height=.1\linewidth]{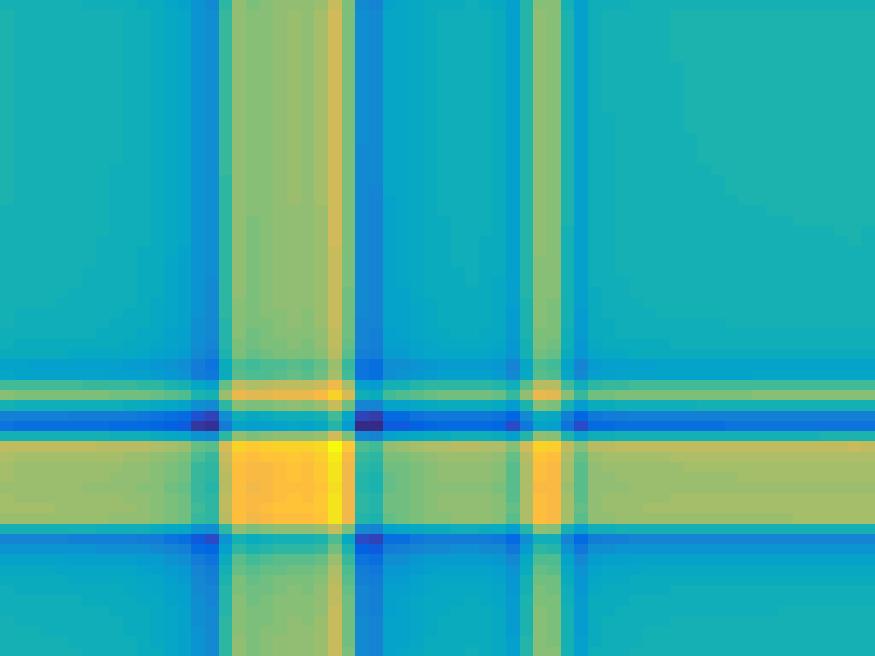}&
\includegraphics[width=.1\linewidth,height=.1\linewidth]{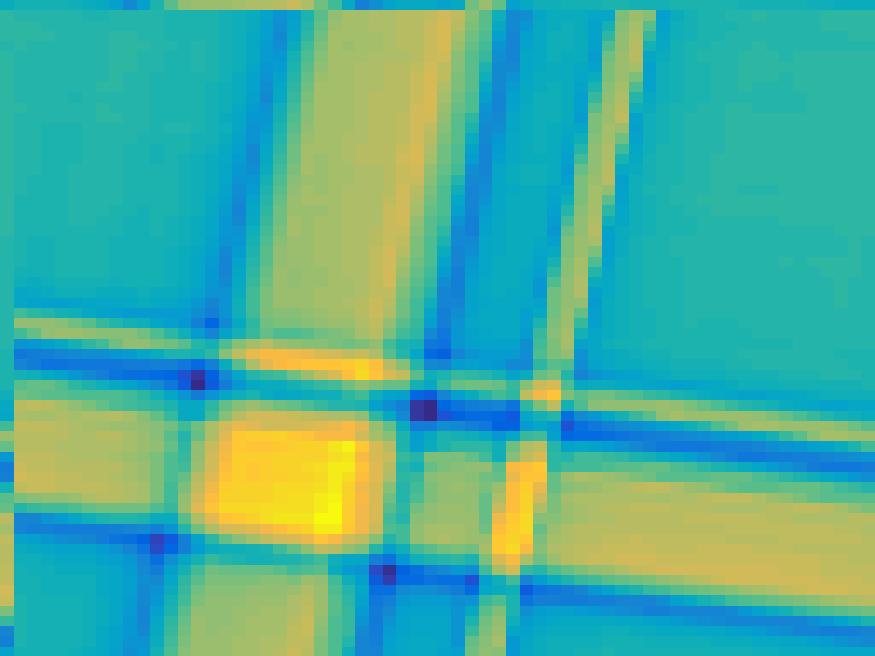}&
\includegraphics[width=.1\linewidth,height=.1\linewidth]{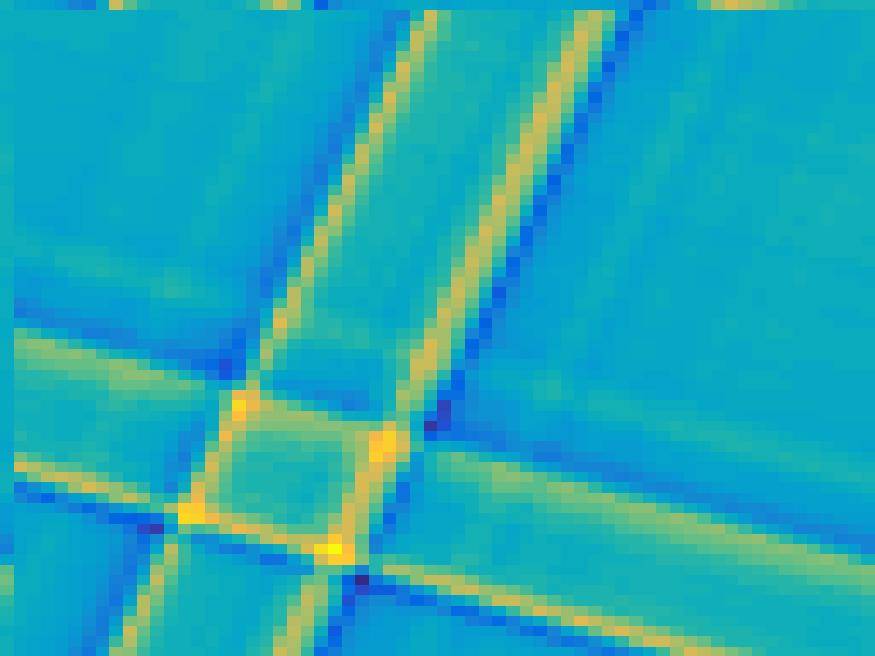}&
\includegraphics[width=.1\linewidth,height=.1\linewidth]{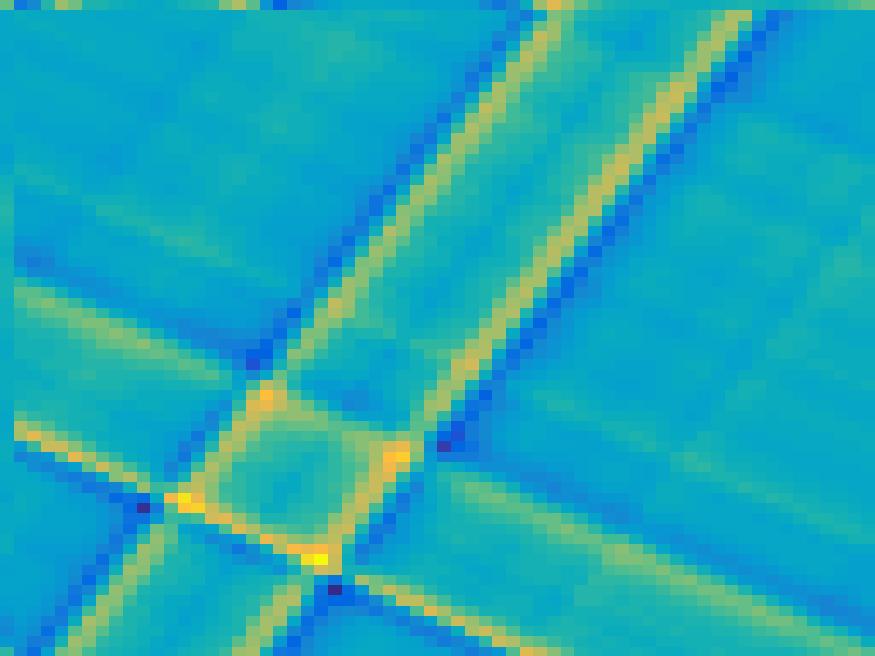}&
\includegraphics[width=.1\linewidth,height=.1\linewidth]{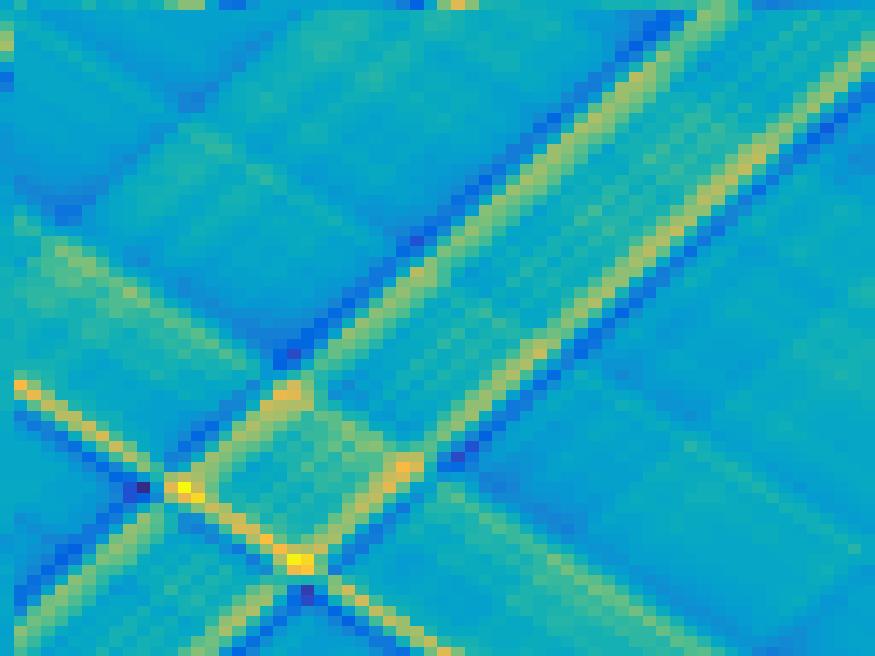}&
\includegraphics[width=.1\linewidth,height=.1\linewidth]{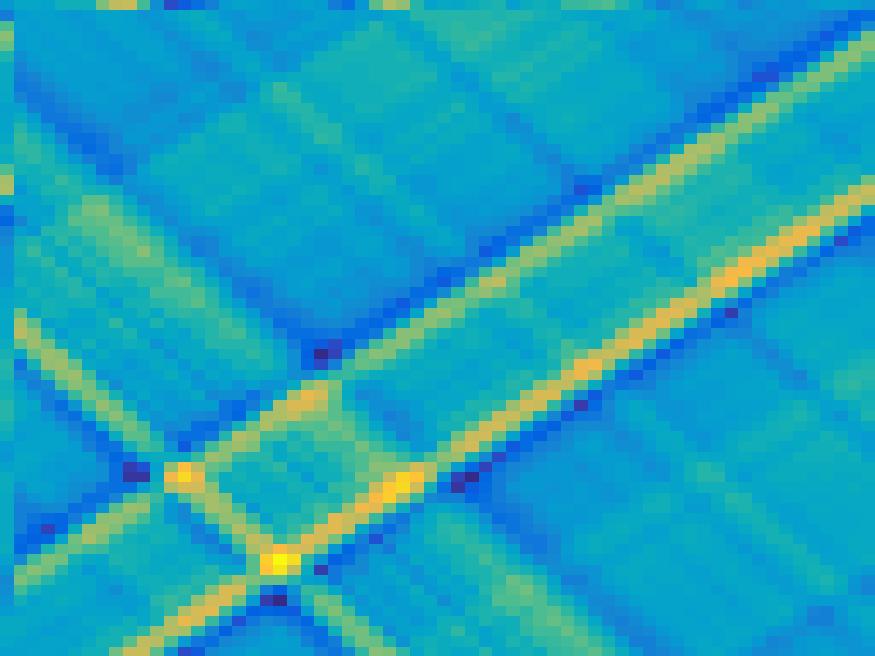}&
\includegraphics[width=.1\linewidth,height=.1\linewidth]{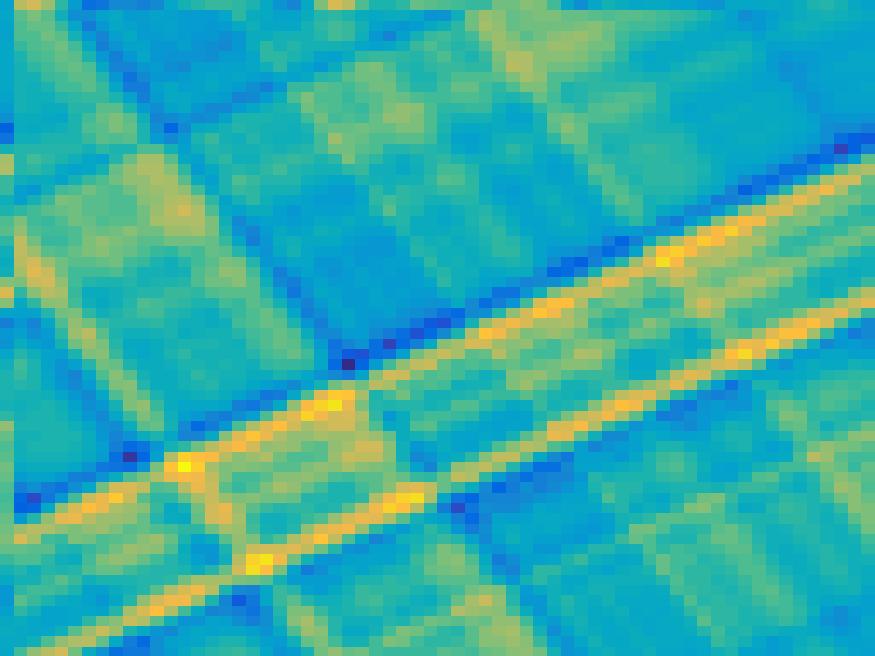}&
\includegraphics[width=.1\linewidth,height=.1\linewidth]{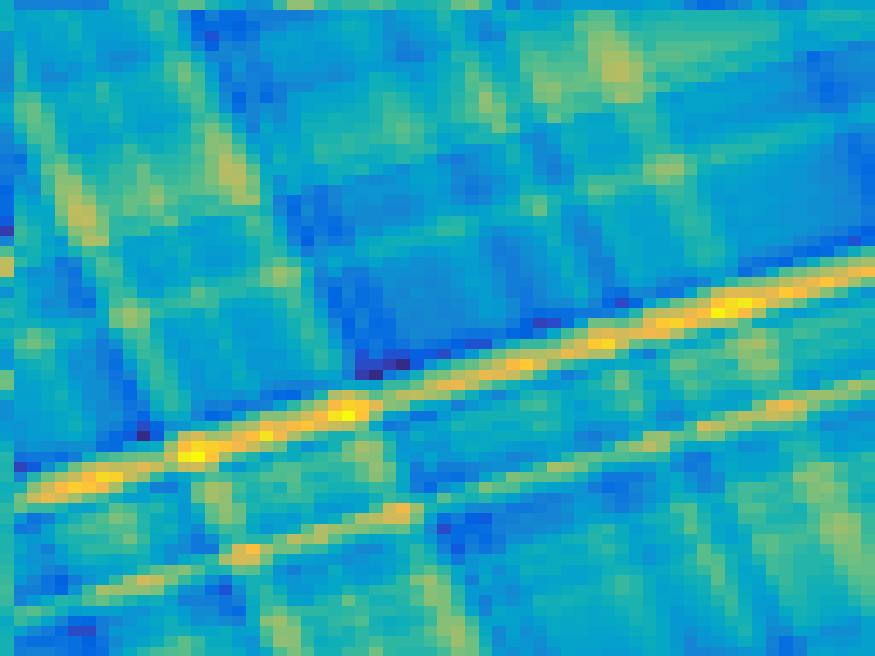}&
\includegraphics[width=.1\linewidth,height=.1\linewidth]{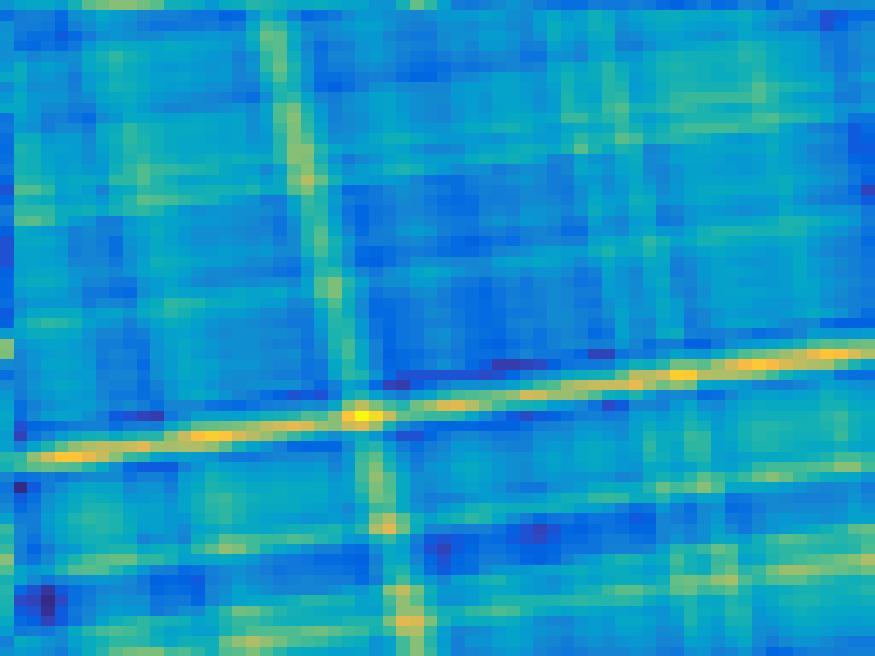}\\
\includegraphics[width=.1\linewidth,height=.1\linewidth]{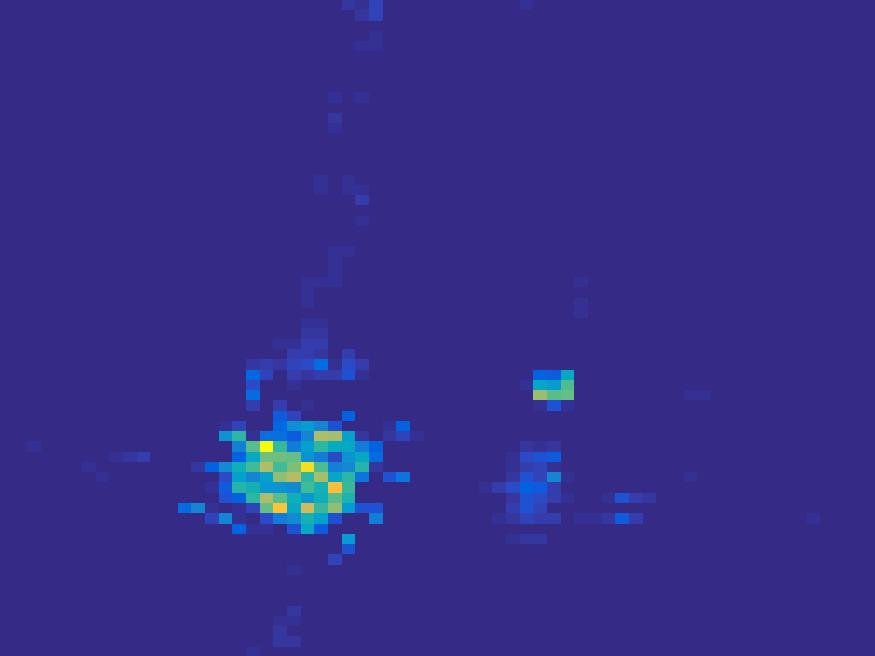}&
\includegraphics[width=.1\linewidth,height=.1\linewidth]{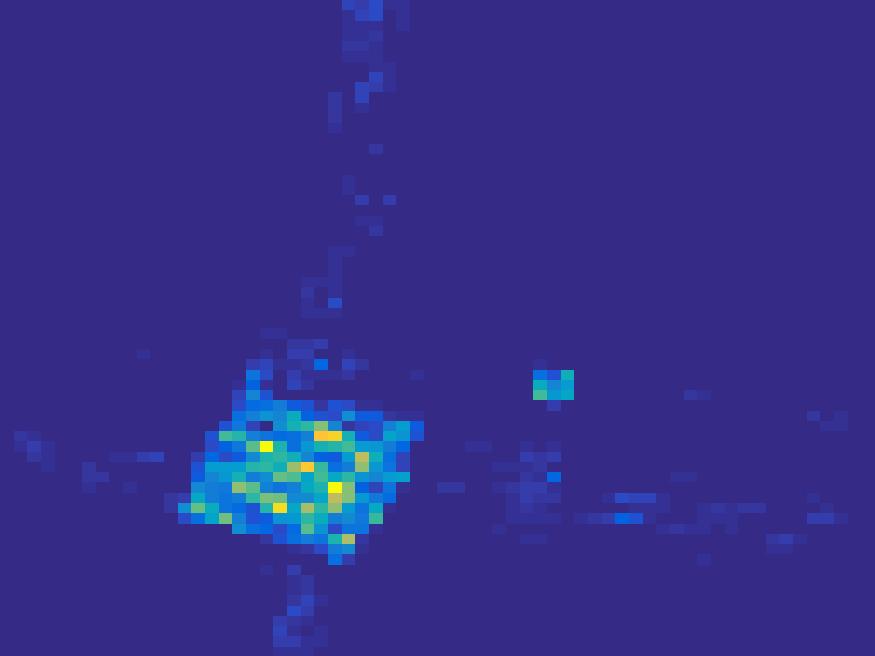}&
\includegraphics[width=.1\linewidth,height=.1\linewidth]{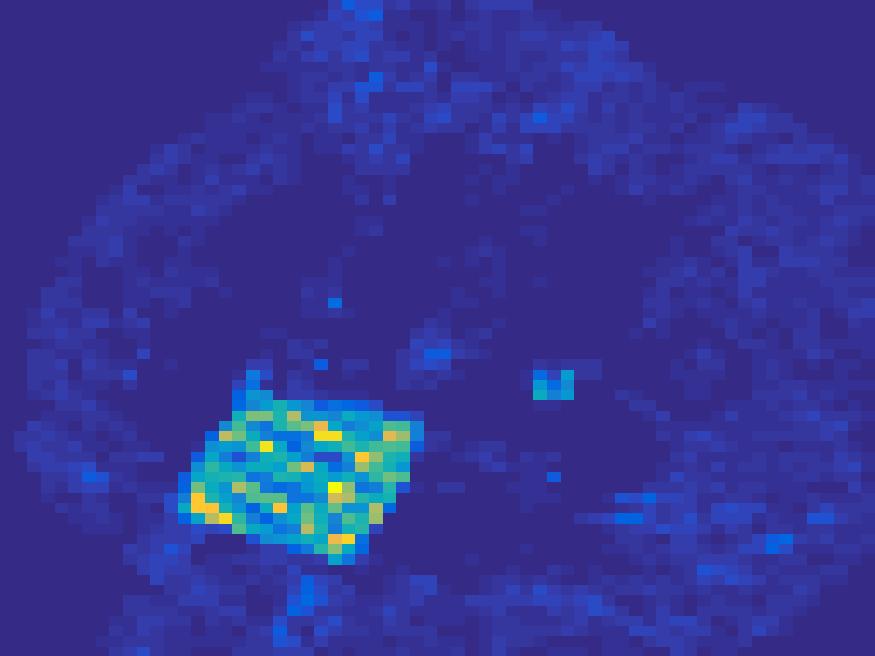}&
\includegraphics[width=.1\linewidth,height=.1\linewidth]{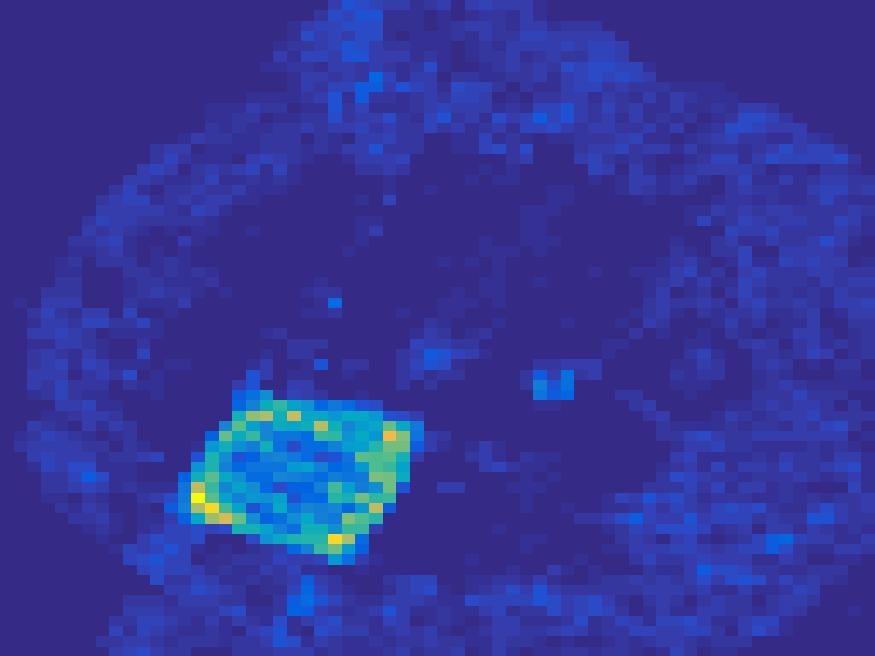}&
\includegraphics[width=.1\linewidth,height=.1\linewidth]{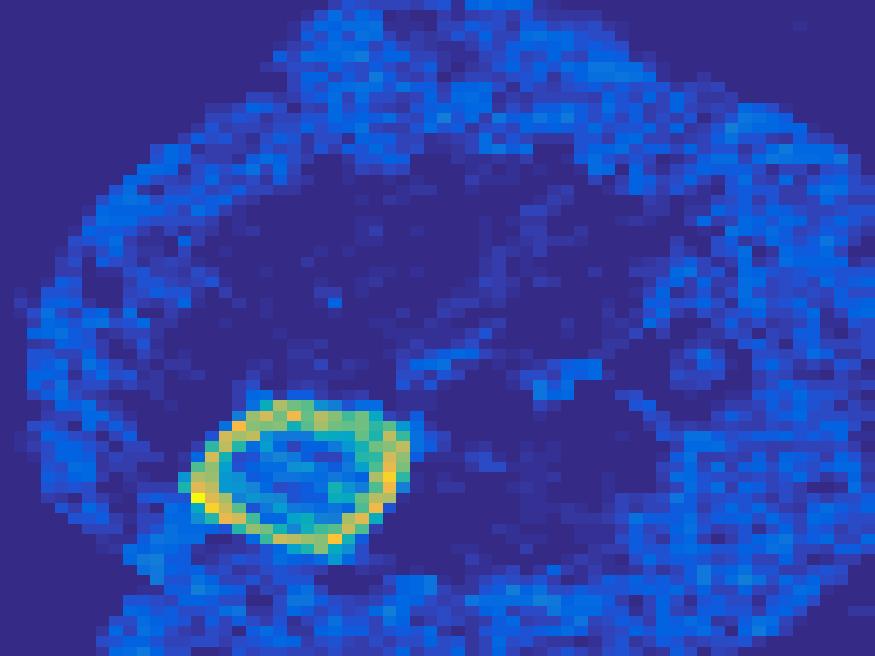}&
\includegraphics[width=.1\linewidth,height=.1\linewidth]{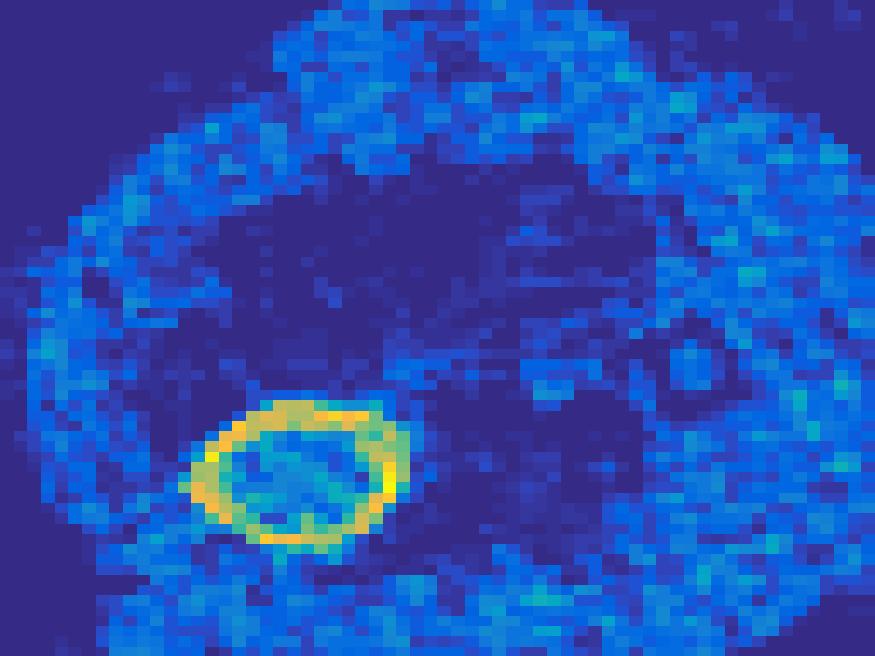}&
\includegraphics[width=.1\linewidth,height=.1\linewidth]{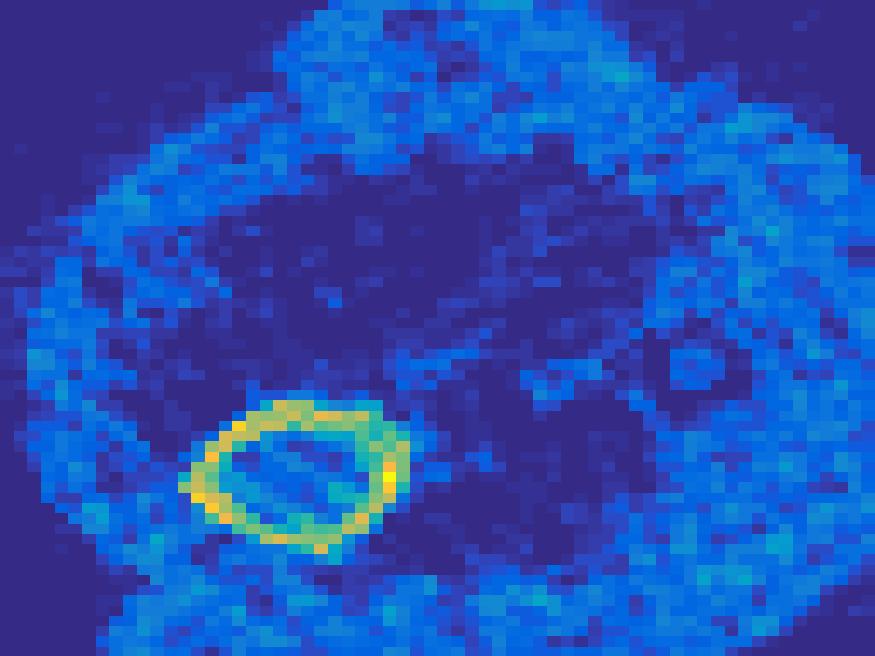}&
\includegraphics[width=.1\linewidth,height=.1\linewidth]{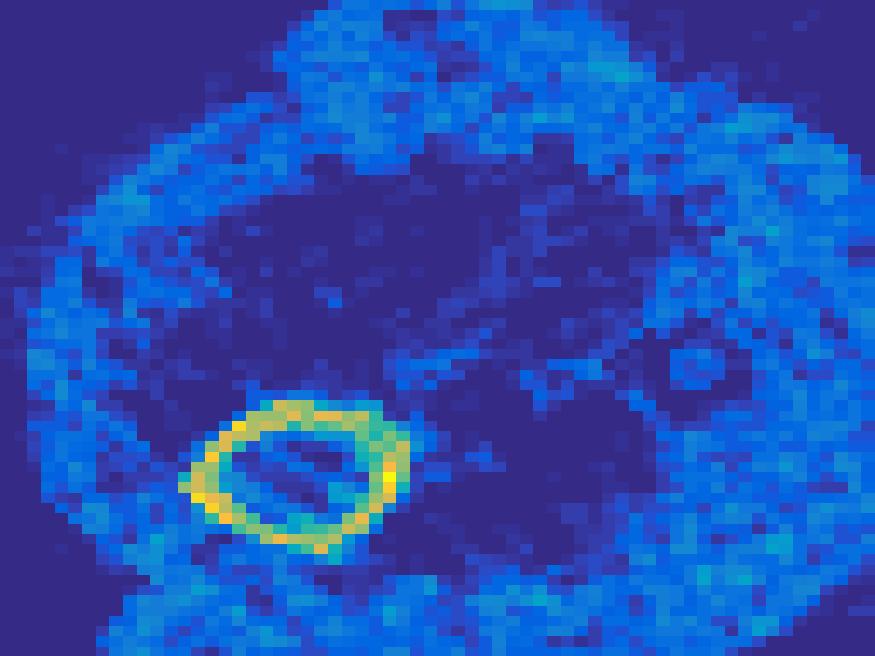}&
\includegraphics[width=.1\linewidth,height=.1\linewidth]{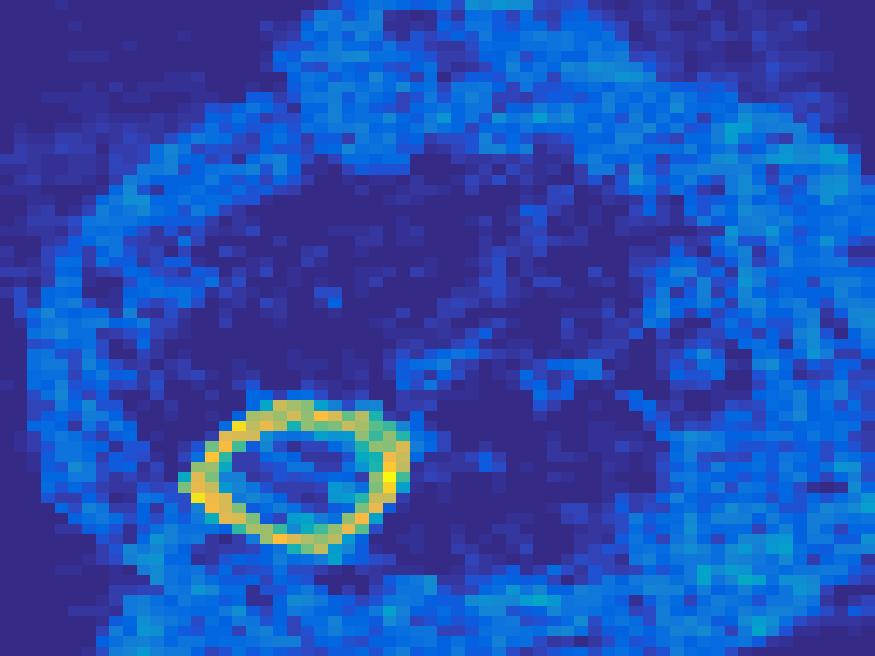}\\
\includegraphics[width=.1\linewidth,height=.1\linewidth]{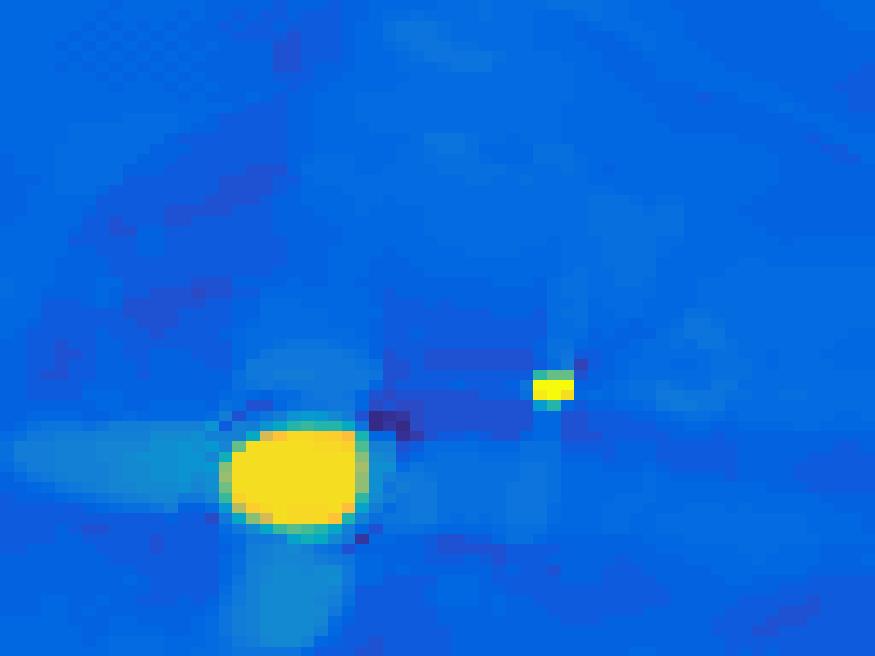}&
\includegraphics[width=.1\linewidth,height=.1\linewidth]{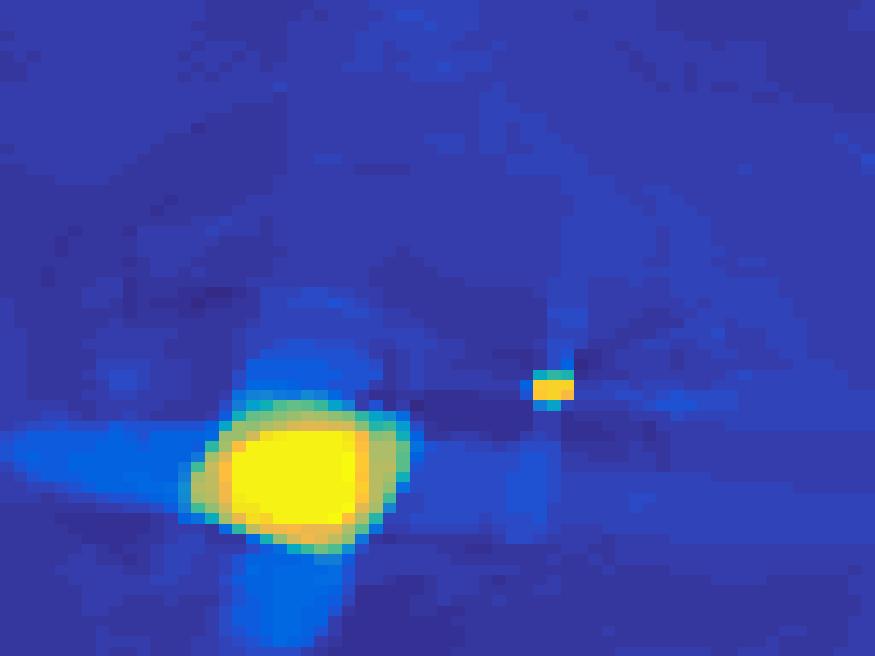}&
\includegraphics[width=.1\linewidth,height=.1\linewidth]{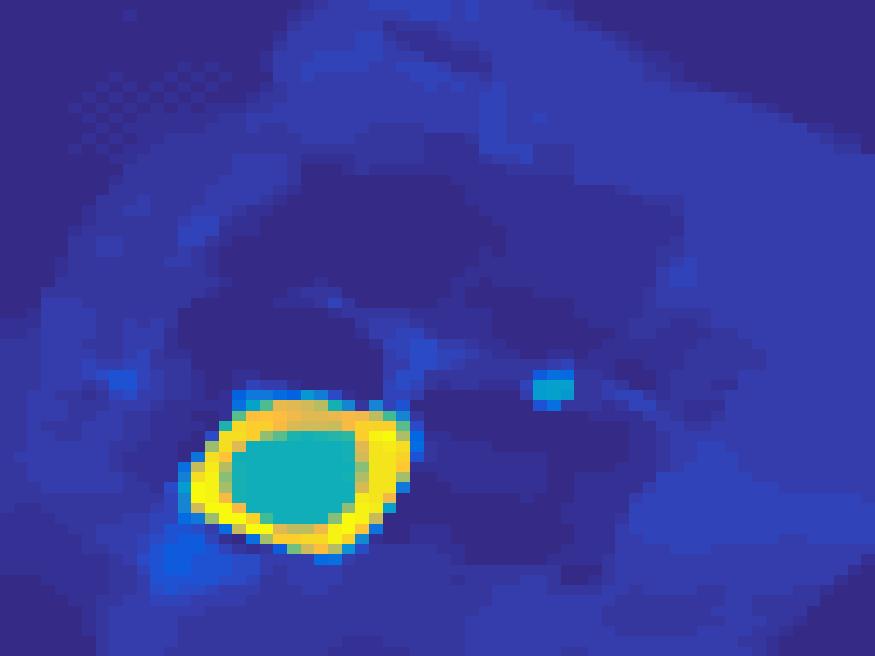}&
\includegraphics[width=.1\linewidth,height=.1\linewidth]{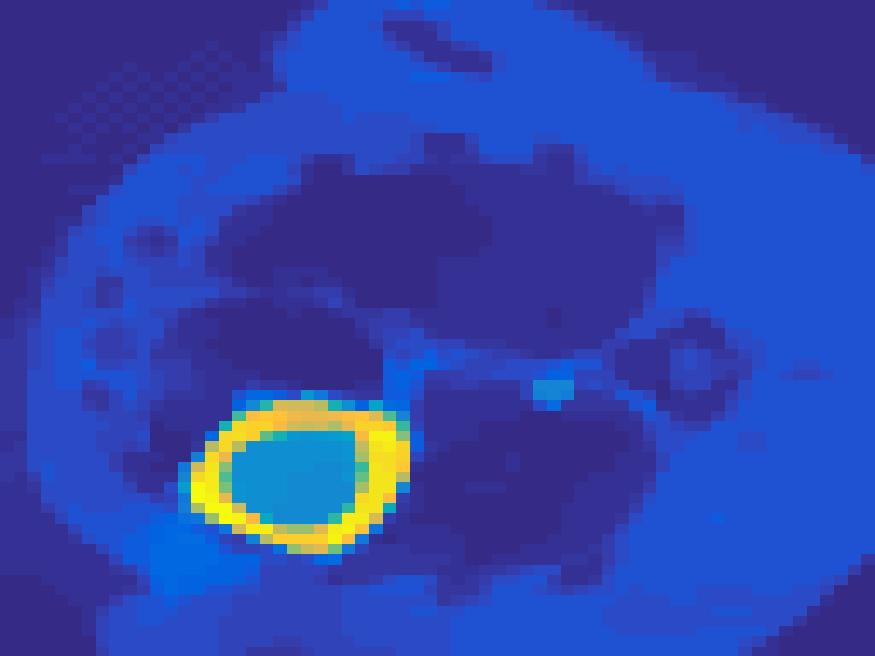}&
\includegraphics[width=.1\linewidth,height=.1\linewidth]{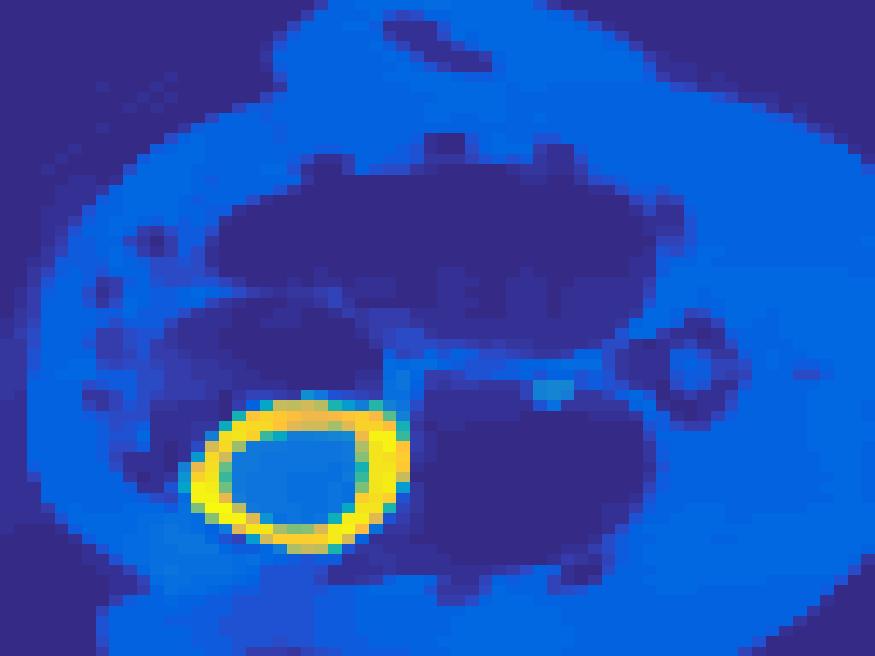}&
\includegraphics[width=.1\linewidth,height=.1\linewidth]{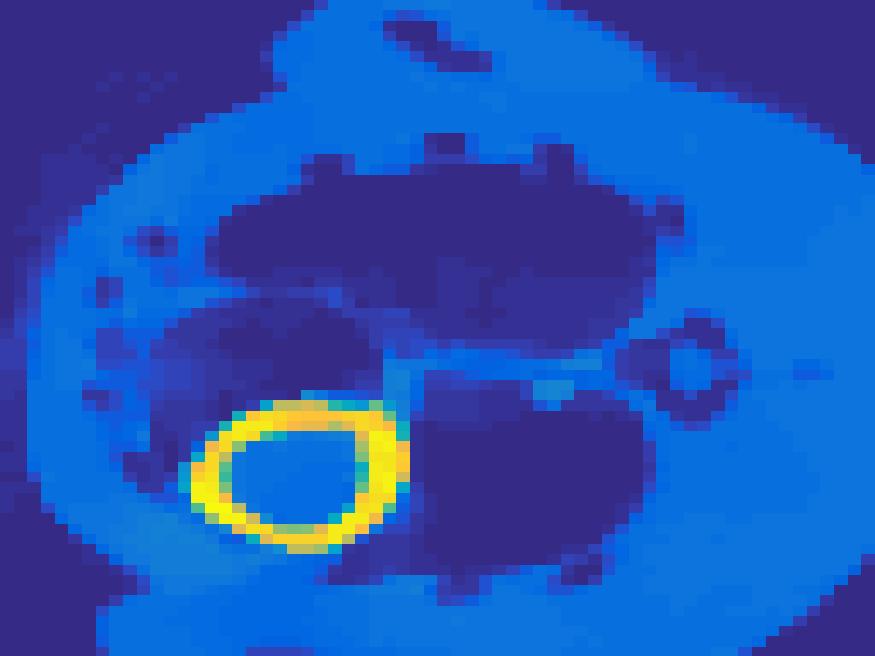}&
\includegraphics[width=.1\linewidth,height=.1\linewidth]{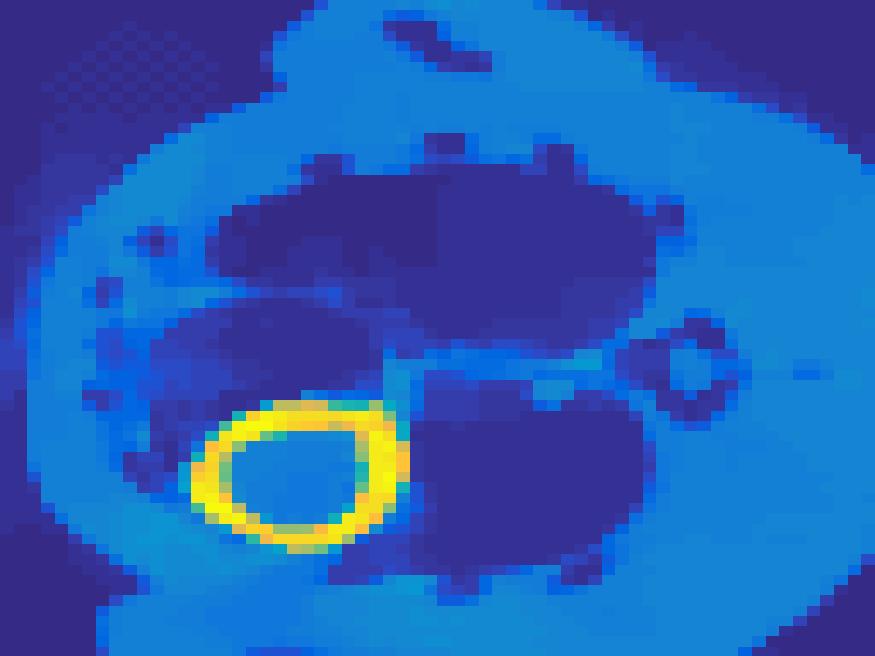}&
\includegraphics[width=.1\linewidth,height=.1\linewidth]{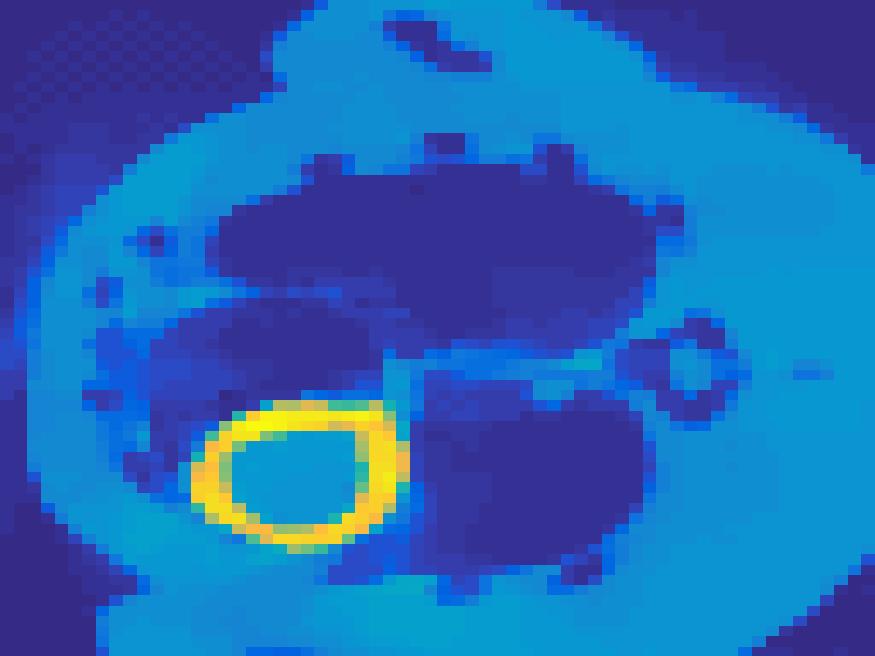}&
\includegraphics[width=.1\linewidth,height=.1\linewidth]{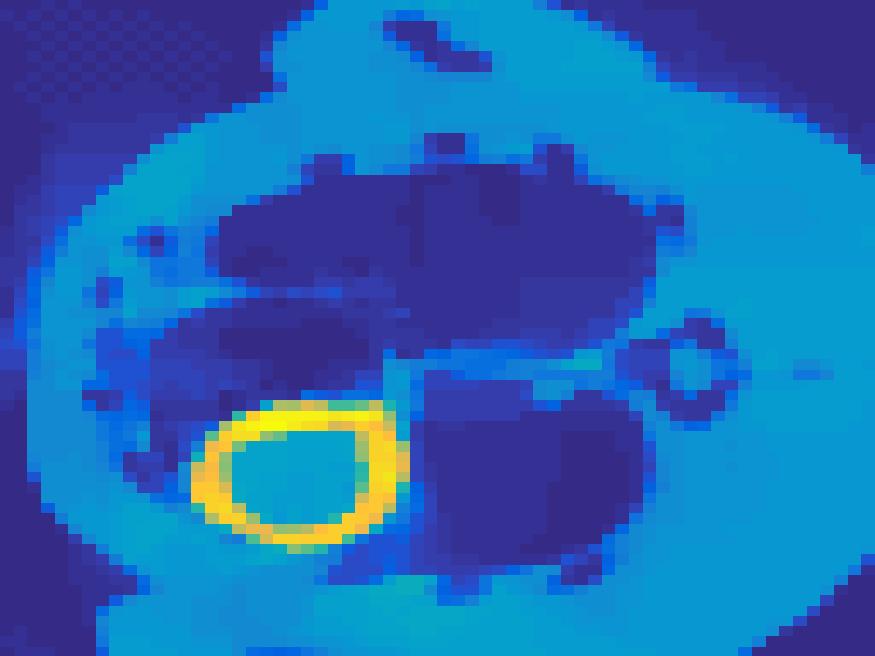}\\
\includegraphics[width=.1\linewidth,height=.1\linewidth]{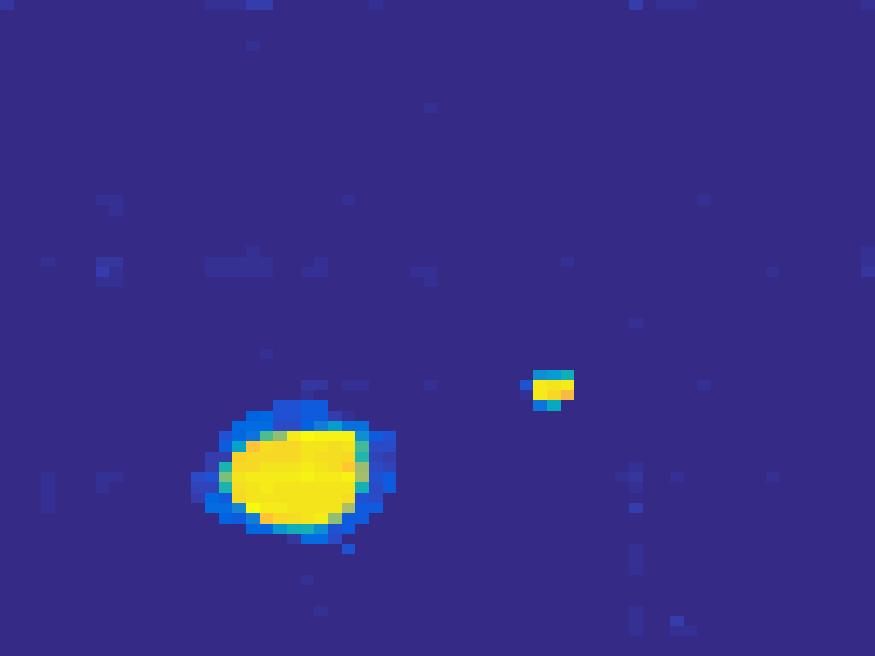}&
\includegraphics[width=.1\linewidth,height=.1\linewidth]{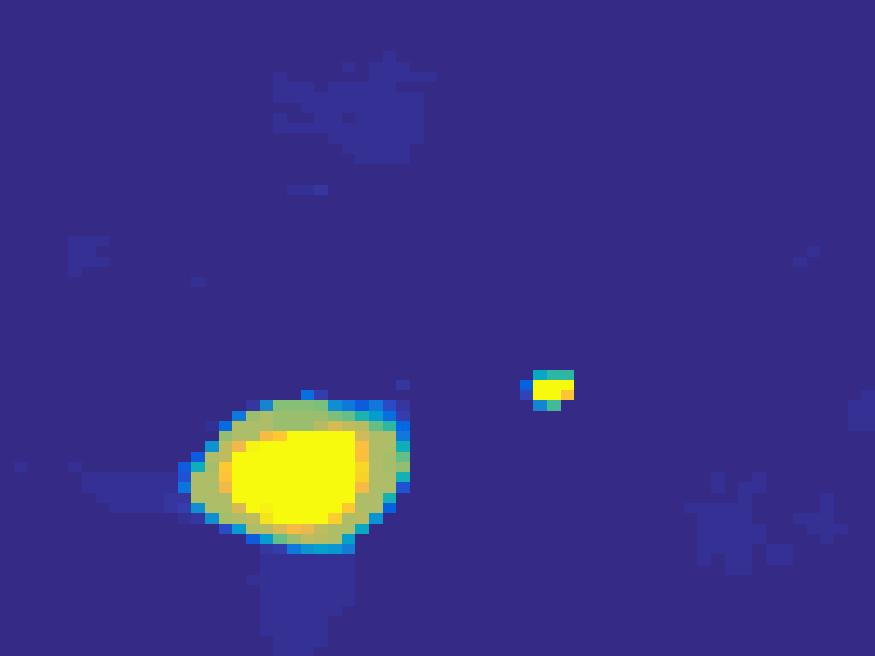}&
\includegraphics[width=.1\linewidth,height=.1\linewidth]{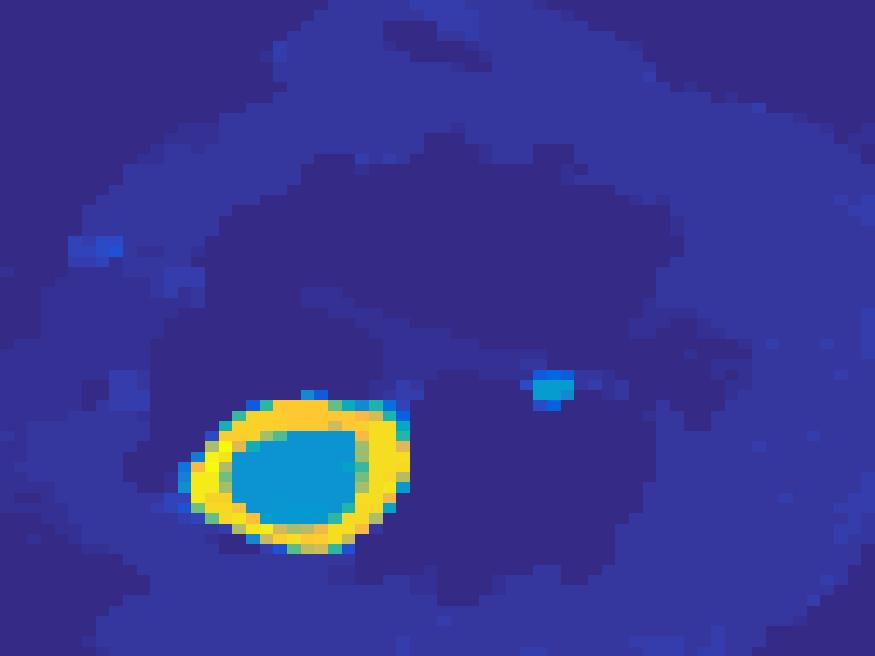}&
\includegraphics[width=.1\linewidth,height=.1\linewidth]{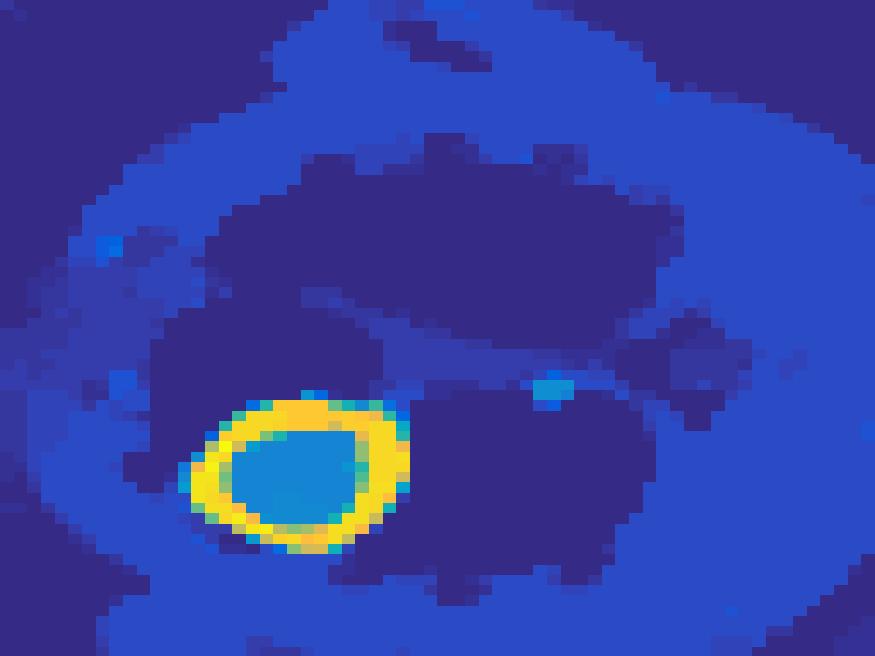}&
\includegraphics[width=.1\linewidth,height=.1\linewidth]{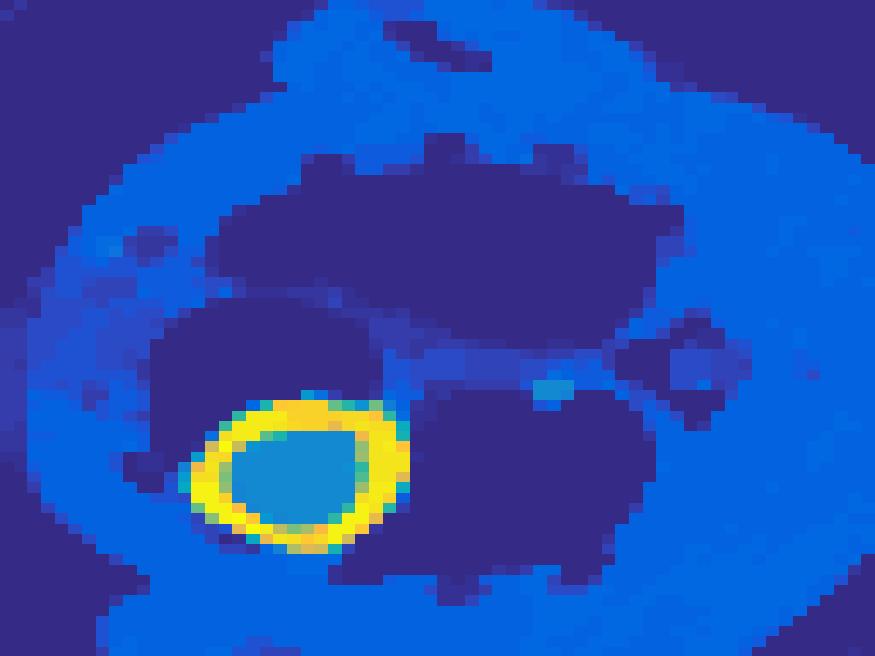}&
\includegraphics[width=.1\linewidth,height=.1\linewidth]{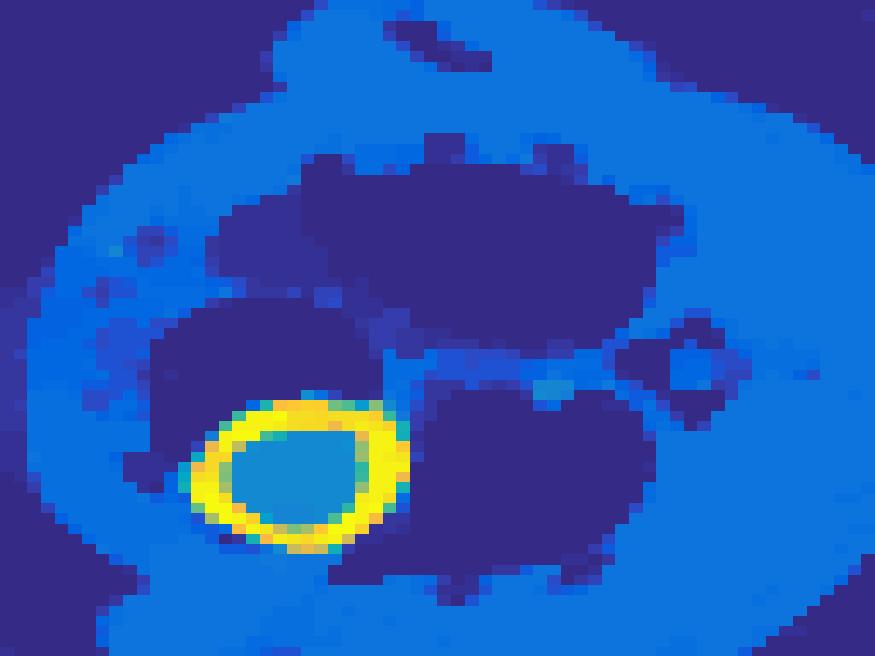}&
\includegraphics[width=.1\linewidth,height=.1\linewidth]{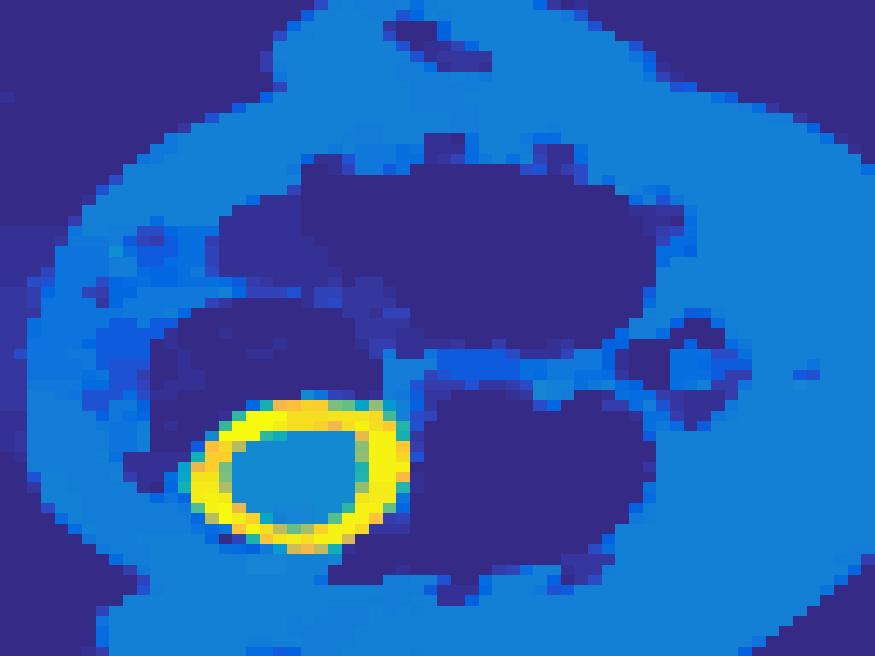}&
\includegraphics[width=.1\linewidth,height=.1\linewidth]{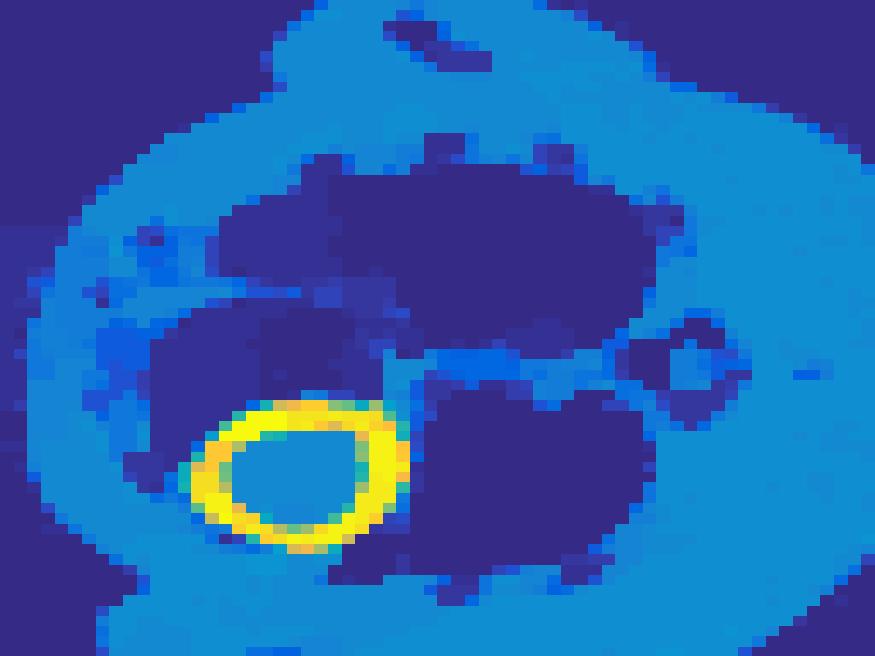}&
\includegraphics[width=.1\linewidth,height=.1\linewidth]{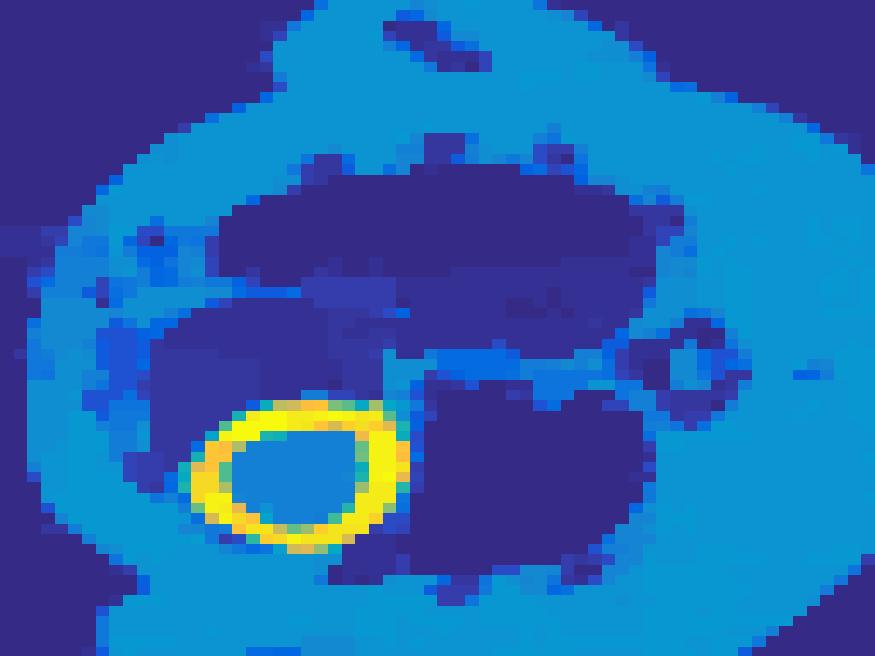}\\
{\footnotesize Frame 1}&
{\footnotesize Frame 11}&
{\footnotesize Frame 21}&
{\footnotesize Frame 31}&
{\footnotesize Frame 41}&
{\footnotesize Frame 51}&
{\footnotesize Frame 61}&
{\footnotesize Frame 71}&
{\footnotesize Frame 81}\\
\end{tabular}
\caption {First row: Ground truth; Second row: FBP; Third  row: least square method;  Forth row: SEMF \cite{ding2015dynamic}; Fifth row: Proposed model.}
\label{fig:ICgaussLiver}
\end{figure}

\begin{figure}[ht]
\begin{center}
\subfigure{
\includegraphics[width=.45\linewidth,height=.3\linewidth]{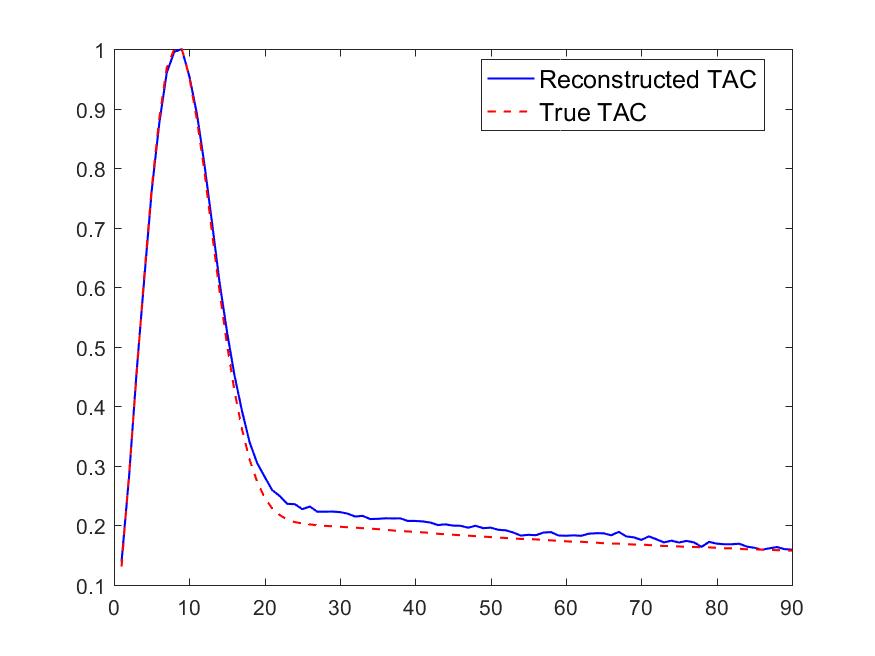}}
\subfigure{
\includegraphics[width=.45\linewidth,height=.3\linewidth]{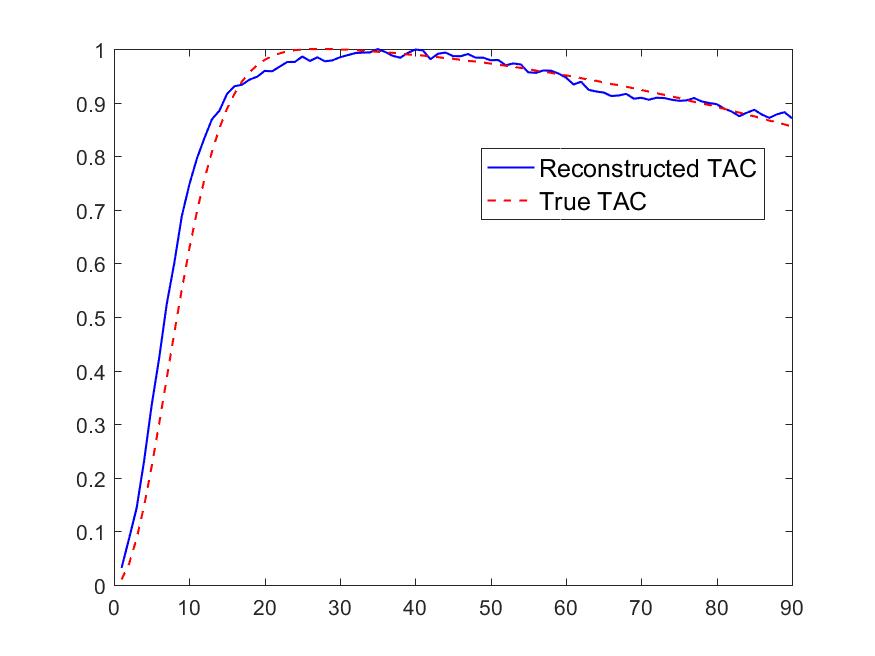}}\\
\subfigure{
\includegraphics[width=.45\linewidth,height=.3\linewidth]{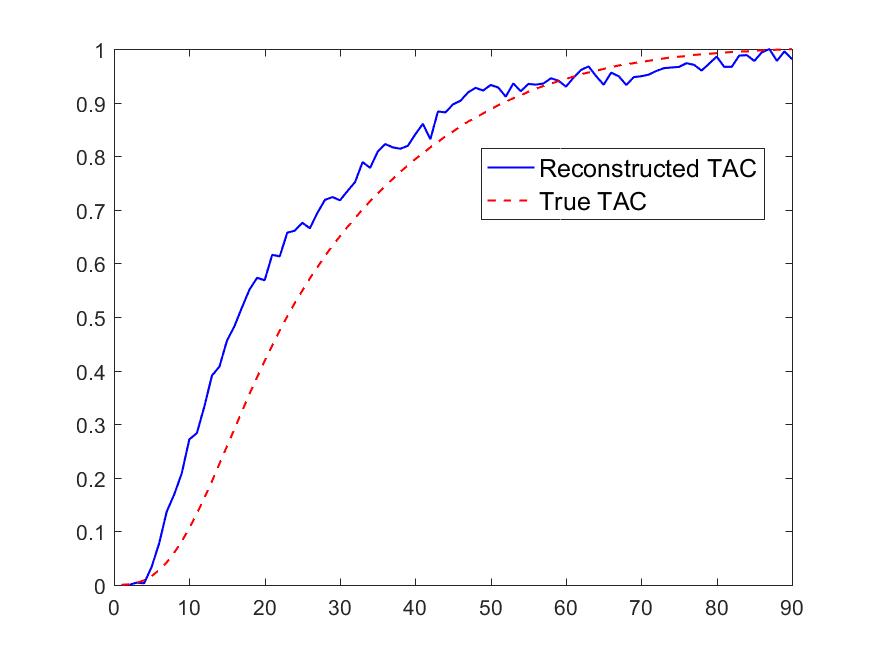}}
\subfigure{
\includegraphics[width=.45\linewidth,height=.3\linewidth]{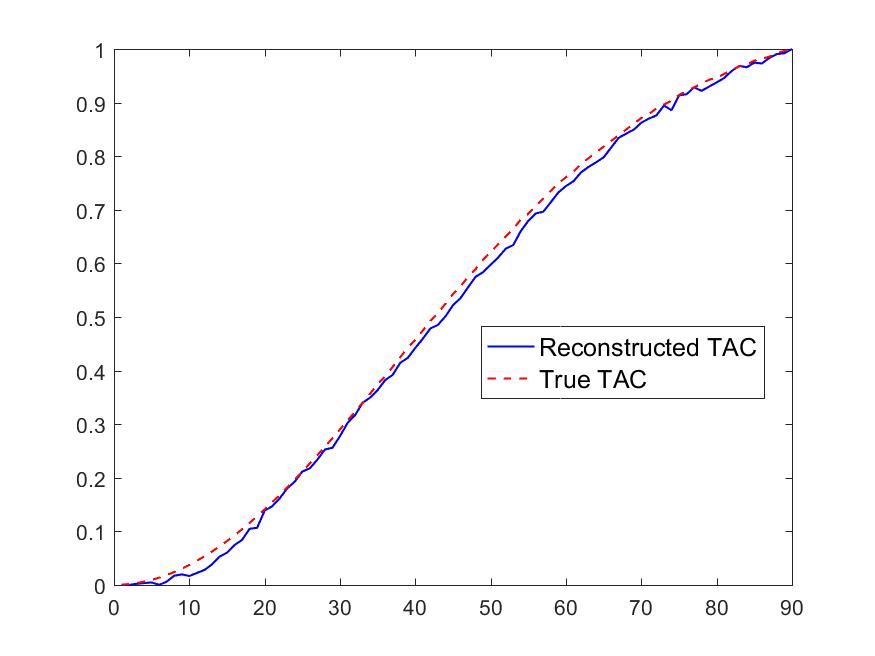}}
\end{center}
\caption{Comparison of the true TACs and the reconstructed TACs by the proposed method.}
\label{fig:ICGauLiverTAC}
\end{figure}

 For the simulated  images of rat's abdomen, the same procedure is applied to generate projection data. Also,  $10\%$  noise was added to the sinogram. Figure \ref{fig:ICgaussLiver} compares the frames reconstructed by different methods. Clearly, the traditional FBP method and least square method cannot reconstruct the dynamic images with very few projections, however the proposed method reconstructs the images quite accurately. Figure \ref{fig:ICGauLiverTAC} illustrates the comparison of the true TACs and those reconstructed
by the proposed method. We can see that they are quite accurate and  present small errors.

%
\begin{figure}[ht]
\begin{center}
\subfigure{
\includegraphics[width=.45\linewidth]{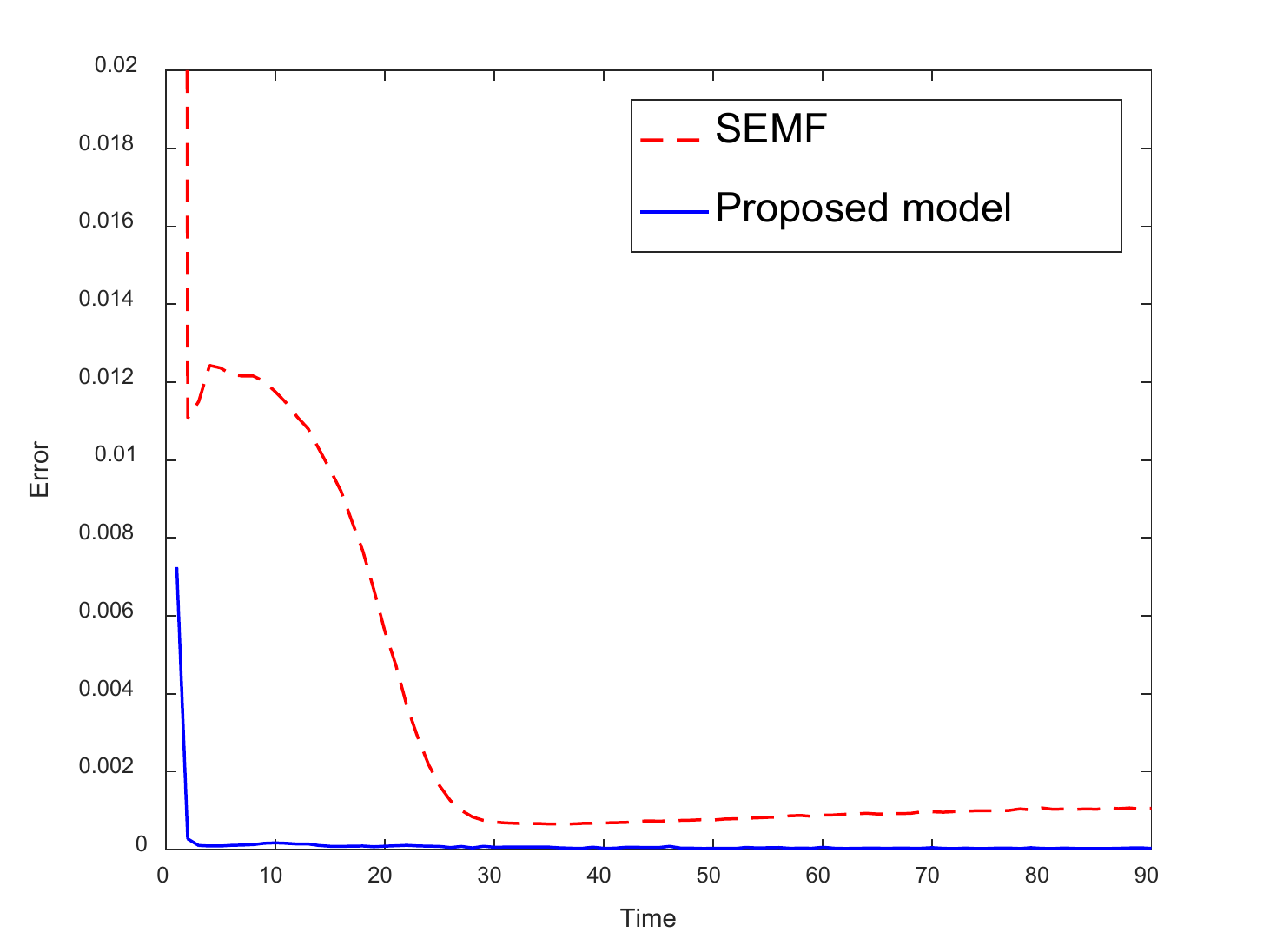}}
\subfigure{
\includegraphics[width=.45\linewidth]{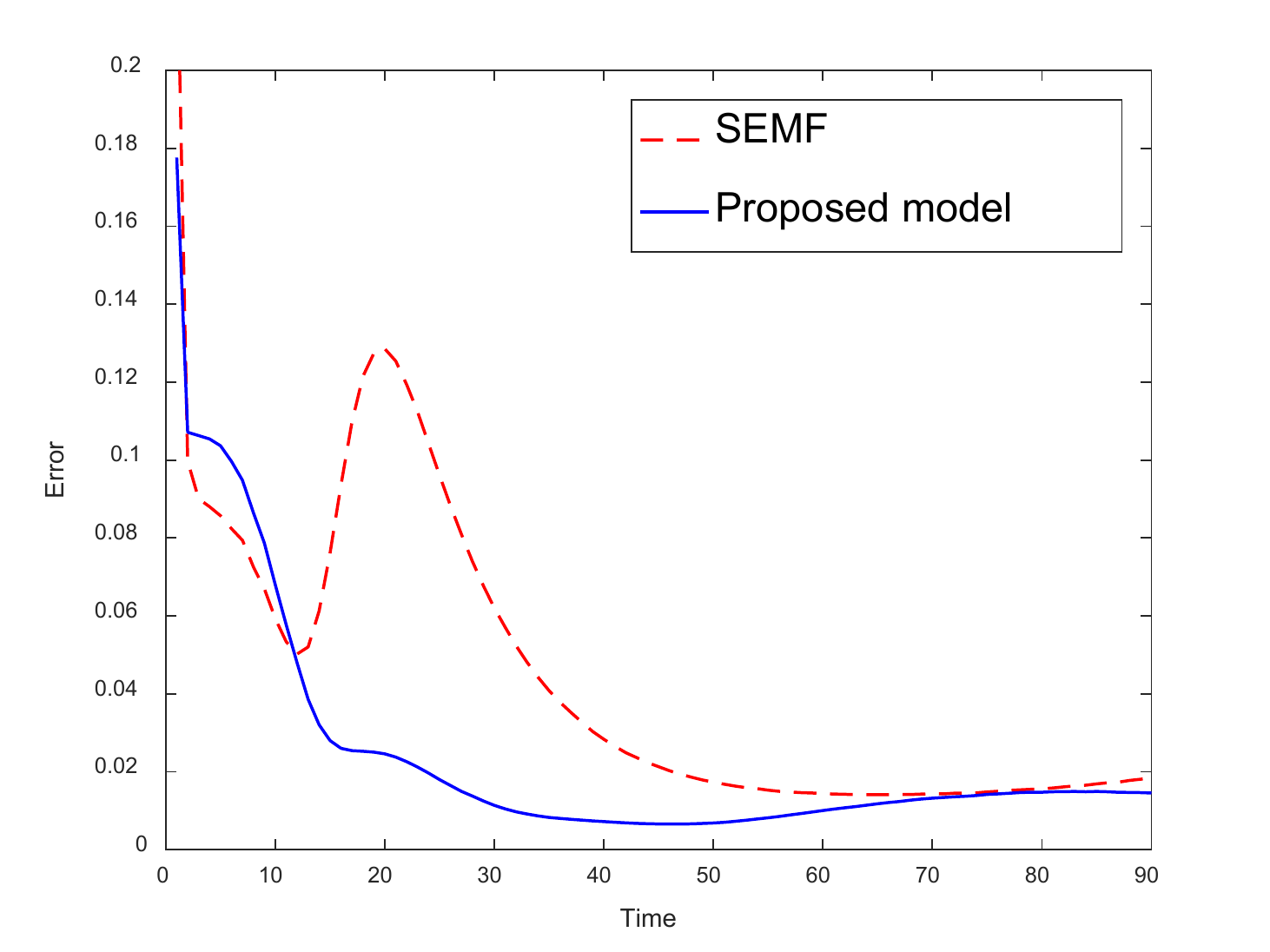}}
\end{center}
\caption{Comparison of the relative error  for $t$-th frames by SEMF and the proposed method. Left: ellipse phantom; Right: rat's abdomen phantom.}
\label{fig:errGau}
\end{figure}

 The relative error of the reconstructed image and the true one for $t$-th frame is defined as  $\frac{\|U(:,t)_{rec}-U(:,t)_{true}\|_{2}^{2}}{\|U(:,t)_{true} \|_{2}^{2}}$, where $U_{rec}$ is the reconstructed frame by the proposed method and $U_{true}$ is the ground truth image.
 Figure \ref{fig:errGau} demonstrates the relative error of $T$ images reconstructed by proposed model and SEMF of two different datasets.
 The solid lines are the errors of the reconstruction images by proposed model and the dash lines are the errors of the reconstruction images by SEMF.
 We can see that the relative error is smaller by proposed model compared with SEMF. This is due to the fact that in the proposed method, we set the former and later images as reference and  the referred images can provide edge information for the image.

\subsection{Poisson noise}

\subsubsection{Simulated Poisson noise} In SPECT/PET,  Poisson noise are usually more common.
To obtain a Poisson corrupted   projection data, we scale the data by the maximum of $f$ and corrupted the data with Poisson noise by using the Matlab command \texttt{ poissrnd}. The reconstructed image is obtained by  applying the proposed model on the rescaling back sinogram .

Again, the proposed method is compared with  FBP,  and alternating applying the EM algorithm and update of the basis  for solving  
$\min_{\alpha,B}D_{KL}(f,A\alpha B^{T}).$
As for the initials $U$, $\alpha$ and $B$ of our methods,  we use  uniformed B-spline, $B\in\mathbb{R}^{90\times20}$ as initial basis $B$  and solve the above model to obtain
 for $\alpha$ and $B $ and  $U=\alpha B^T$.

\begin{figure}
\begin{tabular}{c@{\hspace{2pt}}c@{\hspace{2pt}}c@{\hspace{2pt}}c@{\hspace{2pt}}c@{\hspace{2pt}}c@{\hspace{2pt}}c@{\hspace{2pt}}c@{\hspace{2pt}}c@{\hspace{2pt}}c}
\includegraphics[width=.1\linewidth,height=.1\linewidth]{IellReal1}&
\includegraphics[width=.1\linewidth,height=.1\linewidth]{IellReal11}&
\includegraphics[width=.1\linewidth,height=.1\linewidth]{IellReal21}&
\includegraphics[width=.1\linewidth,height=.1\linewidth]{IellReal31}&
\includegraphics[width=.1\linewidth,height=.1\linewidth]{IellReal41}&
\includegraphics[width=.1\linewidth,height=.1\linewidth]{IellReal51}&
\includegraphics[width=.1\linewidth,height=.1\linewidth]{IellReal61}&
\includegraphics[width=.1\linewidth,height=.1\linewidth]{IellReal71}&
\includegraphics[width=.1\linewidth,height=.1\linewidth]{IellReal81}\\
\includegraphics[width=.1\linewidth,height=.1\linewidth]{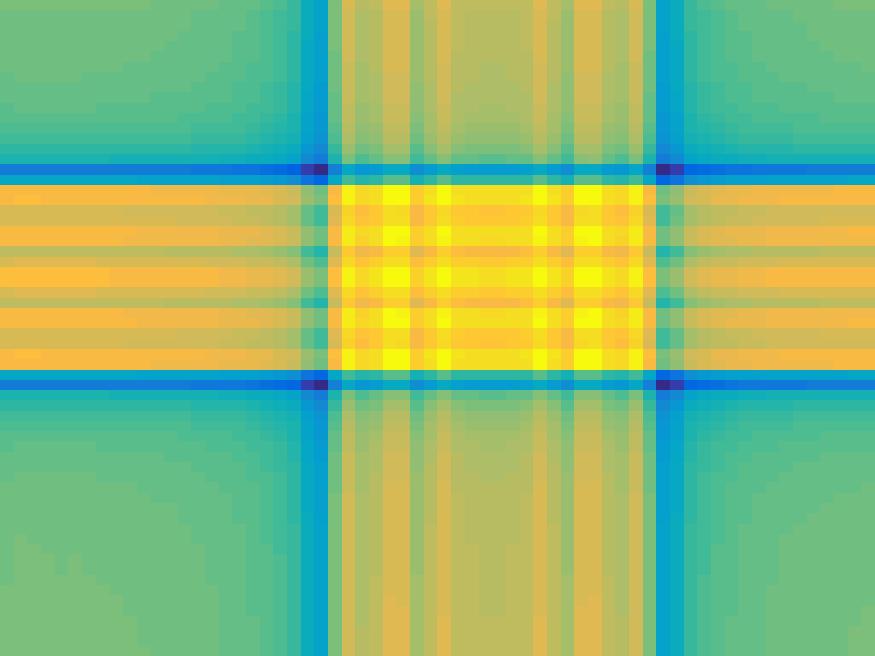}&
\includegraphics[width=.1\linewidth,height=.1\linewidth]{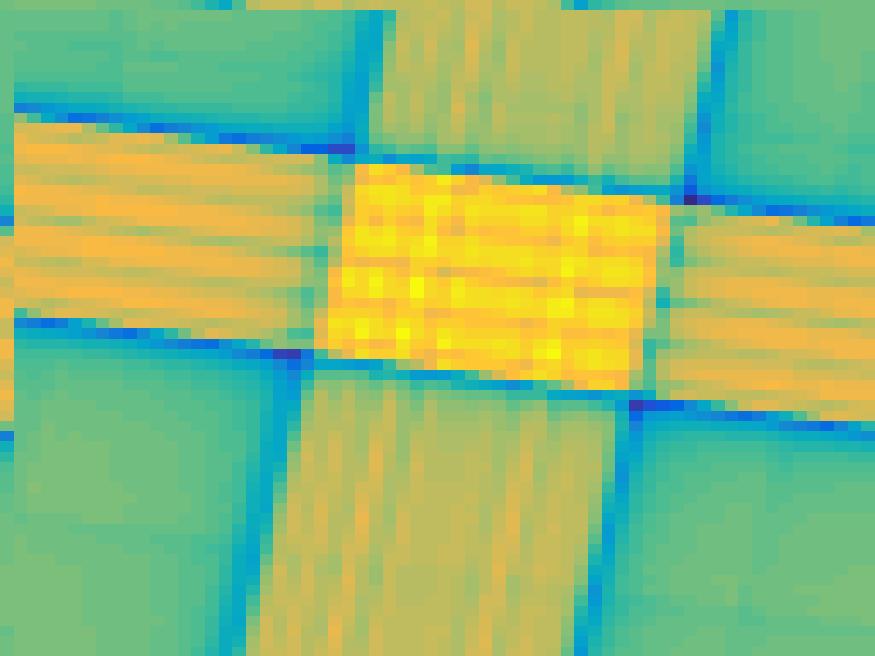}&
\includegraphics[width=.1\linewidth,height=.1\linewidth]{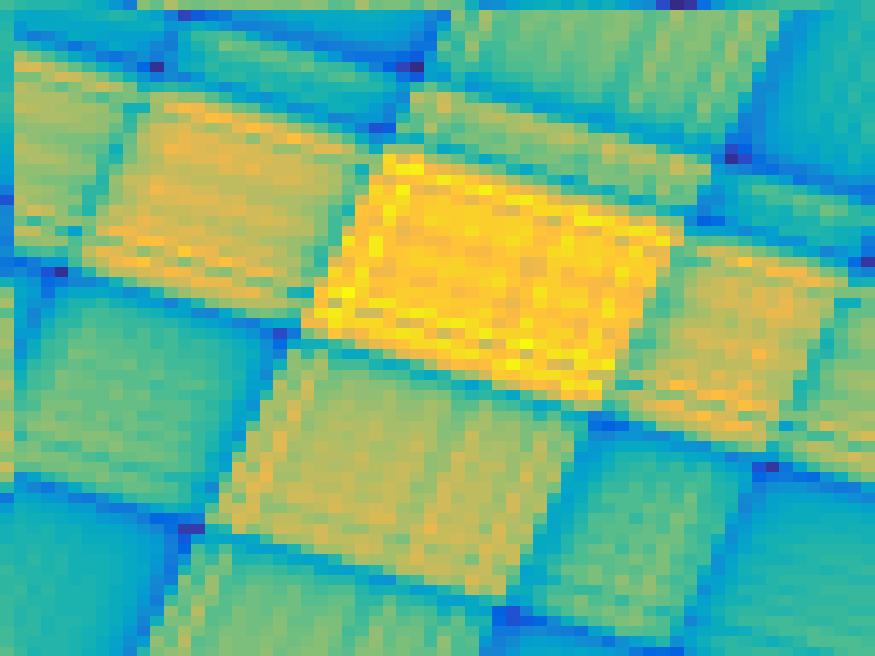}&
\includegraphics[width=.1\linewidth,height=.1\linewidth]{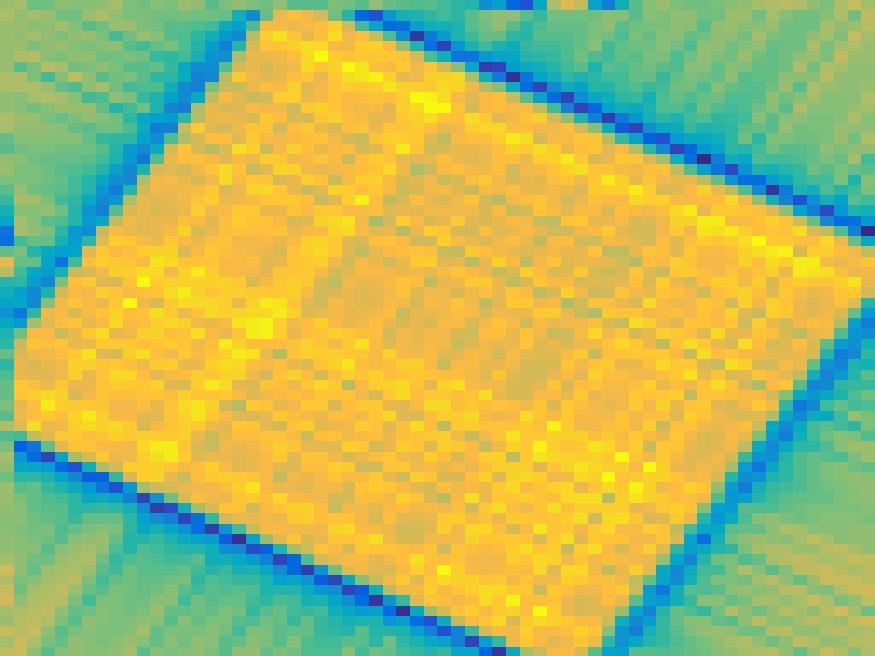}&
\includegraphics[width=.1\linewidth,height=.1\linewidth]{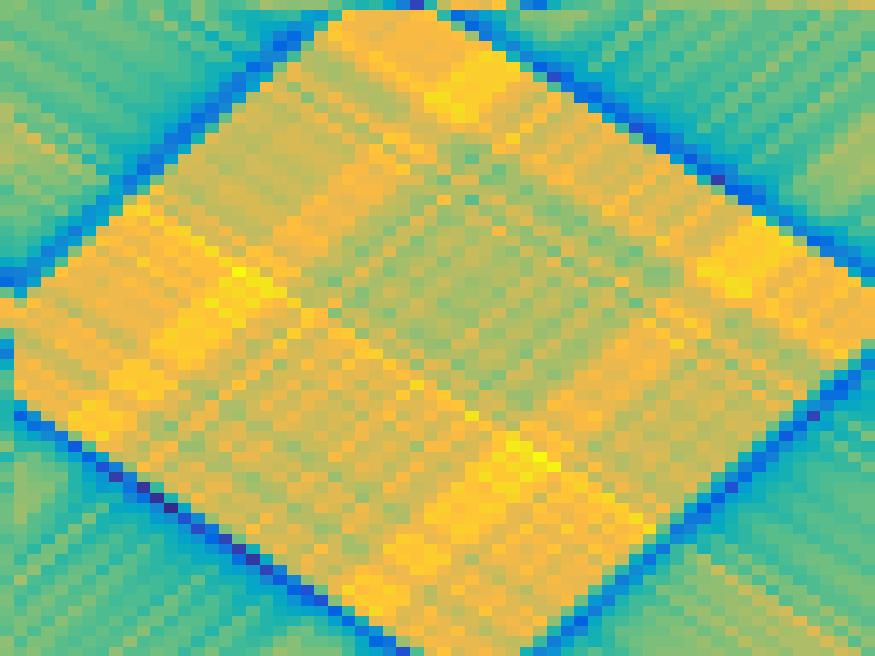}&
\includegraphics[width=.1\linewidth,height=.1\linewidth]{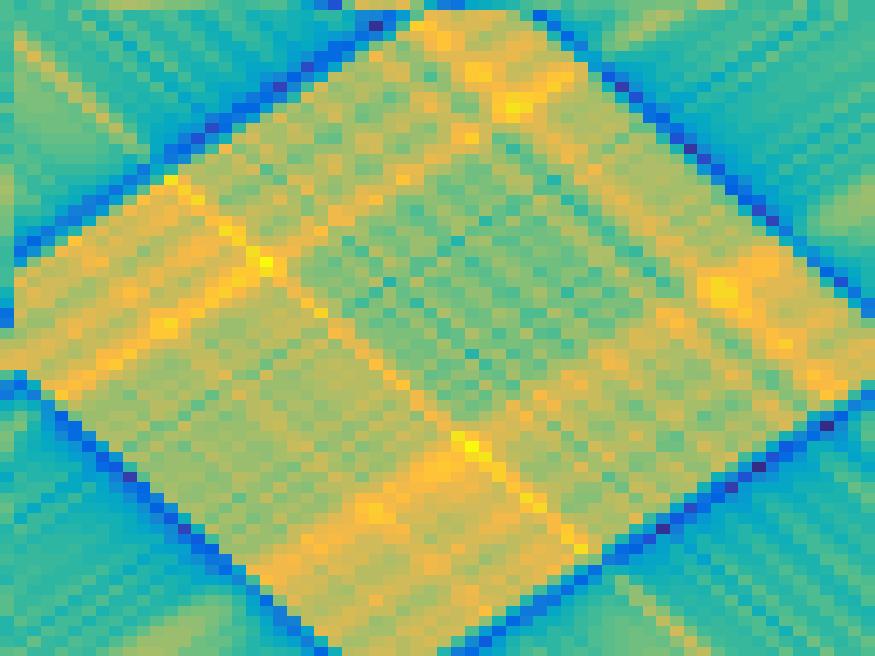}&
\includegraphics[width=.1\linewidth,height=.1\linewidth]{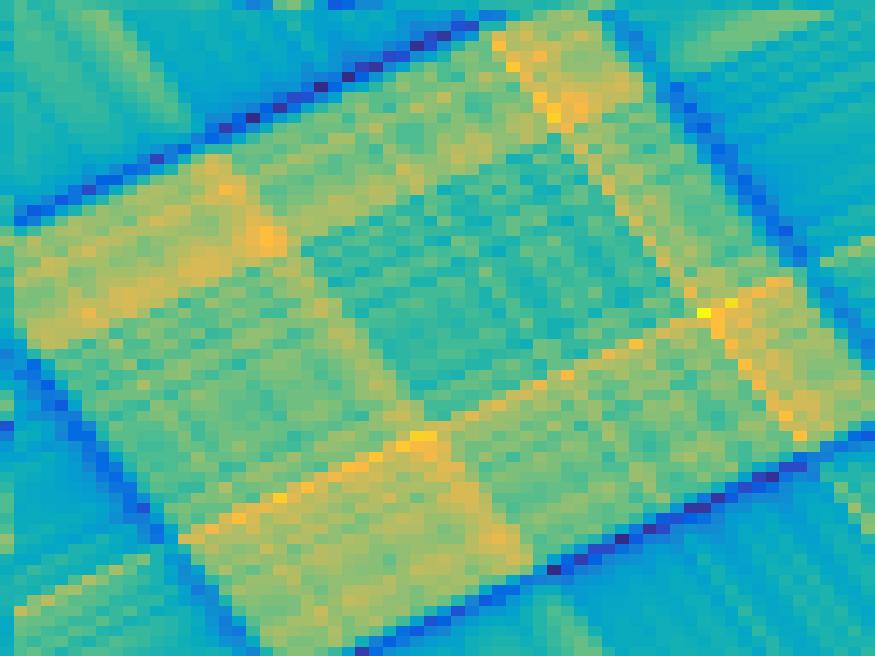}&
\includegraphics[width=.1\linewidth,height=.1\linewidth]{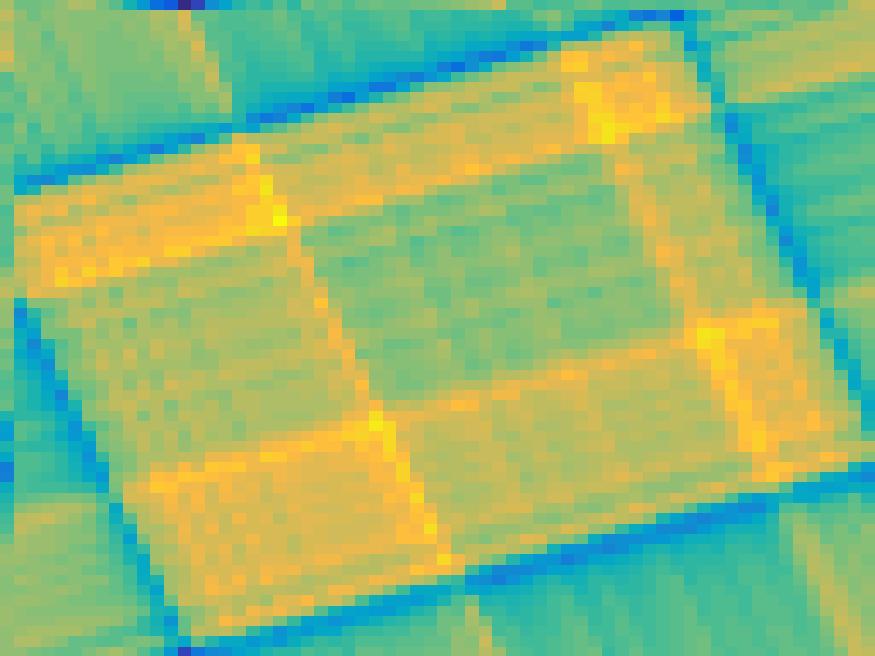}&
\includegraphics[width=.1\linewidth,height=.1\linewidth]{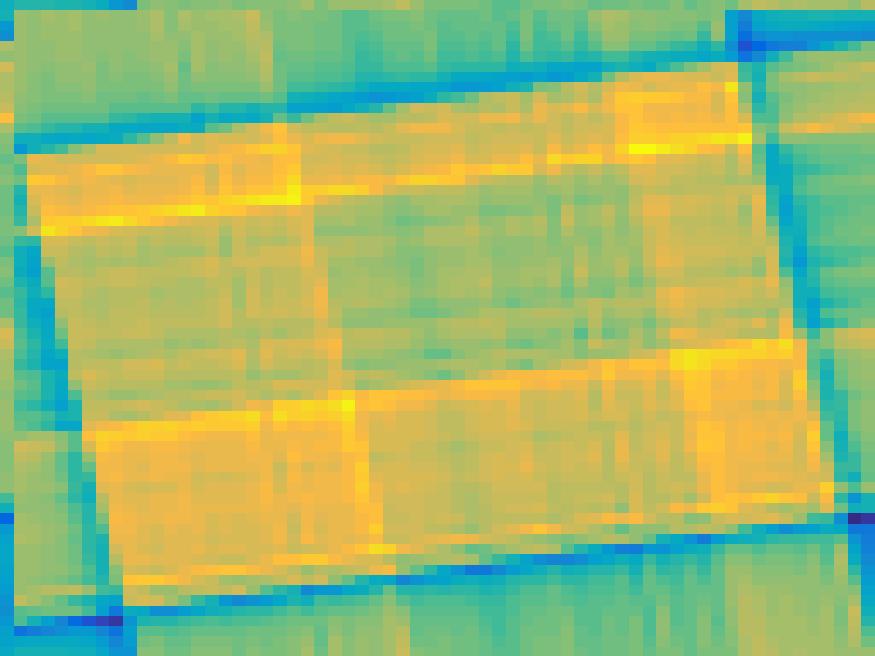}\\
\includegraphics[width=.1\linewidth,height=.1\linewidth]{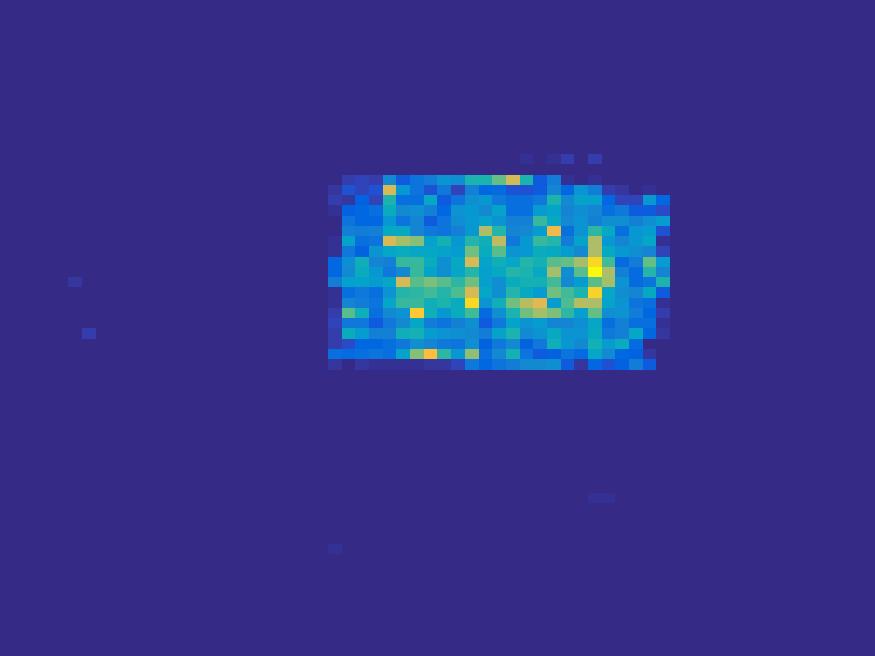}&
\includegraphics[width=.1\linewidth,height=.1\linewidth]{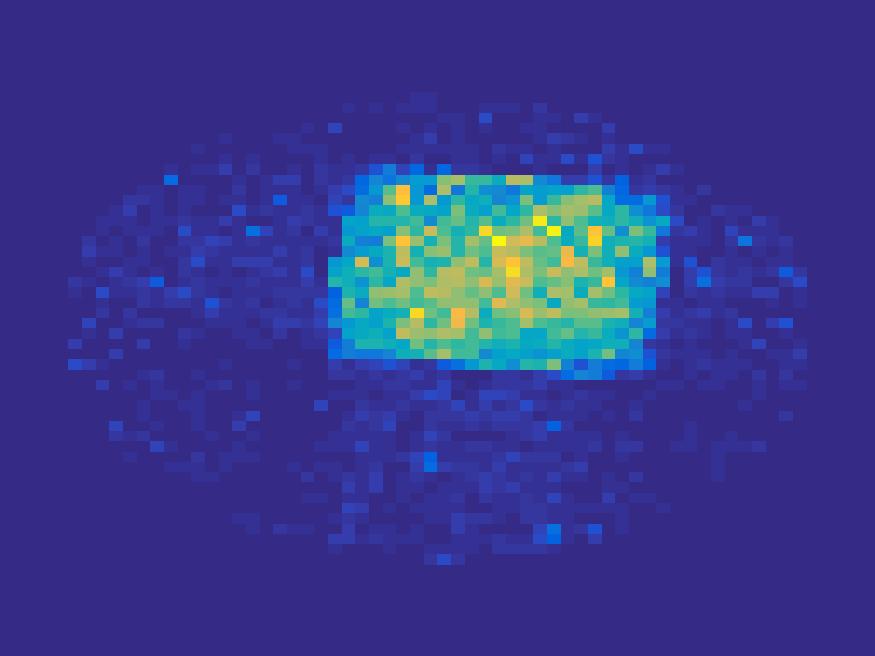}&
\includegraphics[width=.1\linewidth,height=.1\linewidth]{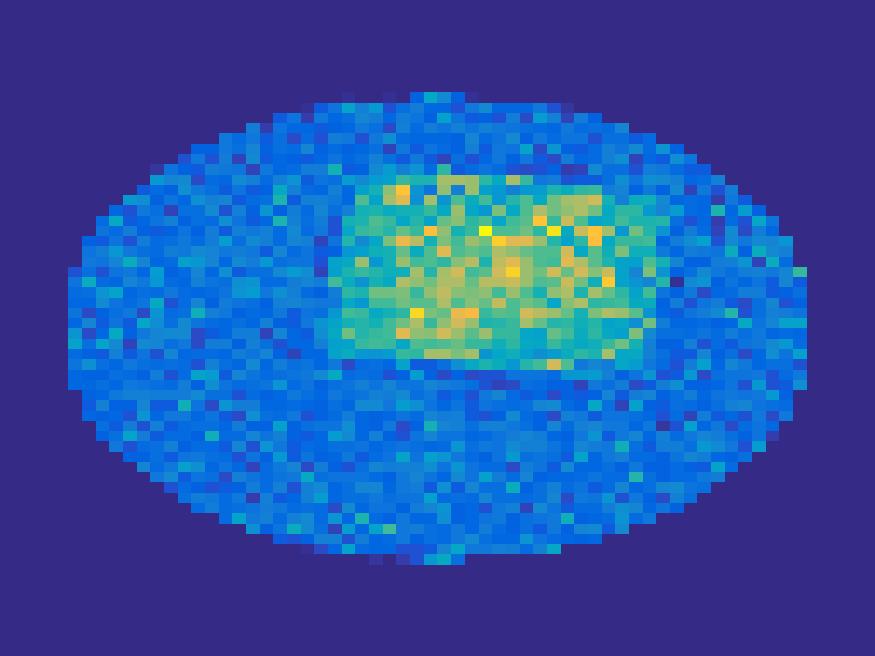}&
\includegraphics[width=.1\linewidth,height=.1\linewidth]{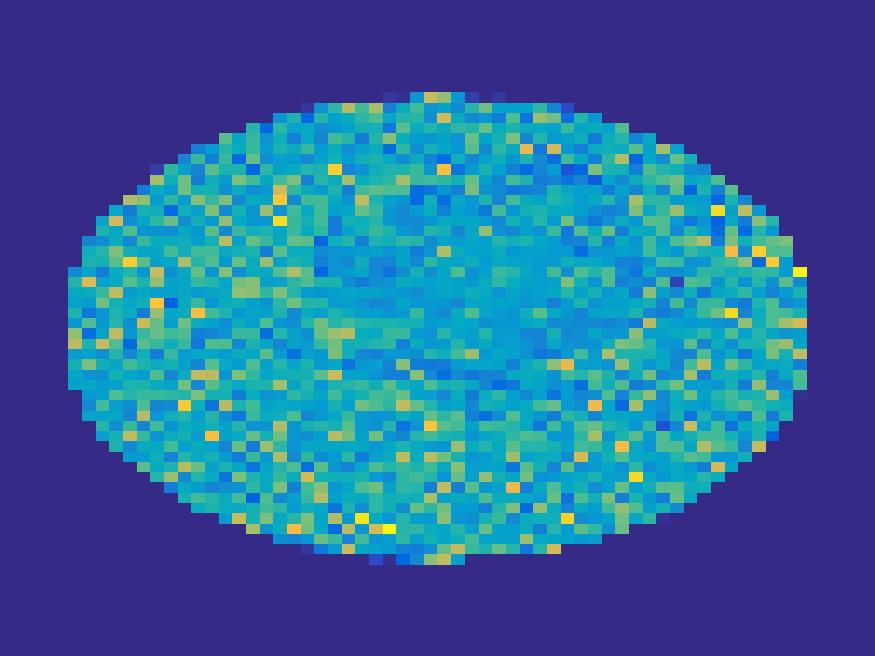}&
\includegraphics[width=.1\linewidth,height=.1\linewidth]{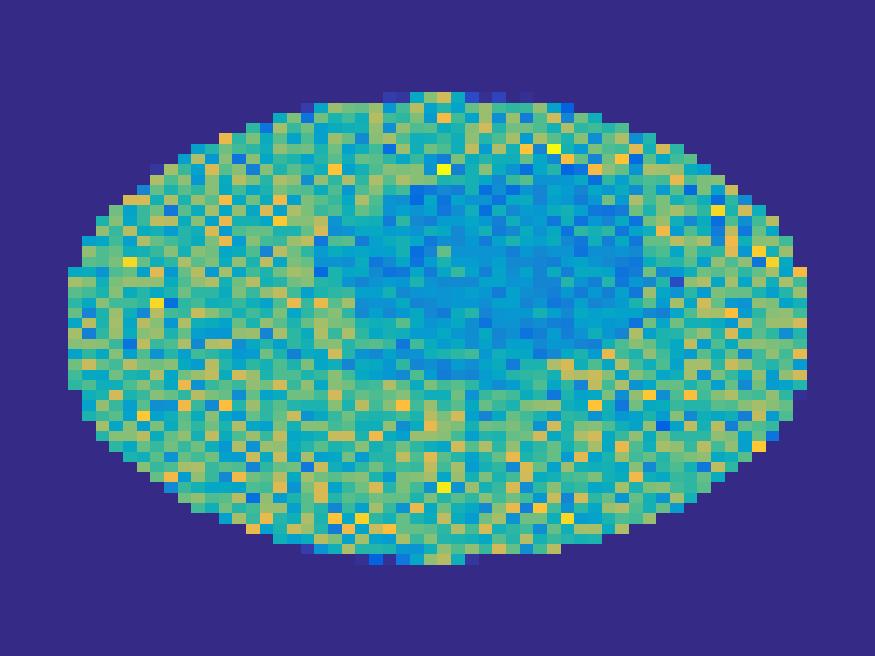}&
\includegraphics[width=.1\linewidth,height=.1\linewidth]{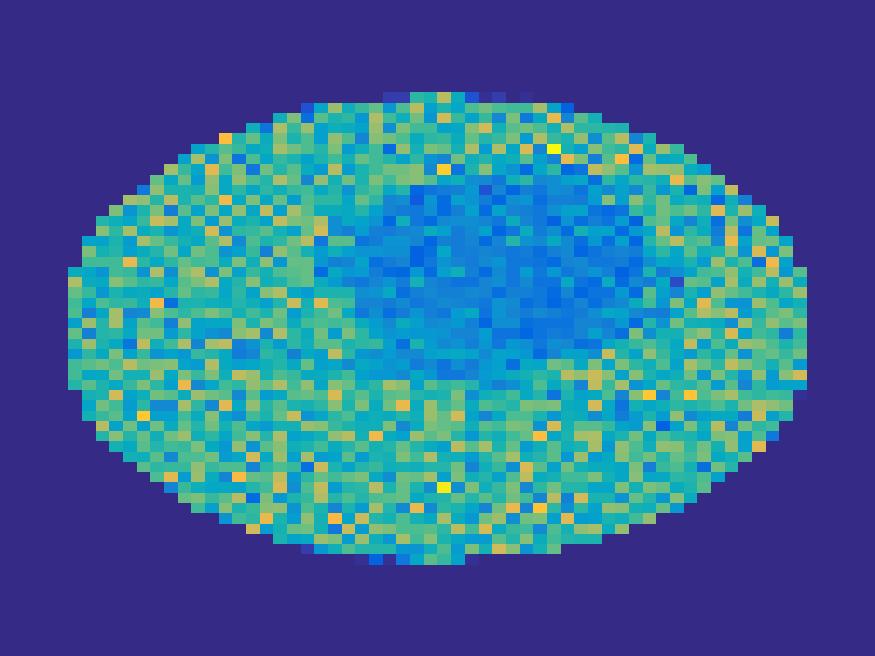}&
\includegraphics[width=.1\linewidth,height=.1\linewidth]{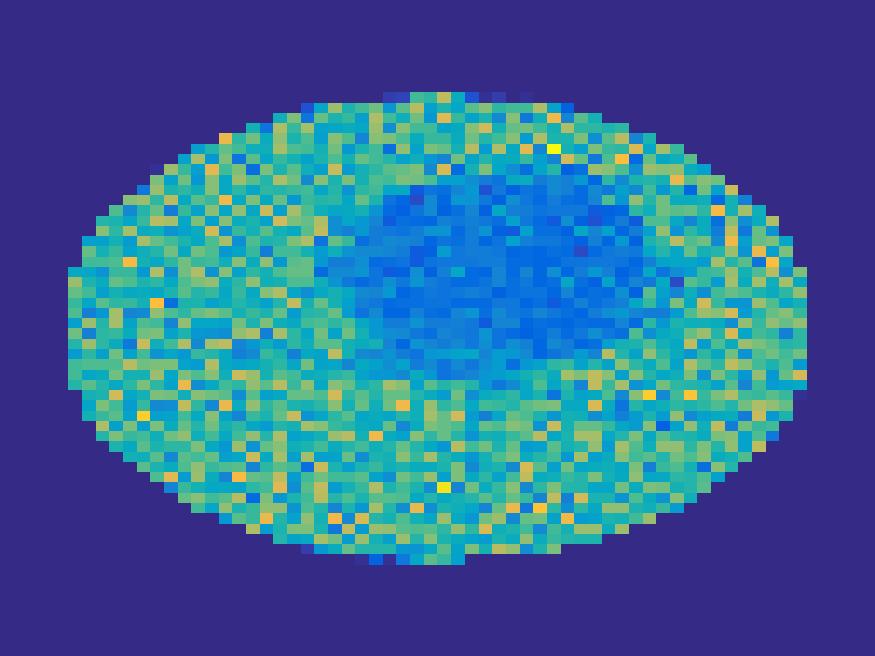}&
\includegraphics[width=.1\linewidth,height=.1\linewidth]{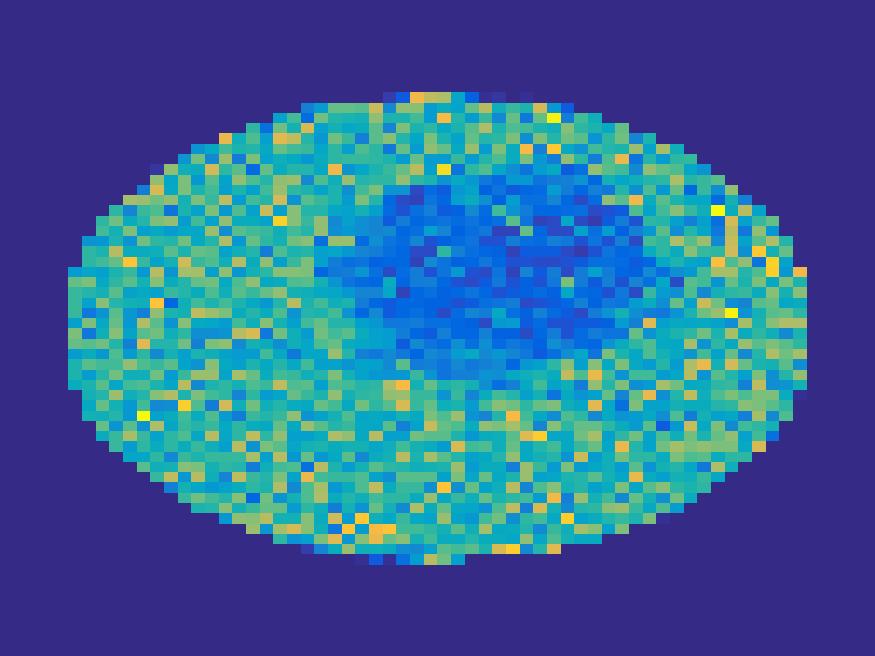}&
\includegraphics[width=.1\linewidth,height=.1\linewidth]{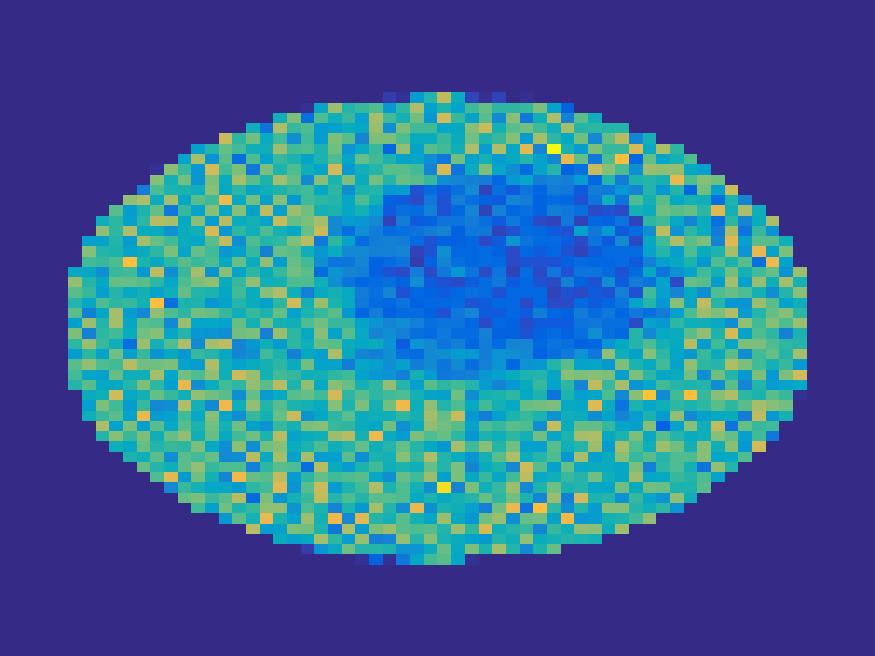}\\
\includegraphics[width=.1\linewidth,height=.1\linewidth]{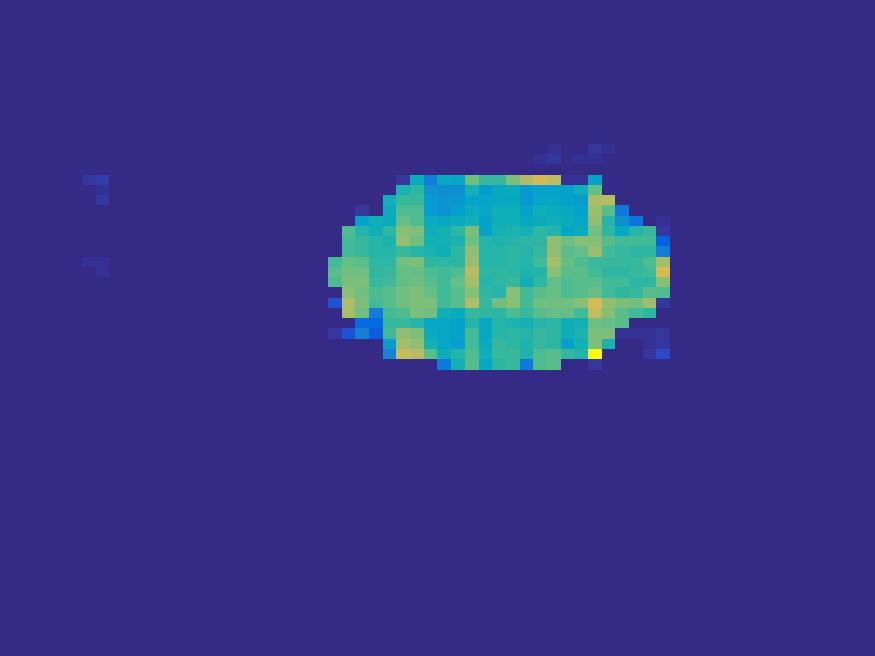}&
\includegraphics[width=.1\linewidth,height=.1\linewidth]{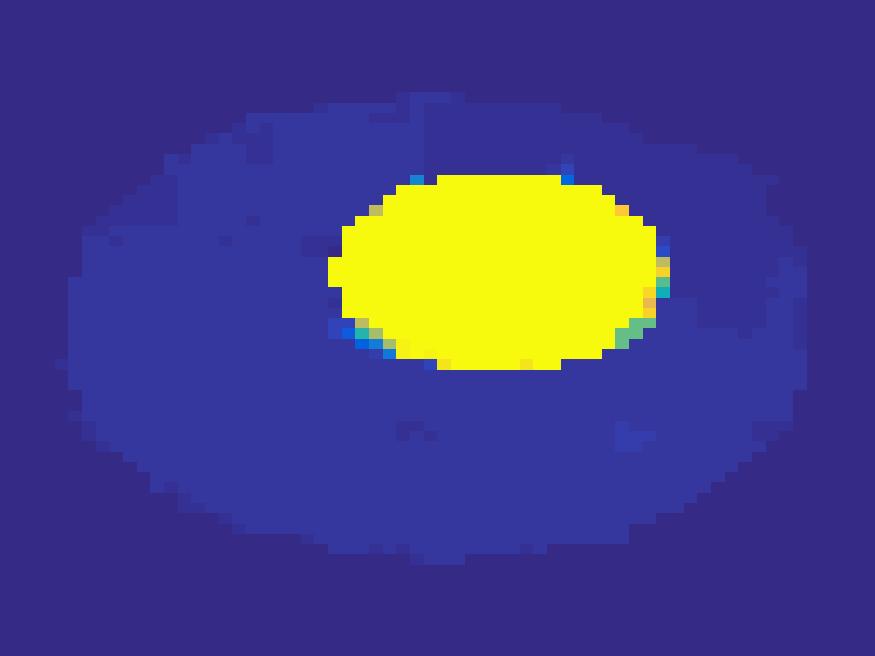}&
\includegraphics[width=.1\linewidth,height=.1\linewidth]{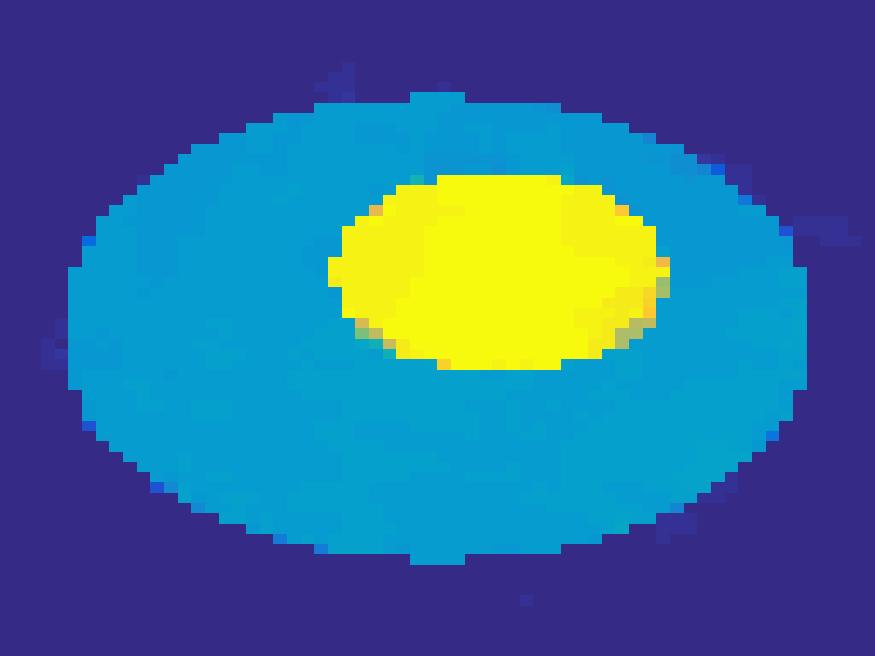}&
\includegraphics[width=.1\linewidth,height=.1\linewidth]{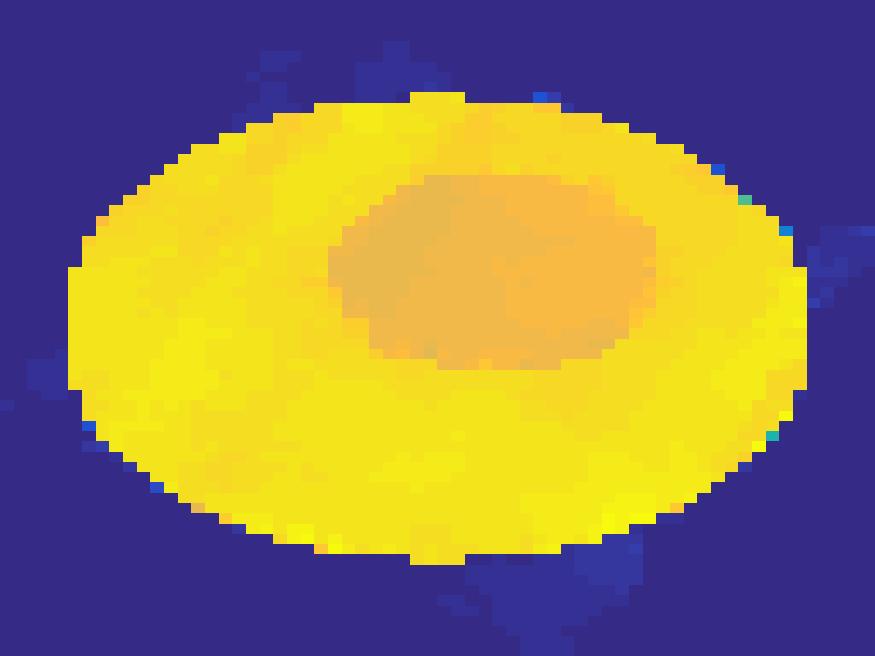}&
\includegraphics[width=.1\linewidth,height=.1\linewidth]{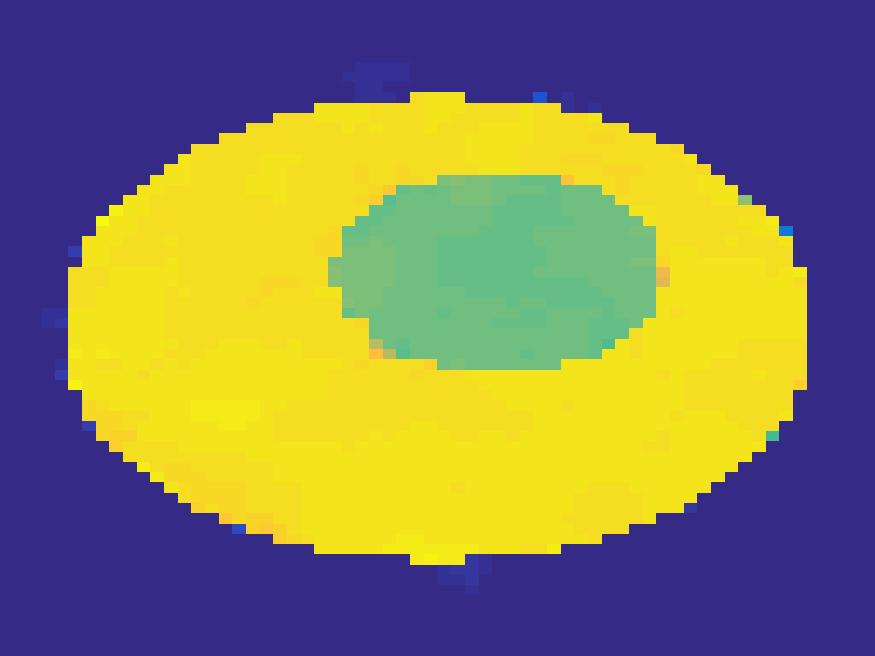}&
\includegraphics[width=.1\linewidth,height=.1\linewidth]{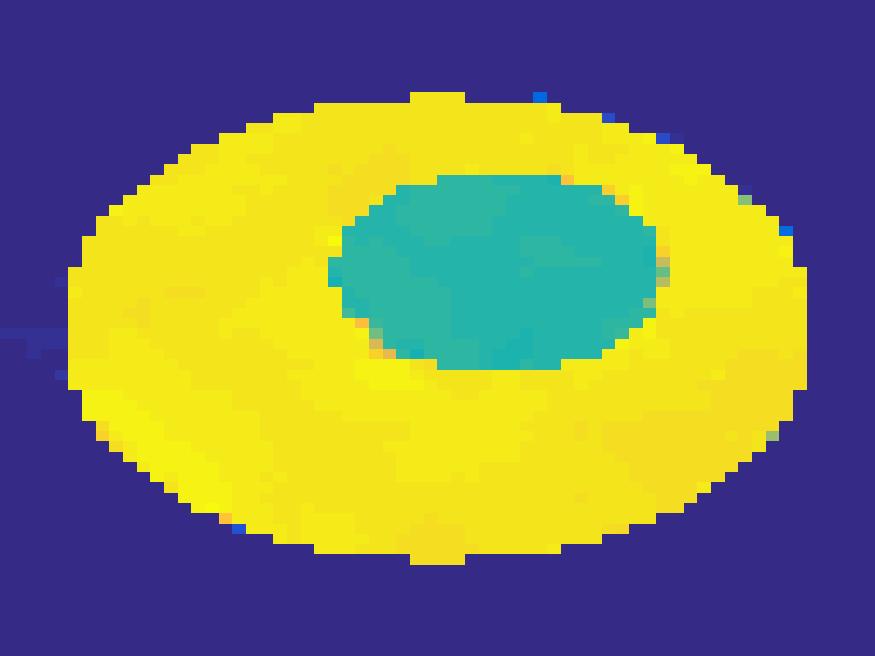}&
\includegraphics[width=.1\linewidth,height=.1\linewidth]{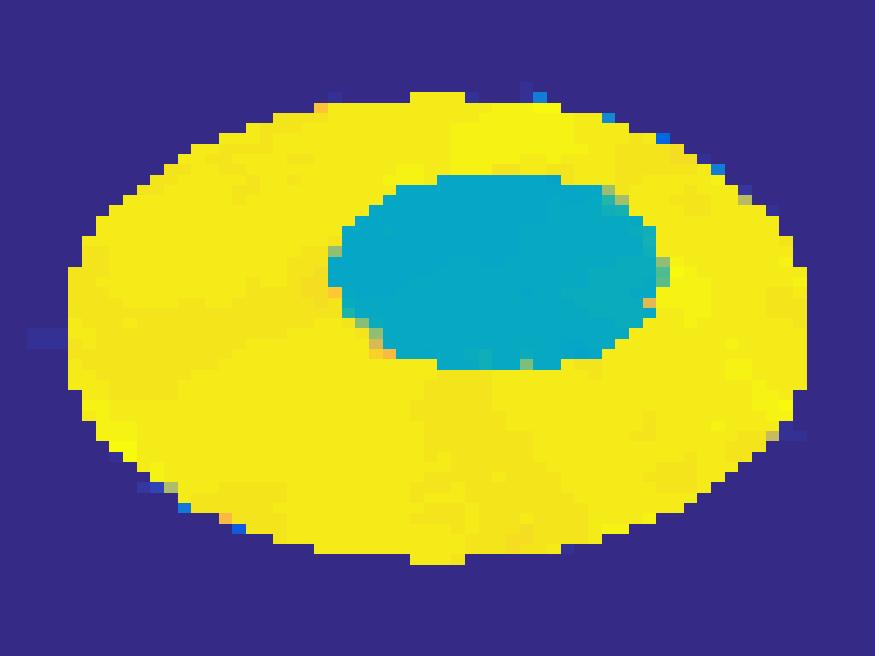}&
\includegraphics[width=.1\linewidth,height=.1\linewidth]{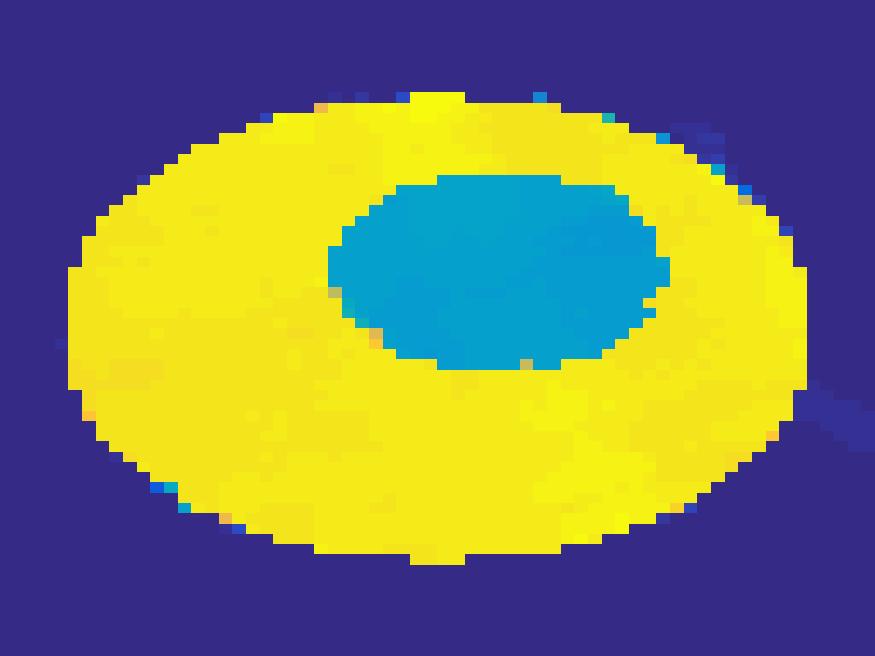}&
\includegraphics[width=.1\linewidth,height=.1\linewidth]{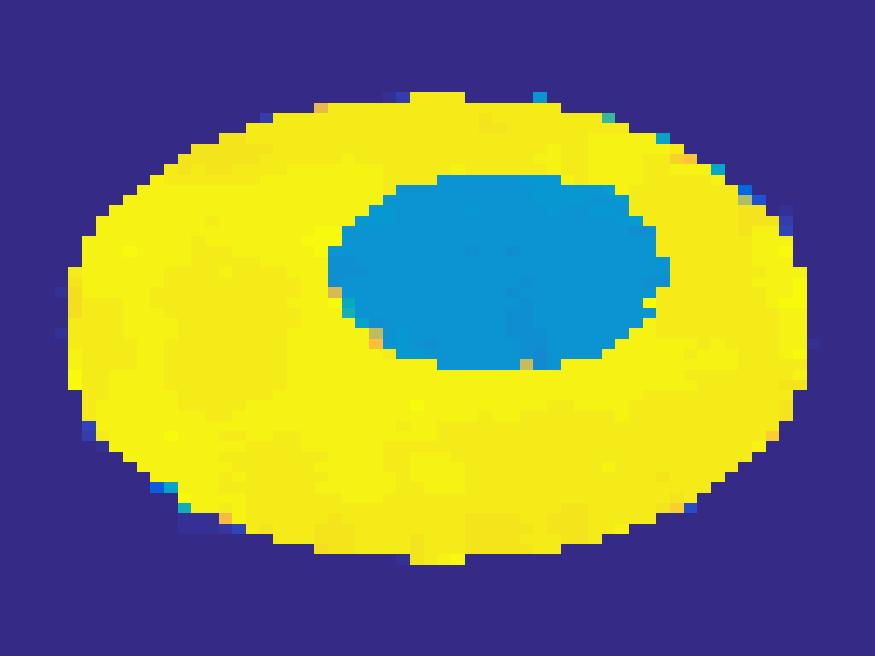}\\
{\footnotesize Frame 1}&
{\footnotesize Frame 11}&
{\footnotesize Frame 21}&
{\footnotesize Frame 31}&
{\footnotesize Frame 41}&
{\footnotesize Frame 51}&
{\footnotesize Frame 61}&
{\footnotesize Frame 71}&
{\footnotesize Frame 81}\\
\end{tabular}
\caption {First row: Ground truth; Second row: FBP; Third  row: EM algorithm with updating $\alpha$ and $B$; Forth row: Proposed method.}
\label{fig:ICPoiEll}
\end{figure}

\begin{figure}
\begin{tabular}{c@{\hspace{2pt}}c@{\hspace{2pt}}c@{\hspace{2pt}}c@{\hspace{2pt}}c@{\hspace{2pt}}c@{\hspace{2pt}}c@{\hspace{2pt}}c@{\hspace{2pt}}c@{\hspace{2pt}}c}
\includegraphics[width=.1\linewidth,height=.1\linewidth]{Liver1}&
\includegraphics[width=.1\linewidth,height=.1\linewidth]{Liver11}&
\includegraphics[width=.1\linewidth,height=.1\linewidth]{Liver21}&
\includegraphics[width=.1\linewidth,height=.1\linewidth]{Liver31}&
\includegraphics[width=.1\linewidth,height=.1\linewidth]{Liver41}&
\includegraphics[width=.1\linewidth,height=.1\linewidth]{Liver51}&
\includegraphics[width=.1\linewidth,height=.1\linewidth]{Liver61}&
\includegraphics[width=.1\linewidth,height=.1\linewidth]{Liver71}&
\includegraphics[width=.1\linewidth,height=.1\linewidth]{Liver81}\\
\includegraphics[width=.1\linewidth,height=.1\linewidth]{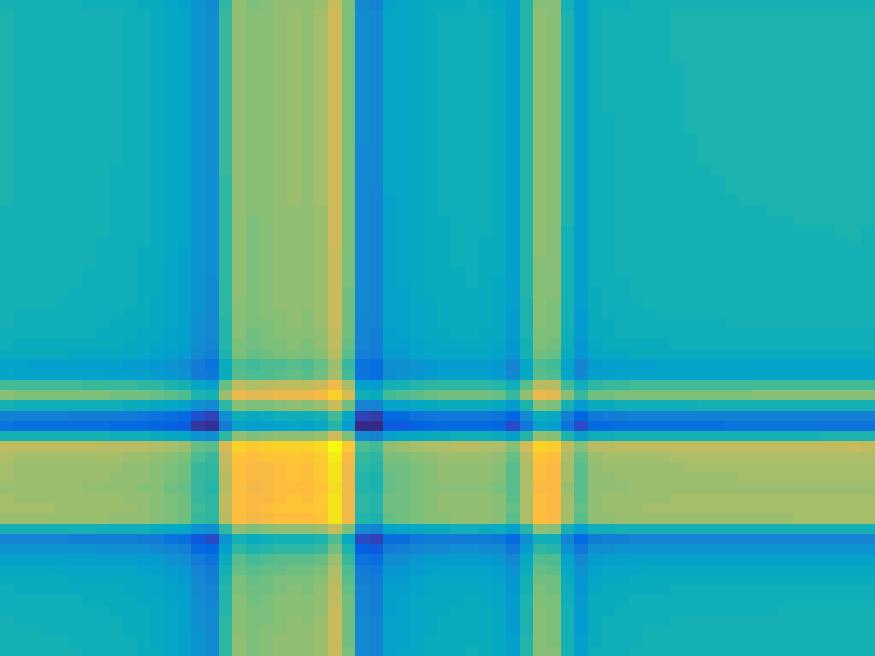}&
\includegraphics[width=.1\linewidth,height=.1\linewidth]{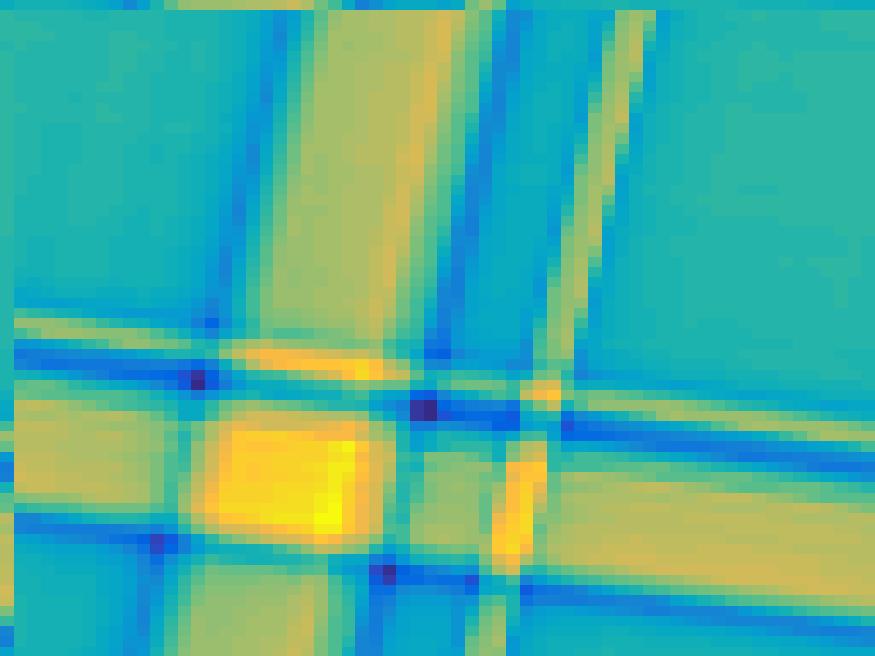}&
\includegraphics[width=.1\linewidth,height=.1\linewidth]{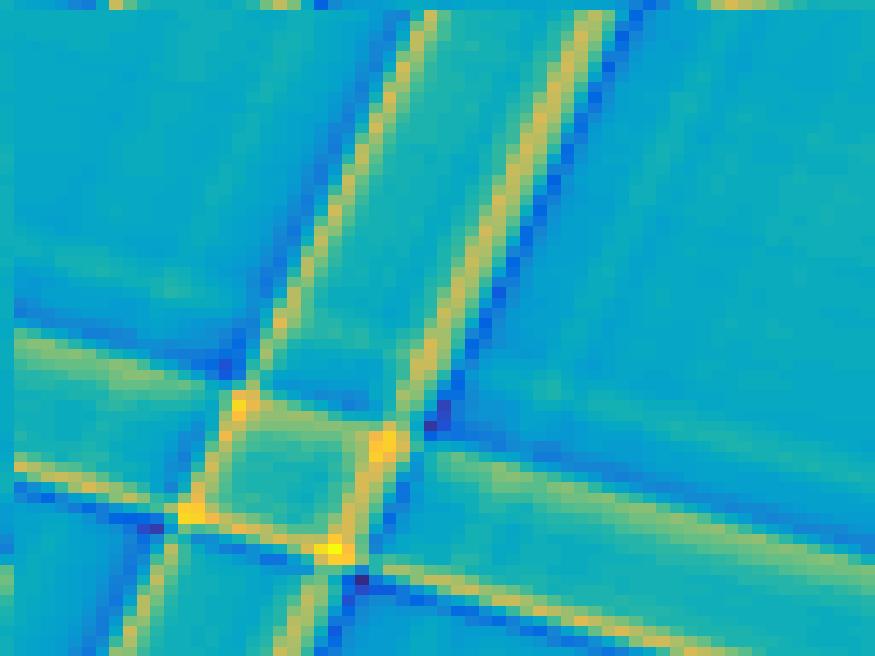}&
\includegraphics[width=.1\linewidth,height=.1\linewidth]{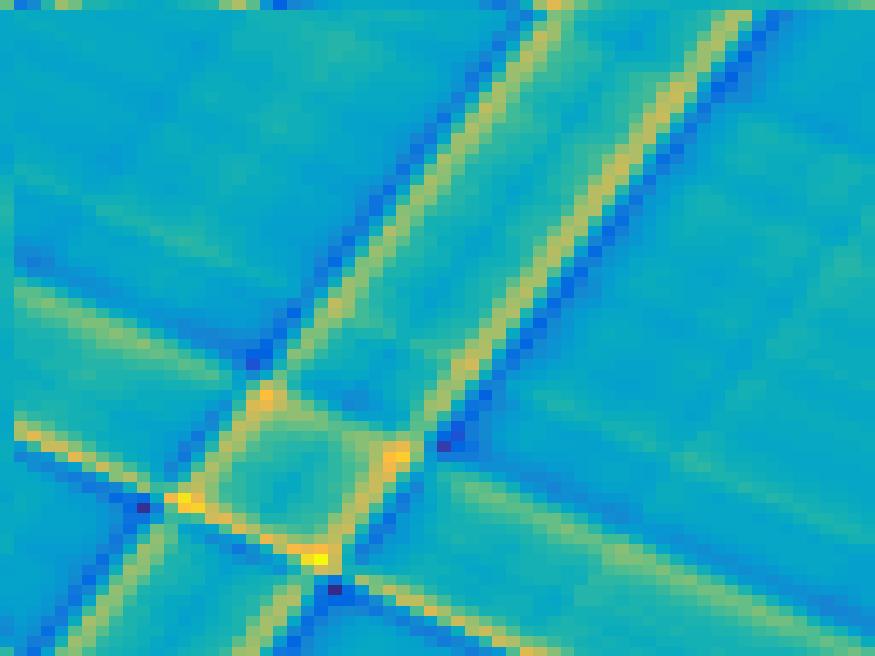}&
\includegraphics[width=.1\linewidth,height=.1\linewidth]{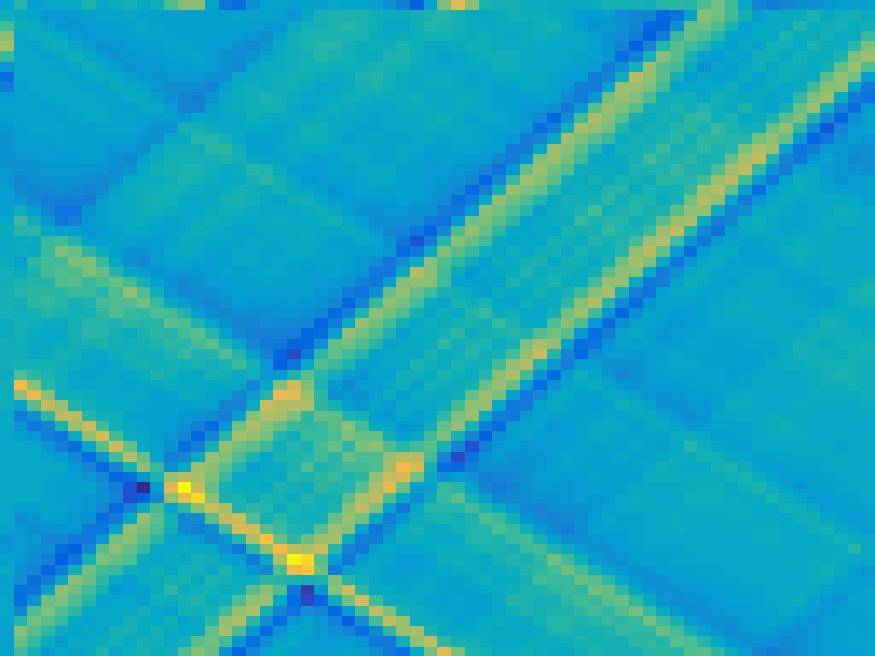}&
\includegraphics[width=.1\linewidth,height=.1\linewidth]{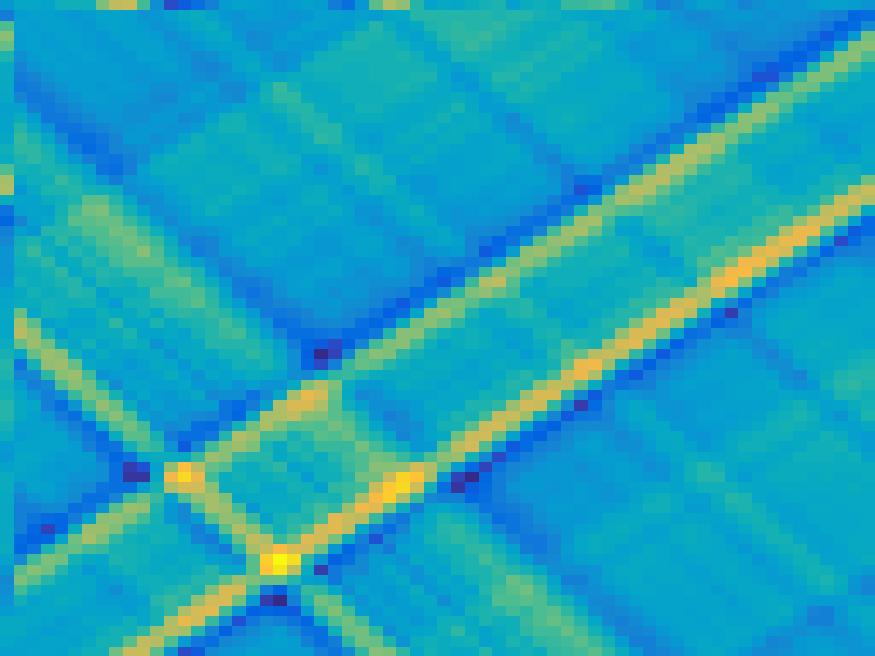}&
\includegraphics[width=.1\linewidth,height=.1\linewidth]{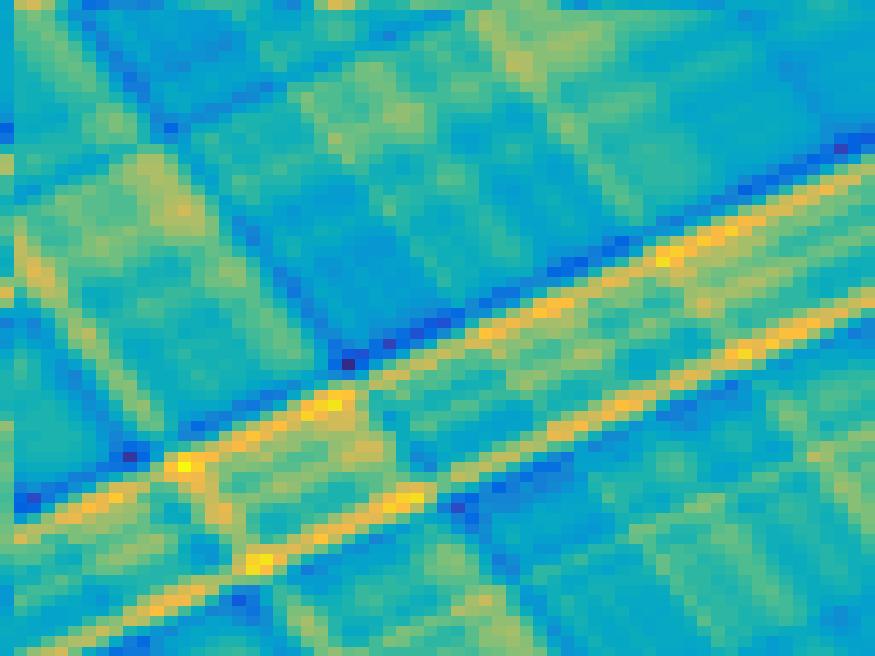}&
\includegraphics[width=.1\linewidth,height=.1\linewidth]{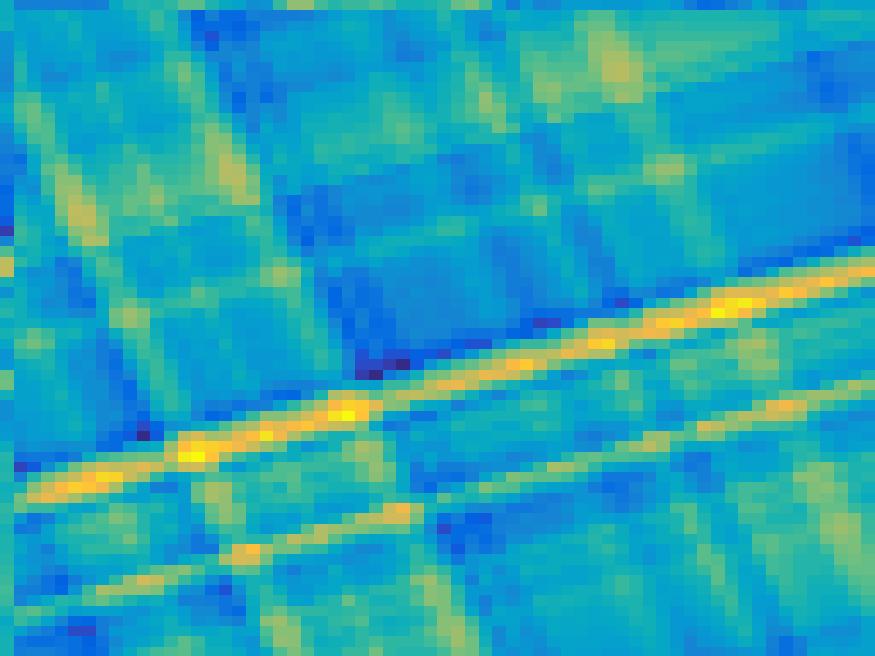}&
\includegraphics[width=.1\linewidth,height=.1\linewidth]{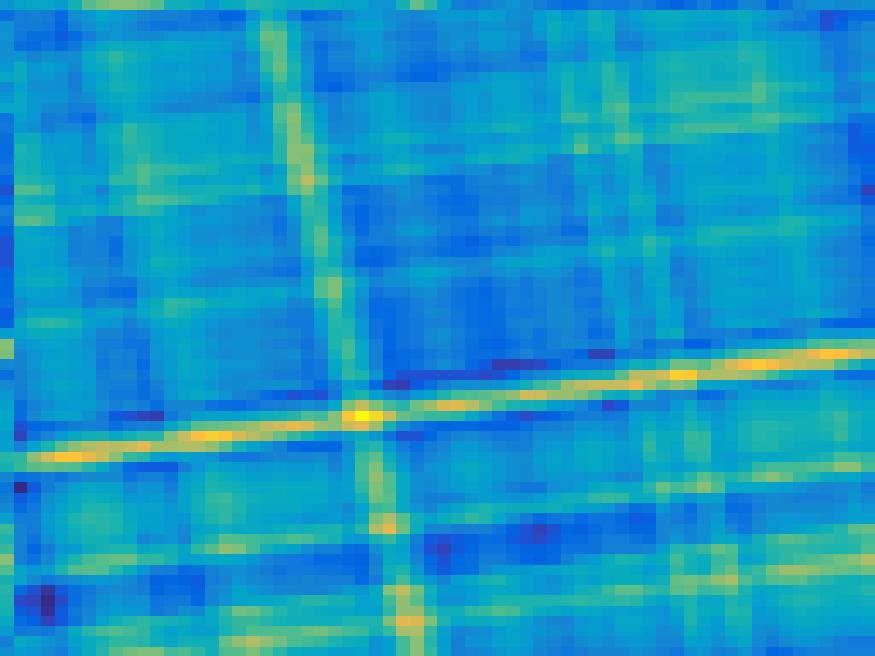}\\
\includegraphics[width=.1\linewidth,height=.1\linewidth]{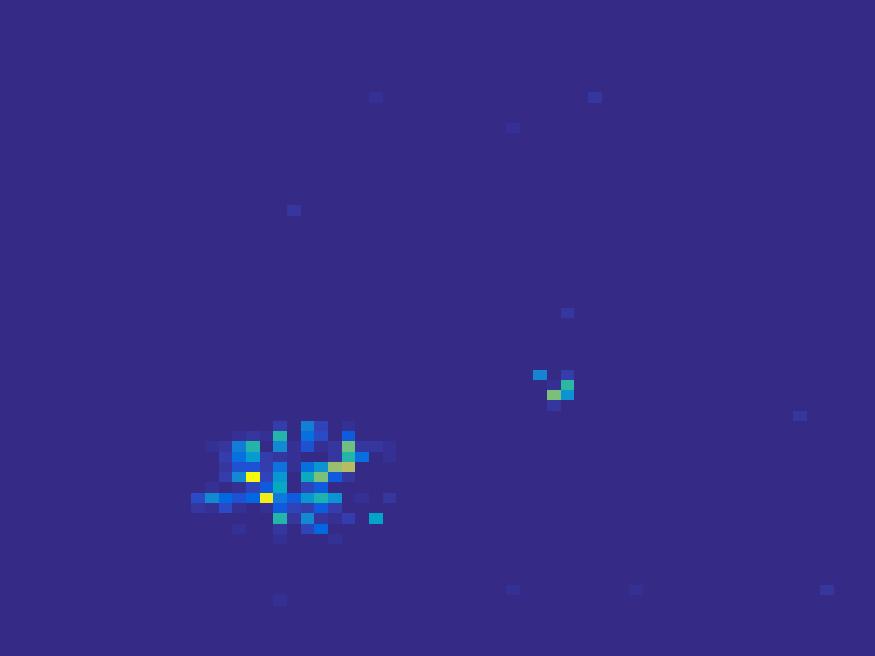}&
\includegraphics[width=.1\linewidth,height=.1\linewidth]{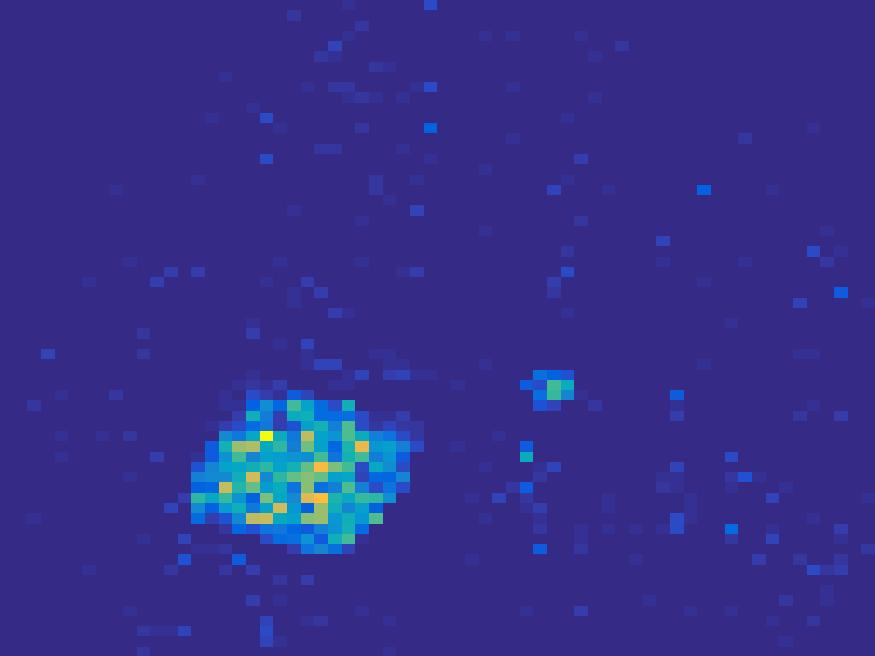}&
\includegraphics[width=.1\linewidth,height=.1\linewidth]{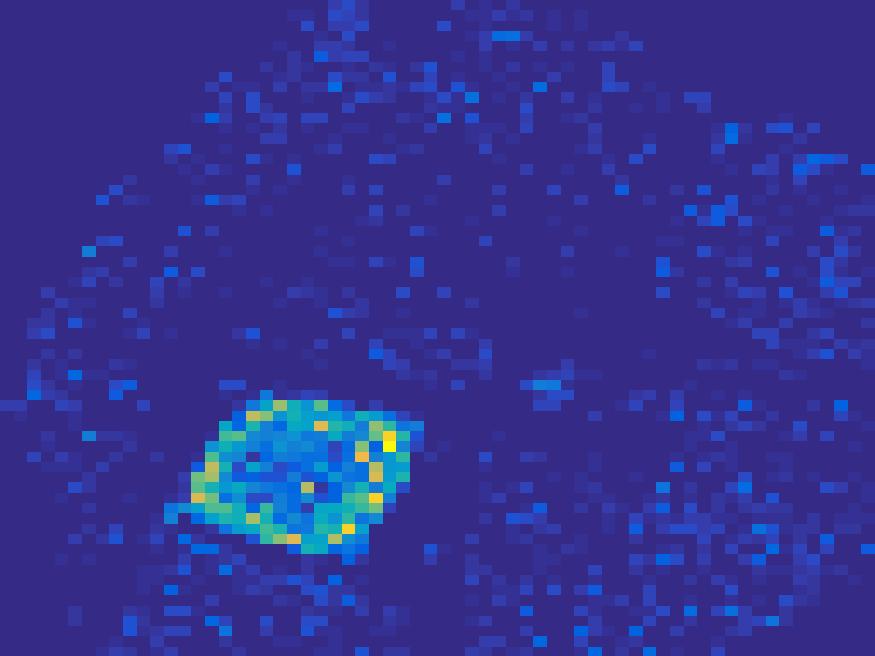}&
\includegraphics[width=.1\linewidth,height=.1\linewidth]{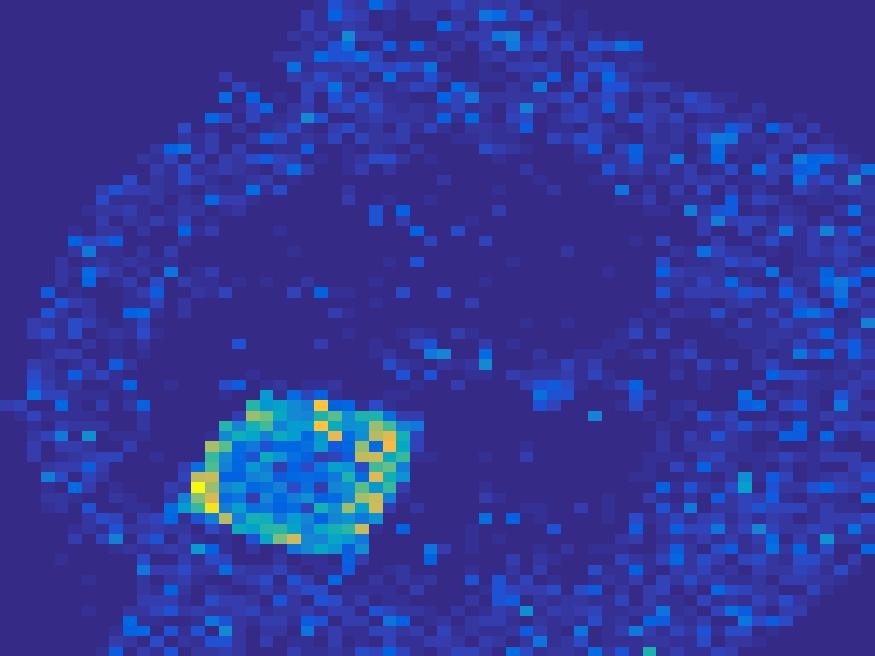}&
\includegraphics[width=.1\linewidth,height=.1\linewidth]{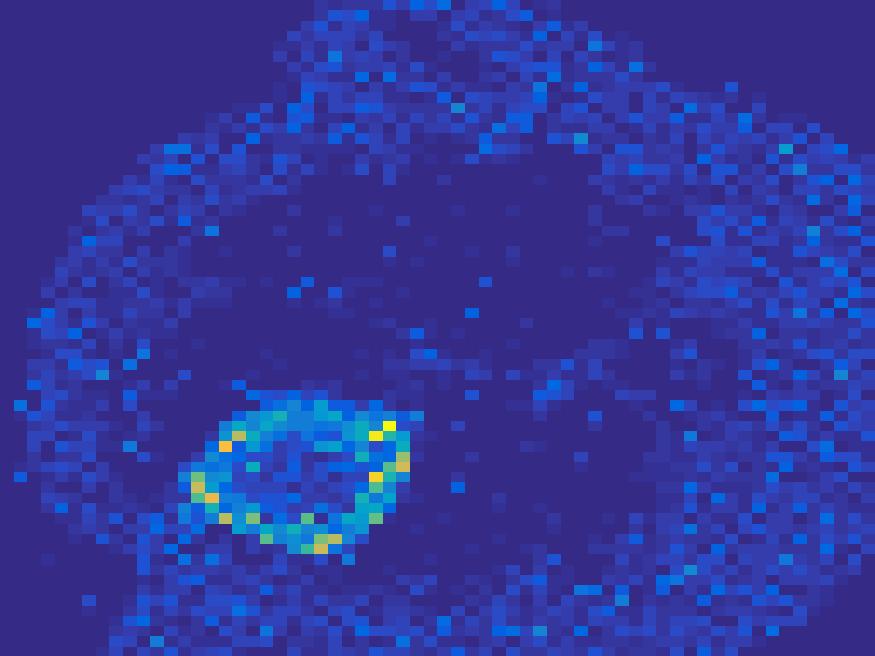}&
\includegraphics[width=.1\linewidth,height=.1\linewidth]{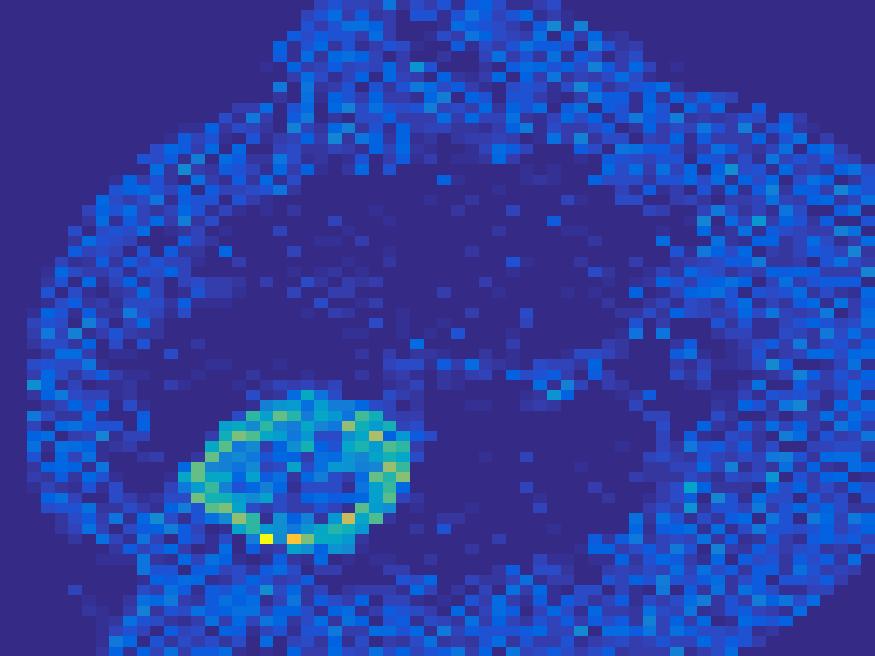}&
\includegraphics[width=.1\linewidth,height=.1\linewidth]{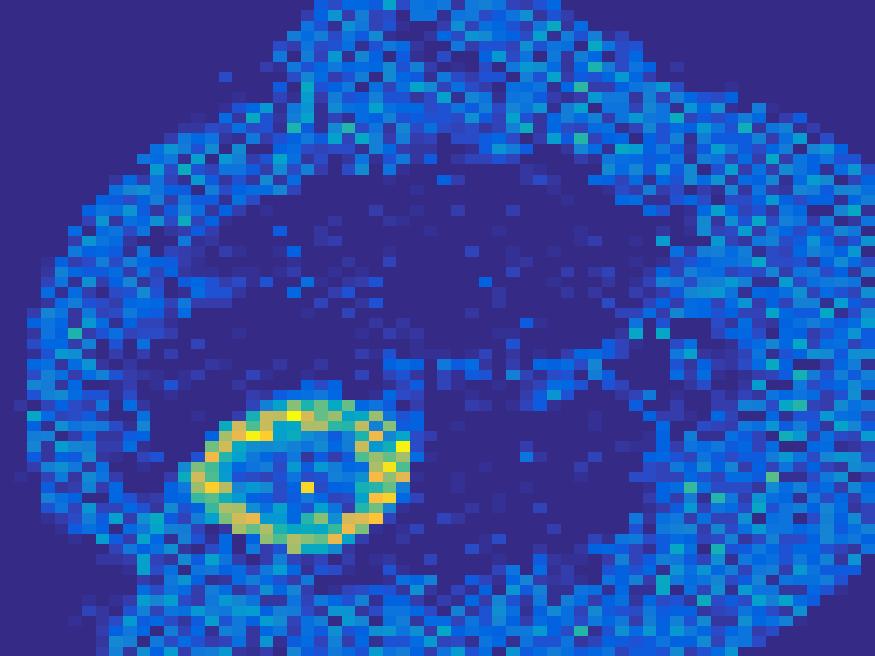}&
\includegraphics[width=.1\linewidth,height=.1\linewidth]{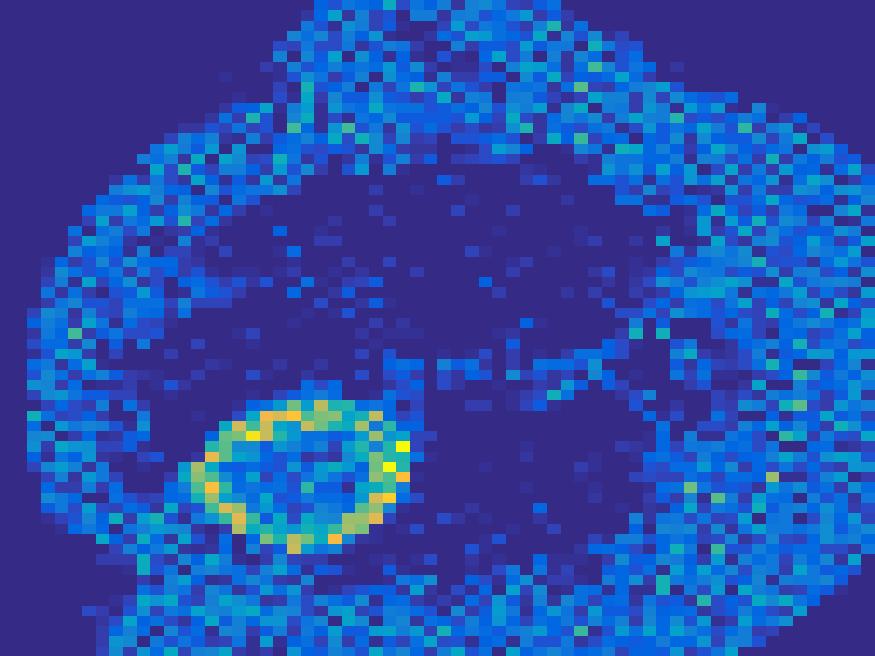}&
\includegraphics[width=.1\linewidth,height=.1\linewidth]{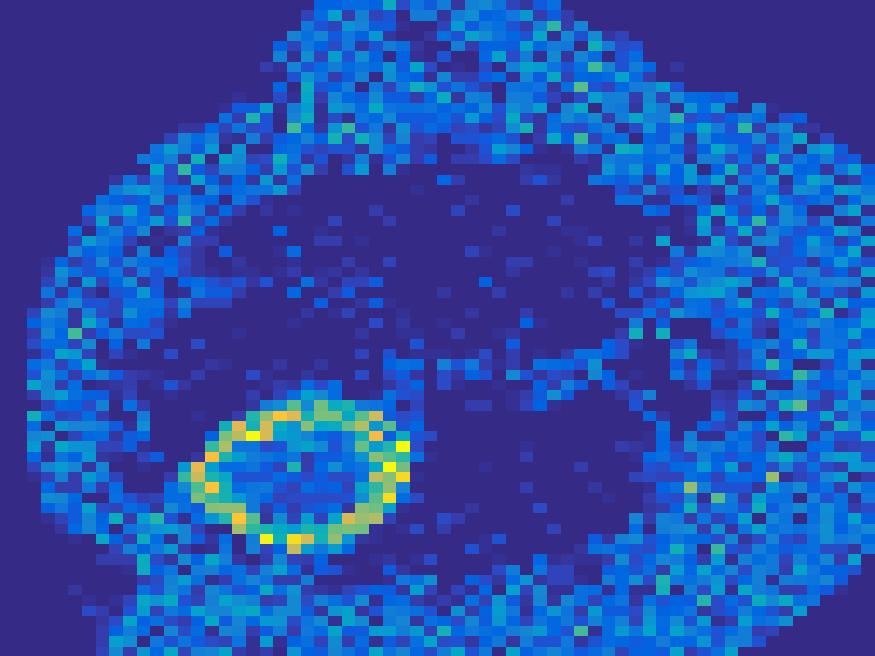}\\
\includegraphics[width=.1\linewidth,height=.1\linewidth]{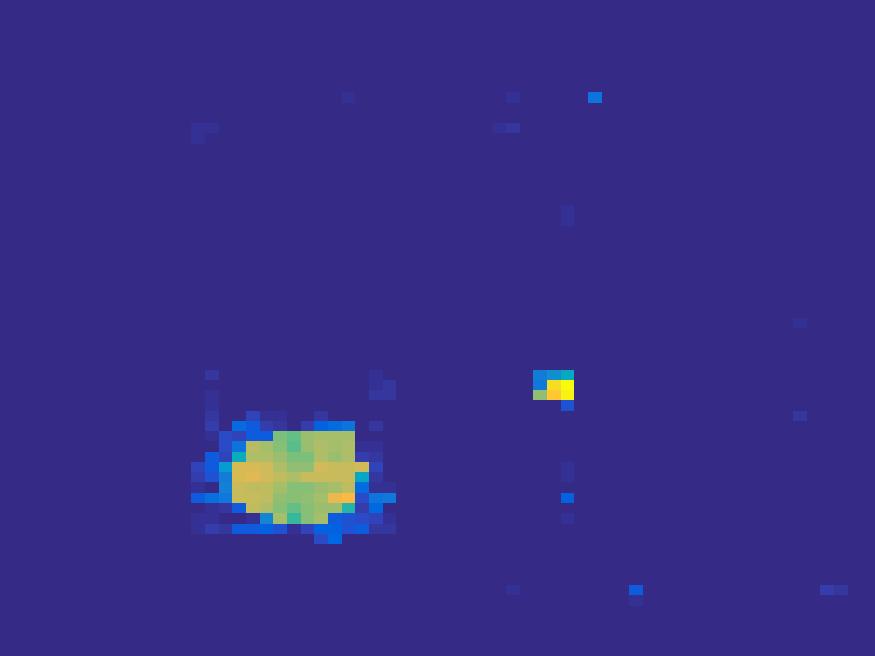}&
\includegraphics[width=.1\linewidth,height=.1\linewidth]{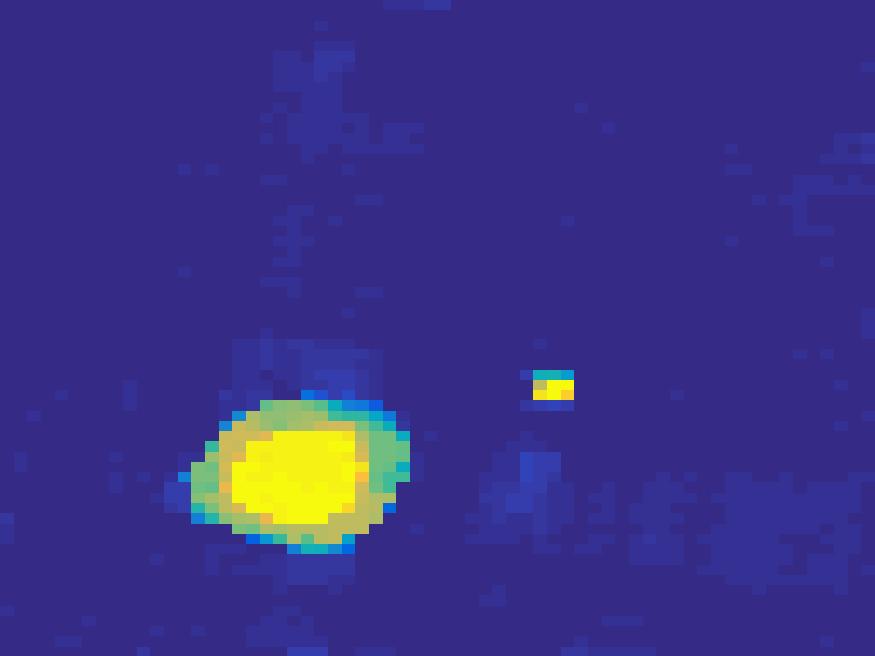}&
\includegraphics[width=.1\linewidth,height=.1\linewidth]{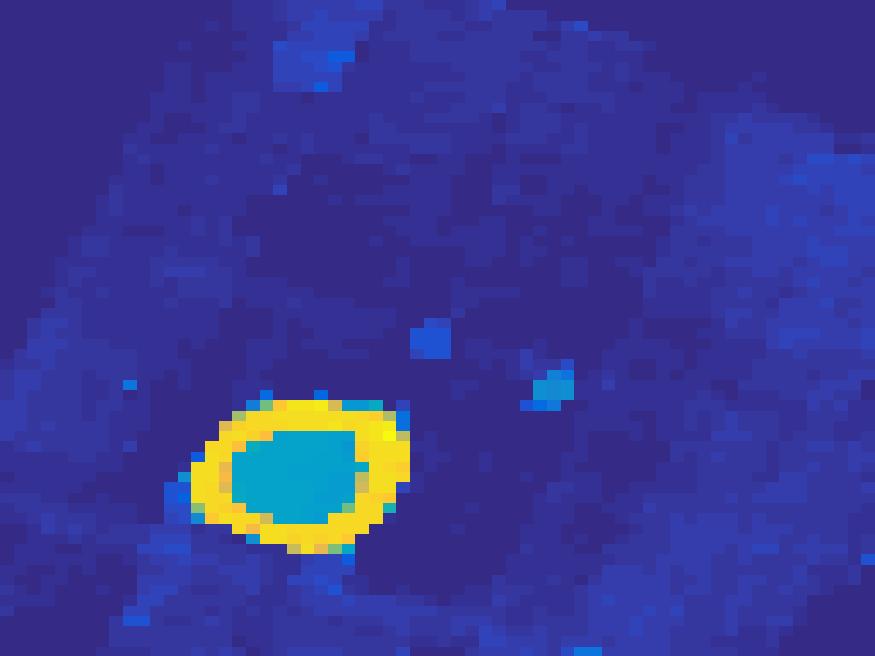}&
\includegraphics[width=.1\linewidth,height=.1\linewidth]{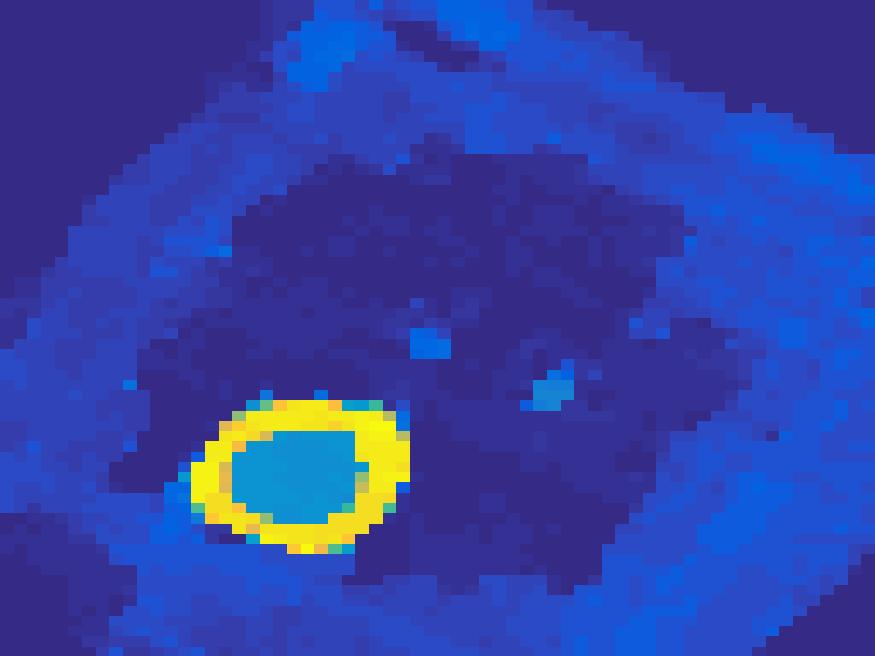}&
\includegraphics[width=.1\linewidth,height=.1\linewidth]{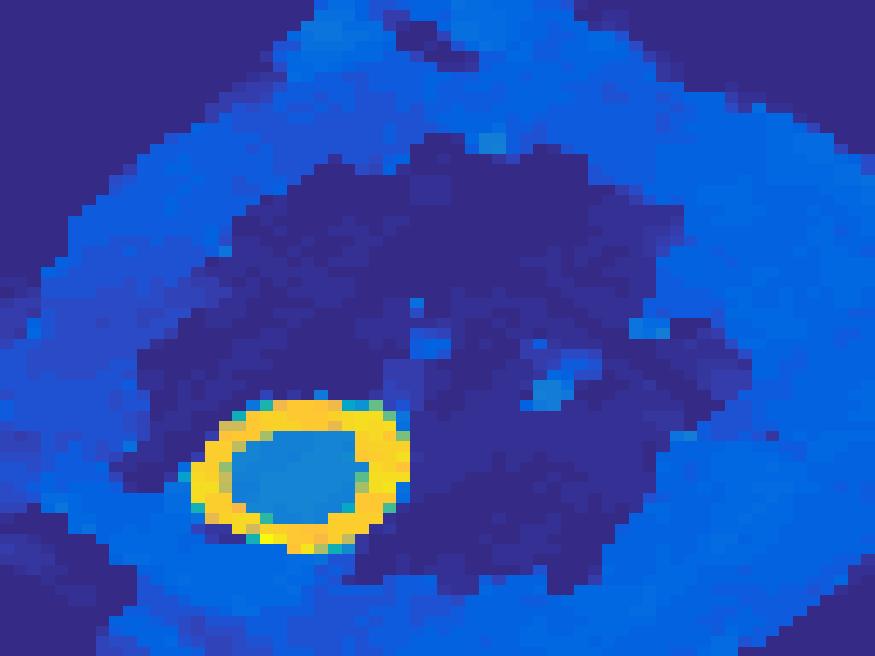}&
\includegraphics[width=.1\linewidth,height=.1\linewidth]{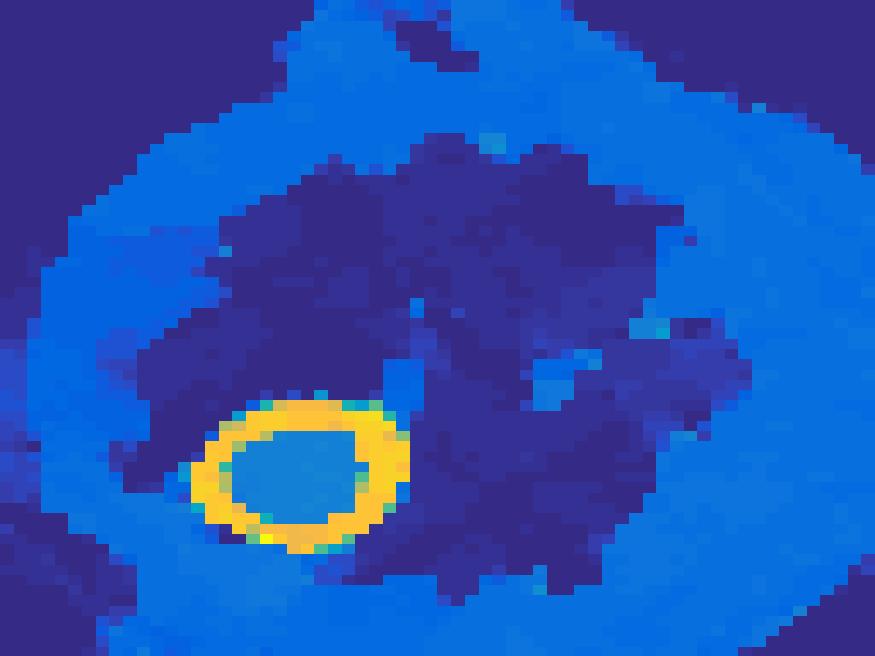}&
\includegraphics[width=.1\linewidth,height=.1\linewidth]{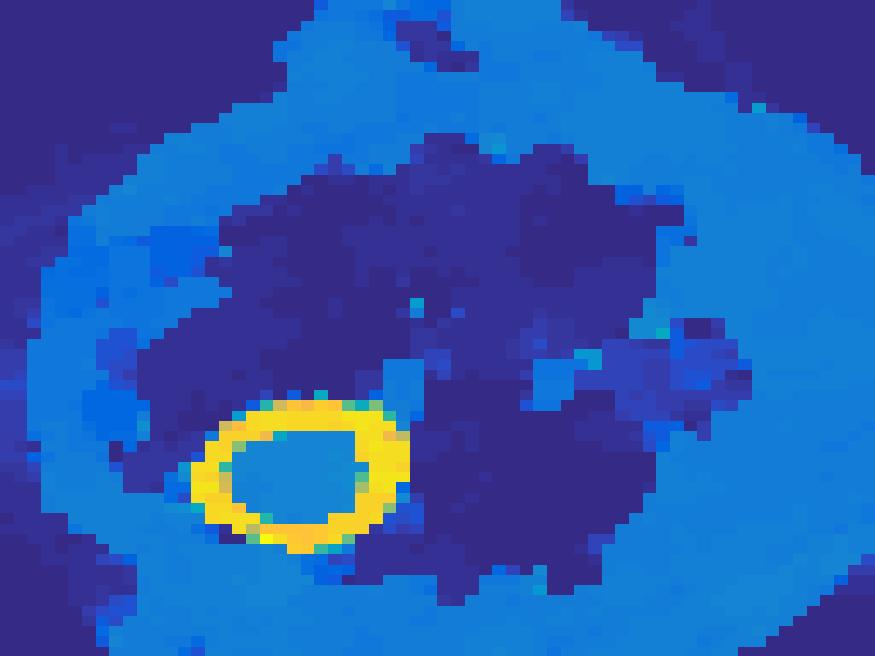}&
\includegraphics[width=.1\linewidth,height=.1\linewidth]{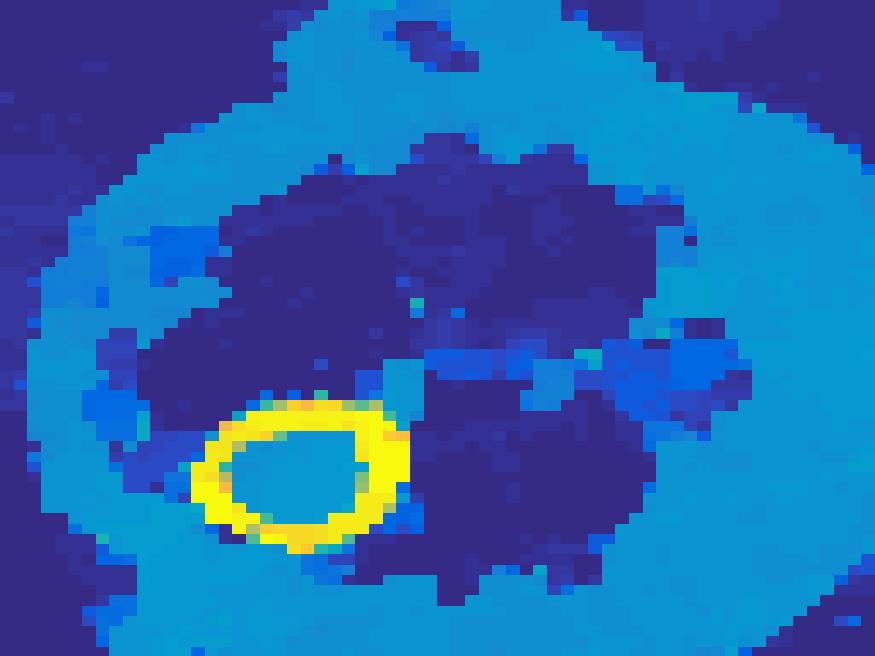}&
\includegraphics[width=.1\linewidth,height=.1\linewidth]{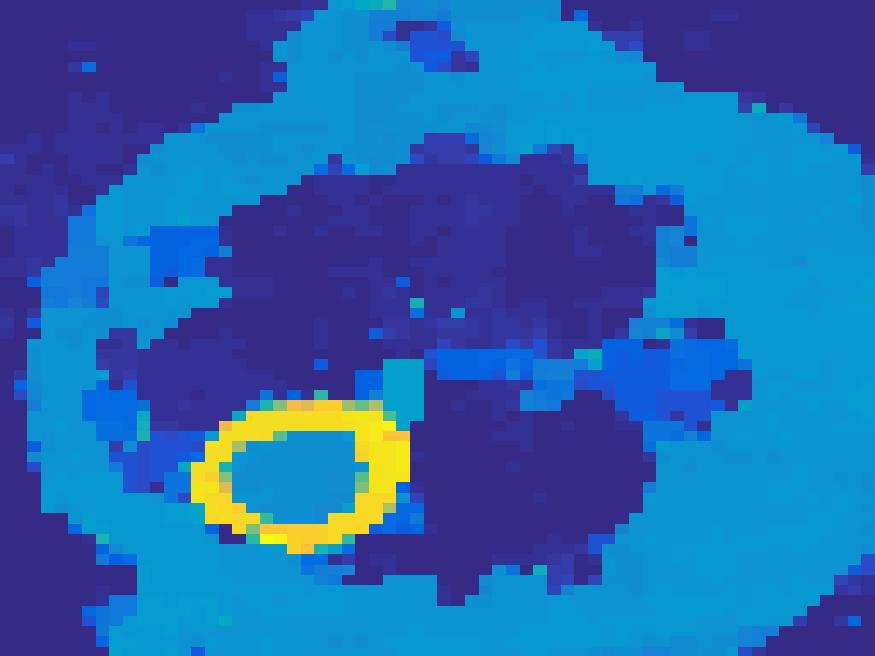}\\
{\footnotesize Frame 1}&
{\footnotesize Frame 11}&
{\footnotesize Frame 21}&
{\footnotesize Frame 31}&
{\footnotesize Frame 41}&
{\footnotesize Frame 51}&
{\footnotesize Frame 61}&
{\footnotesize Frame 71}&
{\footnotesize Frame 81}\\
\end{tabular}
\caption {First row: Ground truth; Second row: FBP; Third  row: EM algorithm with   updating $\alpha$ and $B$; Forth row: Proposed method.}
\label{fig:ICPoiLiver}
\end{figure}

Figure \ref{fig:ICPoiEll} and  \ref{fig:ICPoiLiver}  show the results of the ellipse and rat phantom with Poisson noise. Since the number of projections is very limited  and the corruption by Poisson noise, the reconstruction by both  FBP and EM (with updating basis) are not  satisfactory, while the proposed method is capable to reconstruct the main structure of the images faithfully.

\begin{figure}[ht]
\begin{center}
\subfigure{
\includegraphics[width=.45\linewidth,height=.3\linewidth]{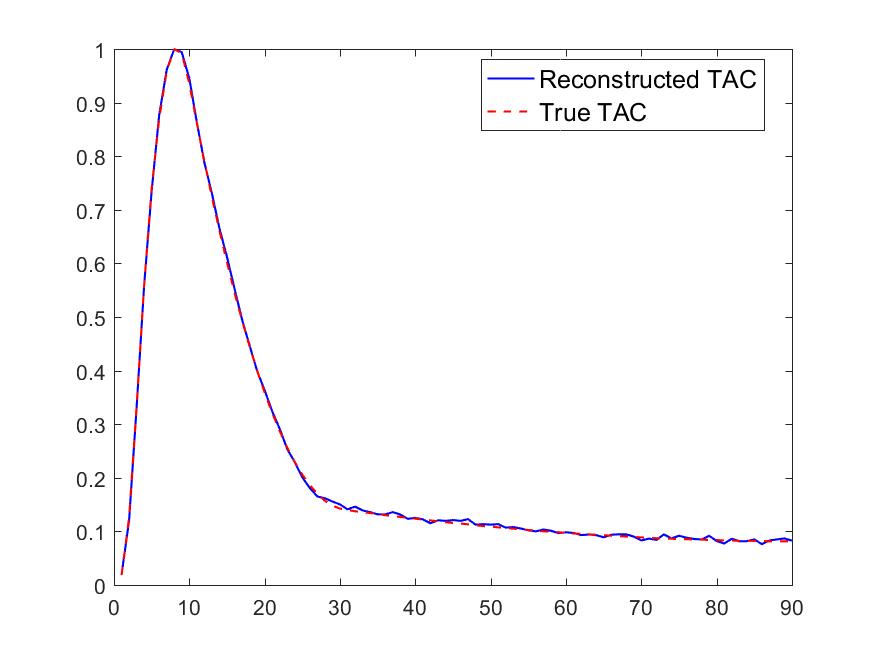}}
\subfigure{
\includegraphics[width=.45\linewidth,height=.3\linewidth]{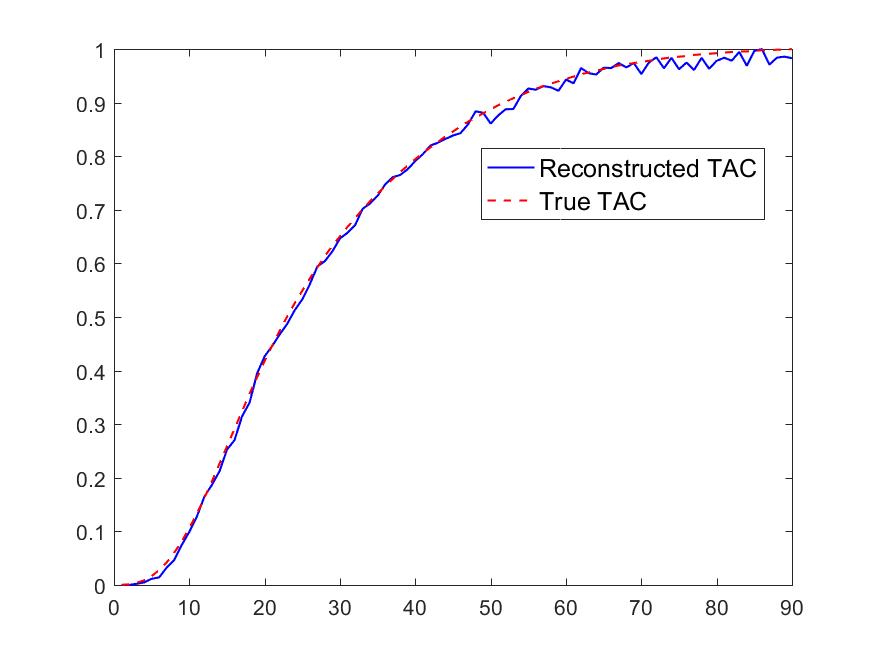}}
\end{center}
\caption{Comparison of the true TACs and the Reconstructed TACs in two regions (Ellipse phantom).}
\label{fig:ICPoiEllTAC}
\end{figure}

\begin{figure}[ht]
\begin{center}
\subfigure{
\includegraphics[width=.45\linewidth,height=.3\linewidth]{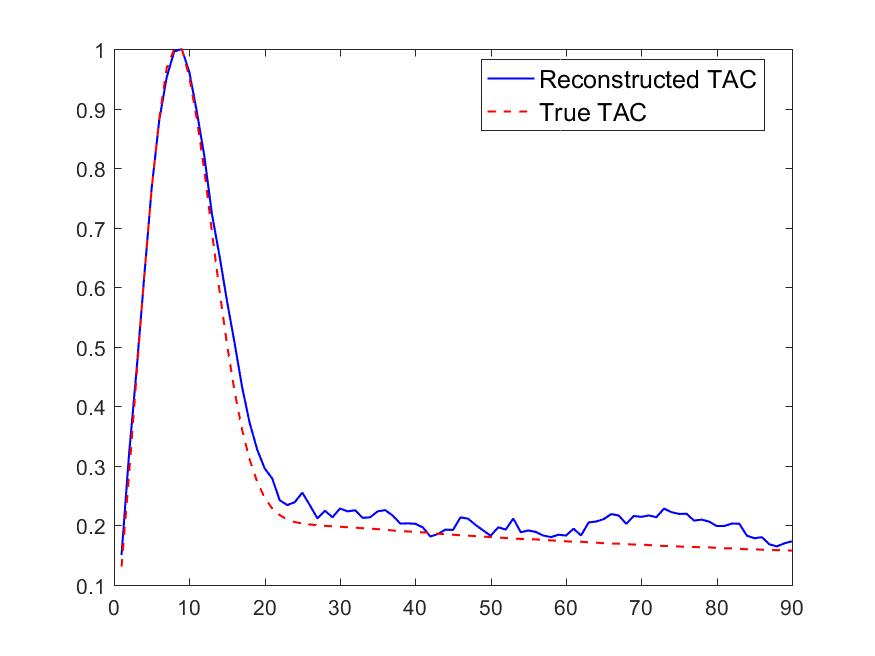}}
\subfigure{
\includegraphics[width=.45\linewidth,height=.3\linewidth]{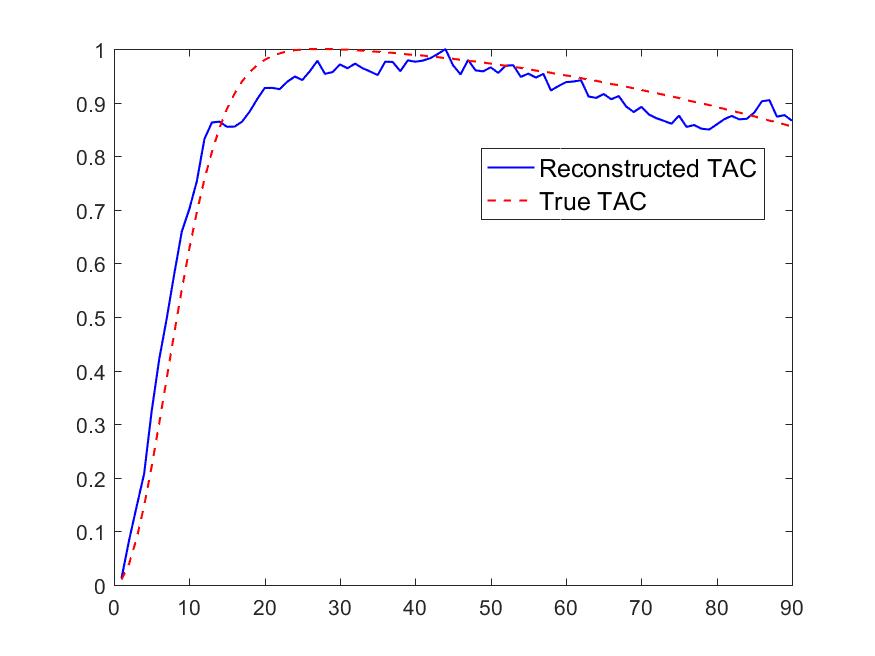}}\\
\subfigure{
\includegraphics[width=.45\linewidth,height=.3\linewidth]{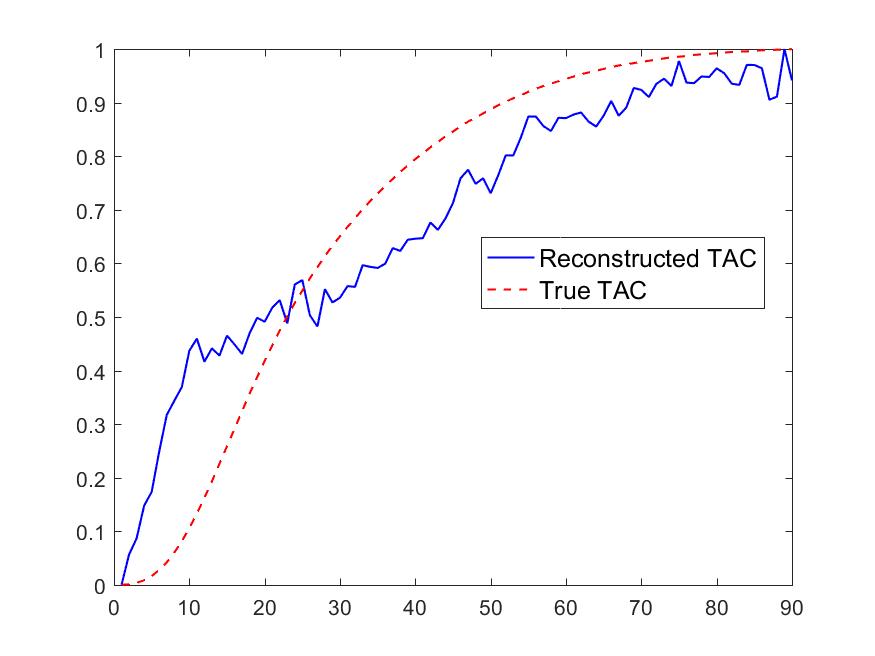}}
\subfigure{
\includegraphics[width=.45\linewidth,height=.3\linewidth]{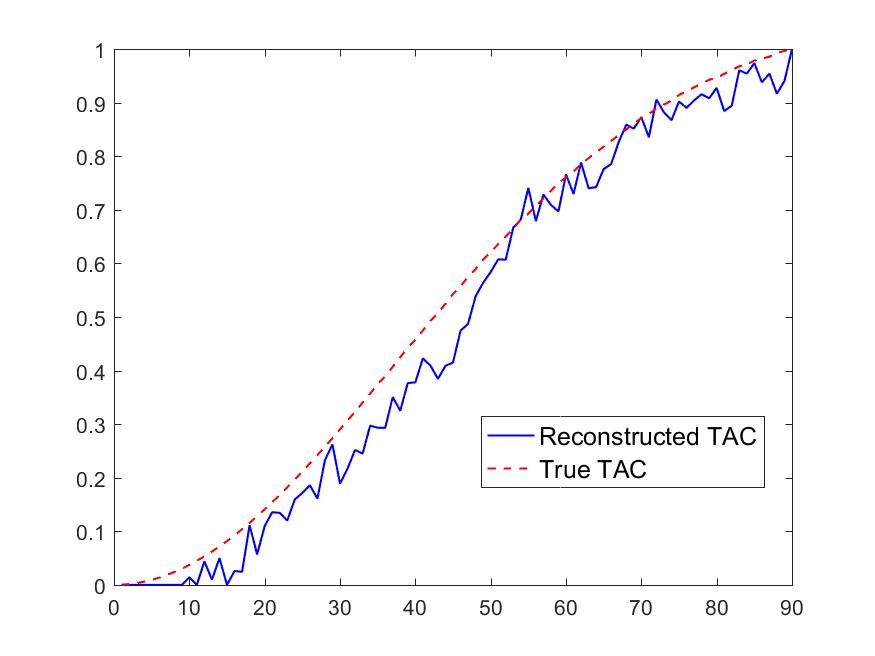}}
\end{center}
\caption{Comparison of true TACs and Reconstructed TACs in four regions (Rat's abdomen phantom).}
\label{fig:ICPoiLiverTAC}
\end{figure}

Figure \ref{fig:ICPoiEllTAC} and \ref{fig:ICPoiLiverTAC} shows the comparison of the TACs of blood and liver for the two  phantoms. The  dash lines are
the normalized true TACs and the solid lines are the normalized one extracted from the reconstruction images by our
method. Even with high level noise and fast change of radioisotope,  the reconstructed one  fit closely to the true one.

\begin{figure}[ht]
\begin{center}
\subfigure{
\includegraphics[width=.45\linewidth]{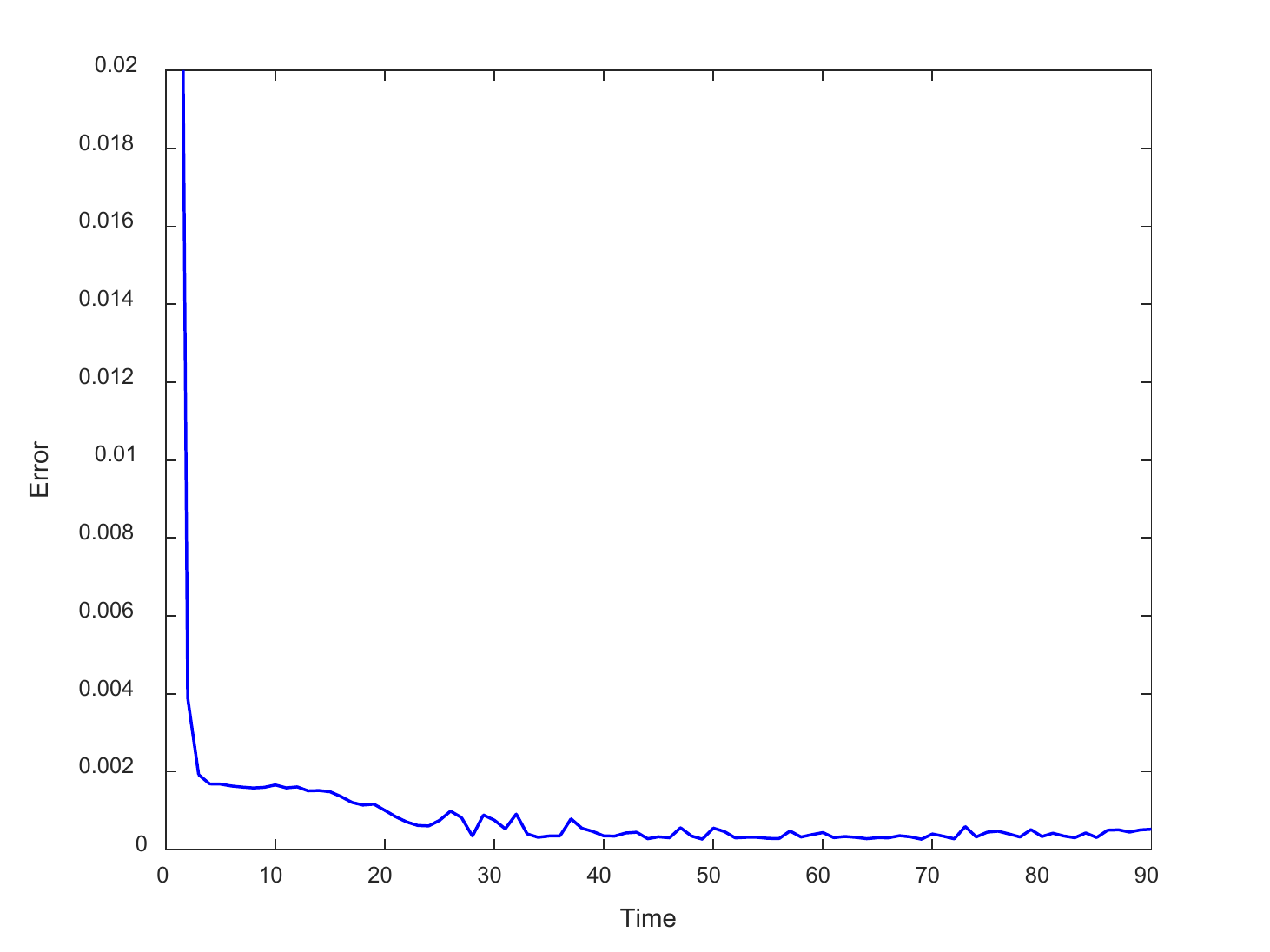}}
\subfigure{
\includegraphics[width=.45\linewidth]{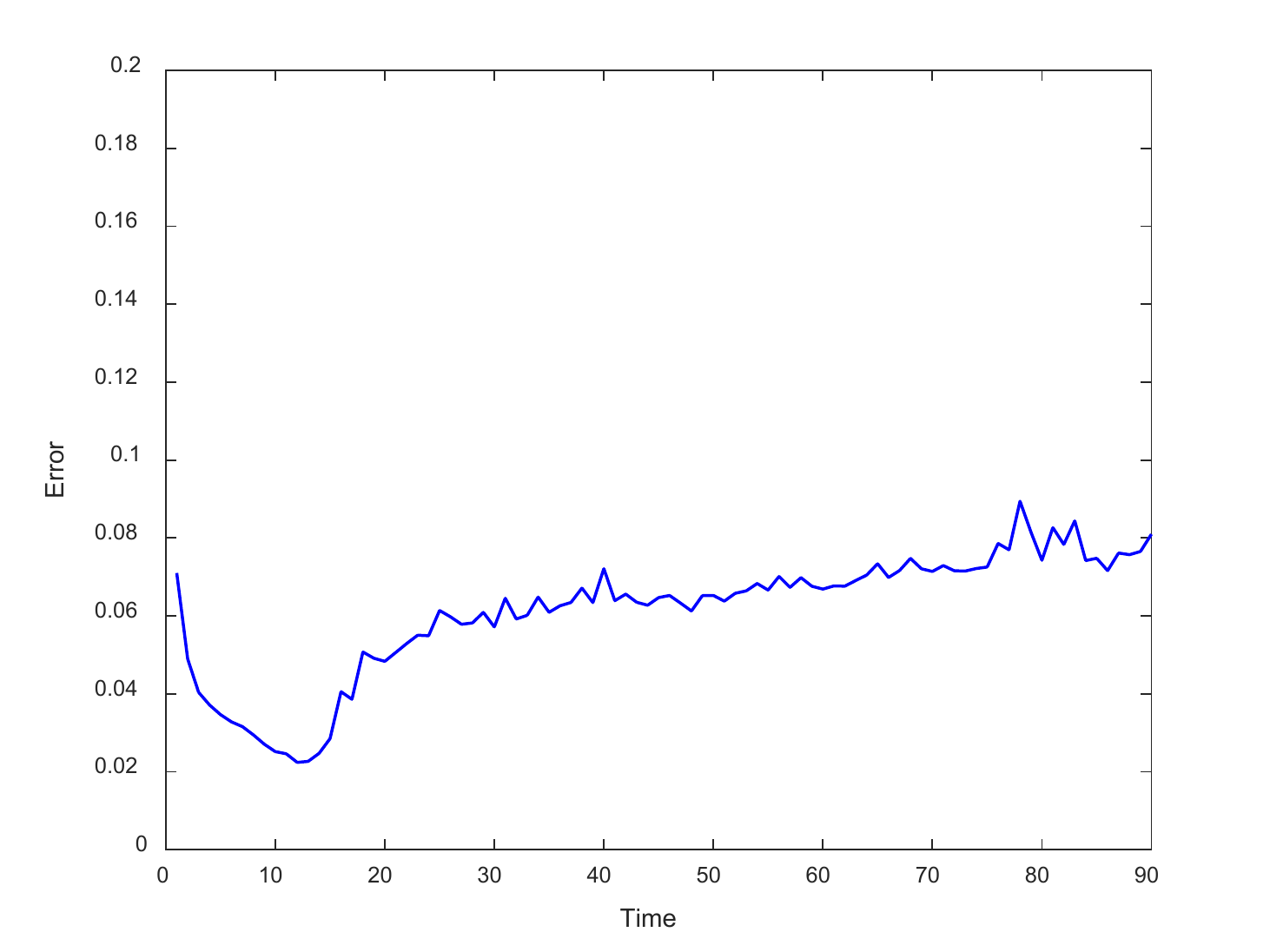}}
\end{center}
\caption{The relative error of the image reconstructed by proposed method with poisson noisy data.
 The left is error of the first numerical experiment  and the right is the second one.}
\label{fig:errPoi}
\end{figure}
 Figure \ref{fig:errPoi} demonstrates the relative error of the image reconstructed by the proposed method with poisson noisy data for the two dataset. By the proposed model, the relative error are small for the proposed method, while for the structure of the second image is complex, the relative error is bigger than the first one.
\subsubsection{Monte Carlo simulation}

\begin{figure}[ht]
\begin{center}
\subfigure[Simulated image Phantom]{
\includegraphics[width=.3\linewidth,height=.3\linewidth]{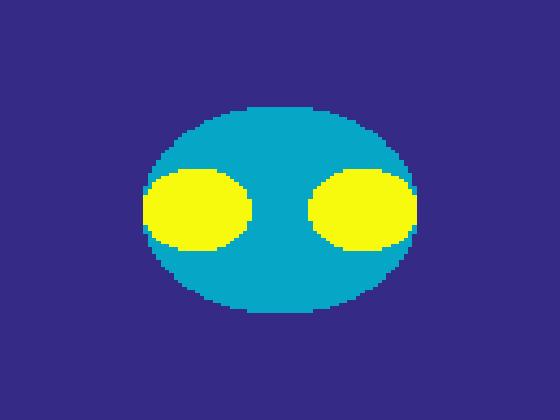}}
\subfigure[Simulated concentration curves]{
\includegraphics[width=.45\linewidth,height=.3\linewidth]{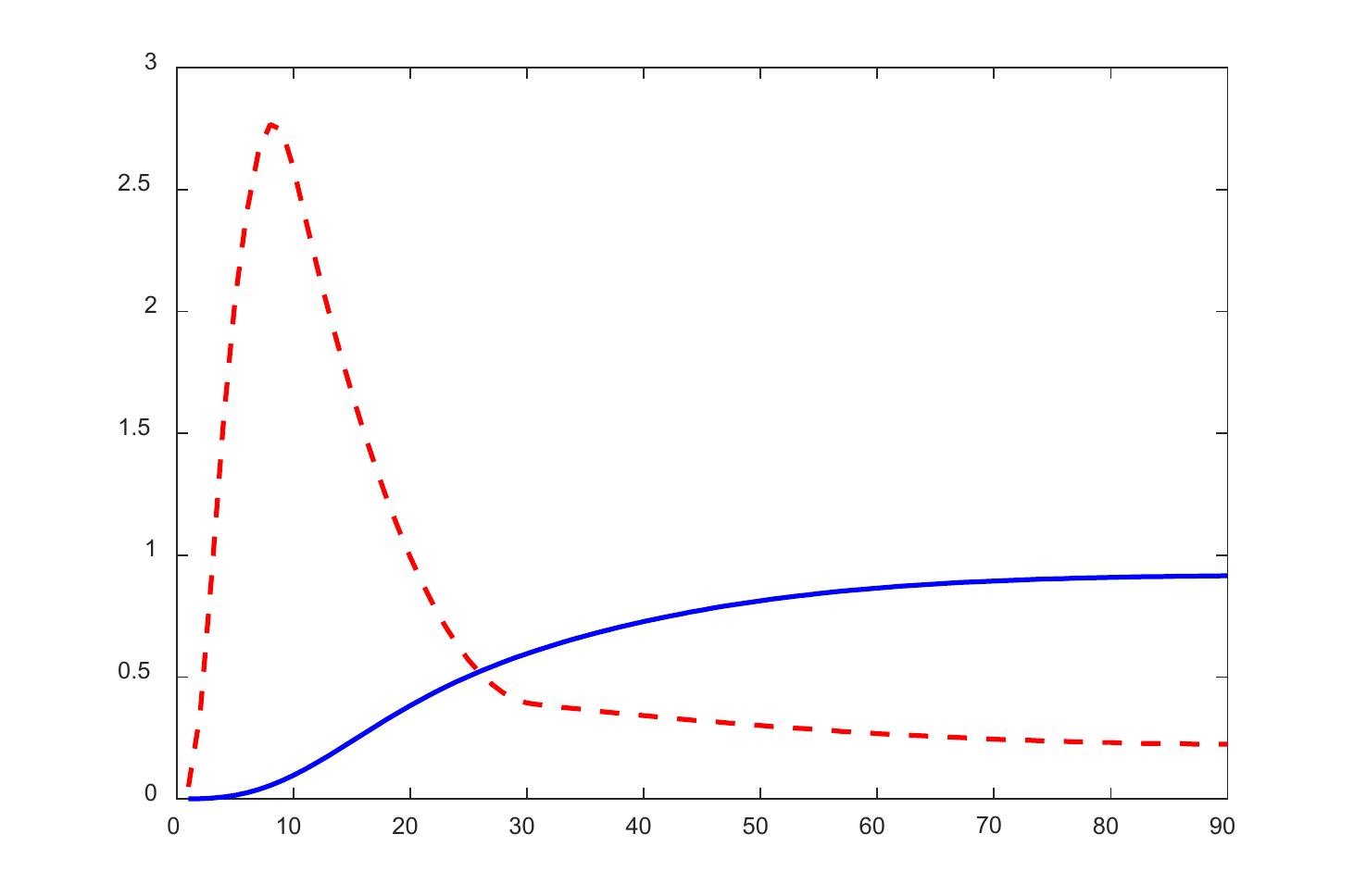}}
\end{center}
\caption{Monte Carlo simulation. Blue solid  line corresponds to the outer circle,  red dash  line corresponds to
the two inner circles.}
\label{fig:MCreal}
\end{figure}

In order to test the performance of the proposed method in a more realistic scenario, we perform a Monte Carlo simulation for dynamic SPECT imaging.
First, we created a $129 \times129$ phantom image consisting of three circles as region of interests, shown in  Figure \ref{fig:MCreal}.
The TAC over  a time period of 90 time steps  of the outer and the two inner circles were displayed in \ref{fig:MCreal}(b).

For each single frame, the photon counts  is a probability proportional to the concentration in every region. The events are detected by a virtual double heads gamma camera rotating around the patient by $1$ degrees per time step, which consists of $374$ detector bins. Every simulated decay event is projected  and counted by the corresponding detector bin.

\begin{figure}[ht]
\begin{center}
\subfigure[sinogram for $\mathrm{events}=20000$]{
\includegraphics[width=.45\linewidth]{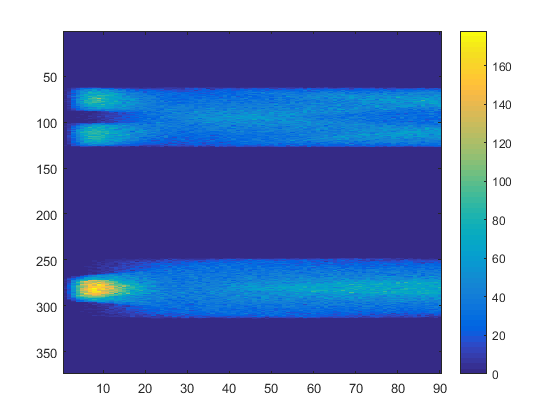}}
\subfigure[sinogram for $\mathrm{events}=200000$]{
\includegraphics[width=.45\linewidth]{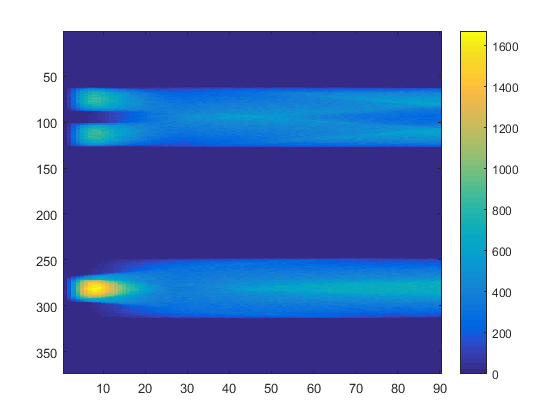}}
\end{center}
\caption{Monte Carlo sinogram data.}
\label{fig:MCSinogram}
\end{figure}

We set the number of events counted by the detector as $\mathrm{events} = 2\times 10^4$ and $\mathrm{events} = 2\times 10^5$) times the average concentration in one pixel of two different tests.
The signogram images the count in each bin of two settings are shown in Figure \ref{fig:MCSinogram}.

\begin{figure}
\begin{tabular}{c@{\hspace{2pt}}c@{\hspace{2pt}}c@{\hspace{2pt}}c@{\hspace{2pt}}c@{\hspace{2pt}}c@{\hspace{2pt}}c@{\hspace{2pt}}c@{\hspace{2pt}}c@{\hspace{2pt}}c}
\includegraphics[width=.1\linewidth,height=.1\linewidth]{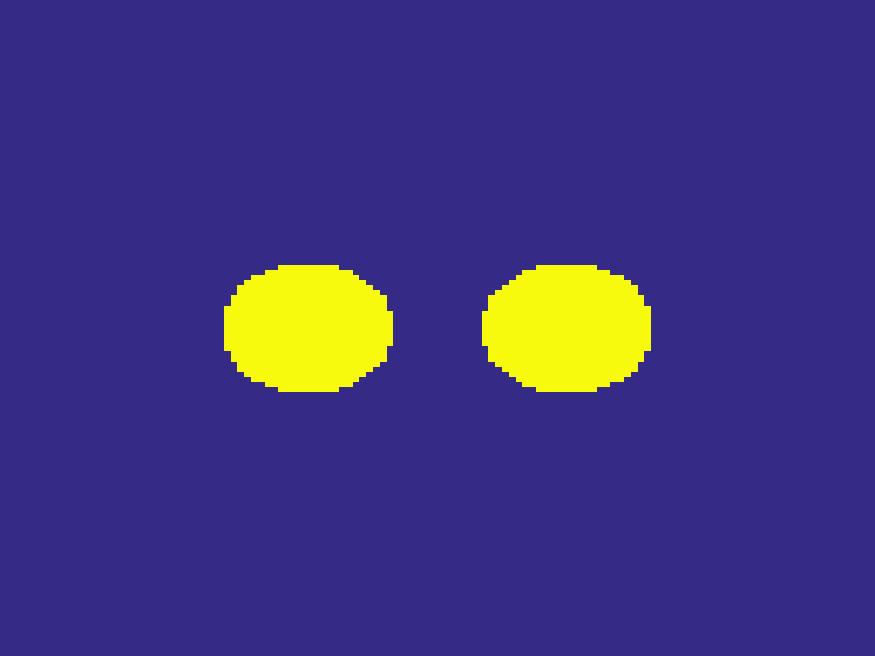}&
\includegraphics[width=.1\linewidth,height=.1\linewidth]{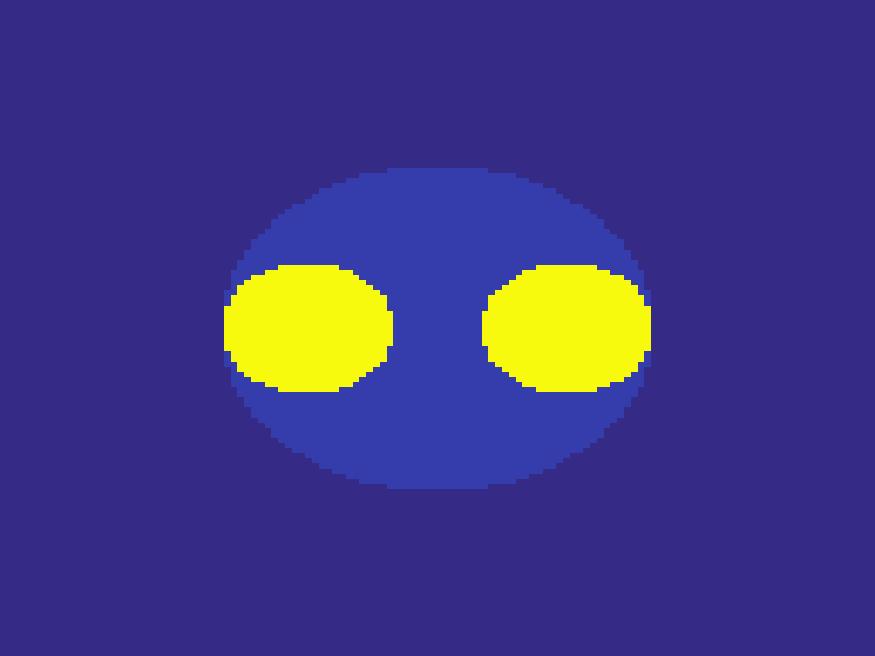}&
\includegraphics[width=.1\linewidth,height=.1\linewidth]{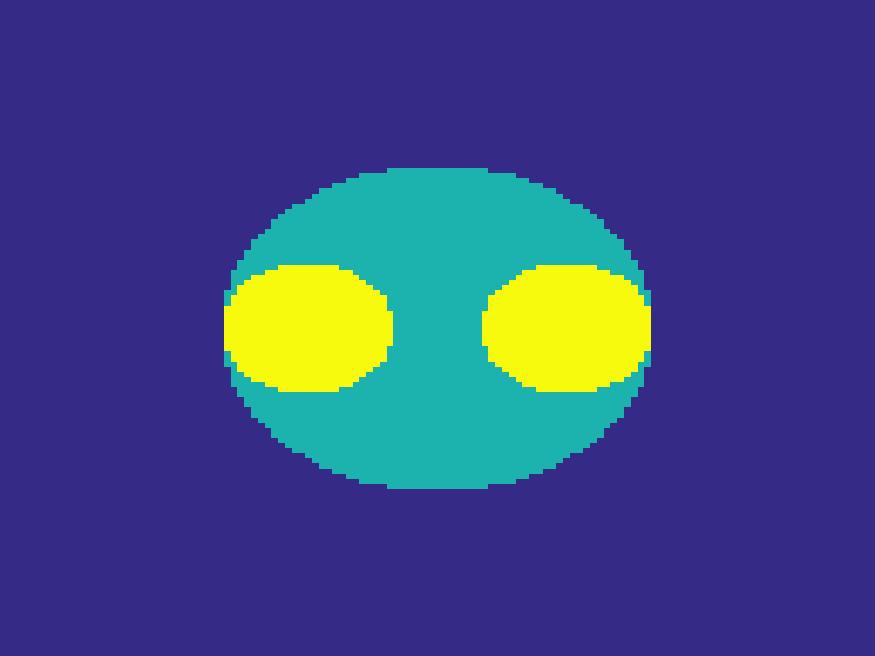}&
\includegraphics[width=.1\linewidth,height=.1\linewidth]{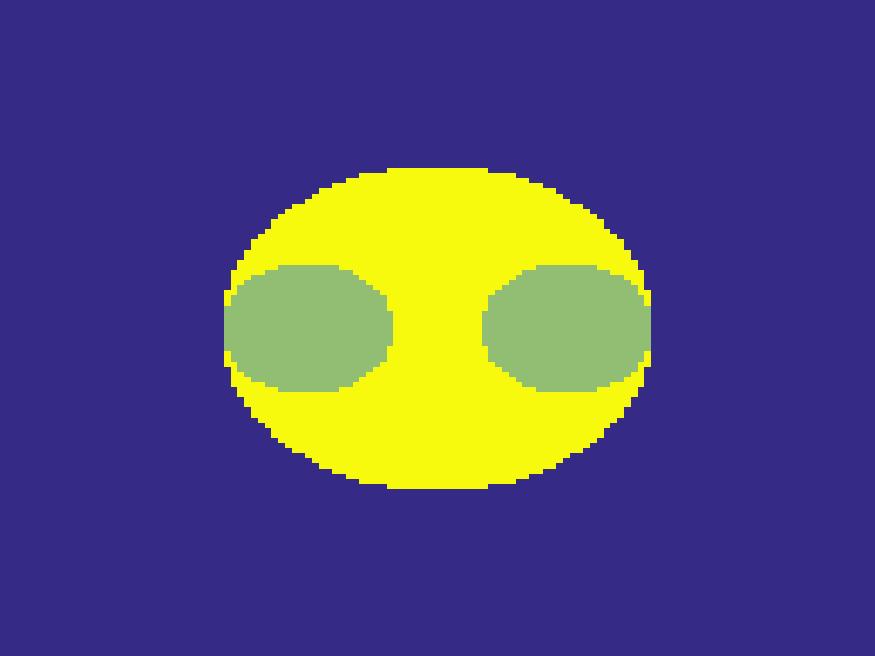}&
\includegraphics[width=.1\linewidth,height=.1\linewidth]{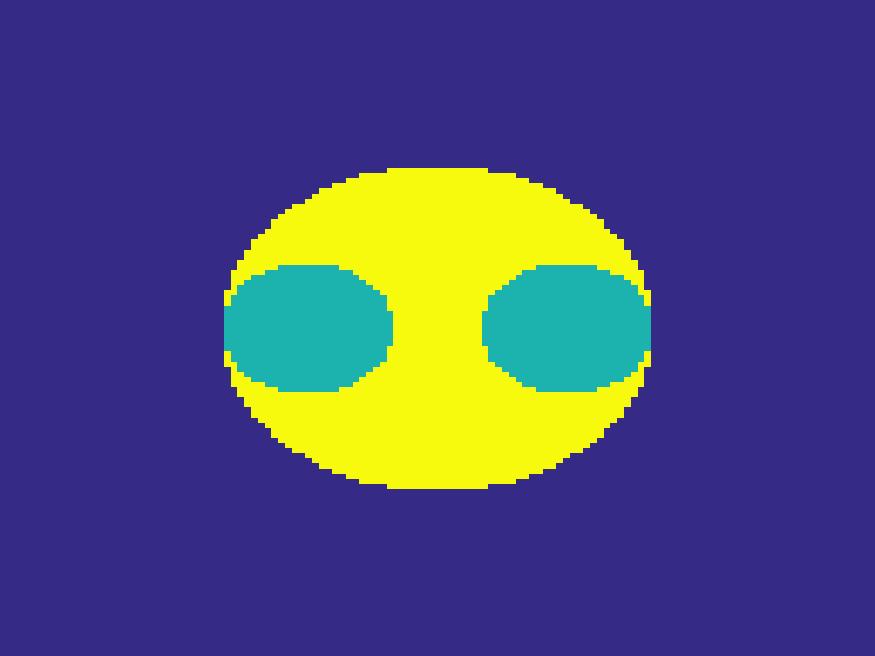}&
\includegraphics[width=.1\linewidth,height=.1\linewidth]{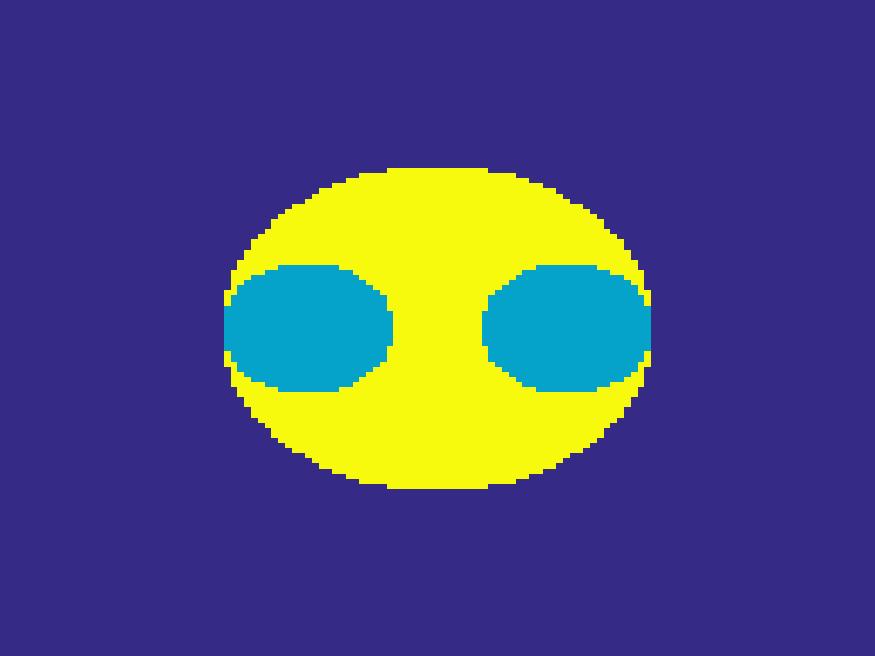}&
\includegraphics[width=.1\linewidth,height=.1\linewidth]{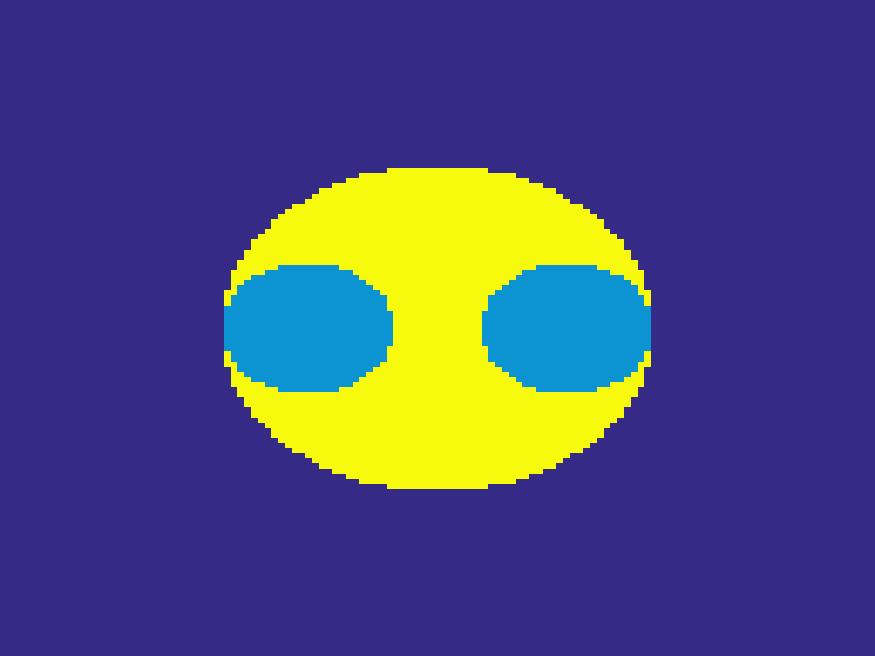}&
\includegraphics[width=.1\linewidth,height=.1\linewidth]{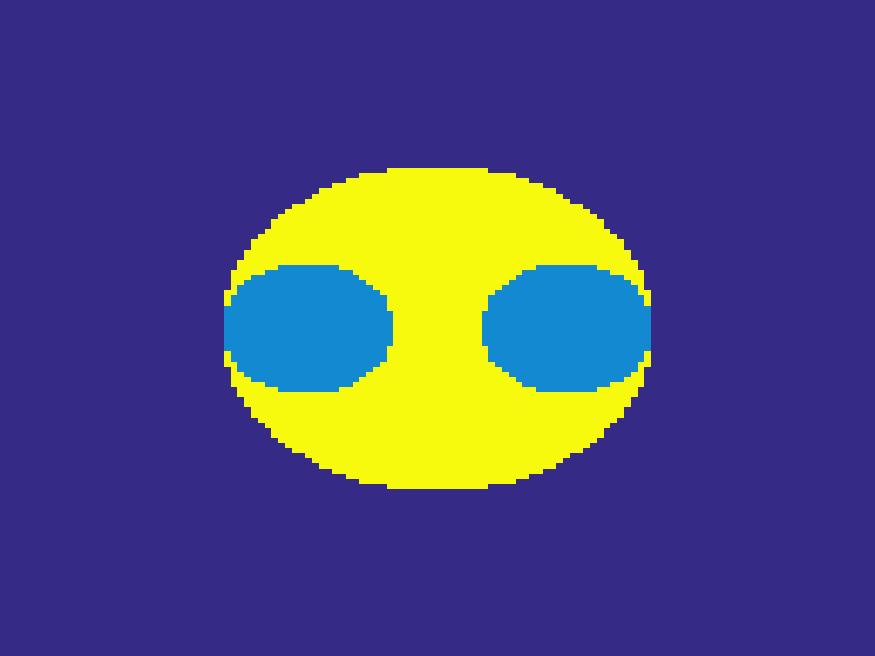}&
\includegraphics[width=.1\linewidth,height=.1\linewidth]{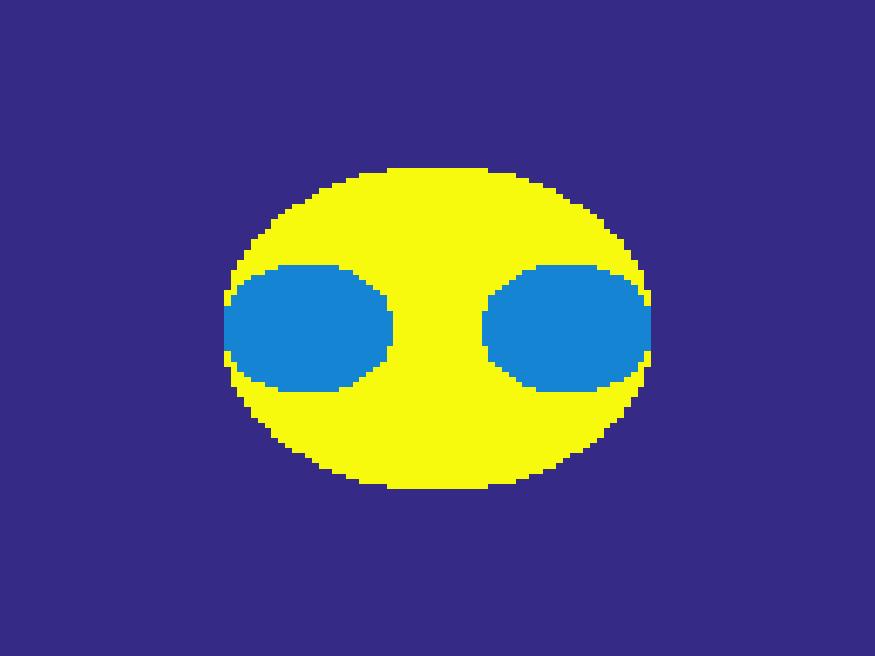}\\
\includegraphics[width=.1\linewidth,height=.1\linewidth]{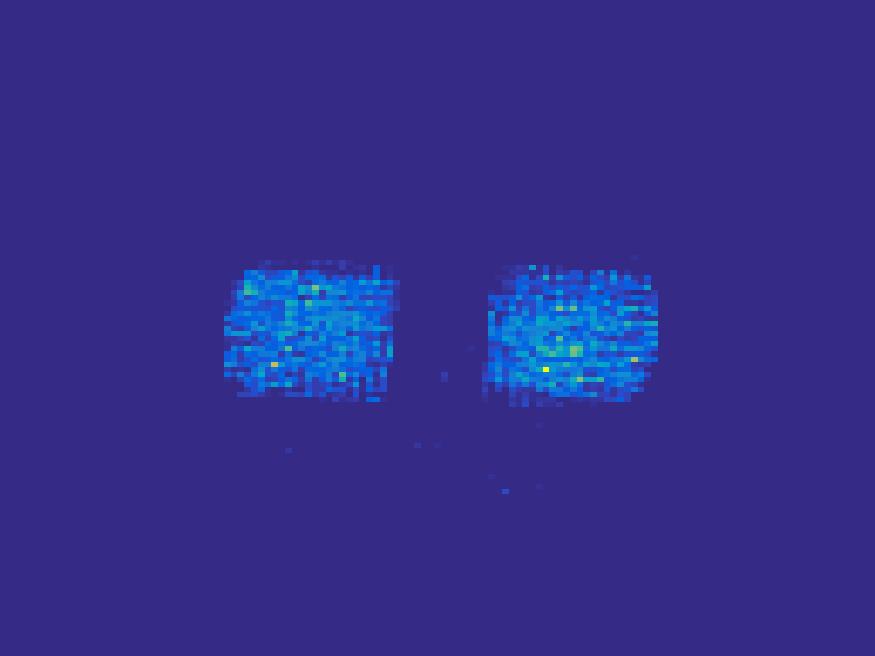}&
\includegraphics[width=.1\linewidth,height=.1\linewidth]{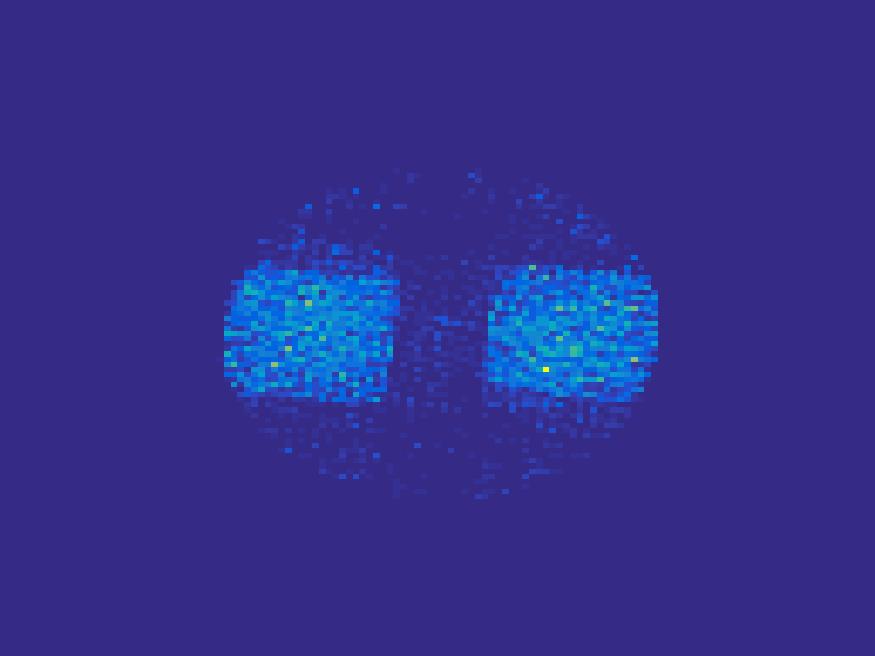}&
\includegraphics[width=.1\linewidth,height=.1\linewidth]{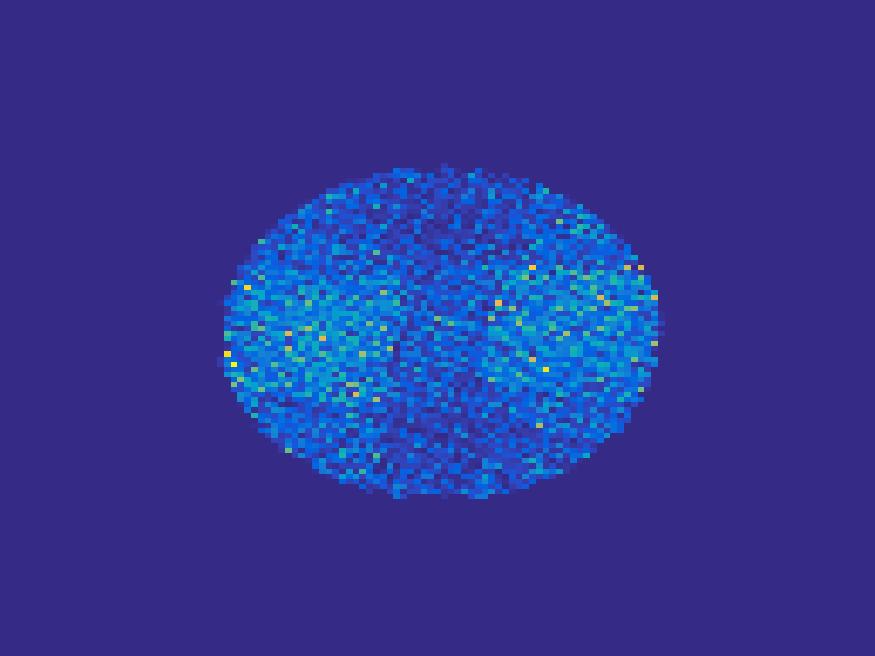}&
\includegraphics[width=.1\linewidth,height=.1\linewidth]{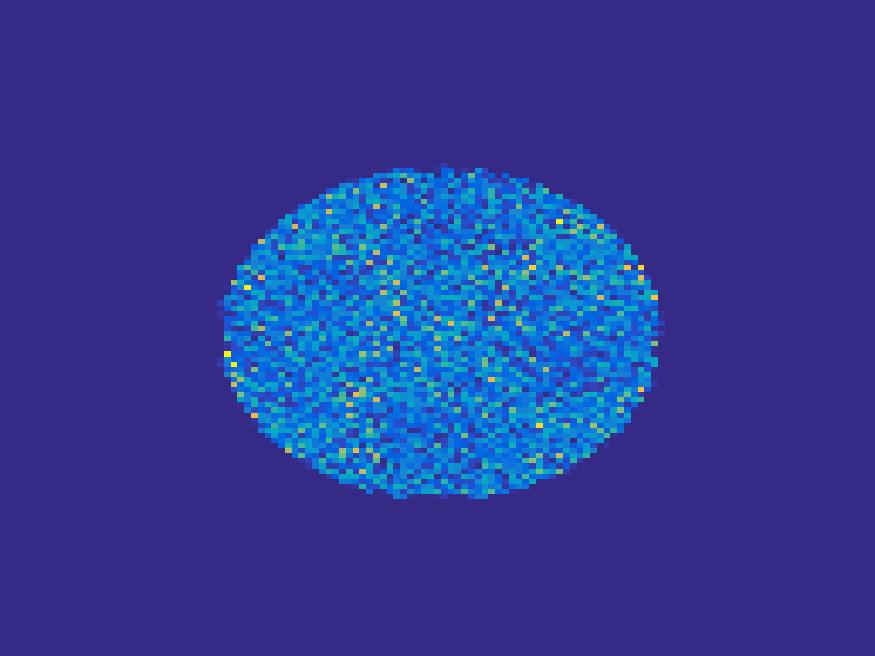}&
\includegraphics[width=.1\linewidth,height=.1\linewidth]{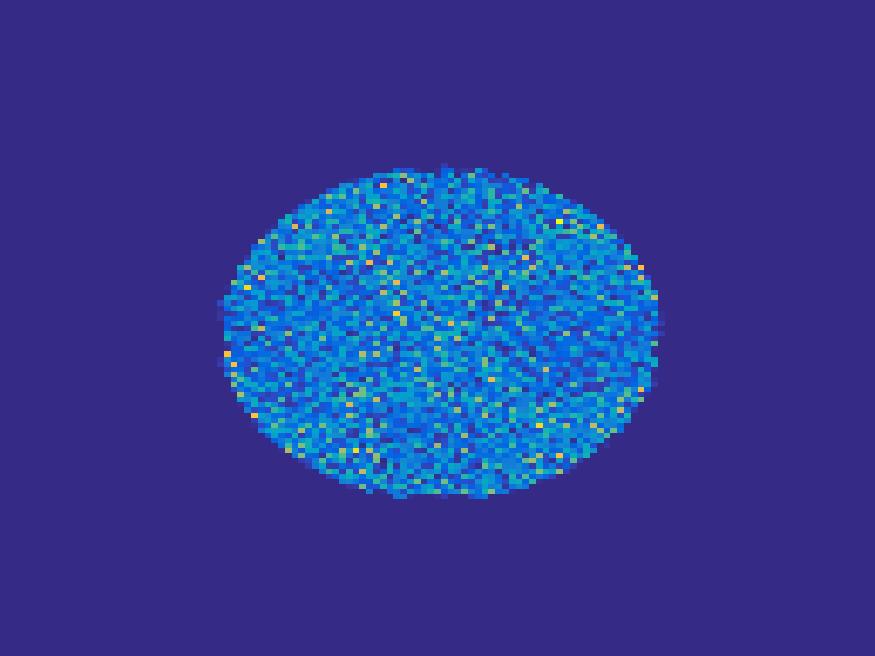}&
\includegraphics[width=.1\linewidth,height=.1\linewidth]{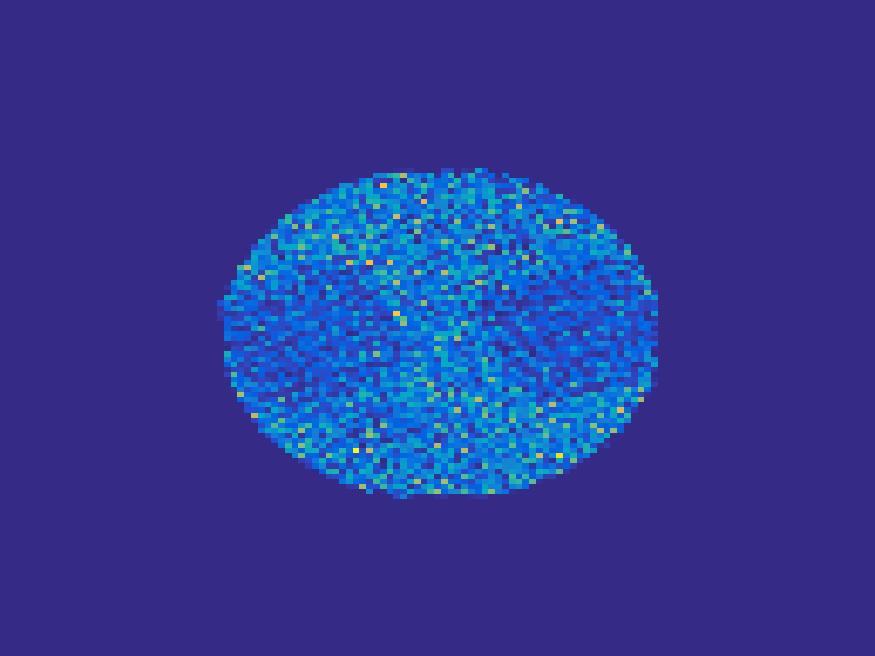}&
\includegraphics[width=.1\linewidth,height=.1\linewidth]{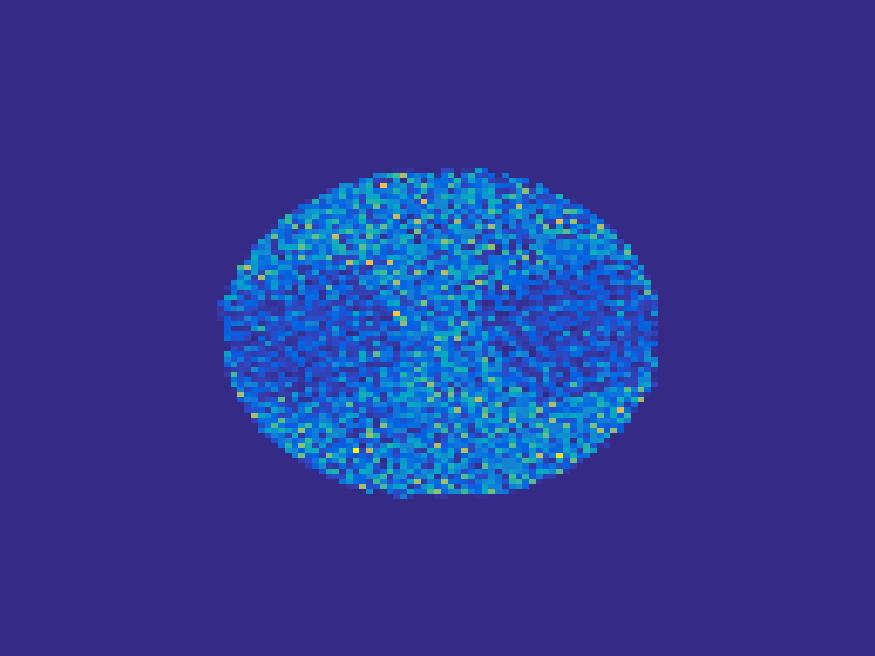}&
\includegraphics[width=.1\linewidth,height=.1\linewidth]{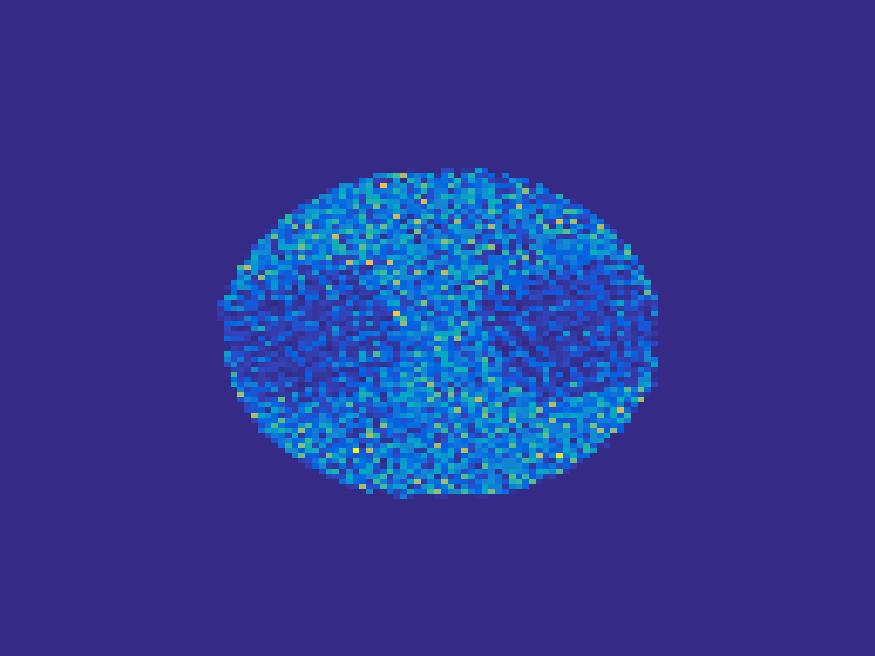}&
\includegraphics[width=.1\linewidth,height=.1\linewidth]{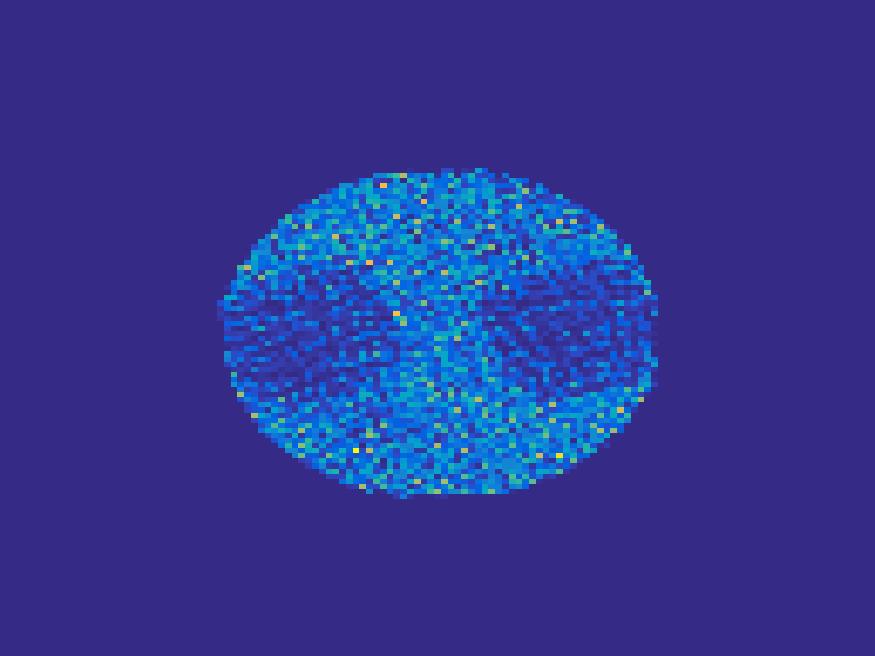}\\
\includegraphics[width=.1\linewidth,height=.1\linewidth]{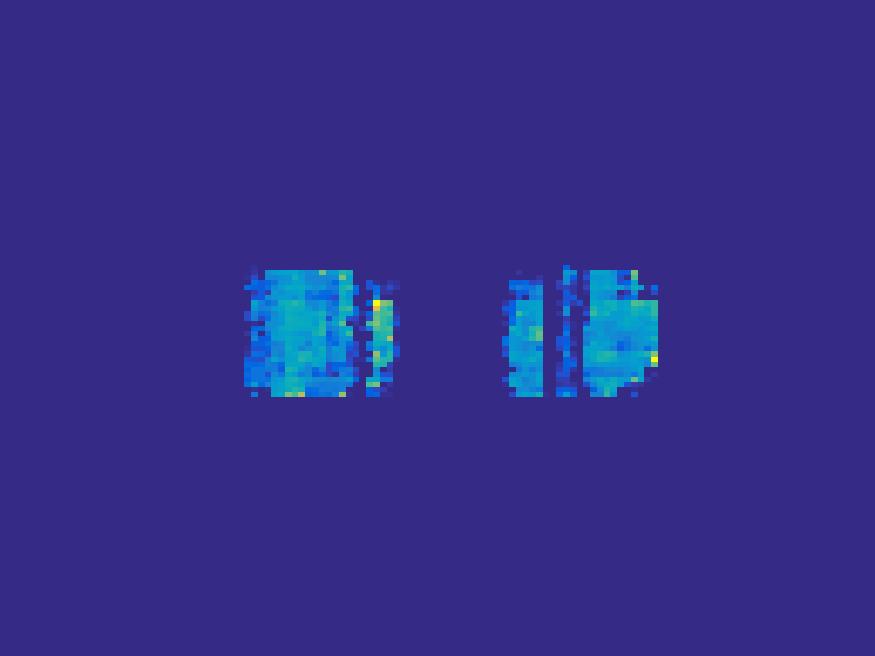}&
\includegraphics[width=.1\linewidth,height=.1\linewidth]{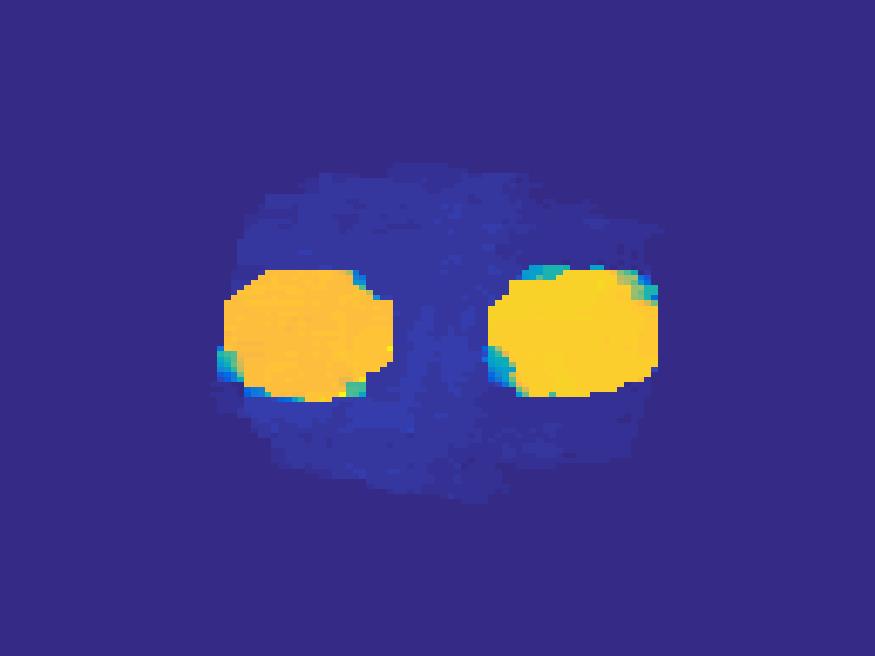}&
\includegraphics[width=.1\linewidth,height=.1\linewidth]{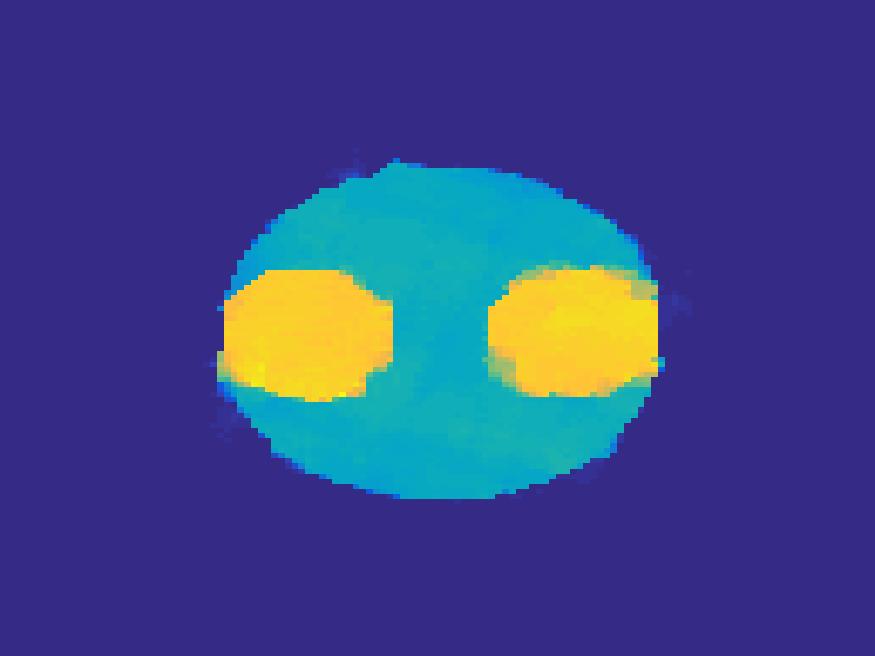}&
\includegraphics[width=.1\linewidth,height=.1\linewidth]{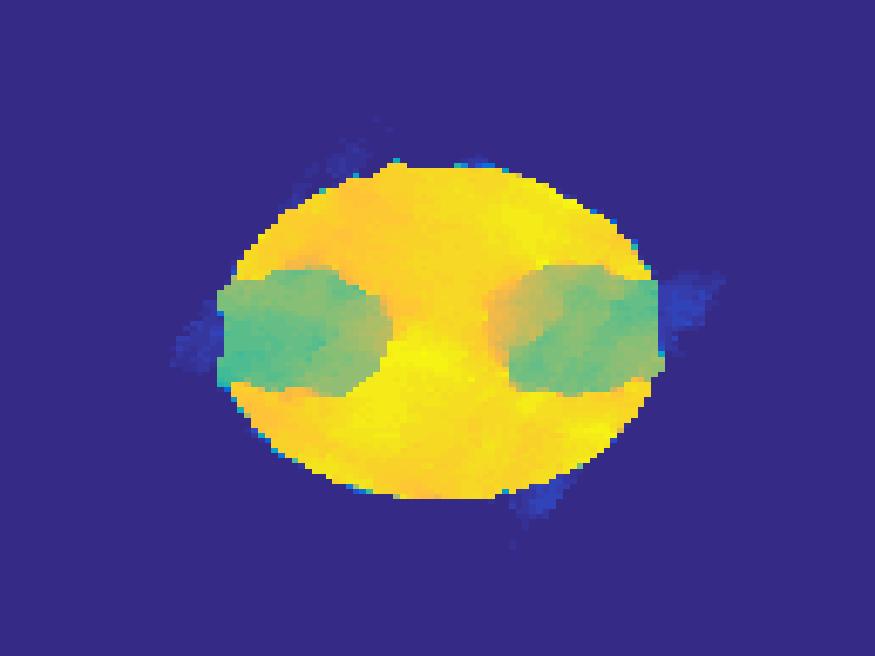}&
\includegraphics[width=.1\linewidth,height=.1\linewidth]{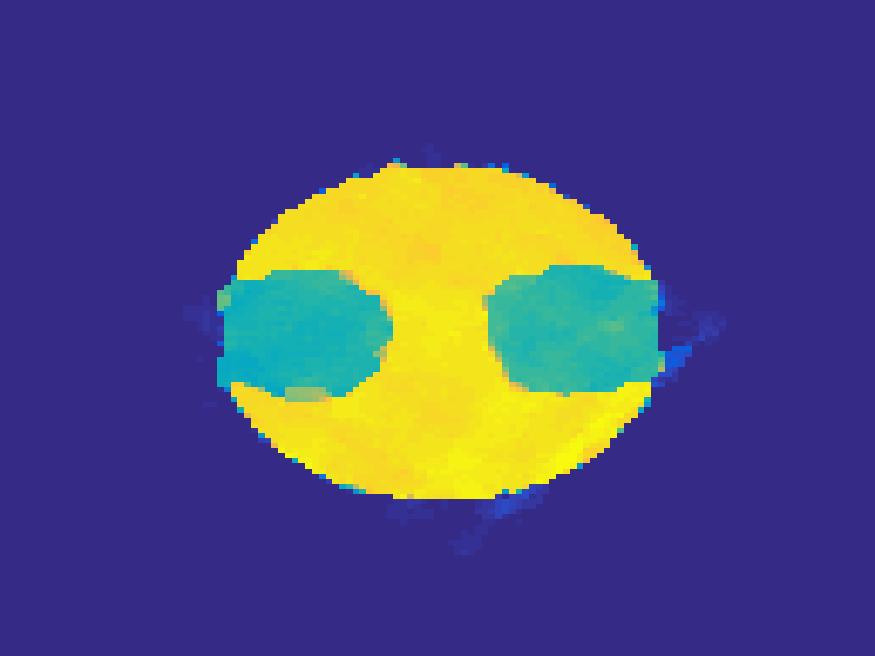}&
\includegraphics[width=.1\linewidth,height=.1\linewidth]{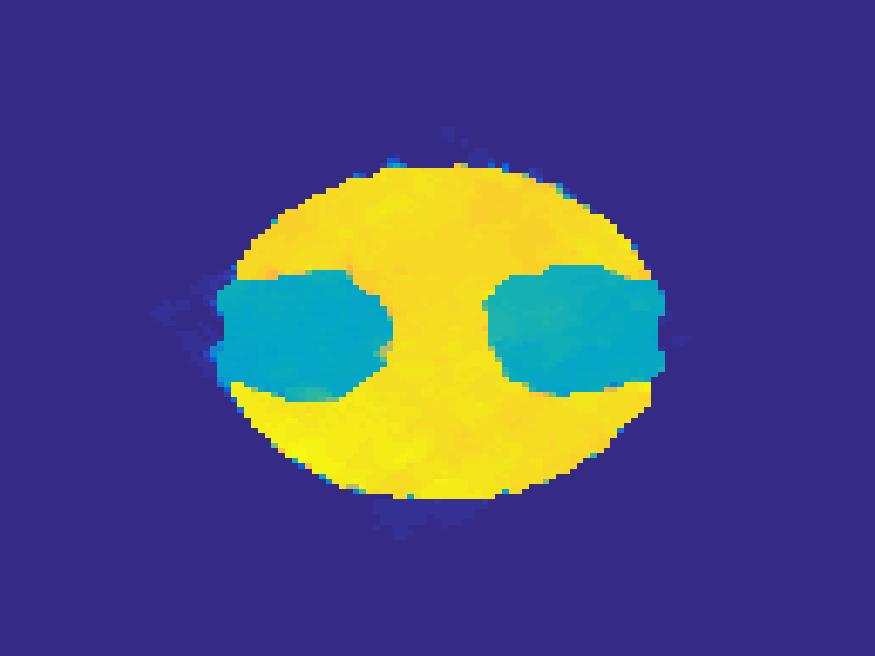}&
\includegraphics[width=.1\linewidth,height=.1\linewidth]{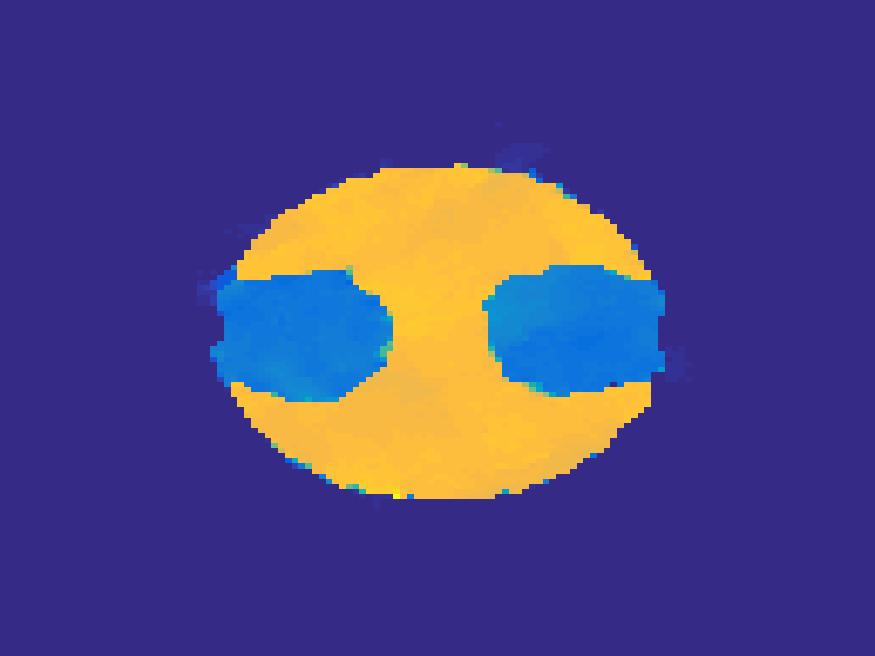}&
\includegraphics[width=.1\linewidth,height=.1\linewidth]{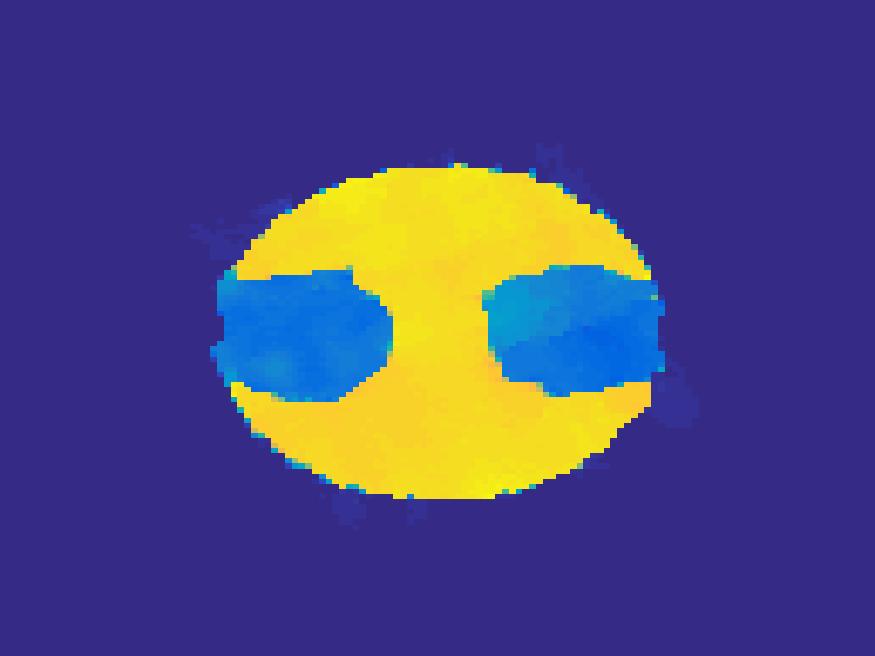}&
\includegraphics[width=.1\linewidth,height=.1\linewidth]{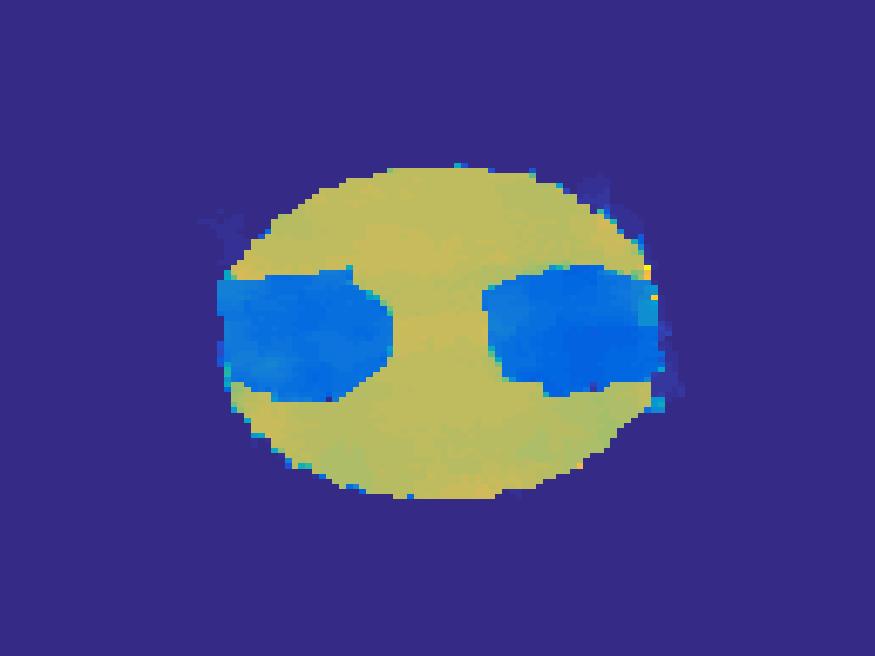}\\
\includegraphics[width=.1\linewidth,height=.1\linewidth]{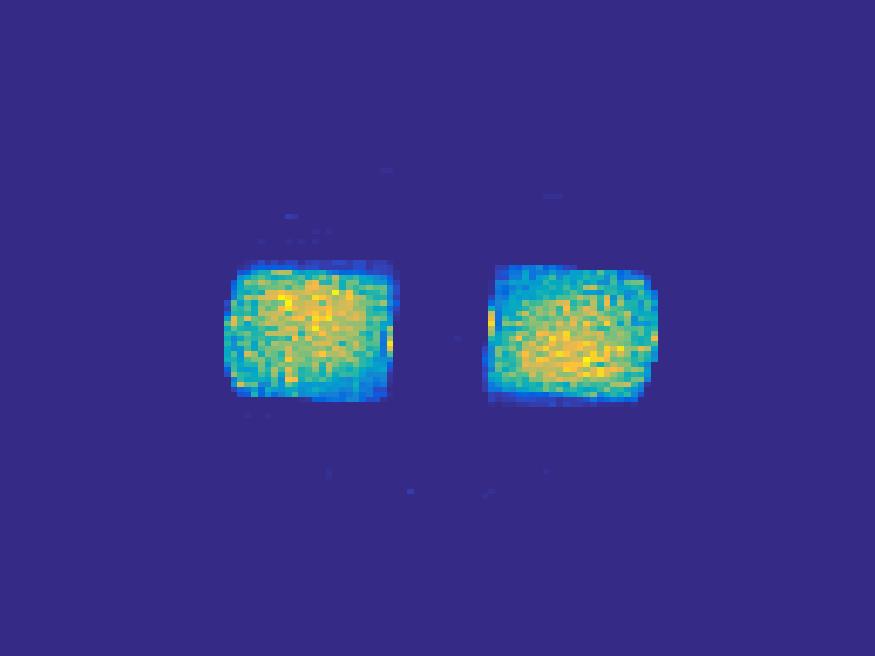}&
\includegraphics[width=.1\linewidth,height=.1\linewidth]{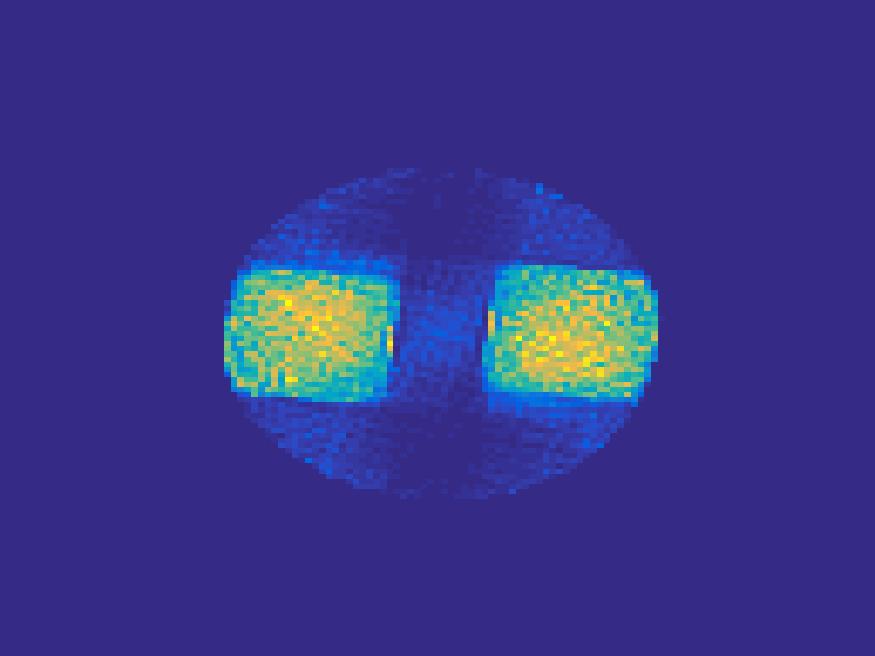}&
\includegraphics[width=.1\linewidth,height=.1\linewidth]{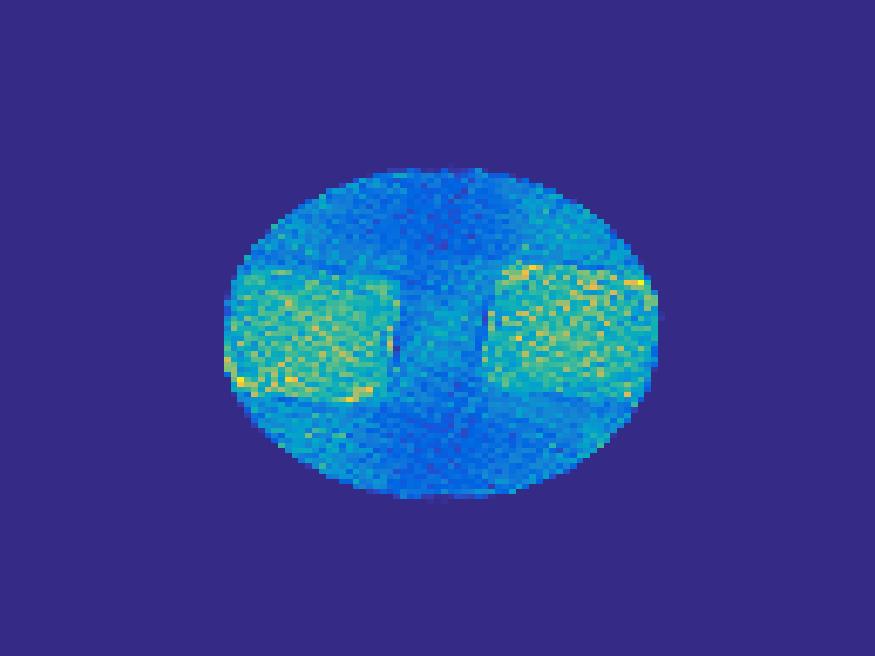}&
\includegraphics[width=.1\linewidth,height=.1\linewidth]{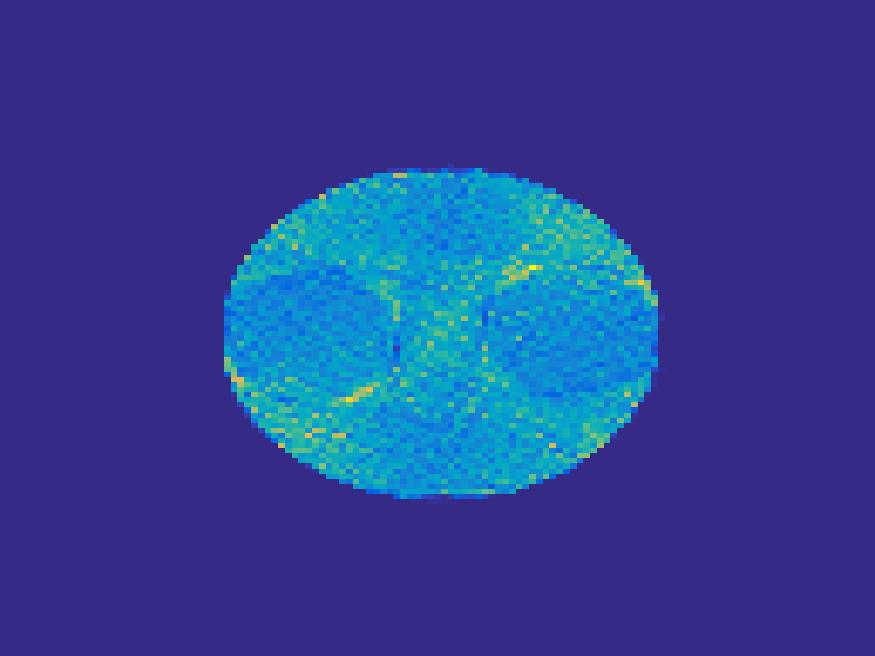}&
\includegraphics[width=.1\linewidth,height=.1\linewidth]{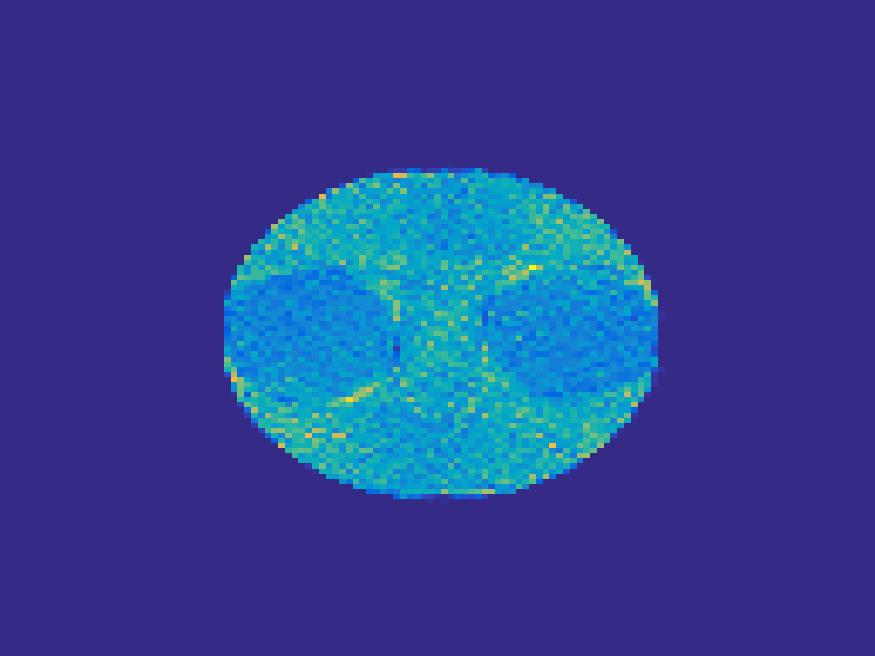}&
\includegraphics[width=.1\linewidth,height=.1\linewidth]{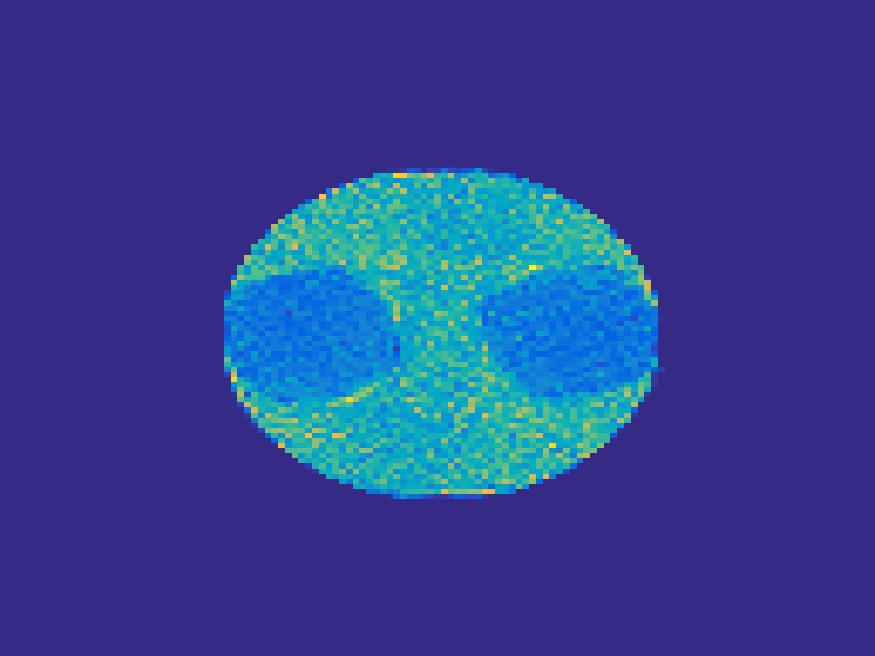}&
\includegraphics[width=.1\linewidth,height=.1\linewidth]{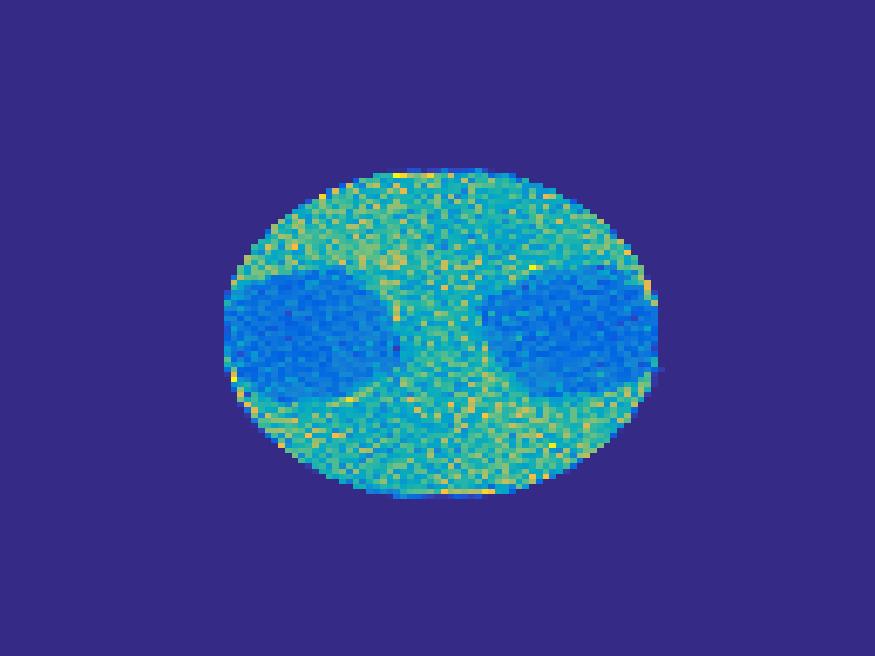}&
\includegraphics[width=.1\linewidth,height=.1\linewidth]{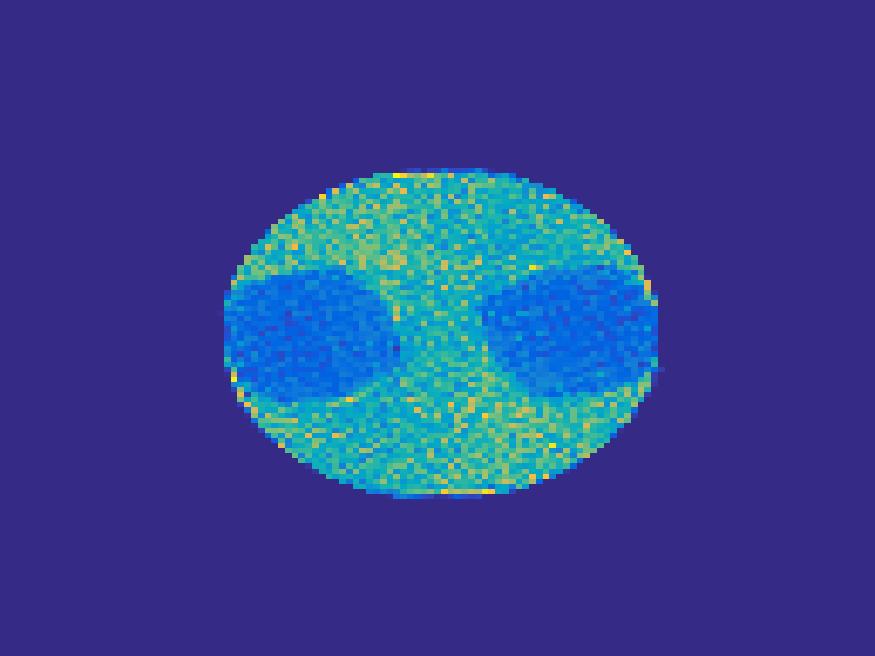}&
\includegraphics[width=.1\linewidth,height=.1\linewidth]{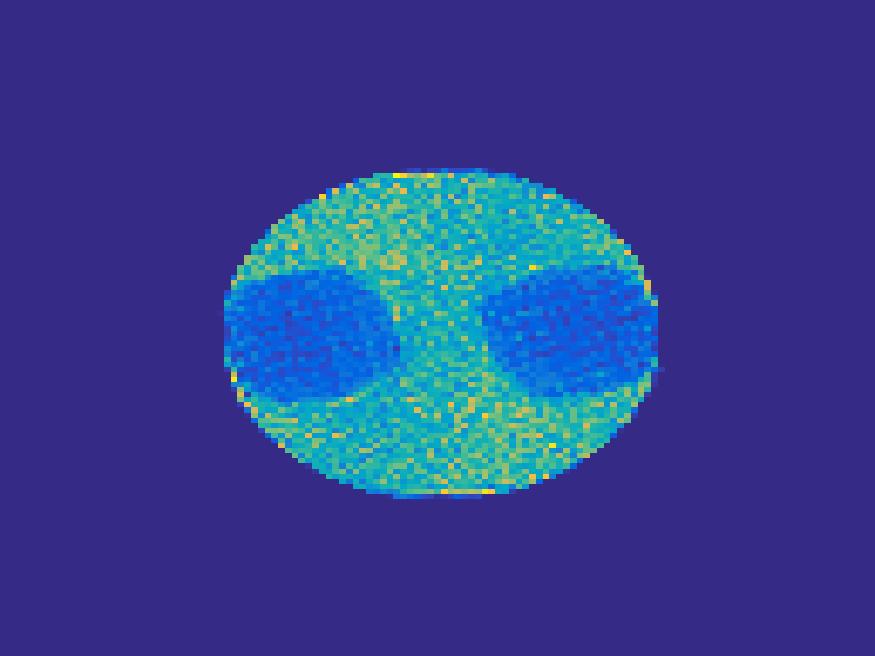}\\
\includegraphics[width=.1\linewidth,height=.1\linewidth]{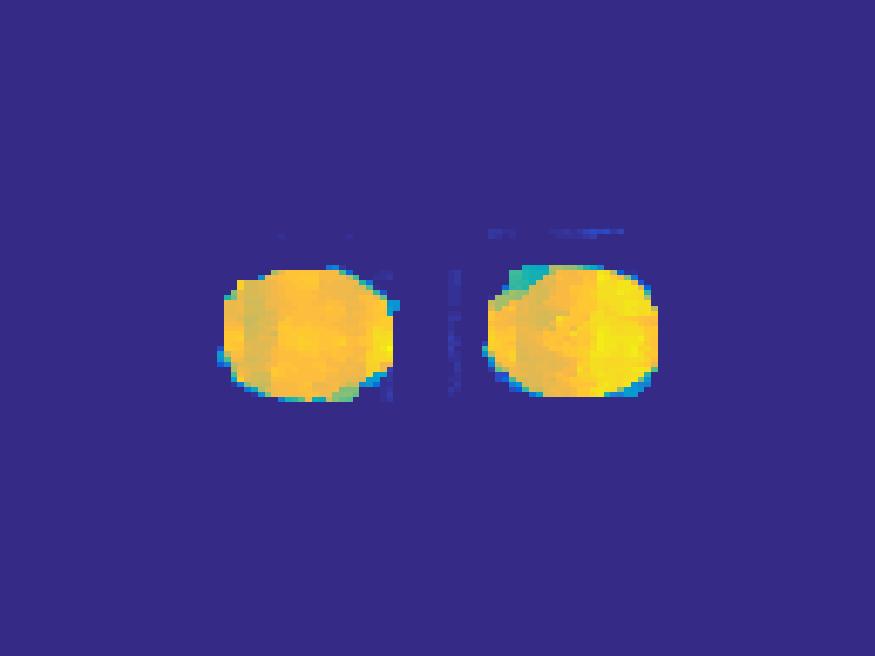}&
\includegraphics[width=.1\linewidth,height=.1\linewidth]{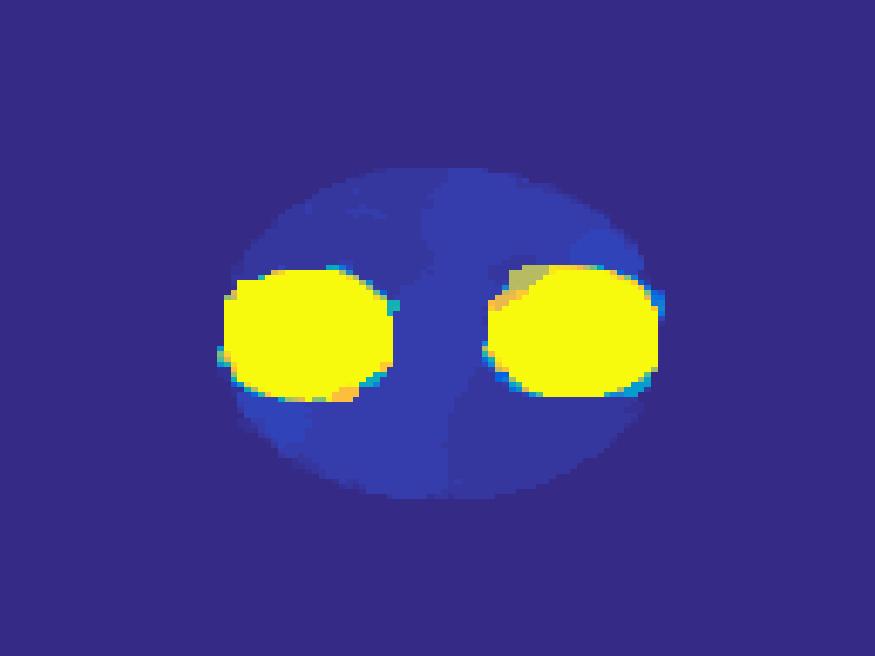}&
\includegraphics[width=.1\linewidth,height=.1\linewidth]{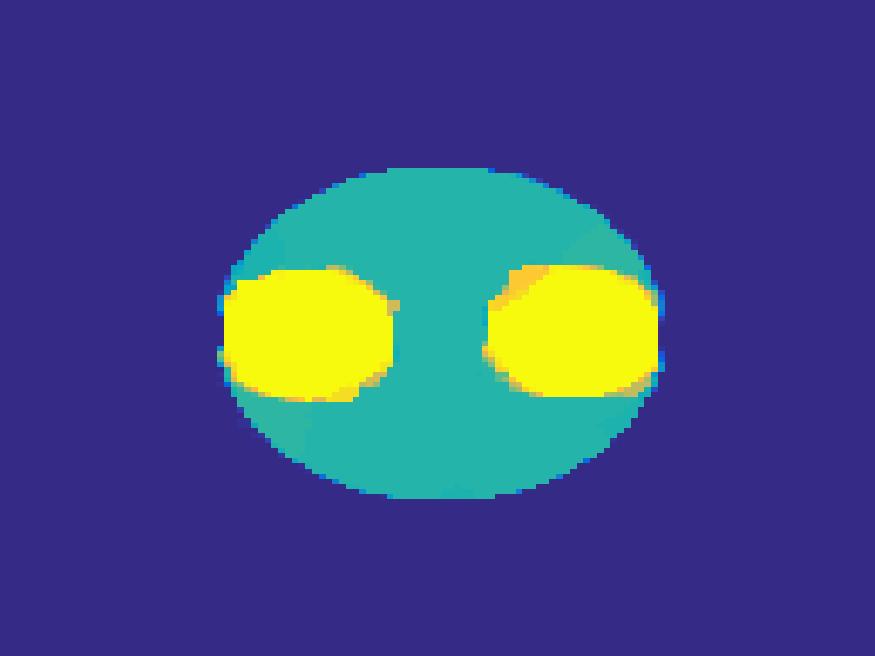}&
\includegraphics[width=.1\linewidth,height=.1\linewidth]{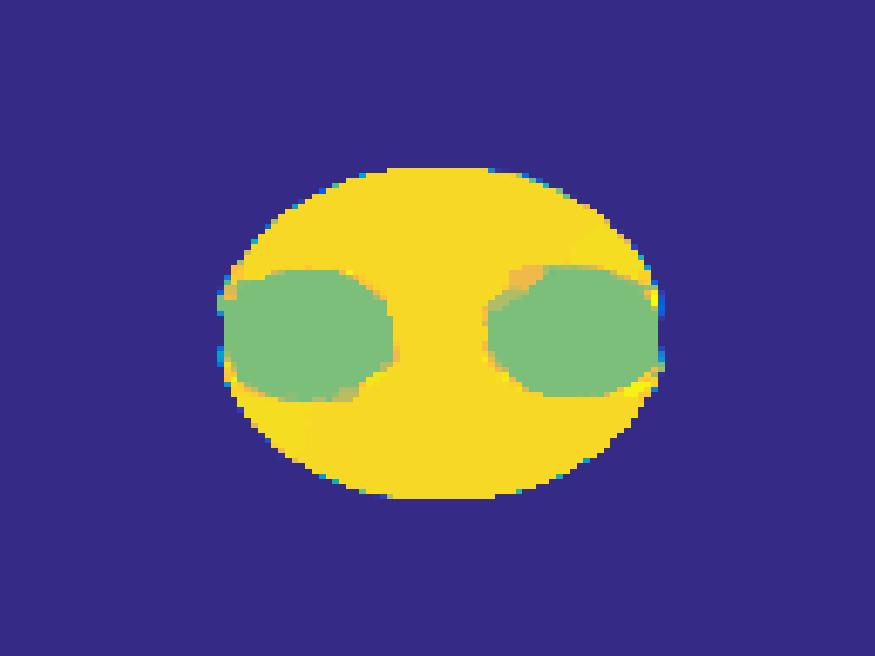}&
\includegraphics[width=.1\linewidth,height=.1\linewidth]{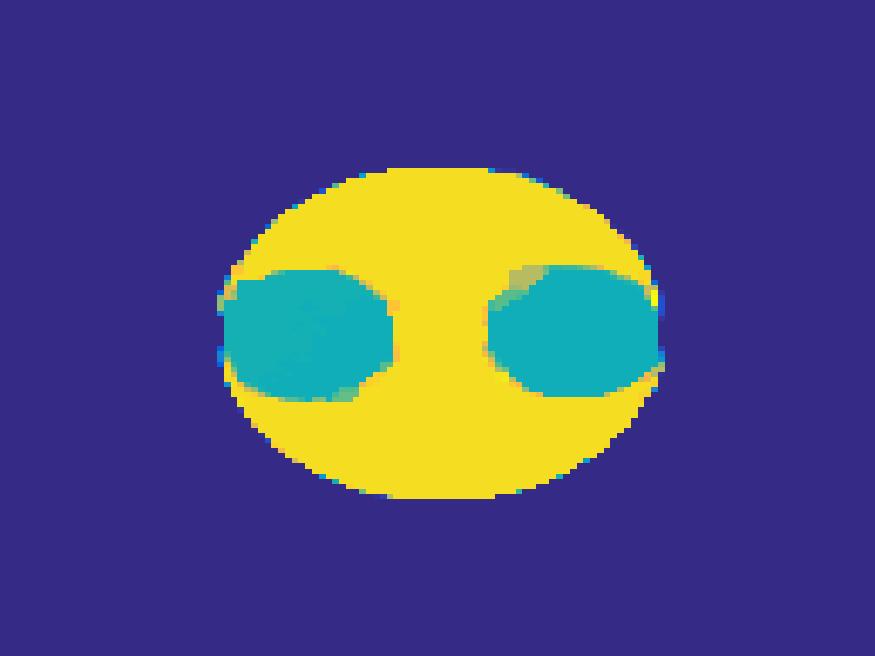}&
\includegraphics[width=.1\linewidth,height=.1\linewidth]{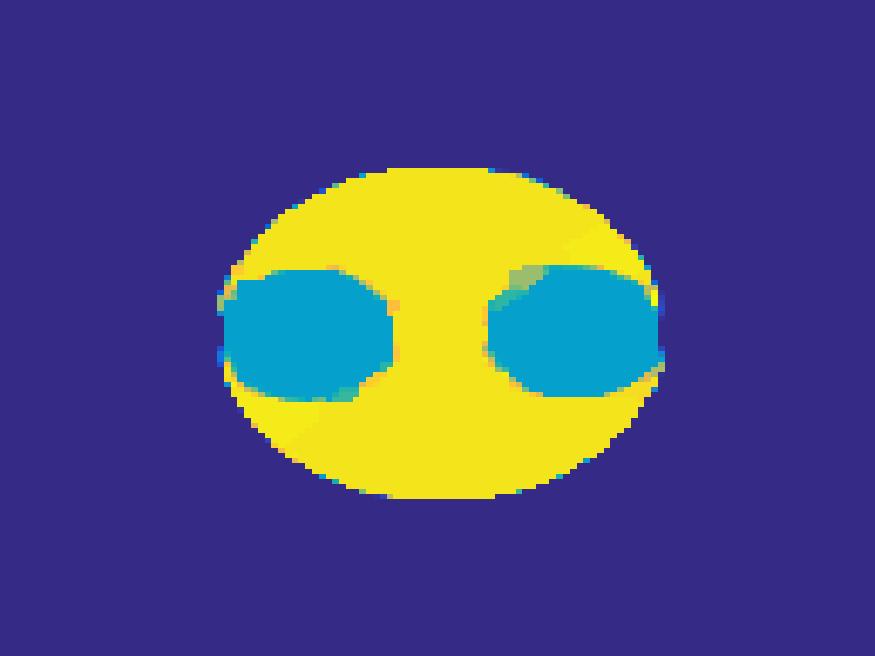}&
\includegraphics[width=.1\linewidth,height=.1\linewidth]{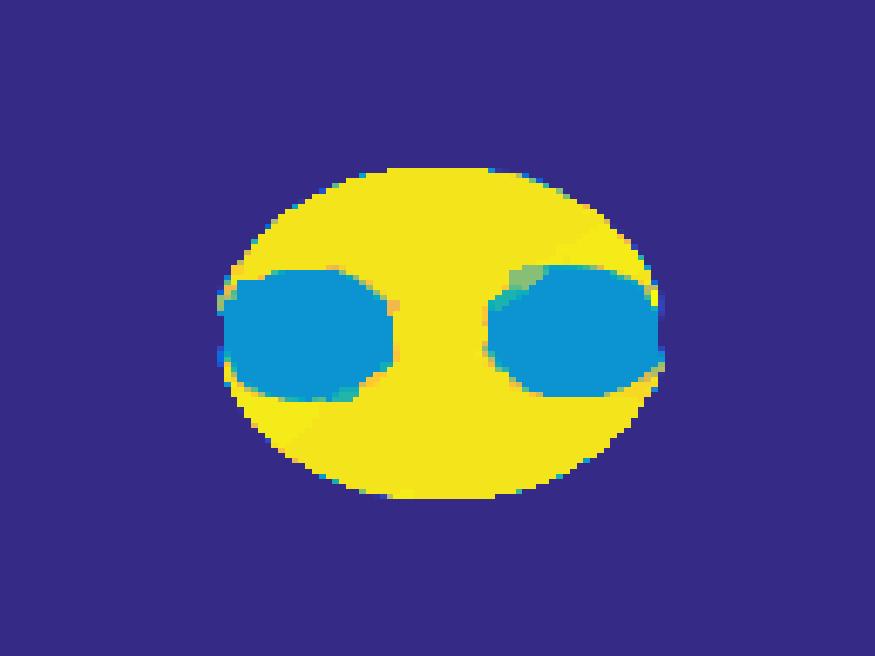}&
\includegraphics[width=.1\linewidth,height=.1\linewidth]{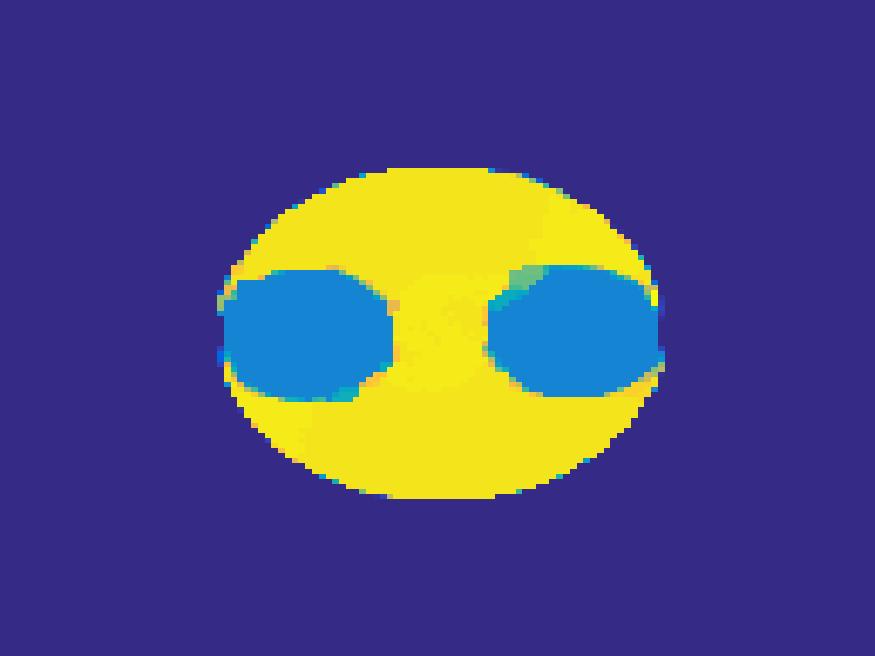}&
\includegraphics[width=.1\linewidth,height=.1\linewidth]{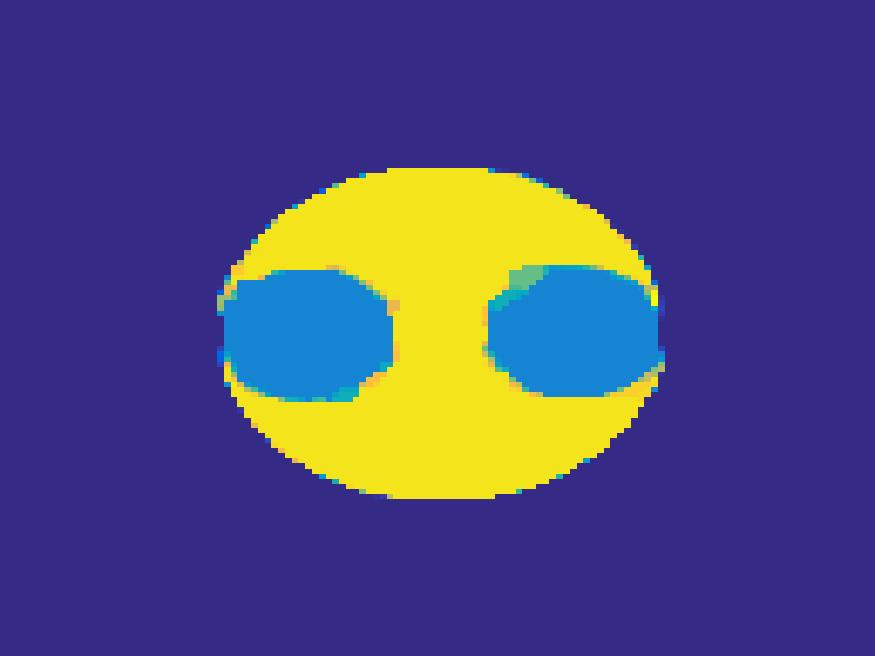}\\
{\footnotesize Frame 1}&
{\footnotesize Frame 11}&
{\footnotesize Frame 21}&
{\footnotesize Frame 31}&
{\footnotesize Frame 41}&
{\footnotesize Frame 51}&
{\footnotesize Frame 61}&
{\footnotesize Frame 71}&
{\footnotesize Frame 81}
\end{tabular}
\caption {Reconstruction with Monte Carlo simulated data. First row: Ground truth. Second and third row are the results with $\mathrm{events}=20000$. Second row : EM with updating $\alpha$ and $B$; third row: proposed method. Forth  and fifth row are the results with $\mathrm{events}=200000$. Forth  row: EM with updating $\alpha$ and $B$; Fifth row: proposed method.}
\label{fig:MCIC}
\end{figure}

Based on the sinogram data, we compare  the proposed method with the alternating EM algorithm. The results for both test cases are shown
in Figure \ref{fig:MCIC}. We can see that for the case of a low count number,  the proposed method is able to reconstruct the regions properly. Within a number of iterations, the algorithm presents a reasonable reconstruction of the region of interest and the corresponding regional tracer concentration curves.

\begin{figure}[ht]
\begin{center}
\subfigure[TAC of $\mathrm{events}=20000$  ]{
\includegraphics[width=.45\linewidth,height=.3\linewidth]{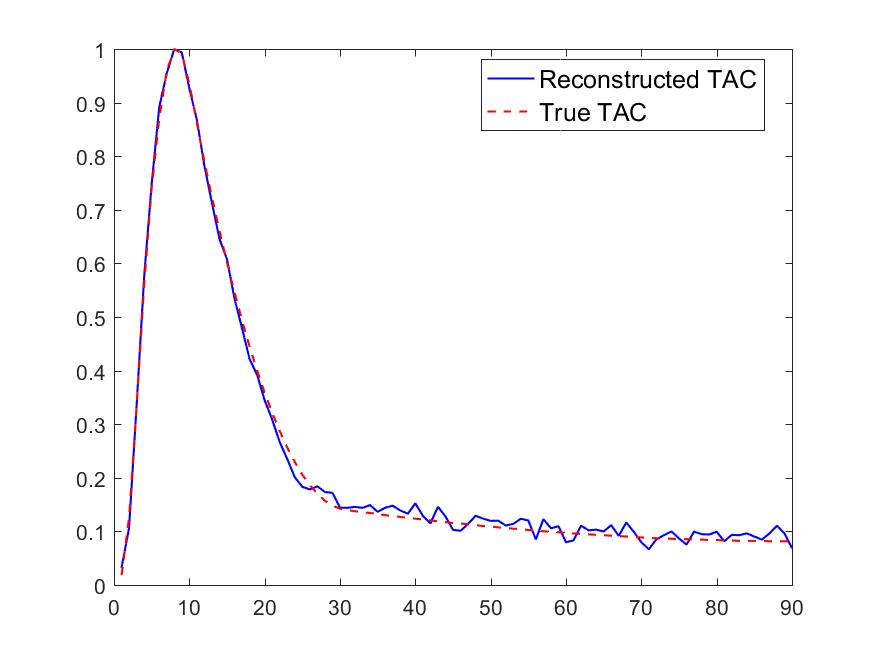}}
\subfigure[ TAC of $\mathrm{events}=20000$ ]{
\includegraphics[width=.45\linewidth,height=.3\linewidth]{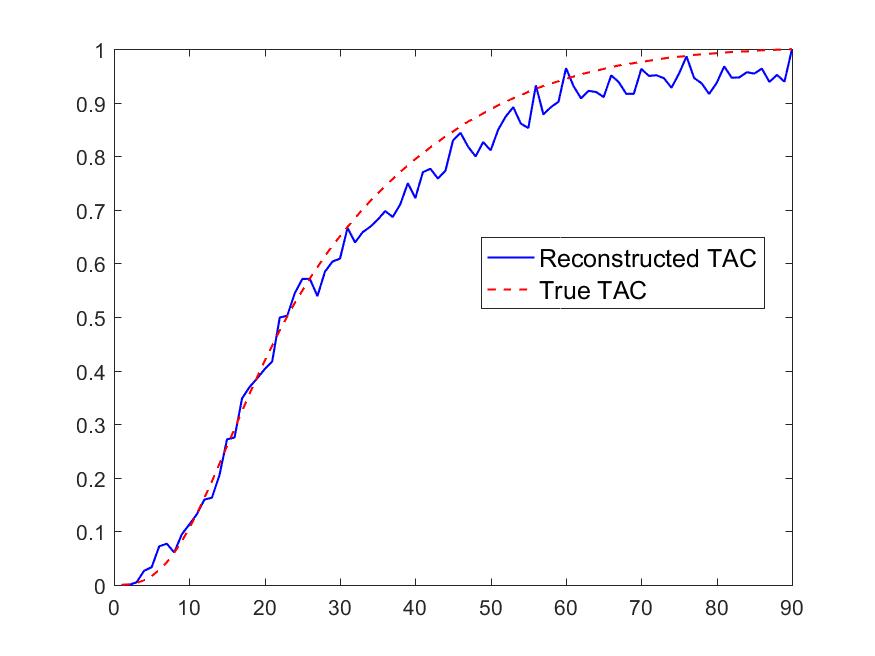}}\\
\subfigure[TAC of $\mathrm{events}=200000$ ]{
\includegraphics[width=.45\linewidth,height=.3\linewidth]{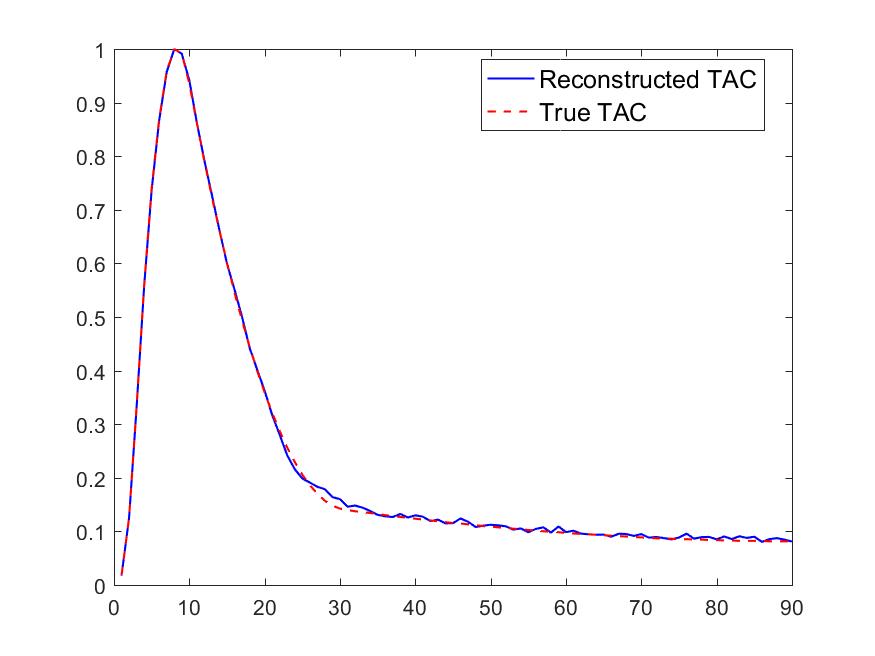}}
\subfigure[TAC of $\mathrm{events}=200000$ ]{
\includegraphics[width=.45\linewidth,height=.3\linewidth]{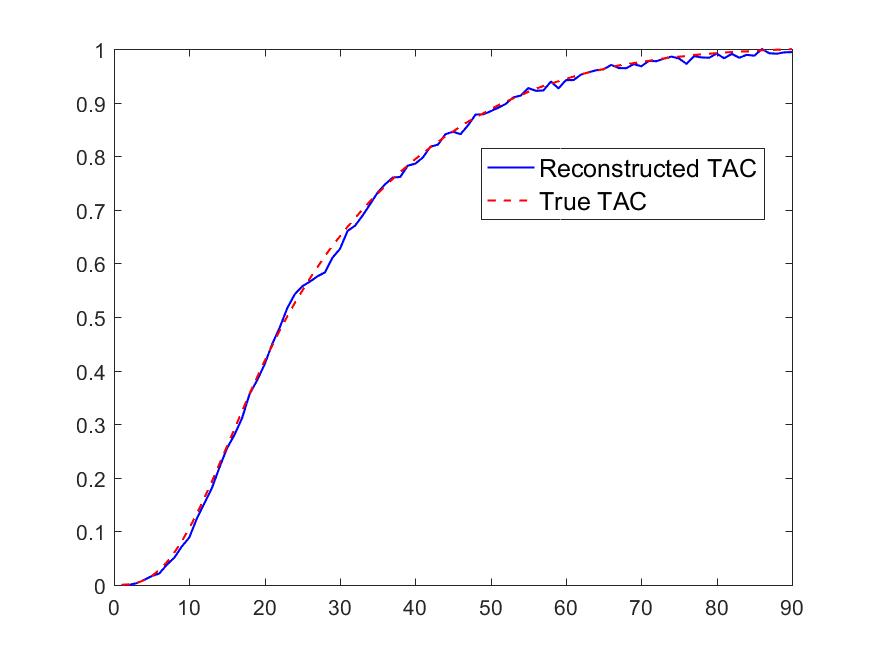}}
\end{center}
\caption{Time activity curve of two regions.}
\label{fig:MCTAC}
\end{figure}

Figure \ref{fig:MCTAC} illustrates the comparison of the TACs of two regions. The dash  lines are
the normalized true TACs and the solid lines are the normalized one extracted from the reconstruction images by our
method. The first row are TACs of the events equal to $20000$ and the second row are TACs of $200000$ events.

\begin{figure}[h]
\begin{center}
\subfigure[Noiseless sinogram]{
\includegraphics[width=.45\linewidth]{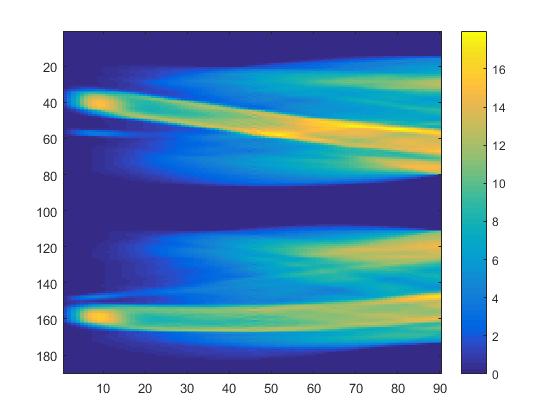}}
\subfigure[sinogram for $\mathrm{events}=200000$]{
\includegraphics[width=.45\linewidth]{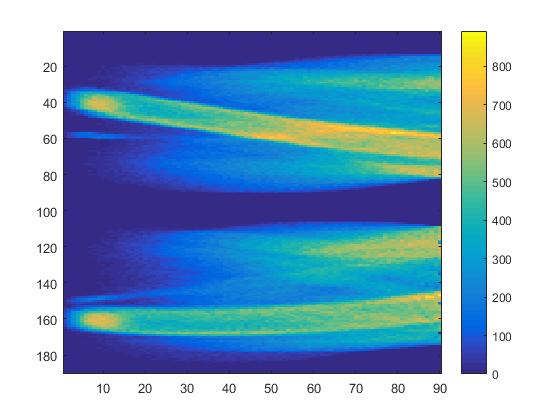}}
\end{center}
\caption{Monte Carlo sinogram data.}
\label{fig:MCLiverSinogram}
\end{figure}

\begin{figure}
\begin{tabular}{c@{\hspace{2pt}}c@{\hspace{2pt}}c@{\hspace{2pt}}c@{\hspace{2pt}}c@{\hspace{2pt}}c@{\hspace{2pt}}c@{\hspace{2pt}}c@{\hspace{2pt}}c@{\hspace{2pt}}c}
\includegraphics[width=.1\linewidth,height=.1\linewidth]{Liver1}&
\includegraphics[width=.1\linewidth,height=.1\linewidth]{Liver11}&
\includegraphics[width=.1\linewidth,height=.1\linewidth]{Liver21}&
\includegraphics[width=.1\linewidth,height=.1\linewidth]{Liver31}&
\includegraphics[width=.1\linewidth,height=.1\linewidth]{Liver41}&
\includegraphics[width=.1\linewidth,height=.1\linewidth]{Liver51}&
\includegraphics[width=.1\linewidth,height=.1\linewidth]{Liver61}&
\includegraphics[width=.1\linewidth,height=.1\linewidth]{Liver71}&
\includegraphics[width=.1\linewidth,height=.1\linewidth]{Liver81}\\
\includegraphics[width=.1\linewidth,height=.1\linewidth]{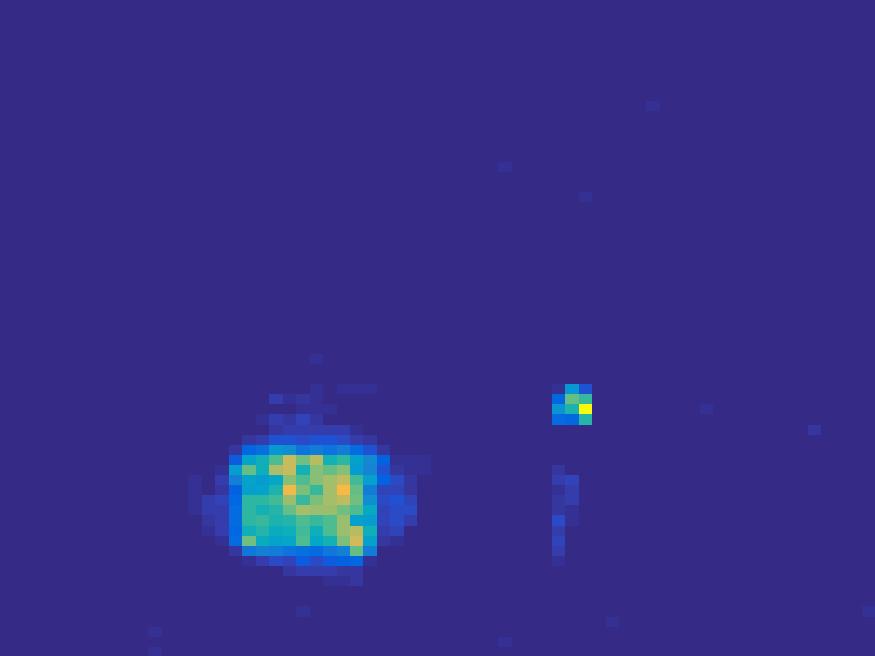}&
\includegraphics[width=.1\linewidth,height=.1\linewidth]{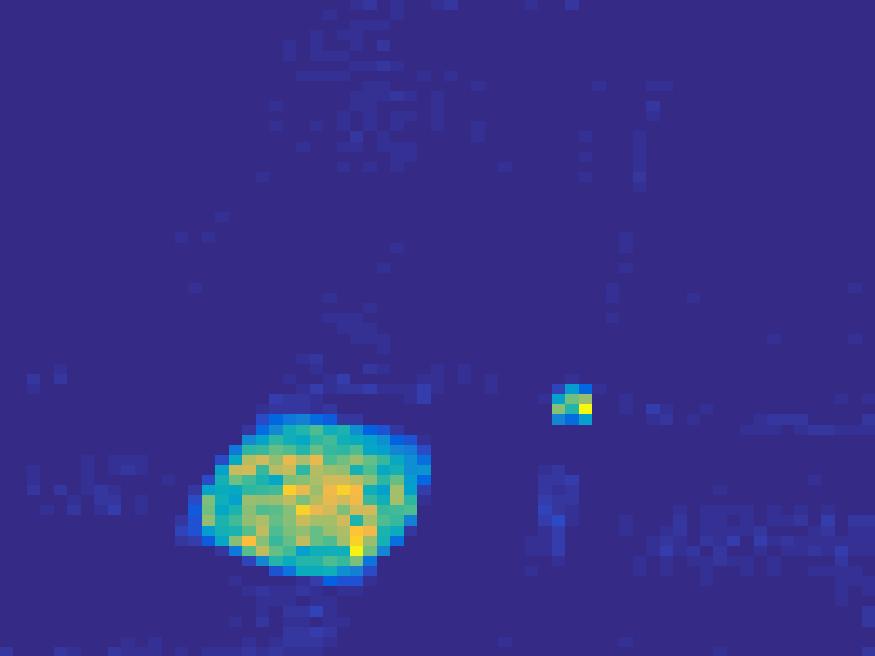}&
\includegraphics[width=.1\linewidth,height=.1\linewidth]{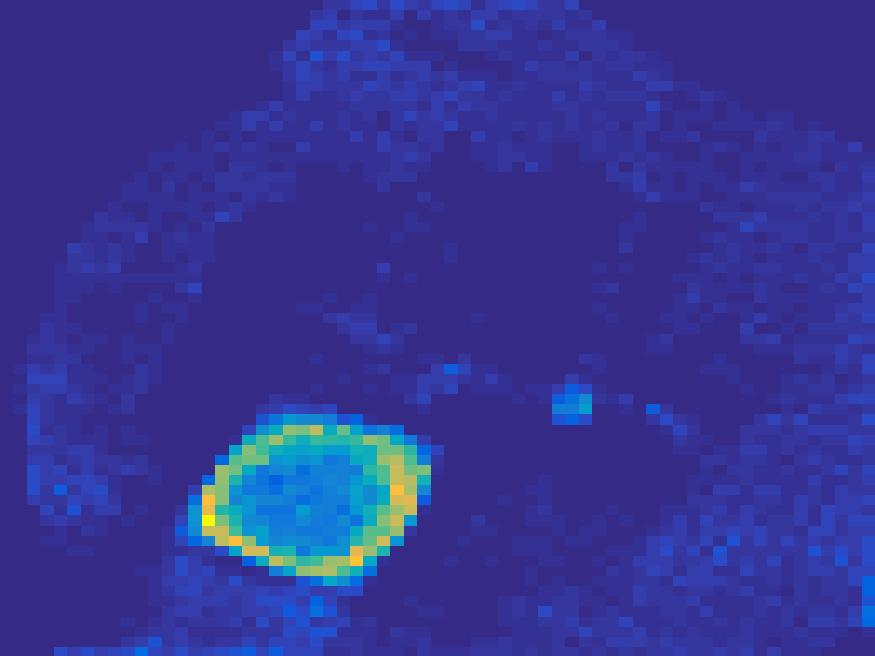}&
\includegraphics[width=.1\linewidth,height=.1\linewidth]{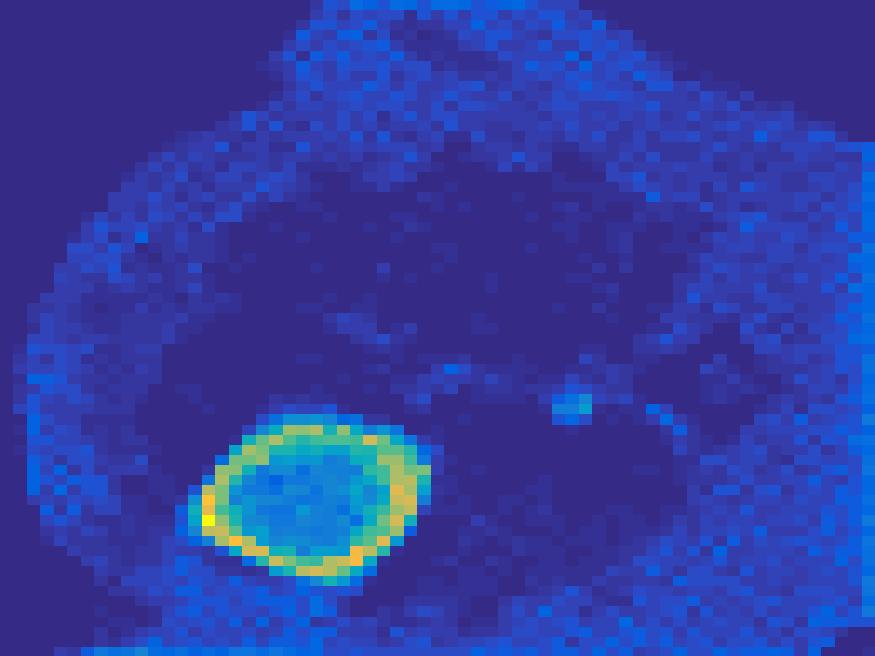}&
\includegraphics[width=.1\linewidth,height=.1\linewidth]{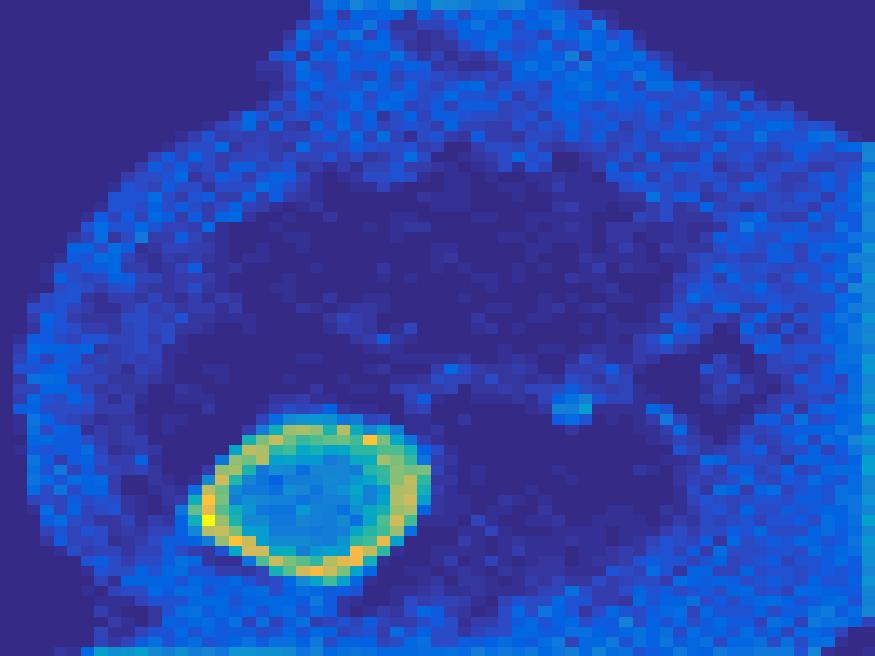}&
\includegraphics[width=.1\linewidth,height=.1\linewidth]{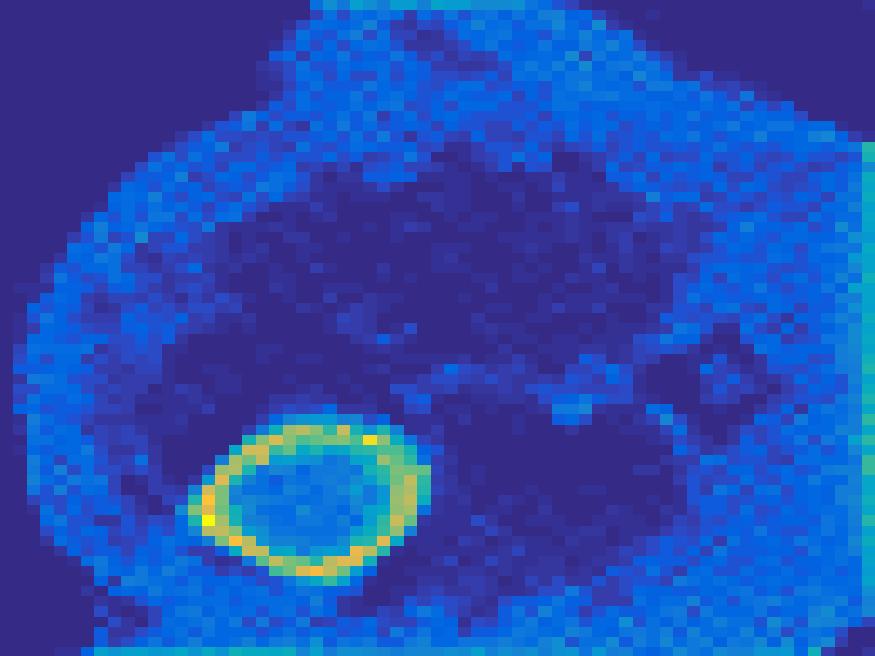}&
\includegraphics[width=.1\linewidth,height=.1\linewidth]{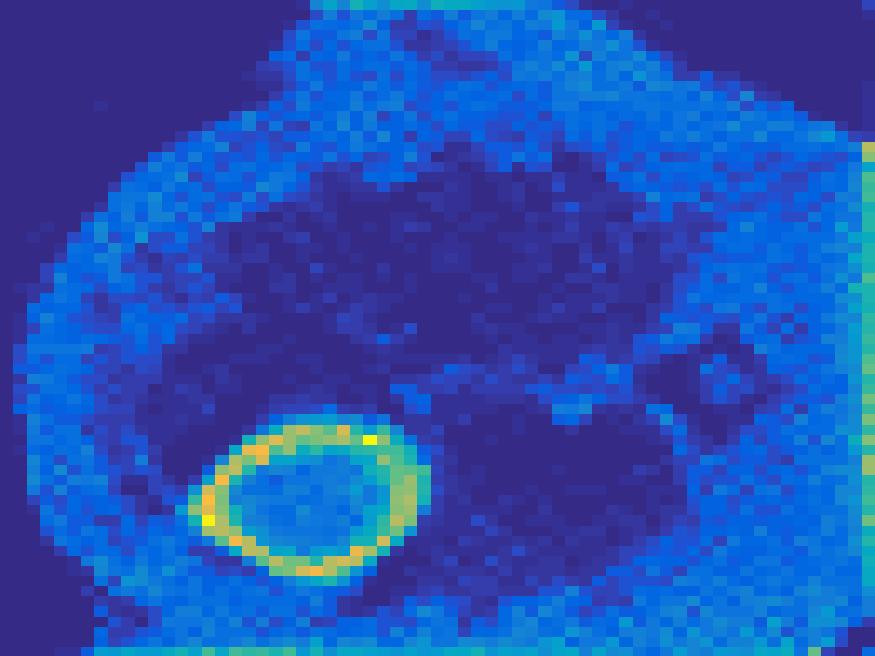}&
\includegraphics[width=.1\linewidth,height=.1\linewidth]{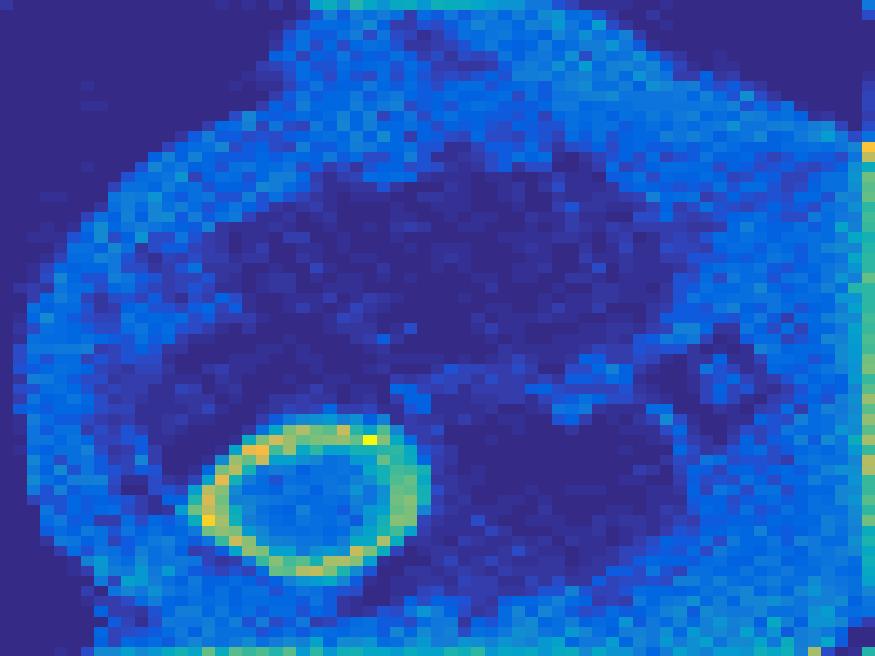}&
\includegraphics[width=.1\linewidth,height=.1\linewidth]{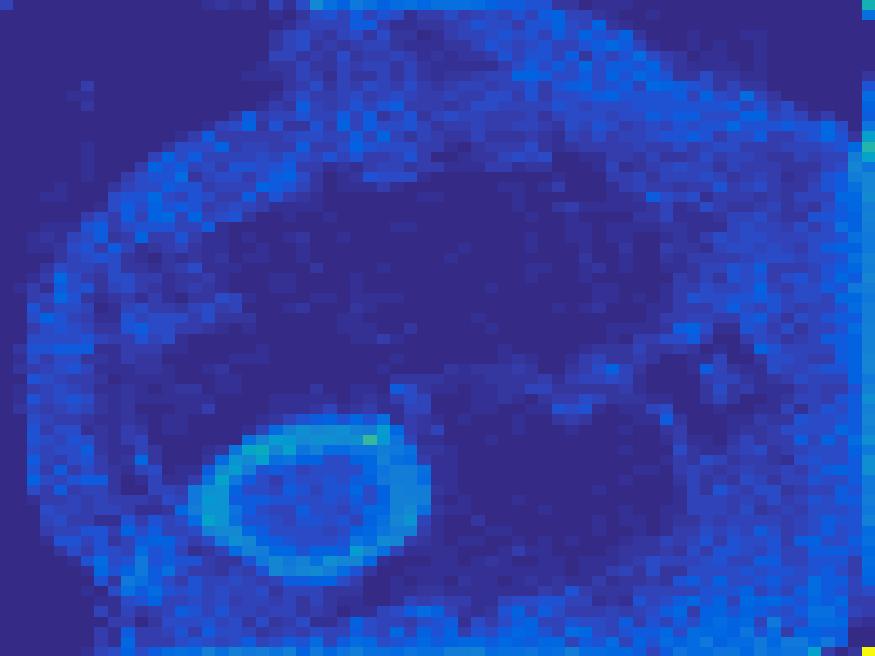}\\
\includegraphics[width=.1\linewidth,height=.1\linewidth]{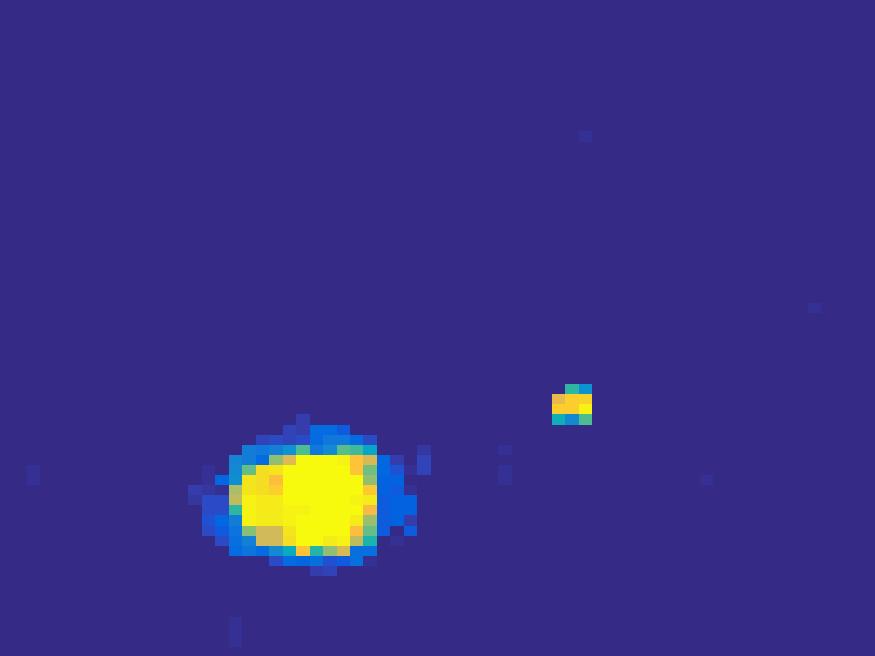}&
\includegraphics[width=.1\linewidth,height=.1\linewidth]{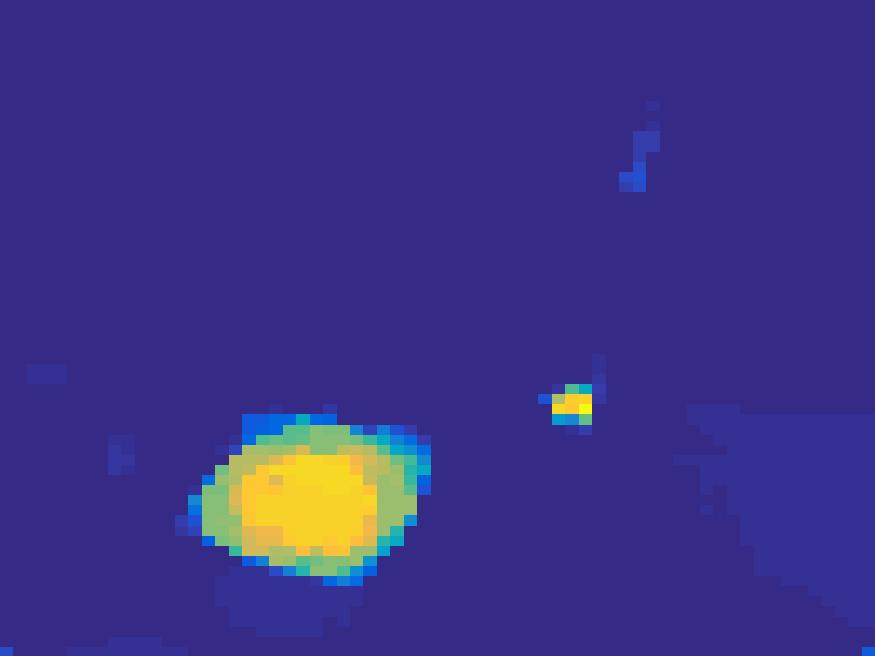}&
\includegraphics[width=.1\linewidth,height=.1\linewidth]{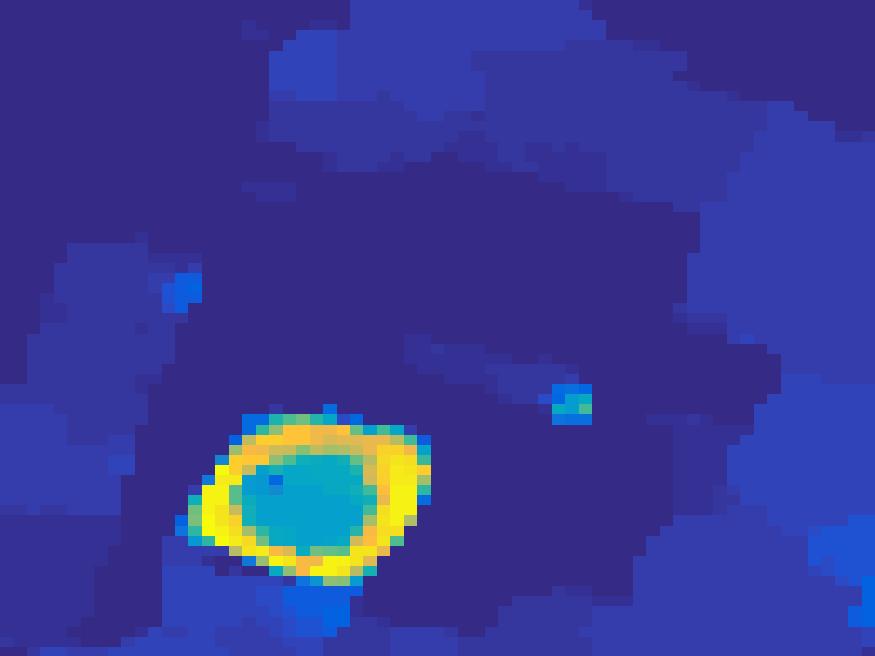}&
\includegraphics[width=.1\linewidth,height=.1\linewidth]{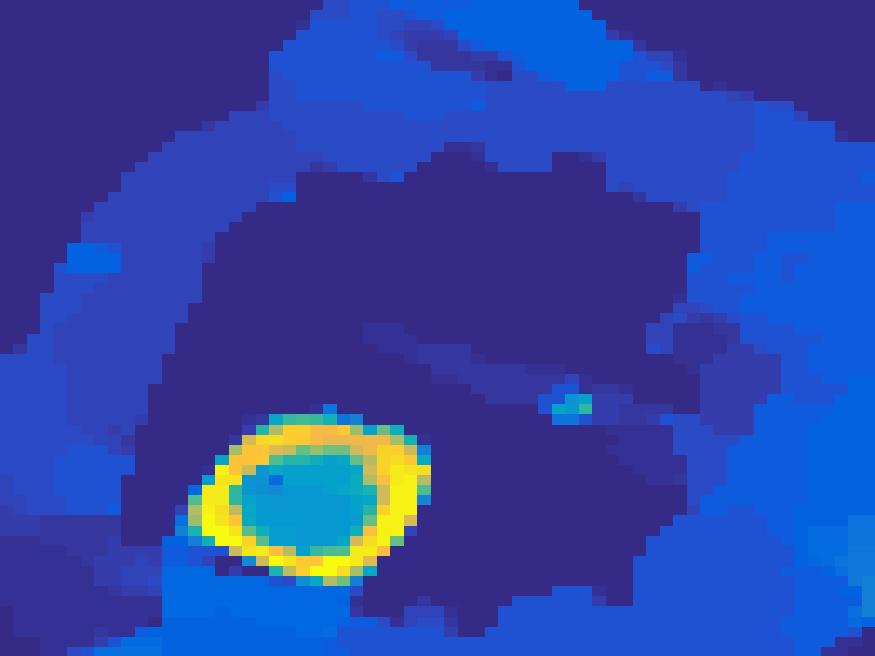}&
\includegraphics[width=.1\linewidth,height=.1\linewidth]{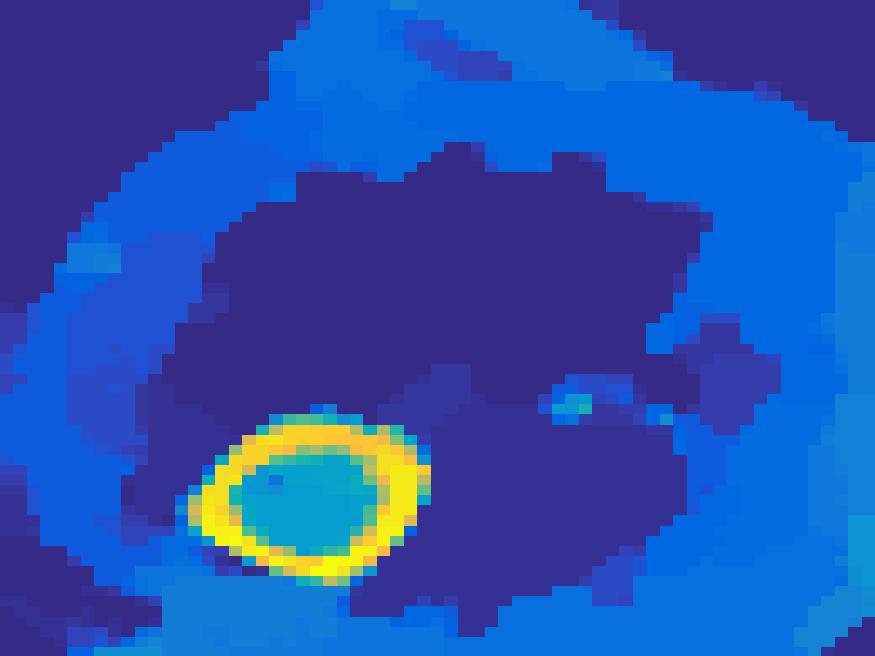}&
\includegraphics[width=.1\linewidth,height=.1\linewidth]{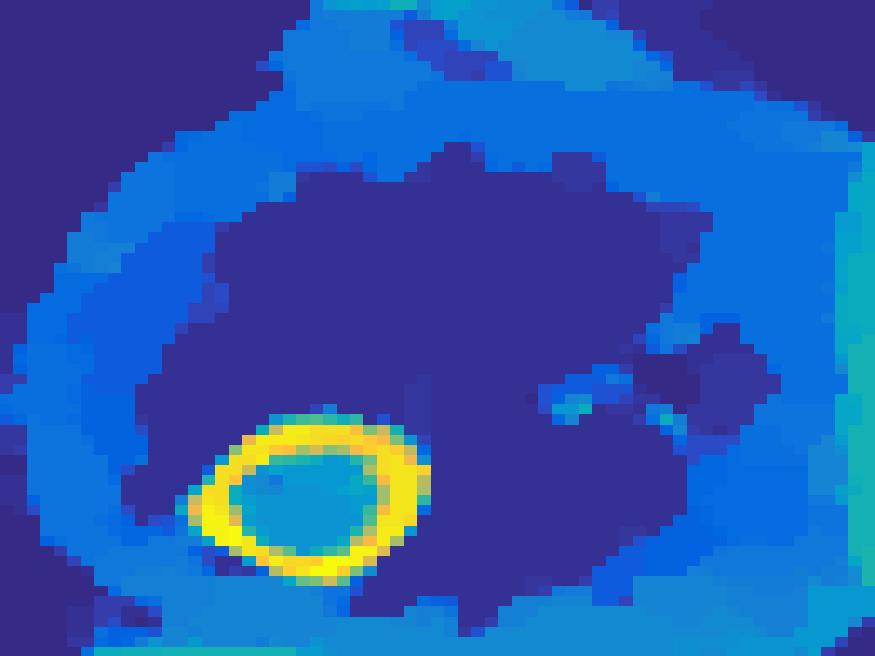}&
\includegraphics[width=.1\linewidth,height=.1\linewidth]{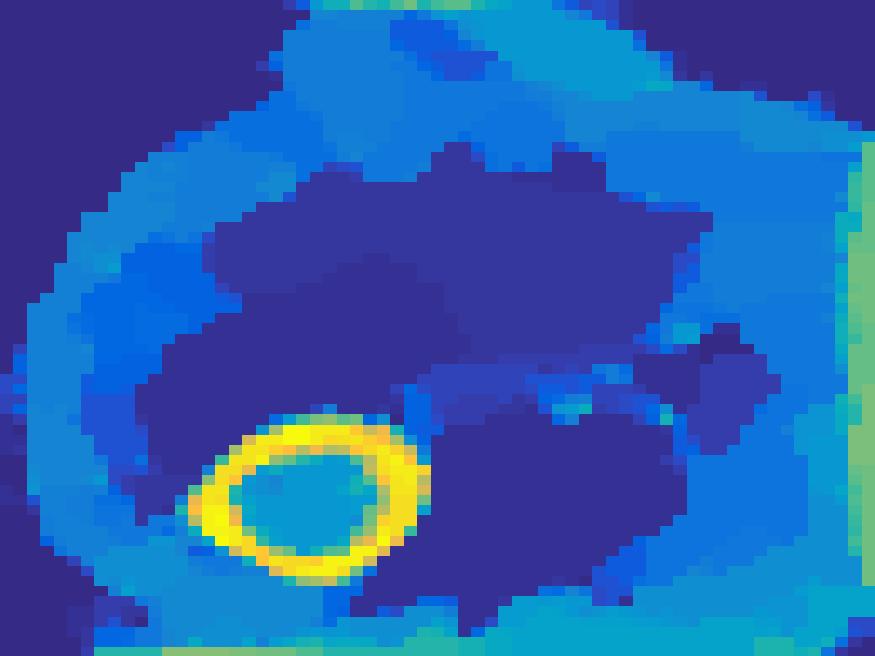}&
\includegraphics[width=.1\linewidth,height=.1\linewidth]{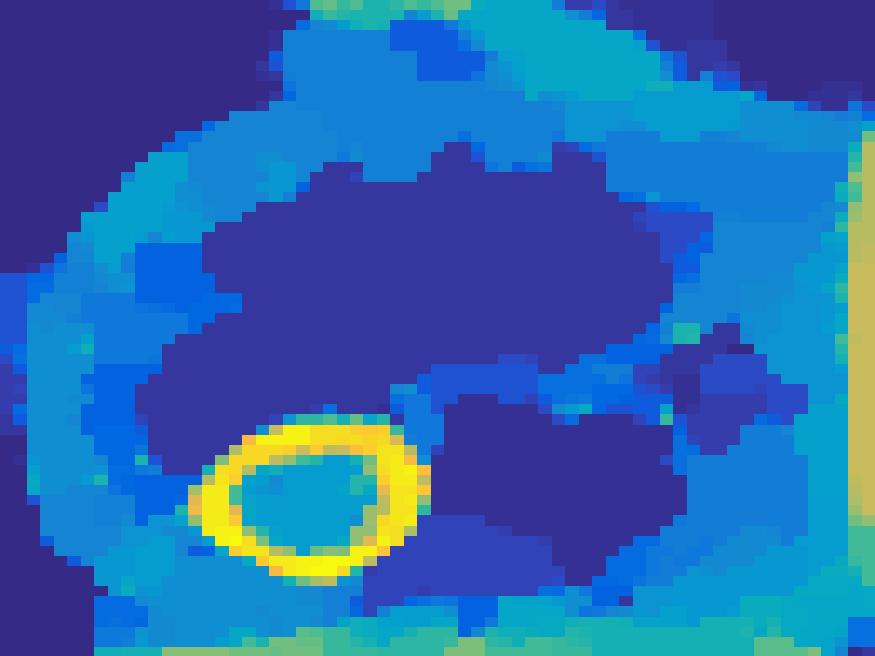}&
\includegraphics[width=.1\linewidth,height=.1\linewidth]{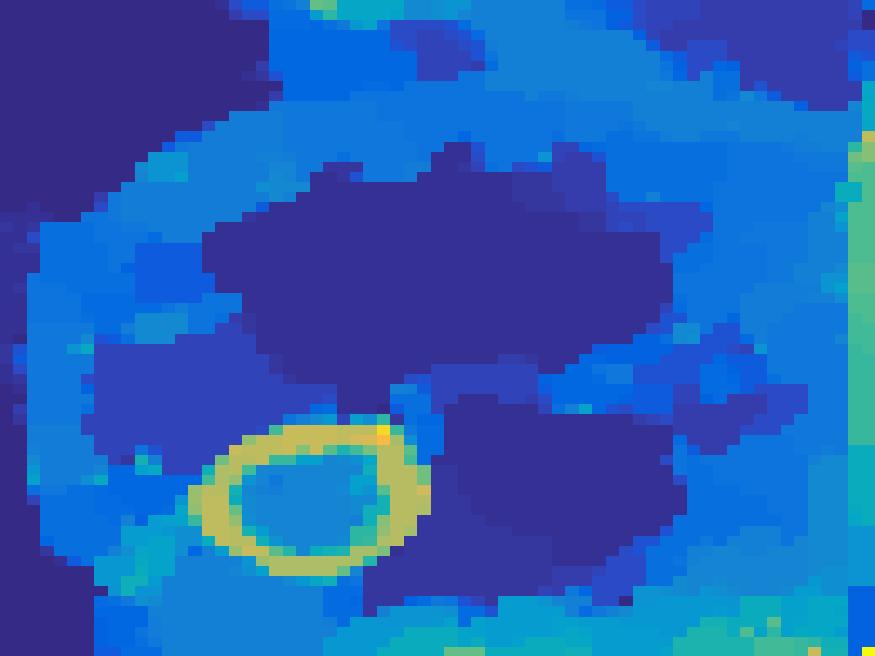}\\
{\footnotesize Frame 1}&
{\footnotesize Frame 11}&
{\footnotesize Frame 21}&
{\footnotesize Frame 31}&
{\footnotesize Frame 41}&
{\footnotesize Frame 51}&
{\footnotesize Frame 61}&
{\footnotesize Frame 71}&
{\footnotesize Frame 81}
\end{tabular}
\caption {Reconstruction of Monte Carlo simulated data. First row: Ground truth; Second row : $\mathrm{events}=200000$, EM with updating $\alpha$ and $B$;  third row: proposed method.
\label{fig:LiverMCIC}
}
\end{figure}

\begin{figure}[ht]
\begin{center}
\subfigure{
\includegraphics[width=.45\linewidth,height=.3\linewidth]{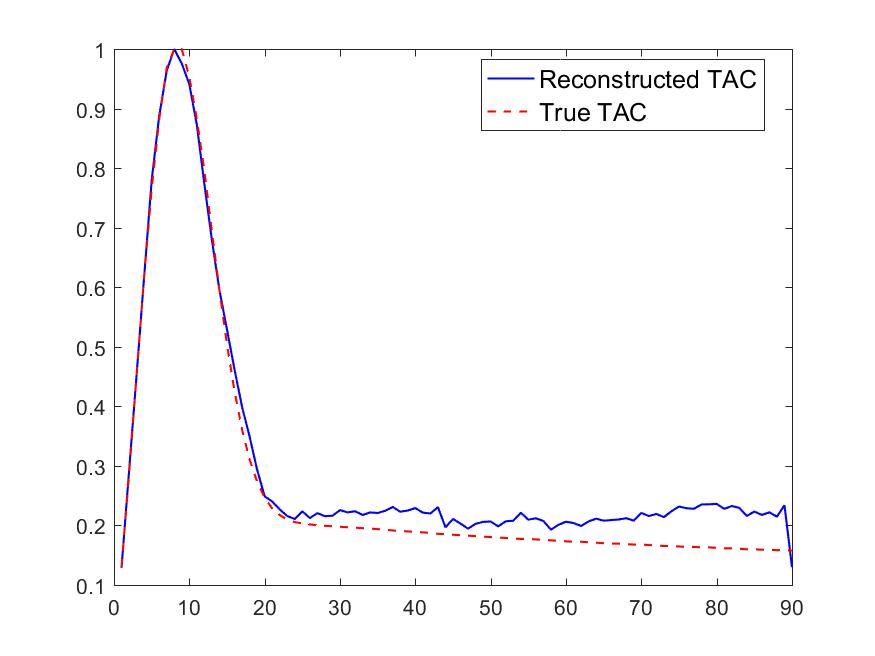}}
\subfigure{
\includegraphics[width=.45\linewidth,height=.3\linewidth]{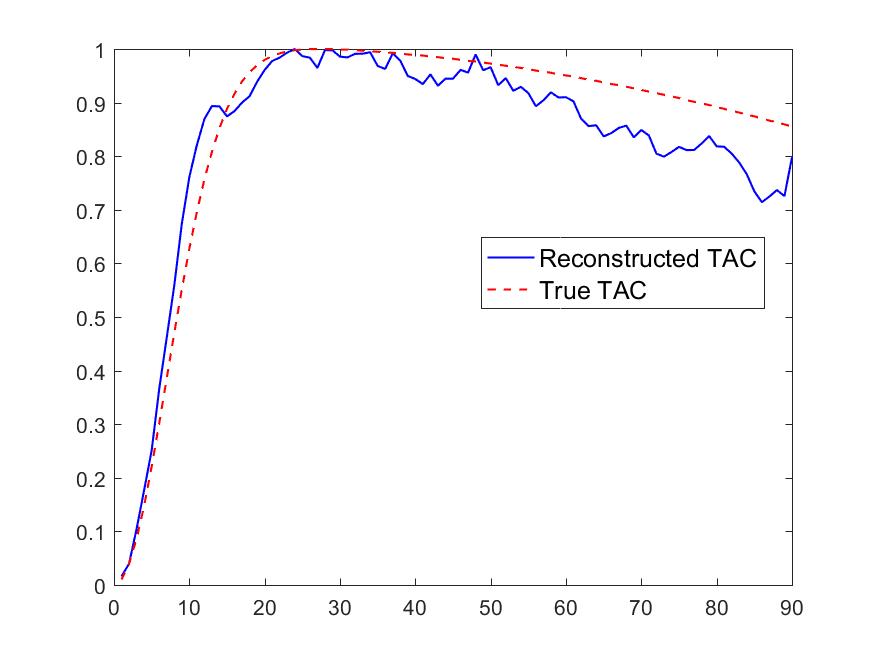}}\\
\subfigure{
\includegraphics[width=.45\linewidth,height=.3\linewidth]{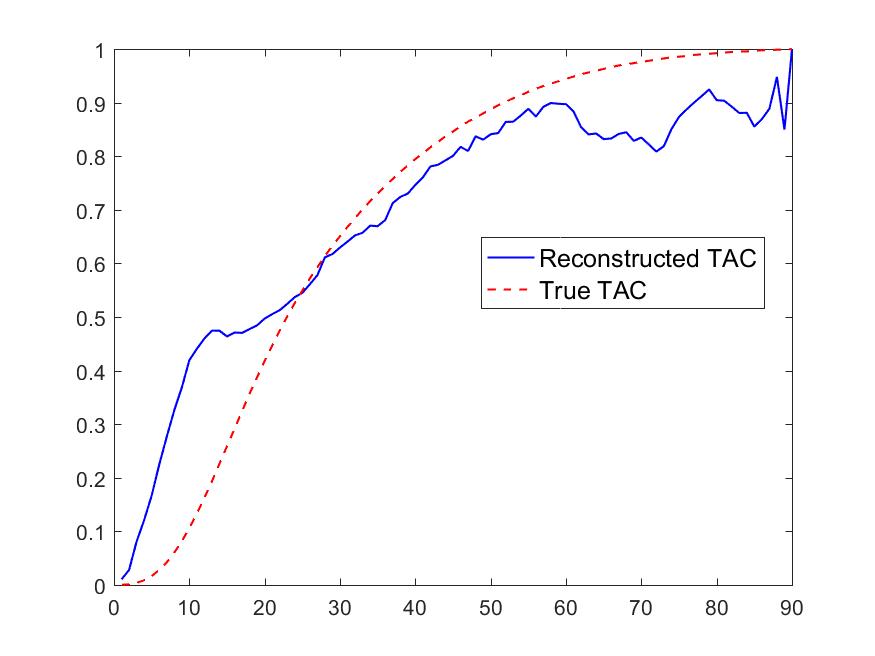}}
\subfigure{
\includegraphics[width=.45\linewidth,height=.3\linewidth]{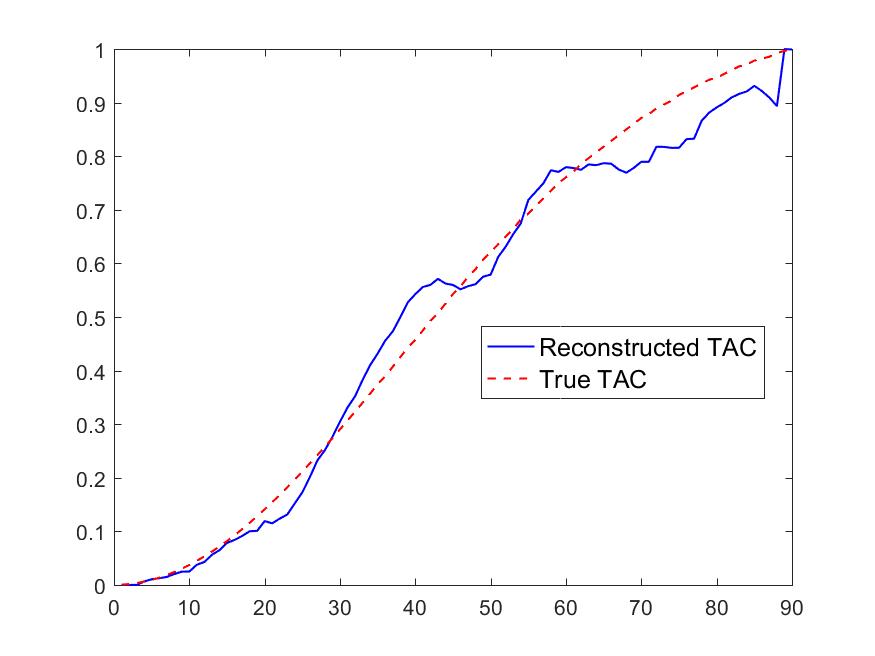}}
\end{center}
\caption{Reconstructed TACs.}
\label{fig:MCLiverTAC}
\end{figure}

We also perform Monte carlo simulation on the more complex images: the rat's abdomen phantom. By  setting $\mathrm{events} = 200000$, the sinogram image is shown in Figure \ref{fig:MCLiverSinogram}.
The images reconstructed can be found in Figure \ref{fig:LiverMCIC} for $\mathrm{events}=200000$.  Figure \ref{fig:MCLiverTAC} illustrates the comparison of the true TACs and those reconstructed by the proposed method. We can see that the proposed method is robust to reconstruct most of the structures present in the images.

\section{Conclusion and Outlook}
\label{sec:Conclusion}
In this paper, we presented a new reconstruction model for dynamic SPECT from few and incomplete projections based on edge correlation.
Both Gaussian noise and Poisson noise are investigated. The proposed  nonconvex model is solved by an alternating scheme. The reconstruction results on two 2D phantoms indicate that our algorithm outperforms the conventional FBP type reconstruction algorithm,  least square/EM method and the former SEMF model. The reconstructed image sequences are very close to the exact ones, especially for those frames with  changed edge directions.
Extensive numerical results show that the choice of the regularization methods as well as the reconstruction approach is effective for a proof of concept study.
Nevertheless, there are still many aspects that needs to be improved in future. Firstly, the method is tested on simulated 2D images with low spatial resolution while real
clinical dynamic SPECT is 3D with higher spatial resolution. Consequently, computation time and acceleration method  should be taken the into account.  Furthermore, the model involves many parameters that needs to be set in a more automatical way. Therefore it is necessary to discuss the parameter choice in a future work.

\section*{Acknowledgements}
XZ and QD are supported in part by Chinese 973 Program (2015CB856000) and National Youth Top-notch Talent program in China. This work has been initiated during a stay of Qiaoqiao Ding in Germany funded by the China Scholarship Council, whose support is gratefully acknowledged. MB acknowledges support by ERC via Grant EU FP 7 - ERC Consolidator Grant 615216 LifeInverse and the German Science Foundation DFG via EXC 1003 Cells in Motion Cluster of Excellence, M\"{u}nster, Germany.

\section*{References}

\end{document}